\documentclass[reqno]{amsart}
\usepackage{amsmath,amssymb}
\usepackage[colorlinks=true,hypertexnames=false,pdftitle={On Minima of 3D Lennard-Jones Lattice Energy and FCC-Tetragonal-Simple Hexagonal Phase Transitions}]{hyperref}
\renewcommand{\qed}{%
\ifmmode\mathqed\else
\leavevmode\unskip\nobreak\hfill\nobreak\hspace{1em}\nobreak\hbox{\qedsymbol}%
\fi}
\let\oldbibliography\thebibliography
\renewcommand{\thebibliography}[1]{%
\oldbibliography{#1}%
\setlength{\itemsep}{0pt}%
}
\numberwithin{equation}{section}
\newtheorem{theorem}{Theorem}[section]
\newtheorem{quotedtheorem}{Theorem}

\newtheorem{conjecture}{Conjecture}
\newtheorem{question}{Question}
\newtheorem{proposition}[theorem]{Proposition}
\newtheorem{lemma}[theorem]{Lemma}
\newtheorem{numericallemma}[theorem]{Numerical lemma}
\newtheorem{corollary}[theorem]{Corollary}
\theoremstyle{remark}
\newtheorem{remark}[theorem]{Remark}
\newcommand{\R}{\mathbb R}
\newcommand{\Z}{\mathbb Z}
\newcommand{\tr}{\operatorname{tr}}
\newcommand{\diag}{\operatorname{diag}}
\newcommand{\dd}{\,\mathrm d}
\newcommand{\FCC}{\mathrm{FCC}}
\newcommand{\BCC}{\mathrm{BCC}}
\newcommand{\SC}{\mathrm{SC}}
\title[3D Lennard--Jones lattice energy and phase transitions]{\bfseries On Minima of 3D Lennard-Jones Lattice Energy and FCC-Tetragonal-Simple Hexagonal Phase Transitions}
\author{Senping Luo}

\author{Juncheng Wei}

\address[S.~Luo]{School of Mathematics and Statistics, Jiangxi Normal University, Nanchang, 330022, China}
\address[J.~Wei]{Department of Mathematics, Chinese University of Hong Kong, Shatin, NT, Hong
Kong}

\email[S.~Luo]{luosp1989@163.com}

\email[J.~Wei]{wei@math.cuhk.edu.hk}

\hypersetup{pdfauthor={Senping Luo and Juncheng Wei}}
\date{}
\begin{document}
\begin{abstract}
Denote the Epstein zeta function by
\begin{equation}\label{eq:Epstein}
E(L,s):=\sum_{v\in L\setminus\{0\}}|v|^{-2s}, \;\;\hbox{where}\;\;L\;\;\hbox{is a lattice};
\qquad s>\frac32.
\end{equation}
We consider the three-dimensional Lennard--Jones lattice energy, which has the form
\[
\mathcal E(L)=\frac12E(L,6)-E(L,3).
\]
We completely classify the minimizers at every covolume $|L|$,
and we prove that
as the covolume increases, the minimizing phase changes from the FCC
to the body-centered tetragonal and then to the simple-hexagonal lattice.
Thereby, we give a positive answer to an open problem raised by Blanc--Lewin
\cite[Section~2.5,~(27), pp.~268--269]{BL}
and B\'etermin--Petrache \cite[Question~1.15]{BP}, and we prove and clarify
the conjecture of B\'etermin--\v{S}amaj--Trav\v{e}nec \cite{BST}.
\end{abstract}
\maketitle

\section{Introduction and statement of main results}

Let $L=g\Z^3\subset\R^3$ be a Bravais lattice of covolume
$|L|=|\det g|$. Lots of physical and mathematical problems can be formulated as the minimization of
the $f$-potential lattice energy as follows:
\begin{equation}\label{EFL}
 E_f(L):=\sum_{v\in L\setminus\{0\}}f(|v|^2),\qquad
 \inf_{|L|=V}E_f(L),\quad V>0.
\end{equation}
Here $|v|$ denotes Euclidean length, $V^{-1}$ is the number density,
and the sum is assumed absolutely convergent. We identify lattices
up to orthogonal equivalence, denoted by $\cong$.

Considering the Riesz potential $f(t)=t^{-s}$ in \eqref{EFL}, we have the Epstein zeta function
\begin{equation}\label{eq:epstein}
 E_f(L)=E(L,s):=\sum_{v\in L\setminus\{0\}}|v|^{-2s},
 \qquad s>\frac32.
\end{equation}

In our previous work \cite[Theorem~1.2]{LW3}, we proved the minima of 3D Epstein zeta functions
conjectured by Sarnak and Str\"ombergsson \cite[Section~5]{SS}. It is stated as follows:
\begin{quotedtheorem}[\cite{LW3}]\label{thm:riesz}
For every Bravais lattice $L\subset\R^3$ of covolume one,
\begin{equation}\label{eq:rieszinput}
 E(L,s)\ge E(\FCC,s),\qquad s>\frac32,
\end{equation}
with equality if and only if $L\cong\FCC$. Here $\FCC$ denotes
the face-centred cubic lattice of covolume one.
\end{quotedtheorem}

In this paper, we consider the widely used Lennard--Jones potential (\cite{Hansen1969,Verlet1967,John1993}) in dimension three:
\begin{equation}\label{LJ}
 f(r^2)=4\varepsilon\left[
 \left(\frac{\sigma}{r}\right)^{12}
 -\left(\frac{\sigma}{r}\right)^6\right],
 \qquad r>0,\quad\varepsilon,\sigma>0.
\end{equation}
The constants $\varepsilon$ and $\sigma$ specify the well depth
and length scale; see .
Measuring lengths in units $2^{1/6}\sigma$ and energies in units
$\varepsilon$ gives $r^{-12}-2r^{-6}$. Under the Lennard-Jones potential, the lattice energy in \eqref{EFL} becomes
\begin{equation}\label{eq:energy}
 \mathcal E(L)=\frac12E(L,6)-E(L,3),\qquad
 m(V)=\inf_{|L|=V}\mathcal E(L).
\end{equation}
After normalization to covolume one, we have
\[
 2V^4m(V)=\inf_{|L|=1}\{E(L,6)-2V^2E(L,3)\},\qquad
 \inf_{a>0}\mathcal E(aL)=-\frac{E(L,3)^2}{2E(L,6)}.
\]

The Lennard--Jones potential is a basic model for intermolecular
forces, clusters and condensed matter, and has a long history; see a recent review by Schwerdtfeger-Wales
\cite{SW100} in
\emph{100 Years of the Lennard-Jones Potential}.
The crystal-energy problem dates back to Jones-Ingham's 1925
study of minimum-energy cubic crystals \cite{JI1925}.
There are tons of numerical simulations of Lennard--Jones energy; see, for example in
\cite{SGKTM,BCS} and the following citing references.

To state our main result, we introduce various important lattices and their notations.
We use the unit-covolume FCC, body-centred cubic (BCC) and simple
cubic lattices
\[
 D_3=\{n\in\Z^3:n_1+n_2+n_3\in2\Z\},\qquad
 F=\FCC=2^{-1/3}D_3,
\]
\[
 \BCC=2^{1/3}\bigl[\Z^3\cup
       (\Z^3+(\tfrac12,\tfrac12,\tfrac12))\bigr],\qquad
 \SC=\Z^3.
\]
See B\'etermin--\v{S}amaj--Trav\v{e}nec \cite{BST}.
The triangular and simple-hexagonal families are denoted by
\[
 T_A=\sqrt{\frac{2A}{\sqrt3}}\,
       \operatorname{span}_{\Z}
       \left\{(1,0),\left(\tfrac12,\tfrac{\sqrt3}2\right)\right\},
 \qquad A>0,
\]
\[
 H_Y=\left(\frac4{3Y}\right)^{1/6}
 \operatorname{span}_{\Z}
 \left\{(1,0,0),\left(\tfrac12,\tfrac{\sqrt3}2,0\right),
                         (0,0,\sqrt Y)\right\},\qquad Y>0.
\]
Then $|T_A|=A$, $|L_X|=|H_Y|=1$, and
\[
 V^{1/3}H_Y=T_A\oplus(V/A)\Z,\qquad
 A=\left(\frac{\sqrt3 V^2}{2Y}\right)^{1/3}.
\]
We shall also introduce a deformation family of the FCC lattice, namely, the unit-covolume body-centered tetragonal family:
\[
 L_X=(4X)^{-1/6}\diag(\sqrt X,1,1)D_3,\qquad X>0;
 \qquad L_1=F,\quad L_{1/2}\cong\BCC.
\]

The three-dimensional Lennard-Jones lattice problem \eqref{eq:energy} attracts lots of attentions over the years.
B\'etermin--\v{S}amaj--Trav\v{e}nec \cite{BST} made a conjecture for general Lennard-Jones potential by numerical simulations.
\begin{conjecture}[B\'etermin--\v{S}amaj--Trav\v{e}nec {\cite[Conjecture~2]{BST}}]\label{conj:BST}
Assume that $p>q>3$. Then there exist $V_{\mathrm{pq}}>0$ such that
\[
\operatorname*{argmin}_{|L|=V}\frac{qE(L,p/2)-pE(L,q/2)}{2(p-q)}=
 \begin{cases}
 \{F\},&0<V<V_{\mathrm{pq}},\\
 \{H_{Y(V)}\},&V>V_{\mathrm{pq}}.
 \end{cases}
\]

\end{conjecture}
For the classical $(12,6)$ potential, B\'etermin--\v{S}amaj--Trav\v{e}nec \cite{BST}
predicted a phase transition from the FCC to the simple-hexagonal lattice as the
inverse density increases:
\begin{equation}\label{eq:FCC-SH}
\mathrm{FCC}\ \longrightarrow\ \mathrm{SH}\qquad (V\uparrow).
\end{equation}

The three-dimensional Lennard-Jones lattice problem \eqref{eq:energy} is also listed as an open problem by by Blanc-Lewin
\cite[Section~2.5, (27), pp.~268--269]{BL}.
It is also listed as an open problem as follows:
\begin{question}[B\'etermin--Petrache {\cite[Question~1.15]{BP}}]\label{que:BP}
For the classical potential $f(r^2)=r^{-12}-r^{-6}$, which
unit-covolume Bravais lattice shapes minimize $E_f(\lambda L)$
for each $\lambda>0$?
\end{question}

The change of scale between this question and \eqref{eq:energy} is
\[
 E_f(\lambda L)=\frac12\mathcal E(2^{-1/6}\lambda L),\qquad
 |L|=1,\quad V=\frac{\lambda^3}{\sqrt2}.
\]

We completely classify the minima of problem \eqref{eq:energy}, thereby
providing an positive answer to Question~\ref{que:BP}, on open problem of
Blanc--Lewin on 3D Lennard-Jones lattice energy (\cite[Section~2.5, (27), pp.~268--269]{BL}), and identifying the intervening tetragonal phase
 beyond Conjecture \ref{conj:BST}.

\begin{theorem}[Minima of 3D Lennard-Jones lattice energy]\label{thm:main}
Let $|L|=V$. Up to orthogonal equivalence, there exist two thresholds $0<V_{\mathrm{FT}}<V_{\mathrm{TH}}<\infty$ such that:
\[
 \operatorname*{argmin}_{|L|=V}\Big(\frac12E(L,6)-E(L,3)\Big)=
 \begin{cases}
 \{FCC\},&0<V<V_{\mathrm{FT}},\\
 \{FCC,L_{X(V)}\},&V=V_{\mathrm{FT}},\\
 \{L_{X(V)}\},&V_{\mathrm{FT}}<V<V_{\mathrm{TH}},\\
 \{L_{X(V)},H_{Y(V)}\},&V=V_{\mathrm{TH}},\\
 \{H_{Y(V)}\},&V>V_{\mathrm{TH}}.
 \end{cases}
\]
Numerically,
\[
 V_{\mathrm{FT}}=1.0648235668\cdots,\qquad
 V_{\mathrm{TH}}=1.0759906200\cdots.
\]
Analytically,
the thresholds $V_{\mathrm{FT}}$ and $V_{\mathrm{TH}}$ are the
$v$-components of the unique solutions $(v,x)$ and $(v,x,y)$ of,
respectively,
\begin{equation}\label{eq:FTthreshold}
 \begin{cases}
 \displaystyle\frac{\dd}{\dd x}
       [E(L_x,6)-2v^2E(L_x,3)]=0,\\[3pt]
 E(L_x,6)-E(F,6)=2v^2[E(L_x,3)-E(F,3)],
 \end{cases}
\end{equation}
\begin{equation}\label{eq:THthreshold}
 \begin{cases}
 \displaystyle\frac{\dd}{\dd x}
       [E(L_x,6)-2v^2E(L_x,3)]=0,\\[3pt]
 \displaystyle\frac{\dd}{\dd y}
       [E(H_y,6)-2v^2E(H_y,3)]=0,\\[3pt]
 E(H_y,6)-E(L_x,6)=2v^2[E(H_y,3)-E(L_x,3)],
 \end{cases}
\end{equation}
in the ranges
\[
 \frac{21}{20}<v<\frac{11}{10},\qquad
 \frac{23}{100}<x<\frac{307}{1000},\qquad
 \frac{103}{100}<y<\frac{103}{50}.
\]
Alternatively,
the first transition can be located as follows:
\begin{equation}\label{eq:beta}
 V_{\mathrm{FT}}^2=\frac12
 \inf_{\substack{|L|=1\\L\not\cong F}}
       \frac{E(L,6)-E(F,6)}{E(L,3)-E(F,3)}.
\end{equation}

Furthermore, the parameters $X(V)$ and $Y(V)$ can be located by solving equations.
\begin{enumerate}
\item[\textup{(i)}]
For $21/20\le V\le11/10$, the parameters $X(V)$ and $Y(V)$ are
uniquely determined by
\begin{align}
 \frac{\dd E(L_X,6)}{\dd\log X}
 &=2V^2\frac{\dd E(L_X,3)}{\dd\log X},\label{eq:Tbranch}\\
 \frac{\dd E(H_Y,6)}{\dd\log Y}
 &=2V^2\frac{\dd E(H_Y,3)}{\dd\log Y},\label{eq:Hbranch}
\end{align}
in $23/100<X<307/1000$ and $103/100<Y<103/50$, respectively,
with $X'(V)<0<Y'(V)$.

\item[\textup{(ii)}]
For $11/10\le V\le89/50$, $Y(V)$ is the unique solution of
\eqref{eq:Hbranch} in $11/10<Y<217/50$,

\item[\textup{(iii)}] For $V\ge89/50$, the simple-hexagonal minimizer is
\begin{equation*}
 L_V=T_{A(V)}\oplus\frac{V}{A(V)}\Z,
 \qquad Y(V)=\frac{\sqrt3 V^2}{2A(V)^3},
\end{equation*}
where $A(V)$ is the unique solution of
\[
 \frac{\dd}{\dd A}\mathcal E\bigl(T_A\oplus(V/A)\Z\bigr)=0,
 \qquad \frac25<A<\frac{89}{100}.
\]
As $V\to\infty$,
\begin{align}
 A(V)&=A_*+\frac{\pi\zeta(4)S_6A_*^4}{3S_3^2}V^{-4}
                    +O(V^{-8}),\label{eq:areaexp}\\
 m(V)&=e_2-\pi\zeta(4)A_*^3V^{-4}
       -\frac{\pi^2\zeta(4)^2S_6A_*^6}{2S_3^2}V^{-8}
       +\frac{\pi\zeta(10)A_*^9}{5}V^{-10}
       +O(V^{-12}).\label{eq:energyexp}
\end{align}
Here,
\begin{equation*}
  S_s=\sum_{\substack{(i,j)\in\Z^2\\(i,j)\ne(0,0)}}
       (i^2+ij+j^2)^{-s}\quad(s=3,6), \;A_*=\frac{\sqrt3}2\left(\frac{S_6}{S_3}\right)^{1/3},\;
 e_2=-\frac{S_3^2}{2S_6}.
\end{equation*}

\item[\textup{(iv)}]
 $m(V)$ in \eqref{eq:energy} increases on $[1,\infty)$.

\end{enumerate}
\end{theorem}

For the Lennard-Jones potential with free parameter $b$ ($r^{-12}-br^{-6}$, $b>0$), the corresponding result follows from
\[
 \mathcal E_b(L)=\frac{b^2}{4}\mathcal E\bigl((b/2)^{1/6}L\bigr),\qquad
 V_{\mathrm{FT}}(b)=\sqrt{2/b}\,V_{\mathrm{FT}},\quad
 V_{\mathrm{TH}}(b)=\sqrt{2/b}\,V_{\mathrm{TH}}.
\]
By Theorem \ref{thm:main}, we provide a complete solution the seemingly simple three-dimensional Lennard-Jones lattice problem \eqref{eq:energy}.

In dimension two, Montgomery \cite{Mont} proved triangular
optimality for every Gaussian scale. For 2D Lennard--Jones
energy, we refer the to B\'etermin and his collaborators  \cite{Bet2016,BZ,B2,BP} for pioneering works; our previous work \cite{LW2}
gives the complete classification. In three dimensions, some local minimal and stability results can be found in
\cite{BCubic,BP,BPdim}; see also \cite{BCS}
for cuboidal BCC instability.

In dimensions 8 and 24,
 Cohn-Kumar-Miller-Radchenko-Viazovska
\cite{CKMRV} proved universal optimality of $E_8$ and the Leech lattice.
Petrache-Serfaty \cite{PetracheSerfatyCK}
transferred universal optimality to Coulomb and Riesz
renormalized and jellium energies ,
obtaining the $E_8$ and Leech minima in dimensions $8$ and $24$.

There are many related crystallization results, and we refer to Blanc-Lewin \cite{BL} and Lewin \cite{LewinRiesz} for reviews.
 Sandier-Serfaty
\cite{SandierSerfaty} proved triangular optimality among lattices
of fixed density for their established Coulomb renormalized energy; Chen-Oshita \cite{Che2007} proved triangular lattice in di-block copolymers. Petrache-Serfaty \cite{PetracheSerfatyRiesz}
identified the renormalized energy governing the microscopic
next-order asymptotics of confined Riesz systems.
 Heitmann-Radin \cite{HR} proved triangular crystallization for sticky disks;
Theil \cite{Theil} proved a triangular crystallization, and E-Li
\cite{ELi} proved a hexagonal crystallization for a three-body potential. Flatley-Theil \cite{FTcryst} proved
FCC crystallization for a class of three-dimensional two- and
three-body energies.

We now discuss some ideas in the proof of Theorem \ref{thm:main}.
Theorem~\ref{thm:riesz} provides the $s=3$ comparison used for the
high-density phase. For the dilute phase, we use B\'etermin's
 theorem \cite[Theorem~1.1]{B2} and its fixed-area
comparison \cite[Proposition~1.2]{B2}: for every unit-area planar
lattice $\Gamma$,
\begin{equation}\label{eq:planarinput}
 E(\Gamma,6)-E(T_1,6)
       \ge2\bigl(E(\Gamma,3)-E(T_1,3)\bigr).
\end{equation}

The proof uses a planar reduction for large covolumes and a profile
comparison on the remaining finite covolume intervals. For large $V$, write
$L=\Gamma+\Z(u+h\nu)$, $A=|\Gamma|$, $Ah=V$. The planar comparison
and Fourier inequalities give
\[
 \mathcal E(L)\ge \mathcal E\bigl(T_A\oplus(V/A)\Z\bigr),\qquad
 \partial_A^2\mathcal E\bigl(T_A\oplus(V/A)\Z\bigr)>0
\]
on the localized area interval. On each remaining covolume interval
$[V_0,V_1]$, let $f(V)$ be the least candidate value of
$E(L,6)-2V^2E(L,3)$. Then
\[
 g(t)=\max_{V_0\le V\le V_1}\{f(V)+2V^2t\},\qquad
 E(L,6)-g(E(L,3))\ge0\quad(|L|=1).
\]
Every profile minimum minimizes $E(L,6)-2V^2E(L,3)$ at each
supporting $V$. Moment bounds give compactness; Fourier and
stationary estimates determine the equality cases.
Comparison \eqref{eq:rieszinput} extends the result to $V<21/20$.

We collect some preliminaries in
Section~2. Section~3 deals with the dilute phase.
Section~\ref{sec:compact} establishes the endpoint comparisons;
Sections~\ref{sec:transition} and~\ref{sec:longphase} prove the
transition and intervening simple-hexagonal phases. The long but finite computations are collected in the Appendix
.

\section{Preliminaries and stationary branches}

Packing, convexity and Ewald estimates give existence, fixed-systole
comparisons and the stationary branches.

For $1\le\operatorname{rank}L\le3$, put
$\lambda=\lambda_1(L):=\min_{v\in L\setminus\{0\}}|v|>0$.
Disjoint balls of radius $\lambda/2$ give
\[
 \#(L\cap B_R)\le(1+2R/\lambda)^3.
\]
Summation over the annuli \(2^j\lambda\le |v|<2^{j+1}\lambda\) yields
\begin{equation}\label{eq:packing}
 E(L,s)\le\frac{125}{1-2^{3-2s}}\lambda^{-2s}
                  \quad(s>3/2).
\end{equation}
For \(\lambda\ge\rho>0\) and \(R\ge\rho\), similarly,
\begin{equation}\label{eq:uniformtail}
 \sum_{\substack{v\in L\\|v|>R}} |v|^{-2s}
 \le\frac{125}{1-2^{3-2s}}\rho^{-3}R^{3-2s}.
\end{equation}
Integration gives
\[
 E(L,3)=6\int_\lambda^\infty
       \#\{v\in L:0<|v|\le t\}\,t^{-7}\dd t
 \le6\int_\lambda^\infty
       \bigl[(1+2t/\lambda)^3-1\bigr]t^{-7}\dd t
 =\frac{206}{5}\lambda^{-6}.
\]
Hence
\begin{equation}\label{eq:coercive}
 \mathcal E(L)\ge \lambda^{-12}-\frac{206}{5}\lambda^{-6}.
\end{equation}
Thus
\[
 \inf_L\mathcal E(L)>-\infty,\qquad
 \sup_n\mathcal E(L_n)<\infty\ \Longrightarrow\ \inf_n\lambda_1(L_n)>0.
\]

\begin{lemma}[A quantitative separation bound]\label{lem:separation}
For a lattice $L\subset\R^3$,
\[
 \mathcal E(L)\le-3\quad\Longrightarrow\quad
 \min_{v\in L\setminus\{0\}}|v|>\frac7{10}.
\]
\end{lemma}

\begin{proof}
Suppose $0<\lambda\le7/10$ and $|v_0|=\lambda$, $v_0\in L$.
Disjoint radius-$\lambda/2$ balls give
\[
 \#\{v\in L\setminus\{0,\pm v_0\}:|v|\le r\}
 \le (1+2r/\lambda)^3-3\qquad(r\ge1).
\]
For every \(r>0\),
\[
 r^{-12}-2r^{-6}\ge
 -12\int_{\max\{1,r\}}^\infty(t^{-7}-t^{-13})\,dt.
\]
Tonelli's theorem gives
\[
\begin{aligned}
 \mathcal E(L)&\ge\lambda^{-12}-2\lambda^{-6}
 -6\int_1^\infty(t^{-7}-t^{-13})
                  \bigl((1+2t/\lambda)^3-3\bigr)\,dt\\
 &=\lambda^{-12}-2\lambda^{-6}+1
       -\frac{216}{55\lambda}-\frac{54}{5\lambda^2}
       -\frac{32}{3\lambda^3}.
\end{aligned}
\]

For $x=1/\lambda$, the $x^{-5}$-weighted derivative is
\[
 12x^6-12-\frac{216}{55x^5}-\frac{108}{5x^4}-\frac{32}{x^3},
\]
strictly increasing on $(0,\infty)$ and positive at $10/7$.
The energy lower bound at $x=10/7$ is $>-3$; hence
$\mathcal E(L)>-3$, contradicting $\mathcal E(L)\le-3$.
This completes the proof of Lemma~\ref{lem:separation}.
\end{proof}

Mahler's compactness criterion \cite{Mahler} states that, for
$a>0$ and $0<V_1\le V_2<\infty$, the set
\[
 \left\{g\mathbb Z^3:\ g\in \mathrm{GL}_3(\mathbb R),\quad
 V_1\le|\det g|\le V_2,\quad
 \min_{n\in\mathbb Z^3\setminus\{0\}}|gn|\ge a\right\}
\]
is compact in $\mathrm{GL}_3(\mathbb R)/\mathrm{GL}_3(\mathbb Z)$.
Equivalently, every sequence $g_j\mathbb Z^3$ in this set admits
\[
 j_k\uparrow\infty,\qquad A_k\in \mathrm{GL}_3(\mathbb Z),\qquad
 g_{j_k}A_k\longrightarrow g\in \mathrm{GL}_3(\mathbb R).
\]
By \eqref{eq:coercive}, \eqref{eq:uniformtail} and this criterion,
$m(V)$ is attained. The chain and planar minima are
\[
 e_1=-\frac{\zeta(6)^2}{\zeta(12)},\qquad
 e_2=\mathcal E(T_{A_*})=-\frac{S_3^2}{2S_6};
\]
the second identity and planar uniqueness follow from \cite{B2}.
The following fixed-systole comparison extends the prescribed-bond
viewpoint of \cite{Bonds} and excludes large planar areas.

\begin{theorem}[A convex fixed-systole comparison]
\label{thm:convex-fixed-systole}
Fix \(a>0\). Suppose that
\(f:[a^2,\infty)\longrightarrow(0,\infty)\) is continuous,
strictly decreasing and strictly convex, and that for some
\(s>3/2\) and \(C<\infty\),
\[
 f(q)\le Cq^{-s}\qquad(q\ge a^2).
\]
Let \(F_1=2^{-1/2}D_3\) denote the FCC lattice of shortest
length one. Every three-dimensional lattice with
\(\lambda_1(L)\ge a\) satisfies
\[
 \sum_{v\in L\setminus\{0\}}f(|v|^2)
 \le\sum_{v\in aF_1\setminus\{0\}}f(|v|^2),
\]
with equality if and only if \(L\) is an orthogonal image of
\(aF_1\).
\end{theorem}

\begin{proof}
For a discrete additive subgroup \(L\) with separation at least \(a\),
disjoint balls of radius \(a/2\) give, for \(R>a/2\),
\begin{equation}\label{eq:fixedsystole-tail}
 \sum_{\substack{v\in L\\|v|>R}}|v|^{-2s}
 \le T_s(a,R):=\frac{24}{a^3(2s-3)}
       \left(1+\frac a{2R}\right)^{2s}
                         (R-a/2)^{3-2s}.
\end{equation}
Let $M_k$ be the supremum over separated additive subgroups of rank
$\le k$, $0\le k\le3$, with $M_0=0$. Point counts give local
sequential compactness; separation, the group operations and vanishing
$(k+1)$-fold exterior products are closed. The tail
\eqref{eq:fixedsystole-tail} gives continuity and attainment.
Adjoining an orthogonal period of length $a$ to an $M_{k-1}$
maximizer preserves separation and adds positive terms, so
\[
 M_k>M_{k-1},\qquad k=1,2,3,
\]
and every maximizer has rank $k$.

Scaling by $a/\lambda_1(L)<1$ strictly increases the sum. Thus
$\lambda_1(L)=a$, and its Gram matrix $Q$ satisfies
\[
 n^TQn\ge a^2\quad(n\in\Z^3\setminus\{0\}),\qquad
 \operatorname{span}\{nn^T:n^TQn=a^2\}=\operatorname{Sym}_3.
\]
Otherwise choose a nonzero symmetric $D$ orthogonal to the span and
\[
 \epsilon\|Q^{-1/2}DQ^{-1/2}\|_{\rm op}<\frac12,\qquad
 0<\epsilon<\min_{\substack{a^2<n^TQn<2a^2\\n^TDn\ne0}}
       \frac{n^TQn-a^2}{|n^TDn|},
\]
where the minimum of an empty set is $+\infty$.
The second inequality handles $n^TQn<2a^2$; otherwise,
\[
 n^TQn\ge2a^2\Longrightarrow
 n^T(Q\pm\epsilon D)n
 \ge\bigl(1-\epsilon\|Q^{-1/2}DQ^{-1/2}\|_{\rm op}\bigr)n^TQn>a^2.
\]
Thus $Q\pm\epsilon D\succ0$ preserve all constraints, with the active ones fixed.

Since $D\ne0$, some coordinate vector or pairwise sum satisfies
$n^TDn\ne0$. Strict convexity gives
\begin{align*}
 \sum_{n\in\Z^3\setminus\{0\}}f(n^TQn)
 &<\frac12\sum_{n\in\Z^3\setminus\{0\}}f(n^T(Q+\epsilon D)n)\\
 &\quad+\frac12\sum_{n\in\Z^3\setminus\{0\}}f(n^T(Q-\epsilon D)n).
\end{align*}
Convergence follows from the tail hypothesis. One strict term
contradicts maximality, so at least six shortest-vector pairs occur.

For a superbase $v_0+v_1+v_2+v_3=0$, every three being a basis,
and $\{i,j,k,l\}=\{0,1,2,3\}$ with $v_i\cdot v_j>0$,
Selling's superbase move satisfies
\[
 (\widetilde v_i,\widetilde v_j,\widetilde v_k,\widetilde v_l)
 =(-v_i,v_j,v_k+v_i,v_l+v_i),\quad
 \sum_{m=0}^3|\widetilde v_m|^2
 =\sum_{m=0}^3|v_m|^2-2v_i\cdot v_j<\sum_{m=0}^3|v_m|^2.
\]
Finiteness under the squared-length bound gives termination at
$v_i\cdot v_j\le0$.
Scale $a$ to one and put $p_{ij}=-v_i\cdot v_j\ge0$.
The seven cut weights are
\[
 \begin{gathered}
 p_{01}+p_{02}+p_{03},\quad p_{01}+p_{12}+p_{13},\\
 p_{02}+p_{12}+p_{23},\quad p_{03}+p_{13}+p_{23},\\
 p_{02}+p_{03}+p_{12}+p_{13},\quad
 p_{01}+p_{03}+p_{12}+p_{23},\quad
 p_{01}+p_{02}+p_{13}+p_{23}.
 \end{gathered}
\]
Every cut weight is at least one. Integer level sets give
\[
 \left|\sum_{i=0}^3n_iv_i\right|^2
 =\sum_{0\le i<j\le3}p_{ij}(n_i-n_j)^2
 \ge\sum_{k=\min_{0\le i\le3}n_i+1}^{\max_{0\le i\le3}n_i}
       \sum_{\substack{0\le i,j\le3\\n_i\ge k>n_j}}p_{ij}
 \ge\max_{0\le i\le3}n_i-\min_{0\le i\le3}n_i.
\]
Thus shortest vectors are binary cuts, with at least six weights
equal to one. The first four and last three weights both sum to
$2\sum_{0\le i<j\le3}p_{ij}$; hence not all seven equal one. If the missing
equality were a singleton, its weight would be zero: the other
three singleton weights and the last three weights both sum to three.

Relabeling the missing cut as the last one, the other six equations give
\[
 p_{01}=p_{02}=p_{13}=p_{23}=\frac12,\qquad p_{03}=p_{12}=0.
\]
In the basis $(v_0,v_3,v_1)$, the Gram matrix is
\[
 \begin{pmatrix}1&0&-1/2\\0&1&-1/2\\-1/2&-1/2&1\end{pmatrix},
\]
the shortest-length-one FCC Gram matrix. Thus $L\cong aF_1$.
This completes the proof of Theorem~\ref{thm:convex-fixed-systole}.
\end{proof}

\begin{theorem}[Fixed-systole consequences]
\label{thm:fccsystolemaximum}\label{cor:lj-convex-systole}
For $\lambda_1(L)\ge a>0$,
\[
 E(L,s)\le a^{-2s}E(F_1,s)\quad(s>3/2),\qquad
 \mathcal E(L)\ge\frac{E(F_1,6)}{2a^{12}}-\frac{E(F_1,3)}{a^6}
 \quad(a^6\ge7/4).
\]
In either inequality equality holds exactly for $L\cong aF_1$.
\end{theorem}
\begin{proof}
Apply Theorem~\ref{thm:convex-fixed-systole} to $q^{-s}$ and
$f(q)=2q^{-3}-q^{-6}$, respectively. On $q\ge a^2$, $a^6\ge7/4$,
\[
 0<f(q)\le2q^{-3},\qquad f'(q)=6q^{-7}(1-q^3)<0,\qquad
 f''(q)=6q^{-8}(4q^3-7)\ge0.
\]
Thus $f$ is strictly convex also at $a^6=7/4$, and
$\mathcal E(L)=-\frac12\sum_{v\in L\setminus\{0\}}f(|v|^2)$ gives the second
comparison and equality case.
This completes the proof of Theorem~\ref{thm:fccsystolemaximum}.
\end{proof}

\begin{lemma}[An explicit fixed-systole Riesz constant]
\label{lem:fccsystoleconstant}
For every three-dimensional lattice \(L\),
\[
 E(L,3)<\frac{2901}{200}\lambda_1(L)^{-6},\qquad
 \mathcal E(L)>\lambda_1(L)^{-12}
                  -\frac{2901}{200}\lambda_1(L)^{-6}.
\]
\end{lemma}
\begin{proof}
Scale to shortest length one and apply
Theorem~\ref{thm:fccsystolemaximum}. With $2^{-1/6}F=2^{-1/2}D_3$,
rational summation gives
\[
 \sum_{\substack{(i,j,k)\in\{-11,\ldots,11\}^3\\i+j+k\in2\mathbb Z\\
                     0<i^2+j^2+k^2\le128}}
 \frac8{(i^2+j^2+k^2)^3}<\frac{14443}{1000}.
\]
There are $3054$ retained points; the others have radius $>8$.
By \eqref{eq:fixedsystole-tail} with $(a,s,R)=(1,3,8)$,
\[
 E(2^{-1/6}F,3)<\frac{14443}{1000}+T_3(1,8)<\frac{2901}{200}.
\]
Thus $E(L,6)/2\ge\lambda_1(L)^{-12}$.
This completes the proof of Lemma~\ref{lem:fccsystoleconstant}.
\end{proof}

\begin{lemma}[A bounded primitive planar section]\label{lem:shortplane}
Let $\lambda_2(L)$ be the least radius containing two independent
lattice vectors. Every rank-three lattice with \(\mathcal E(L)\le e_2\) has a primitive
rank-two section \(\Gamma\) satisfying
\[
 \frac7{10}<\lambda_1(L)<\frac85,\qquad
 \lambda_2(L)\le\frac{25}{8},\qquad |\Gamma|<5.
\]
Writing \(L=\Gamma+\mathbb Z(u+he_3)\), one has \(h>V/5\).
\end{lemma}
\begin{proof}
\[
 e_2<\sum_{\substack{(i,j)\in\mathbb Z^2\\0<i^2+ij+j^2\le4}}
 \left(\frac1{2(i^2+ij+j^2)^6}-\frac1{(i^2+ij+j^2)^3}\right)
 <-\frac{331}{100}.
\]
Lemma~\ref{lem:separation} gives \(\lambda_1>7/10\).
If \(\lambda_1\ge8/5\), then
\[
 \mathcal E(L)\ge-\frac{206}{5}\lambda_1^{-6}>-\frac{331}{100},
\]
a contradiction. If \(\lambda_2>25/8\), the radius-\(25/8\)
points lie on one primitive chain. Its partial energy is at least
the full energy, which is \(\ge e_1\): omitted chain distances
exceed one and their energies are negative.
Disjoint balls of radius \(7/20\) give
\[
 \sum_{\substack{v\in L\\|v|>25/8}}|v|^{-6}
 \le\frac{8000}{343}\left(\frac{139}{125}\right)^6
       \left(\frac{111}{40}\right)^{-3}.
\]
Together with $e_1>-21/20$, obtained from $\zeta(6)<327/320$, this gives
\[
 \mathcal E(L)\ge e_1-\sum_{\substack{v\in L\\|v|>25/8}}|v|^{-6}>-\frac{331}{100}>e_2,
\]
again a contradiction.

Gauss reduction of
$\Gamma=L\cap\operatorname{span}_{\R}\{v_1,v_2\}$ for shortest
independent $v_1,v_2$ gives basis lengths $\lambda_1$ and
$\le\lambda_2$. Hence
\[
 |\Gamma|\le\lambda_1\lambda_2<5,\qquad
 h=V/|\Gamma|>V/5.
\]
The saturated section $\Gamma$ is primitive and extends to a basis.
This completes the proof of Lemma~\ref{lem:shortplane}.
\end{proof}

\begin{lemma}[The stationary branches]\label{lem:scalarfamilyvalues}
For $21/20\le V\le11/10$, equations \eqref{eq:Tbranch} and
\eqref{eq:Hbranch} have unique real analytic solutions $X(V),Y(V)$
in the stated intervals, with $X'(V)<0<Y'(V)$.
Let
\begin{align*}
 \mathcal D_{TF}(V)&=2V^4
       [\mathcal E(V^{1/3}L_{X(V)})-\mathcal E(V^{1/3}F)],\\
 \mathcal D_{TH}(V)&=2V^4
       [\mathcal E(V^{1/3}L_{X(V)})-\mathcal E(V^{1/3}H_{Y(V)})].
\end{align*}
Then
\[
 \mathcal D_{TF}'(V)<-\frac{16}{25}V,\qquad
 \mathcal D_{TH}'(V)>\frac15V.
\]
\end{lemma}
\begin{proof}
For $s=3,6$, termwise differentiation gives
\begin{align*}
 \frac{\dd}{\dd x}E(L_x,s)
 &=\frac{s(4x)^{s/3}}{3x}
 \sum_{n\in D_3\setminus\{0\}}
 \frac{n_2^2+n_3^2-2xn_1^2}
      {(xn_1^2+n_2^2+n_3^2)^{s+1}},\\
 \frac{\dd}{\dd y}E(H_y,s)
 &=\frac{s(3y/4)^{s/3}}{3y}
 \sum_{(i,j,k)\in\Z^3\setminus\{0\}}
 \frac{i^2+ij+j^2-2yk^2}
      {(i^2+ij+j^2+yk^2)^{s+1}}.
\end{align*}
The series and their fixed-order derivatives converge uniformly on
the stated compact positive intervals.

Numerical Lemma~\ref{num:branchmoments} gives
\[
\begin{aligned}
 \partial_{\log X}E(L_X,3)&<0,\qquad
 \partial_{\log X}\!\left(
  \frac{\partial_{\log X}E(L_X,6)}{\partial_{\log X}E(L_X,3)}
                         \right)<0,\\
 \partial_{\log Y}E(H_Y,3)&>0,\qquad
 \partial_{\log Y}\!\left(
  \frac{\partial_{\log Y}E(H_Y,6)}{\partial_{\log Y}E(H_Y,3)}
                         \right)>0.
\end{aligned}
\]
The endpoint ratios bracket $2V^2$. At either root, with
$L(t)=L_{e^t}$ or $H_{e^t}$,
\[
 \partial_t^2\{E(L(t),6)-2V^2E(L(t),3)\}
 =\partial_tE(L(t),3)\,\partial_t\!\left(
       \frac{\partial_tE(L(t),6)}{\partial_tE(L(t),3)}\right)>0.
\]
Implicit differentiation gives analyticity and $X'(V)<0<Y'(V)$.

Stationarity and homogeneity give
\[
 \frac1{4V}
 \begin{pmatrix}\mathcal D_{TF}'(V)\\\mathcal D_{TH}'(V)\end{pmatrix}
 =\begin{pmatrix}
 E(F,3)-E(L_{X(V)},3)\\
 E(H_{Y(V)},3)-E(L_{X(V)},3)
 \end{pmatrix}.
\]
The certified moment bounds are
\[
 E(F,3)<7.226961<7.39056<E(L_{X(V)},3)
                  <7.88582<7.94245<E(H_{Y(V)},3).
\]
These gaps prove the derivative bounds.
This completes the proof of Lemma~\ref{lem:scalarfamilyvalues}.
\end{proof}

\begin{proposition}[The complete comparison within two lattice families]
\label{prop:familyphases}
For $21/20\le V\le11/10$, the branch equations
\eqref{eq:Tbranch}--\eqref{eq:Hbranch} have unique real analytic
solutions, and the contact systems \eqref{eq:FTthreshold} and
\eqref{eq:THthreshold} have unique solutions in the stated rectangles:
\[
\begin{gathered}
 \frac{23}{100}<X(V)<\frac{307}{1000},\qquad
 \frac{103}{100}<Y(V)<\frac{103}{50},\\
 (v,x)=(V_{\mathrm{FT}},X(V_{\mathrm{FT}})),\qquad
 (v,x,y)=(V_{\mathrm{TH}},X(V_{\mathrm{TH}}),Y(V_{\mathrm{TH}})),\\
 \frac{21}{20}<V_{\mathrm{FT}}<V_{\mathrm{TH}}<\frac{11}{10}.
\end{gathered}
\]
Put $\widetilde m(V)=\min_{X,Y>0}
\{\mathcal E(V^{1/3}L_X),\mathcal E(V^{1/3}H_Y)\}$.
Up to orthogonal equivalence,
\[
 \operatorname*{argmin}_{\substack{\Lambda=L_X,\ X>0\\
                         \text{or }\Lambda=H_Y,\ Y>0}}
 \mathcal E(V^{1/3}\Lambda)=
 \begin{cases}
 \{F\},&V<V_{\mathrm{FT}},\\
 \{F,L_{X(V)}\},&V=V_{\mathrm{FT}},\\
 \{L_{X(V)}\},&V_{\mathrm{FT}}<V<V_{\mathrm{TH}},\\
 \{L_{X(V)},H_{Y(V)}\},&V=V_{\mathrm{TH}},\\
 \{H_{Y(V)}\},&V>V_{\mathrm{TH}},
 \end{cases}
\]
and
\[
 \widetilde m'_+(V)-\widetilde m'_-(V)<0
 \qquad(V\in\{V_{\mathrm{FT}},V_{\mathrm{TH}}\}).
\]
\end{proposition}
\begin{proof}
Lemma~\ref{lem:scalarfamilyvalues} gives the two analytic branches.
For the remaining parameters, use
\begin{align*}
 \mathcal E(V^{1/3}F)&<\frac{31}{20}\left(\frac{20}{21}\right)^4
                  -\frac{36}{5}\left(\frac{10}{11}\right)^2<-4,\\
 \mathcal E(V^{1/3}H_1)&<\frac52\left(\frac{20}{21}\right)^4
                  -\frac{15}{2}\left(\frac{10}{11}\right)^2<-4.
\end{align*}
At a family minimum, Lemma~\ref{lem:separation} gives
$\lambda_1>7/10$, whereas
\begin{align*}
 \lambda_1(V^{1/3}L_X)^2
 &=V^{2/3}(4X)^{-1/3}\min(4X,X+1,2),\\
 \lambda_1(V^{1/3}H_Y)^2
 &=V^{2/3}(4/(3Y))^{1/3}\min(1,Y).
\end{align*}
Hence $77/1000<X<21$ and $27/100<Y<14$.
Lemma~\ref{num:fullfamilysigns} excludes every critical point
outside the two branch windows and the two cubic windows
$12/25\le X\le13/25$, $49/50\le X\le51/50$.
Cubic symmetry fixes $X=1/2,1$; the strict second-derivative signs
make BCC and FCC the unique stationary points in their windows,
a maximum and a minimum, respectively.
Thus $\operatorname*{argmin}_{X>0}\mathcal E(V^{1/3}L_X)\subset\{X(V),1\}$
and $\operatorname*{argmin}_{Y>0}\mathcal E(V^{1/3}H_Y)=\{Y(V)\}$.

In the stated rectangles,
\[
\begin{aligned}
 \eqref{eq:FTthreshold}&\Longleftrightarrow
 x=X(v),\quad\mathcal D_{TF}(v)=0,\\
 \eqref{eq:THthreshold}&\Longleftrightarrow
 (x,y)=(X(v),Y(v)),\quad\mathcal D_{TH}(v)=0.
\end{aligned}
\]
The opposite signs in Numerical Lemma~\ref{num:branchcontacts}
and the derivative bounds in Lemma~\ref{lem:scalarfamilyvalues}
give unique contacts in disjoint brackets and the phase ordering.
At either contact $V$, stationarity and equality of energies give
the derivative jump
\[
 -2V^{-3}\{E(\Lambda_+,3)-E(\Lambda_-,3)\}<0,
\]
where $\Lambda_-$ and $\Lambda_+$ are the adjacent unit-covolume shapes.
This completes the proof of Proposition~\ref{prop:familyphases}.
\end{proof}

\begin{lemma}[Ewald identity and rational kernel bounds]\label{lem:firstloss}
For a positive-definite $3\times3$ Gram matrix $Q$ with
$\det Q=V^2$, write
\[
 E(Q,s)=\sum_{n\in\Z^3\setminus\{0\}}(n^{\mathsf T}Qn)^{-s}
 \quad(s>3/2),\qquad
 \mathcal E(Q)=\tfrac12E(Q,6)-E(Q,3).
\]
Set
\[
 \mathcal G_a(x)=\int_1^\infty e^{-xt}t^{a-1}\,dt.
\]
Splitting Poisson summation at $t=1$ gives, absolutely convergently,
\begin{align*}
 \mathcal E(Q)={}&
 \frac{\pi^6}{240}\left(-\frac16+\frac{2}{9V}\right)
 -\frac{\pi^3}{2}\left(-\frac13+\frac{2}{3V}\right)\\
 &+\sum_{n\in\Z^3\setminus\{0\}}\left\{\frac{\pi^6}{240}\mathcal G_6(\pi n^{\mathsf T}Qn)
                -\frac{\pi^3}{2}\mathcal G_3(\pi n^{\mathsf T}Qn)\right\}\\
 &+\frac1V\sum_{n\in\Z^3\setminus\{0\}}\left\{\frac{\pi^6}{240}\mathcal G_{-9/2}(\pi n^{\mathsf T}Q^{-1}n)
                -\frac{\pi^3}{2}\mathcal G_{-3/2}(\pi n^{\mathsf T}Q^{-1}n)\right\}.
\end{align*}
The following rational enclosures hold for every kernel.
\end{lemma}
\begin{proof}
Use the alternating-series enclosures of
$\pi=16\arctan(1/5)-4\arctan(1/239)$, with $37$ and $11$
terms, respectively.
For $x\ge0$ and the least $m\ge0$ with $x/2^m\le1/2$,
\[
 \left(\sum_{j=0}^{33}\frac{(-x/2^m)^j}{j!}\right)^{2^m}
 \le e^{-x}\le
 \left(\sum_{j=0}^{32}\frac{(-x/2^m)^j}{j!}\right)^{2^m}.
\]
For integers $a\ge1$,
\[
 \mathcal G_a(x)=e^{-x}
 \sum_{j=0}^{a-1}\frac{(a-1)!}{j!x^{a-j}}.
\]
For $x,k>0$ and $M\in\mathbb Z_{\ge1}$, define the symmetric tridiagonal matrix
\[
 (J_{M,k})_{jj}=k+2j\quad(0\le j<M),\qquad
 (J_{M,k})_{j,j+1}^{\,2}=(j+1)(j+k)\quad(0\le j<M-1).
\]
For $e_0=(1,0,\ldots,0)^{\mathsf T}\in\R^M$, its rational resolvent is
\begin{gather*}
 L_{M,k}(x)=e_0^{\mathsf T}(xI+J_{M,k})^{-1}e_0=d_0^{-1},\\
 d_{M-1}=x+k+2(M-1),\qquad
 d_j=x+k+2j-\frac{(j+1)(j+k)}{d_{j+1}}\quad(j=M-2,\ldots,0).
\end{gather*}
For every $k>0$,
\[
 e^{-x}L_{M,k}(x)<\mathcal G_{1-k}(x)
 <\frac{e^{-x}}x\{1-kL_{M,k+1}(x)\}.
\]
Indeed,
\[
 e^x\mathcal G_{1-k}(x)=\frac1{\Gamma(k)}\int_0^\infty
               \frac{t^{k-1}e^{-t}}{x+t}\,dt.
\]
For the monic degree-$M$ gamma-weight orthogonal polynomial $P$,
Gauss quadrature is exact on $(1-P(t)^2/P(-x)^2)/(x+t)$; hence
\[
 e^x\mathcal G_{1-k}(x)-L_{M,k}(x)
 =\frac1{P(-x)^2\Gamma(k)}\int_0^\infty
       \frac{P(t)^2t^{k-1}e^{-t}}{x+t}\dd t>0.
\]
The upper bound follows from
$xe^x\mathcal G_{1-k}(x)=1-ke^x\mathcal G_{-k}(x)$.
For $M=48,32,24,16$ and $x\in[x_-,x_+]$, the bounds for
$e^x\mathcal G_{1-k}(x)$ are
$L_{M,k}(x_+)$ and $\{1-kL_{M,k+1}(x_-)\}/x_-$, respectively.
Use interval arithmetic for signed kernel combinations.

For $N\in\mathbb Z_{\ge0}$ and $b>0$, put
\[
 \Theta_N(b)=\frac{6e^{-\pi b(N+1)^2}(1+b^{-1/2})^2}
                         {1-e^{-\pi b(2N+3)}}.
\]
The coordinate union bound and geometric summation give
\[
 \sum_{\substack{n\in\Z^3\\\max_{1\le i\le3}|n_i|>N}}e^{-\pi b|n|^2}\le\Theta_N(b),\qquad
 \sum_{j\in\Z}e^{-\pi bj^2}\le1+b^{-1/2}.
\]
For the cube $N=5$, take $a\le\lambda_{\min}(Q)$ and
$a^*\le\lambda_{\min}(Q^{-1})$.
The omitted energy is bounded below and above, respectively, by
\[
\begin{gathered}
 -\frac{\pi^3}{2}\left(\frac1{36\pi a}+\frac2{(36\pi a)^2}
            +\frac2{(36\pi a)^3}\right)\Theta_5(a)
       -\frac{\pi^3/2}{V(36\pi a^*)}\Theta_5(a^*),\\
 \begin{aligned}
 &\frac{\pi^6}{240}\biggl(\frac1{36\pi a}+\frac5{(36\pi a)^2}
            +\frac{20}{(36\pi a)^3}+\frac{60}{(36\pi a)^4}\\
 &\qquad+\frac{120}{(36\pi a)^5}+\frac{120}{(36\pi a)^6}\biggr)\Theta_5(a)
       +\frac{\pi^6/240}{V(36\pi a^*)}\Theta_5(a^*).
 \end{aligned}
\end{gathered}
\]

For cube $N$, replace $36,13$ by $(N+1)^2,2N+3$.
Lemma~\ref{lem:finitearithmetic} gives the common arithmetic.
This completes the proof of Lemma~\ref{lem:firstloss}.
\end{proof}
The complete coordinate derivatives and tails are given in
Lemma~\ref{lem:coordinate-ewald-tails}.

\begin{lemma}[One polynomial construction for the stationary basins]
\label{lem:fctuniformbasin}
For a reference lattice $L_0$, use
\begin{equation}\label{eq:common-metric-chart}
 \begin{gathered}
 S=\begin{pmatrix}2u&x&y\\x&-u+v&z\\y&z&-u-v\end{pmatrix},\qquad
 \mathbf s=(u,v,x,y,z)^{\mathsf T},\\
 \|S\|_F^2=6u^2+2v^2+2x^2+2y^2+2z^2,\qquad D=\det(I+S).
 \end{gathered}
\end{equation}
and $L(S)=D^{-1/6}(I+S)^{1/2}L_0$.
For $w\in L_0\setminus\{0\}$ set $q_w=|w|^2$ and
$\ell_w=w^{\mathsf T}Sw/q_w$. Then
\[
 \mathcal E(L(S))=\sum_{w\in L_0\setminus\{0\}}
 \left[\frac{D^2}{2q_w^6(1+\ell_w)^6}
       -\frac{D}{q_w^3(1+\ell_w)^3}\right].
\]
Choose a finite symmetric $\mathcal F\subset L_0\setminus\{0\}$ and
complete bounds
\[
 \sum_{w\in L_0\setminus(\mathcal F\cup\{0\})}|w|^{-2s}\le\tau_s,
 \qquad E(L_0,s)\le E^+(L_0,s)\quad(s=3,6).
\]
Put $p_s(z)=\sum_{k=0}^6(-1)^k\binom{s+k-1}{k}z^k$ and
\[
 P(S)=\sum_{w\in\mathcal F}
 \left[\frac{D^2}{2q_w^6}p_6(\ell_w)-\frac D{q_w^3}p_3(\ell_w)\right].
\]
Retain the determinant factors; $\deg P\le12$.
For $0\le R<6/5$ and $j=0,1,2$, define
\[
 T_{s,j}=\left(\frac56\right)^jE^+(L_0,s)\left\{
 \frac{(s+j-1)!}{(s-1)!(1-5R/6)^{s+j}}
 -\sum_{k=j}^6\frac{k!}{(k-j)!}\binom{s+k-1}{k}(5R/6)^{k-j}
 \right\}.
\]
On $\|S\|_F\le R$, the complete Hessian remainder is at most
\begin{align}
 \varepsilon(R)={}&T_{3,2}+\frac{2R}{1-5R/6}T_{3,1}
                  +\frac{R^2+1}{(1-5R/6)^2}T_{3,0}\notag\\
 &+\frac12\left(T_{6,2}+\frac{4R}{1-5R/6}T_{6,1}
                  +\frac{4R^2+2}{(1-5R/6)^2}T_{6,0}\right)\notag\\
 &+\frac{R^2+5R+28/3}{(1-5R/6)^5}\tau_3
 +\frac{4R^2+20R+187/6}{2(1-5R/6)^8}\tau_6.
 \label{eq:common-basin-remainder}
\end{align}
Thus $\|\nabla_S^2[\mathcal E(L(S))-P(S)]\|_{\rm op}\le\varepsilon(R)$,
where the operator norm is induced by $\|\cdot\|_F$ on trace-free
symmetric matrices.
For $\mathcal E_\theta(L)=\tfrac12E(L,6)-\theta E(L,3)$,
multiply the two $s=3$ contributions by an upper bound for $|\theta|$.
\end{lemma}
\begin{proof}
For trace-free symmetric $H$,
$\|H\|_{\rm op}\le\sqrt{2/3}\|H\|_F\le\frac56\|H\|_F$.
The positive binomial tails of orders $j=0,1,2$ are bounded by
$T_{s,j}$. For $s=3,6$ and $\|H\|_F=1$,
\begin{align*}
 \frac{\left.\partial_t\det(I+S+tH)^{s/3}\right|_{t=0}}
      {\det(I+S)^{s/3}}
 &=\frac{s}{3}\tr\bigl((I+S)^{-1}H\bigr),\\
 \frac{\left.\partial_t^2\det(I+S+tH)^{s/3}\right|_{t=0}}
      {\det(I+S)^{s/3}}
 &=\frac{s^2}{9}\bigl[\tr((I+S)^{-1}H)\bigr]^2\\
 &\quad-\frac{s}{3}\tr\bigl(((I+S)^{-1}H)^2\bigr).
\end{align*}
Since $D\le1$, these identities give
\[
 \left|\left.\partial_t^j\det(I+S+tH)^{s/3}\right|_{t=0}\right|
 \le\begin{cases}
 \dfrac{sR}{3(1-5R/6)},&j=1,\\[4pt]
 \dfrac{s^2R^2/9+s/3}{(1-5R/6)^2},&j=2.
 \end{cases}
\]
The product rule gives
\eqref{eq:common-basin-remainder}, including all omitted vectors
and binomial terms. The finite matrix checks use
\begin{align*}
 D={}&1-3u^2-v^2-x^2-y^2-z^2+2u^3-2uv^2-2uz^2\\
    &+ux^2+vx^2+uy^2-vy^2+2xyz,\\
 \max_{\|S\|_F\le R}|\mathbf s^\alpha|
   & = R^{|\alpha|}\left(\frac{\alpha_1}{6|\alpha|}\right)^{\alpha_1/2}
       \prod_{i=2}^5\left(\frac{\alpha_i}{2|\alpha|}\right)^{\alpha_i/2},
 \quad |\alpha|>0.
\end{align*}
Here $\alpha\in\mathbb Z_{\ge0}^5$ and $|\alpha|=\sum_{i=1}^5\alpha_i$;
zero exponents contribute one. Apply them to $\nabla_{\mathbf s}^2P$,
retaining the nonnegative signs of even monomials.
For the homogeneous decomposition $P=\sum_{d=0}^{12}P_d$, put
\[
 B_{ij}(S)=\sum_{d=2}^{12}\frac1{d-1}
          \frac{\partial^2P_d}{\partial s_i\partial s_j}(S),
 \qquad1\le i,j\le5.
\]
Euler's identity and integration give
\begin{gather*}
 \mathbf s\cdot\nabla_{\mathbf s}P(S)
 =\mathbf s\cdot\nabla_{\mathbf s}P(0)+\mathbf s^{\mathsf T}B(S)\mathbf s,\\
 \begin{aligned}
 &\left|\left.\partial_t[\mathcal E(L(tS))-P(tS)]\right|_{t=0}^{t=1}\right|\\
 &\qquad=\left|\int_0^1\partial_t^2[\mathcal E(L(tS))-P(tS)]\dd t\right|
 \le\varepsilon(R)\|S\|_F^2.
 \end{aligned}
\end{gather*}
The matrix and radial bounds give inner convexity and outer
stationary-point exclusion.
This completes the proof of Lemma~\ref{lem:fctuniformbasin}.
\end{proof}

\begin{lemma}[Coordinate Ewald formulas and uniform truncation bounds]
\label{lem:coordinate-ewald-tails}
Fix $V>0$, $N\in\mathbb Z_{\ge0}$ and a closed box
$B=\prod_{i=1}^5[p_i^-,p_i^+]\subset\{a,\delta>0\}$,
where $p=(a,\delta,x,y,z)$. Put
\[
 U(p)=\begin{pmatrix}a&ax&ay\\0&\delta&\delta z\\0&0&h\end{pmatrix},
 \qquad h=V/(a\delta),
\]
and let $h_->0$ be a lower bound for $h$ on $B$.
Set $U_i=\partial_{p_i}U$, $U_{ij}=\partial_{p_i}\partial_{p_j}U$.
The identities \eqref{eq:criticalKZprimal}--\eqref{eq:criticalderivatives}
and the complete value, gradient and Hessian truncation bounds below
hold uniformly on $B$, with no reduction hypothesis.
\end{lemma}
\begin{proof}
The squared lengths are
\begin{align}
 q_p(i,j,k)&=a^2(i+xj+yk)^2+\delta^2(j+zk)^2+h^2k^2,
                       \label{eq:criticalKZprimal}\\
 q_p^*(i,j,k)&=a^{-2}i^2+\delta^{-2}(j-xi)^2
             +h^{-2}\{k-zj+(xz-y)i\}^2.
                       \label{eq:criticalKZdual}
\end{align}
Put $f(p)=\tfrac12E(U(p)\Z^3,6)-E(U(p)\Z^3,3)$ and set
\[
 g(q)=\frac{\pi^6}{240}\mathcal G_6(\pi q)-\frac{\pi^3}{2}\mathcal G_3(\pi q),\qquad
 \widetilde g(q)=\frac{\frac{\pi^6}{240}\mathcal G_{-9/2}(\pi q)
                         -\frac{\pi^3}{2}\mathcal G_{-3/2}(\pi q)}{V}.
\]
The exact identity of Lemma~\ref{lem:firstloss} is
\begin{align}
 f(p)={}&\frac{\pi^6}{240}\left(-\frac16+\frac{2}{9V}\right)
       -\frac{\pi^3}{2}\left(-\frac13+\frac{2}{3V}\right)\notag\\
 &+\sum_{(i,j,k)\in\Z^3\setminus\{0\}}
  \{g(q_p(i,j,k))+\widetilde g(q_p^*(i,j,k))\}.
 \label{eq:criticalEwald}
\end{align}
For $r=0,1,2$, with $g^{(r)}=\dd^rg/\dd q^r$ and
$\widetilde g^{(r)}=\dd^r\widetilde g/\dd q^r$,
\[
\begin{aligned}
 g^{(r)}(q)&=(-\pi)^r
  \{\frac{\pi^6}{240}\mathcal G_{6+r}(\pi q)-\frac{\pi^3}{2}\mathcal G_{3+r}(\pi q)\},\\
 \widetilde g^{(r)}(q)&=\frac{(-\pi)^r}{V}
  \{\frac{\pi^6}{240}\mathcal G_{-9/2+r}(\pi q)
                 -\frac{\pi^3}{2}\mathcal G_{-3/2+r}(\pi q)\}.
\end{aligned}
\]
For $q=q_p(i_1,i_2,i_3)$, $q^*=q_p^*(i_1,i_2,i_3)$ and $1\le i,j\le5$,
\[
 (f_i,q_i,q_i^*)=\partial_{p_i}(f,q,q^*),\qquad
 (f_{ij},q_{ij},q_{ij}^*)=\partial_{p_i}\partial_{p_j}(f,q,q^*).
\]
Lemma~\ref{lem:firstloss} encloses the kernels; differentiation gives
\begin{align}
 f_i&=\sum_{(i_1,i_2,i_3)\in\Z^3\setminus\{(0,0,0)\}}
  \{g'(q)q_i+\widetilde g'(q^*)q_i^*\},\notag\\
 f_{ij}&=\sum_{(i_1,i_2,i_3)\in\Z^3\setminus\{(0,0,0)\}}
  \{g''(q)q_iq_j+g'(q)q_{ij}
       +\widetilde g''(q^*)q_i^*q_j^*
       +\widetilde g'(q^*)q_{ij}^*\}.
 \label{eq:criticalderivatives}
\end{align}
\paragraph{Complete tails.}
For $p\in B$, $n\in\Z^3\setminus\{0\}$, independently of reduction,
\[
\begin{gathered}
 0<\alpha\le\bigl(\sup_{p\in B}\|U(p)^{-1}\|_{\rm op}^2\bigr)^{-1},\quad
 0<\beta\le\bigl(\sup_{p\in B}\|U(p)\|_{\rm op}^2\bigr)^{-1},\\
 q_p(n)\ge\alpha|n|^2,\qquad q_p^*(n)\ge\beta|n|^2,\qquad
 \|A\|_{\rm op}\le\min\{\|A\|_F,\sqrt{\|A\|_1\|A\|_\infty}\}.
\end{gathered}
\]
Use $\Theta_N$ from Lemma~\ref{lem:firstloss}; for $s\in\{3,6\}$,
$r\in\{0,1,2\}$, set
\[
 A_{s,r}(b)=\sum_{j=0}^{s+r-1}\frac{(s+r-1)!}{j!}
 \begin{cases}
 \{\pi b(N+1)^2\}^{j-s}\Theta_N(b),&j\le s,\\
 \Theta_N(b/2),&j=s+1,
 \end{cases}
\]
Indeed,
\[
 x^r\mathcal G_{s+r}(x)
 =e^{-x}\sum_{j=0}^{s+r-1}\frac{(s+r-1)!}{j!}x^{j-s},
 \qquad xe^{-x}\le e^{-x/2},
\]
and $\mathcal G_v(x)\le e^{-x}/x$ for $v\le1$. Thus, uniformly on $B$,
\[
\begin{aligned}
 \sum_{\substack{n\in\Z^3\\\|n\|_\infty>N}}
 q_p(n)^r|g^{(r)}(q_p(n))|
 &\le P_r=\frac{\pi^3}{2}A_{3,r}(\alpha)+\frac{\pi^6}{240}A_{6,r}(\alpha),\\
 \sum_{\substack{n\in\Z^3\\\|n\|_\infty>N}}
 q_p^*(n)^r|\widetilde g^{(r)}(q_p^*(n))|
 &\le D_r=\frac{\frac{\pi^3}{2}+\frac{\pi^6}{240}}{V}
 \begin{cases}
 \Theta_N(\beta)/\{\pi\beta(N+1)^2\},&r=0,\\
 \Theta_N(\beta),&r=1,\\
 \Theta_N(\beta/2),&r=2.
 \end{cases}
\end{aligned}
\]
The positive attraction/repulsion tails are
\[
 \frac{\pi^3}{2}\left\{A_{3,0}(\alpha)
       +\frac{\Theta_N(\beta)}{V\pi\beta(N+1)^2}\right\},\qquad
 \frac{\pi^6}{240}\left\{A_{6,0}(\alpha)
       +\frac{\Theta_N(\beta)}{V\pi\beta(N+1)^2}\right\}.
\]

For box endpoints $-,+$, use
\[
 \begin{gathered}
 (\gamma_a,\gamma_\delta,\gamma_y,\gamma_z)=
 \left(\frac2{a_-},\frac2{\delta_-},\frac{a_+}{h_-},
                                  \frac{\delta_+}{h_-}\right),\\
 \gamma_x=\sqrt{\left(\frac{a_+}{\delta_-}\right)^2+
       \left(\frac{a_+\max\{|z_-|,|z_+|\}}{h_-}\right)^2}.
 \end{gathered}
\]
For $p\in B$, $n\ne0$ and $1\le i,j\le5$, differentiation gives
\[
\begin{gathered}
 |q_i|\le\gamma_iq,\quad |q_i^*|\le\gamma_iq^*,\qquad
 t_i\ge\sup_{p\in B}\|U_iU^{-1}\|_{\rm op},\quad
 t_{ij}\ge\sup_{p\in B}\|U_{ij}U^{-1}\|_{\rm op},\\
 |q_{ij}|\le2(t_it_j+t_{ij})q,\qquad
 |q_{ij}^*|\le2(3t_it_j+t_{ij})q^*.
\end{gathered}
\]
Use the preceding norm bounds and $2|uv|\le u^2+v^2$. Consequently,
\[
 \sum_{\substack{n\in\Z^3\\\|n\|_\infty>N}}
 \bigl(|g'(q)q_i|+|\widetilde g'(q^*)q_i^*|\bigr)
 \le\gamma_i(P_1+D_1),
\]
\begin{equation}
\begin{aligned}
 &\sum_{\substack{n\in\Z^3\\\|n\|_\infty>N}}
 \bigl(|g''(q)q_iq_j|+|g'(q)q_{ij}|\\
 &\hspace{38mm}+|\widetilde g''(q^*)q_i^*q_j^*|
                      +|\widetilde g'(q^*)q_{ij}^*|\bigr)\\
 &\quad\le\gamma_i\gamma_j(P_2+D_2)
 +2(t_it_j+t_{ij})P_1+2(3t_it_j+t_{ij})D_1.
\end{aligned}
 \label{eq:criticalHtail}
\end{equation}

This completes the proof of Lemma~\ref{lem:coordinate-ewald-tails}.
\end{proof}

\section{The simple-hexagonal phase}

For $V\ge89/50$, planar reduction identifies the simple-hexagonal
minimizer, area convexity gives uniqueness, and Fourier expansion
yields its dilute asymptotics. Put
$F_V(A)=\mathcal E(T_A\oplus(V/A)\Z)$. Since $A_*<1$,
\[
 m(V)\le F_V(A_*)<e_2.
\]
The reduction is
\[
 L\in\operatorname*{argmin}_{|L|=V}\mathcal E(L)
 \ \Longrightarrow\
 L\cong T_A\oplus(V/A)\Z,\quad
 \frac25<A<\frac{89}{100},\quad F_V'(A)=0,
\]
with uniqueness from $F_V''>34$.

For $y>0$, set
\[
 \Gamma_{x,y}=y^{-1/2}\operatorname{span}_{\Z}\{(1,0),(x,y)\}.
\]
Thus $\Gamma_{1/2,\sqrt3/2}=T_1$, and for $s=3,6$,
\[
 E(\Gamma_{x,y},s)=\sum_{\substack{(i,j)\in\Z^2\\(i,j)\ne(0,0)}}
 \frac{y^s}{((i+jx)^2+j^2y^2)^s}.
\]

\begin{lemma}[Watson {\cite{Watson}}]\label{lem:besselproperties}
For $z>0$ and $\nu\ge0$,
\begin{equation}\label{eq:fouriertubebesselintegral}
\begin{aligned}
 K_\nu(z)&=\int_0^\infty e^{-z\cosh a}\cosh(\nu a)\,\dd a\\
 &=\sqrt{\frac\pi{2z}}\frac{e^{-z}}{\Gamma(\nu+1/2)}
 \int_0^\infty e^{-b}b^{\nu-1/2}
                 \left(1+\frac{b}{2z}\right)^{\nu-1/2}\dd b.
\end{aligned}
\end{equation}
For $j\in\mathbb Z_{\ge0}$ and $0\le\mu\le\nu$,
\[
\begin{gathered}
 (-1)^jK_\nu^{(j)}(z)
 =\int_0^\infty e^{-z\cosh a}(\cosh a)^j\cosh(\nu a)\,\dd a>0,\\
 K_{-\nu}(z)=K_\nu(z),\qquad K_\mu(z)\le K_\nu(z).
\end{gathered}
\]
The differential and recurrence identities are
\[
\begin{gathered}
 z^2K_\nu''+zK_\nu'-(z^2+\nu^2)K_\nu=0,\qquad K_0'=-K_1,\\
 K_{\nu+1}=K_{\nu-1}+\frac{2\nu}{z}K_\nu,\qquad
 \frac{\dd}{\dd z}\{z^\nu K_\nu(z)\}=-z^\nu K_{\nu-1}(z).
\end{gathered}
\]
For $n\in\mathbb Z_{\ge0}$,
\[
 P_n(z)=\sum_{\ell=0}^n
 \frac{(n+\ell)!}{\ell!(n-\ell)!(2z)^\ell},\qquad
 K_{n+1/2}(z)=\sqrt{\frac\pi{2z}}e^{-z}P_n(z).
\]
For $z\ge1$, $\nu\ge0$ and $j\in\mathbb Z_{\ge0}$,
\begin{equation}\label{eq:besselbounds}
 z^{-1/2}e^{-z-1}\le K_\nu(z),\qquad
 |K_\nu^{(j)}(z)|\le\sqrt{\frac{2\pi}{z}}
             e^{-z+(j+\nu)^2/(2z)}.
\end{equation}
\end{lemma}

\begin{lemma}[Watson {\cite{Watson}}]\label{lem:criticalcusp}
For $z>0$, set
\[
\begin{aligned}
 \mathcal K_2(z)&=\sqrt{\frac{\pi}{2z}}e^{-z}
       \left(1+\frac{15}{8z}+\frac{105}{128z^2}\right),\\
 \mathcal K_5(z)&=\sqrt{\frac{\pi}{2z}}e^{-z}
       \left(1+\frac{99}{8z}+\frac{9009}{128z^2}
       +\frac{225225}{1024z^3}
       +\frac{11486475}{32768z^4}
       +\frac{43648605}{262144z^5}\right).
\end{aligned}
\]
They satisfy $K_\nu(z)\le\mathcal K_\nu(z)$ for
$\nu=2,5$ and $z>0$.
\end{lemma}

\begin{lemma}[An explicit planar Epstein gap]\label{lem:planargap}
For $0\le x\le1/2$ and $y\ge\sqrt3/2$,
\begin{align*}
 E(\Gamma_{x,y},3)-E(T_1,3)
 &\ge \frac95\pi^3e^{-2\pi y}(1+\cos(2\pi x))
       +3\sqrt3\left(y-\frac{\sqrt3}{2}\right)^2\\
 &\ge \frac14\left[\left(x-\frac12\right)^2
                   +\left(y-\frac{\sqrt3}{2}\right)^2\right].
\end{align*}
If \(0<\delta\le1/10\), then
\[
 \left(x-\frac12\right)^2+
 \left(y-\frac{\sqrt3}{2}\right)^2\ge\delta^2
 \quad\Longrightarrow\quad
 E(\Gamma_{x,y},3)-E(T_1,3)\ge\frac{\delta^2}{2}.
\]
\end{lemma}
\begin{proof}
Poisson summation and Lemma~\ref{lem:besselproperties} with $n=2$ give
\begin{align*}
 E(\Gamma_{x,y},3)&=2\zeta(6)y^3+\frac{3\pi}{4}\zeta(5)y^{-2}\\
 &\quad+2\pi^3\sum_{n=1}^{\infty}
 n^2\left(\sum_{\substack{d\in\mathbb Z_{\ge1}\\d\mid n}}d^{-5}\right)e^{-2\pi ny}
 \left(1+\frac3{2\pi ny}+\frac3{(2\pi ny)^2}\right)
 \cos(2\pi nx).
\end{align*}
The Fourier series and its derivatives converge absolutely and
uniformly for $y\ge\sqrt3/2$. Here $e^{-2\pi y}<1/220<1/200$, and
\[
 \sum_{\substack{d\in\mathbb Z_{\ge1}\\d\mid n}}d^{-5}\le\zeta(5)<\frac{26}{25},\qquad
 \frac{26}{25}\sum_{n=2}^{\infty}n^4(1/200)^{n-1}<\frac1{10}.
\]
The inequalities $|\sin(n\theta)|\le n\sin\theta$
($0\le\theta\le\pi$) and $(1+3/t+3/t^2)'<0$ ($t>0$) give
\[
 -\partial_xE(\Gamma_{x,y},3)
 \ge\frac{18}{5}\pi^4e^{-2\pi y}\sin(2\pi x).
\]
Integration from \(x\) to \(1/2\) yields
\[
 E(\Gamma_{x,y},3)-E(\Gamma_{1/2,y},3)
 \ge\frac95\pi^3e^{-2\pi y}(1+\cos(2\pi x)).
\]

Twice differentiating and using $\cos(2\pi nx)\ge-1$ gives
\begin{align*}
 \partial_{yy}E(\Gamma_{x,y},3)
 &\ge12\zeta(6)y+\frac{9\pi}{2}\zeta(5)y^{-4}\\
 &\quad-8\pi^5\zeta(5)\sum_{n=1}^{\infty}n^4e^{-2\pi ny}
 \left(1+\frac3{2\pi ny}+\frac9{(2\pi ny)^2}
 +\frac{18}{(2\pi ny)^3}+\frac{18}{(2\pi ny)^4}\right).
\end{align*}
Since
\[
 \frac{\dd}{\dd t}\{e^{-t}(t^4+3t^3+9t^2+18t+18)\}
 =-e^{-t}t^3(t-1)\le0\qquad(t\ge1),
\]
the negative Fourier part, divided by $(9\pi/2)\zeta(5)y^{-4}$,
has absolute value at most
\[
 \frac19\sum_{n=1}^{\infty}220^{-n}
 \left[\left(\frac{27n}{5}\right)^4
 +3\left(\frac{27n}{5}\right)^3
 +9\left(\frac{27n}{5}\right)^2
 +18\left(\frac{27n}{5}\right)+18\right]<1.
\]
Thus $\partial_{yy}E(\Gamma_{x,y},3)\ge12y\ge6\sqrt3$.
Modular symmetry gives
$\left.\partial_yE(\Gamma_{1/2,y},3)\right|_{y=\sqrt3/2}=0$;
integrating twice from $\sqrt3/2$ yields
\[
 E(\Gamma_{1/2,y},3)-E(T_1,3)
 \ge3\sqrt3\left(y-\frac{\sqrt3}{2}\right)^2.
\]
For $\sqrt3/2\le y\le1$, $1+\cos(2\pi x)\ge8(x-1/2)^2$ gives
\[
 \frac95\pi^3e^{-2\pi y}(1+\cos(2\pi x))
 \ge\frac12(x-1/2)^2.
\]
Hence the lower bound is
$\ge\tfrac12[(x-1/2)^2+(y-\sqrt3/2)^2]$. For $y\ge1$,
$y-\sqrt3/2>1/8$ and $(x-1/2)^2\le1/4$, so
\[
 3\sqrt3(y-\sqrt3/2)^2
 >5(y-\sqrt3/2)^2
 \ge\frac14\bigl[(x-1/2)^2+(y-\sqrt3/2)^2\bigr].
\]
This exceeds $5/64>\delta^2/2$ for $\delta\le1/10$.
This completes the proof of Lemma~\ref{lem:planargap}.
\end{proof}

\begin{lemma}[Planar constants]\label{lem:planarconstants}
With the notation of the planar triangular minimum,
\[
 S_3<\frac{13}{2},\qquad
 \frac{9\sqrt3}{26}<A_*^3<\frac{611}{972}
 <\left(\frac{43}{50}\right)^3,
 \qquad S_3^2/S_6>6,\qquad E(T_1,3)\ge\frac{9\sqrt3}{4}.
\]
\end{lemma}
\begin{proof}
For $q=i^2+ij+j^2=1,3,4,7,9$, the multiplicities are
$6,6,6,12,6$, respectively; disjoint radius-$1/2$ disks give
\[
 S_3\le6+\frac6{27}+\frac6{64}+\frac{12}{343}+\frac6{729}
 +2\left(\frac76\right)^6\left(\frac52\right)^{-4}
 <\frac{13}{2}.
\]
Since $S_6\ge6$ and the six radius-$\sqrt3$ vectors give
$S_3-S_6\ge52/243$,
\[
 A_*^3=\frac{3\sqrt3}{8}\frac{S_6}{S_3}
 >\frac{9\sqrt3}{26},\qquad
 \frac{S_6}{S_3}<\frac{3055}{3159}.
\]
Hence $A_*^3<611/972<(43/50)^3$, while
$S_3^2/S_6>S_3>6$ and $E(T_1,3)=(3\sqrt3/8)S_3\ge9\sqrt3/4$.
This completes the proof of Lemma~\ref{lem:planarconstants}.
\end{proof}

Along triangular planes,
\begin{equation}\label{eq:planararea}
 U(A):=\mathcal E(T_A)
 =\frac{S_6}{2}\left(\frac{\sqrt3}{2A}\right)^6
       -S_3\left(\frac{\sqrt3}{2A}\right)^3
 =e_2+\frac{S_3^2}{2S_6}\left[(A_*/A)^3-1\right]^2.
\end{equation}
Since $6<S_6<S_3<13/2$, we have $6<S_3^2/S_6<11$.

\begin{lemma}[Compression to planar area one]\label{lem:compression}
For a planar lattice $\Gamma\subset\mathbb R^2$,
\[
\begin{gathered}
 T_1=\sqrt{\frac2{\sqrt3}}
       \begin{pmatrix}1&1/2\\0&\sqrt3/2\end{pmatrix}\mathbb Z^2,
 \qquad |T_1|=1,\\
 A=|\Gamma|\ge1\quad\Longrightarrow\quad
 \mathcal E(\Gamma)\ge\mathcal E(T_1).
\end{gathered}
\]
\end{lemma}
\begin{proof}
Choose a basis $(a,0),(b,h)$ with
$a=\min_{v\in\Gamma\setminus\{0\}}|v|$ and $ah=A$.
For $a\le1$, decrease $h$ to $1/a$; changed vectors retain
vertical component $\ge1/a\ge1$, and
\[
 \frac{\dd}{\dd r}(r^{-12}-2r^{-6})
 =12r^{-13}(r^6-1)\ge0\qquad(r\ge1).
\]
For $a>1$, scale down until $A=1$ or $a=1$, then compress
vertically if needed. Changed lengths remain $\ge1$.
At area one, \eqref{eq:planarinput} gives
$\mathcal E(\Gamma)\ge\mathcal E(T_1)$.
This completes the proof of Lemma~\ref{lem:compression}.
\end{proof}

Write the layered lattice as
\[
 L=\Gamma+\Z(u+he_3),\qquad A=|\Gamma|,\quad h=V/A,
 \quad u\in\R^2.
\]
As in \eqref{eq:energy},
$\mathcal E(\Gamma)=\tfrac12E(\Gamma,6)-E(\Gamma,3)$.
Its reciprocal lattice is
\[
 \Gamma^*=\{w\in\R^2:w\cdot v\in\Z\text{ for every }v\in\Gamma\}.
\]
For \(d>0\) define
\[
 f_d(x)=(|x|^2+d^2)^{-6}-2(|x|^2+d^2)^{-3},
 \qquad x\in\R^2.
\]
We use the Fourier transform with kernel \(e^{-2\pi i w\cdot x}\).
The Gaussian Laplace integral gives, for \(s>1\),
\begin{equation}\label{eq:fourier}
 \widehat{(|x|^2+d^2)^{-s}}(w)
 =\frac{\pi}{\Gamma(s)}
   \int_0^\infty t^{s-2}e^{-d^2t-\pi^2|w|^2/t}\dd t .
\end{equation}
At $w=0$ the value is $\pi d^{2-2s}/(s-1)$; for $r=|w|>0$,
Lemma~\ref{lem:besselproperties} gives
\[
 \widehat{(|x|^2+d^2)^{-s}}(w)
 =\frac{2\pi^s}{\Gamma(s)}(r/d)^{s-1}
                         K_{s-1}(2\pi dr).
\]
The layers \(m\) and \(-m\) combine to give
\begin{equation}\label{eq:poisson}
 \mathcal E(L)=\mathcal E(\Gamma)
   +\frac1A\sum_{m=1}^\infty\sum_{w\in\Gamma^*}
             \widehat f_{mh}(w)\cos(2\pi m w\cdot u).
\end{equation}
The zero-frequency terms are
\begin{equation}\label{eq:zeromode}
 -\frac{\pi\zeta(4)}{Ah^4}
     +\frac{\pi\zeta(10)}{5Ah^{10}}
 =-\pi\zeta(4)\frac{A^3}{V^4}
     +\frac{\pi\zeta(10)}5\frac{A^9}{V^{10}}.
\end{equation}
For $\Gamma$ in a compact neighborhood of $T_{A_*}$, there is
$\rho>0$ such that
\[
 \lambda_1(\Gamma),\lambda_1(\Gamma^*)\ge\rho,\qquad
 \#(\Gamma\cap B_R)+\#(\Gamma^*\cap B_R)\le C(1+R)^2.
\]
Choose a smooth basis and lattice coordinates for $u$.
For any fixed mixed partial derivative $\partial^\beta$ in the area and shape coordinates,
\eqref{eq:fourier} and \eqref{eq:besselbounds} give, for $h\ge1$,
\[
 \left|\partial^\beta\left\{A^{-1}\widehat f_{mh}(w)
                   \cos(2\pi m w\cdot u)\right\}\right|
 \le C_\beta(1+mh+|w|)^{N_\beta}e^{-2\pi mh|w|}.
\]
Quadratic counting gives a uniform $O(e^{-ch})$ sum over
$m\ge1$, $w\ne0$, with any fixed number of derivatives
(after decreasing $c>0$).

\begin{lemma}[Planar and cross-layer energy bounds]\label{lem:crossenvelope}
For a lattice $L\subset\mathbb R^3$, choose
\[
\begin{gathered}
 v_1\in\operatorname*{argmin}_{v\in L\setminus\{0\}}|v|,\qquad
 v_2\in\operatorname*{argmin}_{v\in L\setminus\mathbb Rv_1}|v|,\\
 \Gamma=L\cap\operatorname{span}_{\mathbb R}\{v_1,v_2\},\\
 a=|v_1|,\quad b=|v_2|,\quad A=|\Gamma|,\quad h=V/A,\quad V=|L|.
\end{gathered}
\] If $\mathcal E(L)\le e_2$, the short-plane bounds give
$7/10<a<8/5$ and $A<5$.
Moreover
\[
 h^2\ge\frac34b^2-\frac14a^2
 >\frac{75}{256}A^2-\frac{32}{25}.
\]
The planar energy satisfies
\[
 \mathcal E(\Gamma)\ge a^{-12}-\frac{11}{5}a^{-6}
 -\frac{11}{5}\frac{a^6}{A^6}
 -\frac{363}{140}\frac{a^4}{A^5},\qquad
 \mathcal E(\Gamma)>-\frac{21}{20}-\frac7{A^3}.
\]
The full cross-layer energy, with any layer shift, satisfies
\[
 \mathcal E(L)-\mathcal E(\Gamma)\ge
 -\frac{11}{5h^6}
 -\frac{363}{140}(1/a+a/A)h^{-5}
 -\frac{121}{35A}h^{-4}\ge-J(A,h),
\]
where
\[
 J(A,h)=\frac{11}{5h^6}
 +\frac{363}{140}\left(\frac{10}{7}+\frac8{5A}\right)h^{-5}
 +\frac{121}{35A}h^{-4}.
\]
\end{lemma}
\begin{proof}
Horizontal reduction gives
\[
 [\Gamma:\mathbb Zv_1+\mathbb Zv_2]>1
 \ \Longrightarrow\
 \exists v\in\Gamma\setminus\mathbb Rv_1:
 |v|^2\le\frac{a^2+b^2}{4}<b^2,
\]
contradicting the choice of $v_2$. Choose its sign so that
$0\le t=v_1\cdot v_2\le a^2/2$. The two nonobtuse triangles
in the reduced cell have common circumradius
\[
 R^2=\frac{a^2b^2(a^2+b^2-2t)}{4(a^2b^2-t^2)}
 \le\frac{a^2+b^2}{4}.
\]
Their radius-$R$ vertex disks cover the cell; reducing a height-$h$
vector modulo $\Gamma$ gives $b^2\le h^2+R^2$. Thus
\[
 h^2\ge b^2-R^2\ge\frac34b^2-\frac14a^2
 \ge\frac{3A^2}{4a^2}-\frac{a^2}{4}
 >\frac{75}{256}A^2-\frac{32}{25}
 \qquad(b\ge A/a,\ a<8/5).
\]

For any $c,z>0$ and any shift $t$, positivity of the Fourier
coefficients and monotonicity on the positive half-line give
\begin{align*}
 \sum_{k\in\mathbb Z}((kc+t)^2+z^2)^{-3}
 &\le z^{-6}+\frac{3\pi}{8c}z^{-5},\\
 \sum_{k\in\mathbb Z}((kc+t)^2+z^2)^{-5/2}
 &\le z^{-5}+\frac4{3c}z^{-4}.
\end{align*}
Applying these along $(a,0)$ and the spacing $A/a$ gives
\[
 \sum_{p\in\Gamma}(|p+t|^2+z^2)^{-3}
 \le z^{-6}+\frac{3\pi}{8}(1/a+a/A)z^{-5}
                         +\frac\pi{2A}z^{-4}.
\]
Sum the nonzero layers and omit repulsion;
$\pi<22/7$ and $\zeta(4),\zeta(5),\zeta(6)<11/10$ give $J$.

Separating the plane itself into horizontal chains gives
\[
 \mathcal E(\Gamma)\ge \mathcal E(a\Z)
 -2\zeta(6)(a/A)^6-\frac{3\pi}{4}\zeta(5)a^4/A^5.
\]
Use $\mathcal E(a\Z)\ge a^{-12}-(11/5)a^{-6}$, or
$\mathcal E(a\Z)\ge e_1>-21/20$ with $a^2\le A/(\sqrt3/2)$,
respectively;
\[
 \frac{2\zeta(6)}{(\sqrt3/2)^3}+\pi\zeta(5)<7.
\]
This completes the proof of Lemma~\ref{lem:crossenvelope}.
\end{proof}

\begin{samepage}
\begin{lemma}[A fixed-systole planar maximum]\label{lem:fixedminplanarmax}
Let $\Gamma\subset\mathbb R^2$ be a planar lattice with
\[\min_{v\in\Gamma\setminus\{0\}}|v|=\rho>0,\] and let
$f:(0,\infty)\to(0,\infty)$ be decreasing and convex.
If the following sums converge, then
\[
 \sum_{v\in\Gamma\setminus\{0\}}f(|v|^2)
 \le\sum_{(m,n)\in\mathbb Z^2\setminus\{(0,0)\}}
       f\!\left(\rho^2(m^2+mn+n^2)\right).
\]
\end{lemma}
\end{samepage}
\begin{proof}
Gauss reduction gives the Gram matrix
\[
 \rho^2\begin{pmatrix}1&b\\b&c\end{pmatrix},
 \qquad |b|\le\frac12,\quad c\ge1.
\]
Since $f$ decreases,
\[
 \sum_{(m,n)\in\mathbb Z^2\setminus\{(0,0)\}}f(\rho^2[m^2+2bmn+cn^2])
 \le\sum_{(m,n)\in\mathbb Z^2\setminus\{(0,0)\}}f(\rho^2[m^2+2bmn+n^2]).
\]
Evenness and convexity in $b$ give the maximum at
$b=\pm1/2$, the triangular Gram matrix.
This completes the proof of Lemma~\ref{lem:fixedminplanarmax}.
\end{proof}

\begin{lemma}[A triangular majorant for arbitrary shifted layers]\label{lem:dualtriangle}
Let $\Gamma\subset\mathbb R^2$ be a planar lattice, with
$A=|\Gamma|$ and $a=\min_{v\in\Gamma\setminus\{0\}}|v|$.
Put $\rho=a/A$ and $q_{i,j}=i^2+ij+j^2$. For any shift $u$
and any $d>0$,
\[
 \sum_{p\in\Gamma}(|p+u|^2+d^2)^{-3} \le\frac\pi{2Ad^4} +\frac{\pi^3\rho^2}{Ad^2} \sum_{(i,j)\in\mathbb Z^2\setminus\{(0,0)\}}q_{i,j} K_2(2\pi d\rho\sqrt{q_{i,j}}).
\]
The cross-layer energy is therefore $\ge-J_\triangle(a,A,h)$, where
\[
 J_\triangle(a,A,h)=\frac{\pi\zeta(4)}{Ah^4} +\frac{2\pi^3\rho^2}{Ah^2} \sum_{m=1}^{\infty}\frac1{m^2} \sum_{(i,j)\in\mathbb Z^2\setminus\{(0,0)\}}q_{i,j} K_2(2\pi mh\rho\sqrt{q_{i,j}}).
\]
\end{lemma}
\begin{proof}
Positive Fourier coefficients give the maximum at $u=0$, and
$\Gamma^*\cong A^{-1}\Gamma$ gives $\lambda_1(\Gamma^*)=\rho$.
The Fourier kernel
\[
 F_3(q,d)=\frac\pi2\int_0^\infty
 t e^{-d^2t-\pi^2q/t}\,dt
\]
satisfies $\partial_qF_3<0<\partial_{qq}F_3$ for $q>0$, and
$F_3(r^2,d)=\pi^3r^2d^{-2}K_2(2\pi dr)$.
Lemma~\ref{lem:fixedminplanarmax}, with zero mode $\pi/(2d^4)$,
gives the first bound. Summing over $d=mh$ and dropping repulsion
gives the second.
This completes the proof of Lemma~\ref{lem:dualtriangle}.
\end{proof}

\begin{remark}[An explicit positive series bound]
Lemma~\ref{lem:besselproperties}, with $n=2$, bounds the
nonzero-frequency term by
\[
 \frac{\pi^3\rho^{3/2}}{Ah^{5/2}}
 \sum_{m=1}^{\infty}\sum_{(i,j)\in\mathbb Z^2\setminus\{(0,0)\}}
 \frac{q_{i,j}^{3/4}}{m^{5/2}}
 e^{-2\pi mh\rho\sqrt{q_{i,j}}}
 \left(1+\frac3{2\pi mh\rho\sqrt{q_{i,j}}}
       +\frac3{(2\pi mh\rho\sqrt{q_{i,j}})^2}\right).
\]
Positivity and monotonicity in $h,\rho$ give endpoint bounds.
Reduced-plane geometry gives
\[
 h^2\ge\frac{3A^2}{4a^2}-\frac{a^2}{4},\qquad
 a^2\le\frac{2A}{\sqrt3}
 \quad\Longrightarrow\quad
 h\rho\ge\sqrt{\frac5{12}}>\frac{16}{25}.
\]
This gives uniform tails for finite $m,q_{i,j}$ sums.
\end{remark}

\begin{lemma}[An improved elementary planar lower bound]\label{lem:improvedplanebound}
For a planar lattice $\Gamma\subset\mathbb R^2$, set
$A=|\Gamma|$ and $a=\min_{v\in\Gamma\setminus\{0\}}|v|$. Then
\[
 \mathcal E(\Gamma)\ge a^{-12}-\frac{11}{5}a^{-6}
                         -\frac{8349}{2800}\frac{a^4}{A^5}.
\]
If $A\ge1$, two further lower bounds are
\[
 \mathcal E(\Gamma)>-\frac{947}{320},\qquad
 \mathcal E(\Gamma)\ge a^{-12}-\frac{13}{2}a^{-6}.
\]
\end{lemma}
\begin{proof}
Choose a Gauss-reduced basis $(a,0),(b,c)$, with $c=A/a$ and
$c/a\ge\sqrt3/2$. One-dimensional Poisson summation gives, for
$s=2\pi z/a$,
\[
 \sum_{k\in\mathbb Z}((ka+t)^2+z^2)^{-3}
 \le\frac{3\pi}{8az^5}
 \left[1+2\sum_{n=1}^{\infty}e^{-ns}
             \left(1+ns+\frac{n^2s^2}{3}\right)\right].
\]
Each summand in the correction decreases with $s>0$.
If $z\ge c$, then $s\ge\pi\sqrt3>27/5$ and
$e^{27/5}>220$. Thus the correction is at most
\[
 2\sum_{n=1}^{\infty}220^{-n}
       \left(1+\frac{27n}{5}+\frac{243n^2}{25}\right)
 <\frac3{20}.
\]
Summing over the nonzero planar rows gives
\[
 E(\Gamma,3)\le2\zeta(6)a^{-6}
                  +\frac{69\pi\zeta(5)}{80}\frac{a^4}{A^5}
 <\frac{11}{5}a^{-6}+\frac{8349}{2800}\frac{a^4}{A^5}.
\]
Since $\tfrac12E(\Gamma,6)\ge a^{-12}$, the first bound follows.
For $A\ge1$,
Lemmas~\ref{lem:compression} and~\ref{lem:planarconstants} give
$\mathcal E(\Gamma)\ge U(1)>-947/320$; the fixed-systole bound
$E(\Gamma,3)\le S_3a^{-6}<(13/2)a^{-6}$ gives the last inequality.
This completes the proof of Lemma~\ref{lem:improvedplanebound}.
\end{proof}

\begin{lemma}[A rational planar benchmark]\label{lem:planarbenchmark337}
The planar minimum is
\[
 e_2=\min_{A>0}\mathcal E(T_A)=\mathcal E(T_{A_*})
 =-\frac{S_3^2}{2S_6}<-\frac{337}{100}.
\]
\end{lemma}
\begin{proof}
At nearest-neighbor distance $99/100$, every vector omitted from
the following cutoff has length $>1$ and contributes negatively:
\[
 e_2<\sum_{\substack{(i,j)\in\mathbb Z^2\\0<q=i^2+ij+j^2\le16}}
 \left[\frac{(100/99)^{12}}{2q^6}-\frac{(100/99)^6}{q^3}\right]
 <-\frac{337}{100}.
\]
This completes the proof of Lemma~\ref{lem:planarbenchmark337}.
\end{proof}

\begin{lemma}[A common triangular reference bound]
\label{lem:v2-referencecomparison}
For $V>0$,
\[
 \begin{aligned}
 \mathcal E\left(T_{17/20}\oplus\frac{V}{17/20}\Z\right)
 &\le U(17/20)-\frac{\pi\zeta(4)(17/20)^3}{V^4}\\
 &\quad+\frac{\zeta(10)}{5(49/100)^2}
 \left(\frac{17/20}{V}\right)^{10}
 \left(1+\frac{(49/100)(17/20)}{V}\right)^{12}.
 \end{aligned}
\]
\end{lemma}
\begin{proof}
Since $\lambda_1(T_{17/20})>98/100=2(49/100)$, integration over
disjoint radius-$(49/100)$ disks gives
\[
 \sum_{p\in T_{17/20}}(|p|^2+d^2)^{-6}
 \le\frac{d^{-10}}{5(49/100)^2}
 \left(1+\frac{49/100}{d}\right)^{12}.
\]
Sum over $d=mV/(17/20)$, $m\ge1$, and omit the nonpositive
attractive Fourier terms, retaining $-\pi\zeta(4)(17/20)^3/V^4$.
This completes the proof of Lemma~\ref{lem:v2-referencecomparison}.
\end{proof}

\begin{lemma}[The uniform triangular-dual lower envelope]
\label{num:v2-area}
Suppose
\[
 \frac7{10}\le a\le\frac85,\qquad
 1\le A\le\frac72,\qquad V\ge2,\qquad
 A\ge\frac{433}{500}a^2,
 \qquad \left(\frac VA\right)^2
       \ge\frac{3A^2}{4a^2}-\frac{a^2}{4}.
\]
Put $h=V/A$ and define
\[
 P(a,A)=\max\Bigg\{a^{-12}-\frac{11}{5}a^{-6}
 -\frac{8349}{2800}\frac{a^4}{A^5},\quad
 a^{-12}-\frac{13}{2}a^{-6},\ -\frac{947}{320}\Bigg\}.
\]
With $J_\triangle$ as in Lemma~\ref{lem:dualtriangle},
\[
 \max\left\{P(a,A)-J_\triangle(a,A,h),\,
             a^{-12}-\frac{2901}{200}a^{-6}\right\}
 >-\frac{10}{3}>-\frac{337}{100}.
\]
\end{lemma}

\begin{proof}
Lemma~\ref{lem:besselproperties} gives
\[
 \frac{\pi^3r^2}{d^2}K_2(2\pi dr)<\mathcal U(r,d):=
 \frac{\pi^3}{2}\frac{r^{3/2}}{d^{5/2}}e^{-2\pi dr}
 \left(1+\frac3{2\pi dr}+\frac3{(2\pi dr)^2}\right).
\]
The expression decreases in $r,d$: for $z=2\pi dr$,
\[
 \frac{d}{dz}\{e^{-z}(z^{3/2}+3z^{1/2}+3z^{-1/2})\}
 =-e^{-z}\left(z^{3/2}+\frac32z^{1/2}
                  +\frac32z^{-1/2}+\frac32z^{-3/2}\right)<0.
\]
Monotonicity in $d$ follows termwise.

Put $\rho=a/A$. The hypotheses give
\[
 (h\rho)^2\ge\frac34-\frac{a^4}{4A^2}
 >\left(\frac{16}{25}\right)^2.
\]
Thus $2\pi mh\rho\sqrt q>4m\sqrt q$, using $\pi>25/8$.
The factor $1+3/z+3/z^2$ is less than $2$ when $z>4$.
Furthermore
\[
 \frac{\pi^3\rho^{3/2}}{Ah^{5/2}}
 =\frac{\pi^3a^{3/2}}{V^{5/2}}<12.
\]
Since $\min_{(i,j)\in\mathbb Z^2\setminus\{(0,0)\}}
(i^2+ij+j^2)=1$, disjoint radius-$1/2$ disks give
\[
 \#\{(i,j)\in\mathbb Z^2:n\le\sqrt{i^2+ij+j^2}<n+1\}<25n^2
 \quad(n\in\mathbb Z_{\ge1}).
\]
Here $q=i^2+ij+j^2$. Also
$q^{3/4}<3n^{3/2}$. The complete tail $m\ge6$ or $q>36$ is $<1/1000$:
$\sum_{n=1}^\infty n^4 50^{1-n}<2$, and the two regions have
bounds $4e^{-24}$ and $4\cdot6^4e^{-24}$, respectively; their
overlap may be counted twice. Use the actual $a,A,h$ before rounding
cell endpoints.

Use the rational grid
\[
 a_i=\frac{1792+9i}{2560},\qquad
 A_j=1+\frac{5j}{512},\qquad i,j=0,1,\ldots,256.
\]
Discard only $A_{j+1}<(433/500)a_i^2$: $8562$ cells are discarded
and $56974$ retained. On each retained cell put
\[
 H_{ij}=\frac1{10^6}\left\lfloor10^6
 \sqrt{\max\left\{\frac4{A_{j+1}^2},
        \frac{3A_j^2}{4a_{i+1}^2}-\frac{a_{i+1}^2}{4}\right\}}
                         \right\rfloor,
 \qquad \rho_{ij}=\frac{a_i}{A_{j+1}}.
\]
Then $h\ge H_{ij}$ and $a/A\ge\rho_{ij}$ on the admissible part
of that cell. Define
\begin{align*}
 \overline J_{ij}={}&\frac{\pi^5}{90A_jH_{ij}^4}+\frac1{1000}\\
 &+\frac2{A_j}\sum_{m=1}^{5}\Bigg[
 6\!\!\sum_{q\in\{1,3,4,9,12,16,25,27,36\}}
                 \mathcal U(\rho_{ij}\sqrt q,mH_{ij})\\
 &\hspace{43mm}+
 12\!\!\sum_{q\in\{7,13,19,21,28,31\}}
                 \mathcal U(\rho_{ij}\sqrt q,mH_{ij})\Bigg].
\end{align*}
The multiplicities enumerate $0<i^2+ij+j^2\le36$; monotonicity
and the complete tail give
$J_\triangle(a,A,h)<\overline J_{ij}$.

On each retained cell use the lower bound
\[
\begin{aligned}
 L_{ij}=\max\Bigg\{&
 \max\left\{\begin{aligned}
 &a_{i+1}^{-12}-\frac{11}{5a_i^6}
                   -\frac{8349a_{i+1}^4}{2800A_j^5},\\
 &a_{i+1}^{-12}-\frac{13}{2a_i^6},\quad-\frac{947}{320}
 \end{aligned}\right\}-\overline J_{ij},\\
 &a_{i+1}^{-12}-\frac{2901}{200a_i^6}\Bigg\}.
\end{aligned}
\]
All $56974$ comparisons give $L_{ij}>-10/3$; more precisely,
\[
 -3.319871<\min_{\substack{0\le i,j<256\\A_{j+1}\ge(433/500)a_i^2}} L_{ij}<-3.319870.
\]
The minimum occurs at $(i,j)=(118,18)$, with
\[
 a_i=\frac{1427}{1280},\quad a_{i+1}=\frac{2863}{2560},
 \quad A_j=\frac{301}{256},\quad A_{j+1}=\frac{607}{512},
 \quad H_{ij}=\frac{1686985}{10^6}.
\]

Evaluate the finite sums and complete $1/1000$ tail by interval
addition and multiplication, using Lemma~\ref{lem:finitearithmetic}:
the rational $\pi$ enclosure, exponential bounds of degrees $21,20$,
scaling and checked roots. Compute the floors in $H_{ij}$ by exact
integer square roots.
This completes the proof of Lemma~\ref{num:v2-area}.
\end{proof}

\begin{lemma}[A uniform large-area envelope]
\label{lem:area-envelope-uniform}
With $v_1,v_2,\Gamma$ as in Lemma~\ref{lem:crossenvelope},
\[
 |L|=V\ge2,\qquad\mathcal E(L)=m(V)
 \quad\Longrightarrow\quad A=|\Gamma|<1.
\]
\end{lemma}
\begin{proof}
Since $V/A_*>1$, the competitor $T_{A_*}\oplus(V/A_*)\mathbb Z$
has negative cross-layer energy. Lemmas~\ref{lem:planarbenchmark337},
\ref{lem:separation} and~\ref{lem:shortplane} give
$m(V)<e_2<-337/100<-3$, $7/10<a<8/5$ and $A<5$.
Theorem~\ref{thm:fccsystolemaximum} and
Lemma~\ref{lem:fccsystoleconstant} also give
\[
 \mathcal E(L)\ge a^{-12}-\frac{2901}{200}a^{-6}.
\]
For $A\ge7/2$, Lemma~\ref{lem:crossenvelope} gives $h>3A/7$.
Every term of $7A^{-3}+J(A,3A/7)$ decreases with $A$, and
\[
 \frac7{(7/2)^3}+J(7/2,3/2)<\frac{113}{50}.
\]
The planar and cross-layer bounds give
$\mathcal E(L)>-331/100>e_2$, a contradiction.

For $1\le A\le7/2$, the planar, triangular-dual and FCC bounds
give the maximum in Lemma~\ref{num:v2-area}; its hypotheses follow from
\[
 A\ge\frac{\sqrt3}{2}a^2>\frac{433}{500}a^2,\qquad
 h^2\ge\frac{3A^2}{4a^2}-\frac{a^2}{4}.
\]
Lemma~\ref{num:v2-area} now gives
\[
 \mathcal E(L)>-\frac{10}{3}>-\frac{337}{100}>e_2,
\]
again a contradiction. Hence $A<1$.
This completes the proof of Lemma~\ref{lem:area-envelope-uniform}.
\end{proof}

\begin{lemma}[Exact layer alignment]\label{lem:alignment}
Let $\Gamma\subset\mathbb R^2$ be a planar lattice with
\[
 A=|\Gamma|\le\frac9{10},\qquad
 \min_{v\in\Gamma\setminus\{0\}}|v|\ge\frac7{10},\qquad h\ge2.
\]
Then
\[
 \operatorname*{argmin}_{u\in\mathbb R^2/\Gamma}
 \mathcal E\bigl(\Gamma+\mathbb Z(u+he_3)\bigr)=\{0\}.
\]
\end{lemma}
\begin{proof}
For a reduced reciprocal basis $w_1,w_2$, put $\rho=3/4$ and
$\alpha=4/3$. Then
\[
\begin{gathered}
 |w_j|\le\frac2{\sqrt3\lambda_1(\Gamma)}
 \le\frac{20}{7\sqrt3}<\frac53\quad(j=1,2),\qquad
 \lambda_1(\Gamma^*)\ge\frac79>\rho,\\
 \det(w_i\cdot w_j)_{i,j=1}^2=A^{-2}\ge\frac{100}{81},\qquad
 |w_1|^2+|w_2|^2\le\frac{50}{9},\\
 n^2+m^2\le\frac92|nw_1+mw_2|^2\qquad(n,m\in\mathbb Z).
\end{gathered}
\]
For $r=|nw_1+mw_2|$, $d=kh$, $(n,m)\ne(0,0)$ and $k\ge1$,
\[
 C(r,d)=\frac{\pi^6r^5}{60d^5}K_5(2\pi rd)
             -\frac{2\pi^3r^2}{d^2}K_2(2\pi rd).
\]
Lemma~\ref{lem:besselproperties} gives, for $z=2\pi rd\ge3\pi>47/5$,
\[
 K_5(z)\le\sqrt{\frac\pi{2z}}e^{-z}P_5(z),\qquad
 K_2(z)\ge\sqrt{\frac\pi{2z}}e^{-z},\qquad P_5(z)<\frac92.
\]
For $r=|w_j|$ and $d\ge h\ge2$,
\[
 \frac{\pi^3}{120}(r/d)^3P_5(2\pi rd)<\frac34,\qquad
 -C(|w_j|,h)\ge b(h):=\frac{\pi^3\rho^{3/2}}{4h^{5/2}}
                         e^{-(10\pi/3)h}.
\]
Moreover $\partial_r\{r^3P_5(2\pi dr)\}>0$, since
\[
 3+30/z+105/z^2-945/z^4-1890/z^5>0\qquad(z\ge47/5).
\]
Thus
\[
\begin{gathered}
 r\le\alpha d\ \Longrightarrow\
 \frac{\pi^3}{120}(r/d)^3P_5(2\pi rd)
 \le\frac{\pi^3\alpha^3}{120}P_5(2\pi\alpha d^2)<1,\\
 C(r,d)>0\ \Longrightarrow\ r>\alpha d,\qquad
 C(r,d)_+\le\frac{5\pi^6}{384}r^{9/2}d^{-11/2}e^{-2\pi rd}.
\end{gathered}
\]
Set $\delta_j=1-\cos(2\pi w_j\cdot u)$, $j=1,2$. Then
\[
 1-\cos(2\pi k(nw_1+mw_2)\cdot u)
 \le\frac92k^2r^2(\delta_1+\delta_2).
\]
Equation~\eqref{eq:poisson} yields
\[
 \mathcal E\bigl(\Gamma+\mathbb Z(u+he_3)\bigr)
 -\mathcal E(\Gamma\oplus h\mathbb Z)
 \ge\frac{2b(h)-S(h)}A(\delta_1+\delta_2),
\]
where
\[
\begin{aligned}
 S(h)&:=\frac92\sum_{k=1}^{\infty}k^2
       \sum_{w\in\Gamma^*\setminus\{0\}}|w|^2 C(|w|,kh)_+\\
 &\le\frac{15\pi^6}{256}\sum_{k=1}^{\infty}k^2(kh)^{-11/2}
 \sum_{\substack{w\in\Gamma^*\\|w|>\alpha kh}}
                  |w|^{13/2}e^{-2\pi kh|w|}.
\end{aligned}
\]
For $d\ge2$, $j\in\mathbb Z_{\ge0}$,
\[
 \#\{w\in\Gamma^*:\alpha d+j\rho\le|w|<\alpha d+(j+1)\rho\}
 \le(2\alpha d/\rho+2j+3)^2\le26d^2(j+1)^2.
\]
Since $e^{-2\pi\rho d}<1/10000$ and the radial summand decreases,
geometric comparison gives
\[
\begin{gathered}
 \sum_{\substack{w\in\Gamma^*\\|w|>\alpha d}}
       |w|^{13/2}e^{-2\pi d|w|}
 <52\alpha^{13/2}d^{17/2}e^{-2\pi\alpha d^2},\\
 \sum_{k=1}^{\infty}k^5e^{-2\pi\alpha k^2h^2}
 <2e^{-2\pi\alpha h^2},\\
 \frac{S(h)}{2b(h)}
 <\frac{195\pi^3\alpha^{13/2}h^{11/2}}{16\rho^{3/2}}
       e^{-2\pi h(\alpha h-5/3)}<1\quad(h\ge2).
\end{gathered}
\]
The last majorant decreases for $h\ge2$. Since $w_1,w_2$ form a
reciprocal basis, $\delta_1+\delta_2=0\iff u\in\Gamma$, giving the
stated equality case.
This completes the proof of Lemma~\ref{lem:alignment}.
\end{proof}

\begin{lemma}[Explicit nonzero Fourier bounds]\label{lem:fourierhess}
For a planar lattice \(\Gamma\subset\mathbb R^2\), set \(A=|\Gamma|\) and
\[
 R(\Gamma,h)=\frac1A\sum_{k=1}^{\infty}
 \sum_{w\in\Gamma^*\setminus\{0\}}
 \left[\frac{\pi^6|w|^5}{60(kh)^5}K_5(2\pi kh|w|)
       -\frac{2\pi^3|w|^2}{(kh)^2}K_2(2\pi kh|w|)\right].
\]
For this nonzero aligned-layer Fourier term and \(h\ge2\),
\begin{enumerate}
\item If
\[
 \Gamma=\sqrt A\,\Gamma_{x,y},\quad
 \frac25\le A\le\frac9{10},\quad 0\le x\le\frac12,
 \quad\frac{\sqrt3}{2}\le y\le1,
\]
then, with \(A,h\) fixed,
\[
 \|\nabla^2_{x,y}R(\Gamma,h)\|_{\mathrm{op}}<\frac3{8}.
\]
\item If \(A\le9/10\) and \(\min_{v\in\Gamma\setminus\{0\}}|v|\ge7/10\), then
\[
 |R(\Gamma,h)|<\frac1{80}.
\]
\end{enumerate}
\end{lemma}
\begin{proof}
For $r=|w|$, $d=kh$, $q=r^2$, write
\[
 C(q,d)=\frac{\pi^6r^5}{60d^5}K_5(2\pi rd)
             -\frac{2\pi^3r^2}{d^2}K_2(2\pi rd).
\]
The derivative identity in Lemma~\ref{lem:besselproperties} gives
\[
 C_q=-\frac{\pi^7r^4}{60d^4}K_4(2\pi rd)
 +\frac{2\pi^4r}{d}K_1(2\pi rd),\qquad
 C_{qq}=\frac{\pi^8r^3}{60d^3}K_3(2\pi rd)-2\pi^5K_0(2\pi rd).
\]
In (i), the reciprocal squared lengths satisfy
\[
 q=\frac1A\left(yn^2+\frac{(m+nx)^2}{y}\right),\qquad
 |\nabla_{x,y}q|=q/y,\quad\|\nabla^2_{x,y}q\|_{\mathrm{op}}=2q/y^2.
\]
Thus \(r\ge\rho:=19/20\): \(n=0\) gives
\(q\ge1/(Ay)\ge10/9\), and \(n\ne0\) gives
\(q\ge y/A\ge5\sqrt3/9>361/400\).
In particular \(z=2\pi rd>19\pi/5>119/10\).
Lemma~\ref{lem:besselproperties} gives
\[
 K_\nu(2\pi rd)\le\frac{9}{8\sqrt{rd}}e^{-2\pi rd}
 \quad(\nu\in\{0,1,3,4\}).
\]
Therefore
\begin{align*}
 \|\nabla^2_{x,y}C(q,d)\|_{\mathrm{op}}
 &\le\frac{9e^{-2\pi rd}}{8y^2}
 \left[\frac{\pi^8r^{13/2}}{60d^{7/2}}
 +\frac{2\pi^5r^{7/2}}{d^{1/2}}
 +\frac{\pi^7r^{11/2}}{30d^{9/2}}
 +\frac{4\pi^4r^{5/2}}{d^{3/2}}\right].
\end{align*}
Disjoint reciprocal disks of radius \(\rho/2\) show that the
annulus \(j\rho\le r<(j+1)\rho\) contains at most
\(25j^2\) points. The radial summands decrease in this range.
Put \(q_0=e^{-2\pi\rho h}\). For every positive exponent \(a\)
and each radial power \(p\le13/2\) used here,
\[
 \sum_{k=1}^{\infty}(kh)^{-a}\sum_{w\in\Gamma^*\setminus\{0\}}|w|^pe^{-2\pi kh|w|}
 \le\frac{25\rho^ph^{-a}q_0}{(1-512q_0)(1-q_0)}
 <26\rho^ph^{-a}q_0\qquad(q_0\le1/15000).
\]
The annulus majorant \(\sum_{j=1}^{\infty}j^9q_0^{kj}\) has
successive-term ratio \(\le512q_0^k\); the layer series is
geometric, with $q_0<1/140000$.
Since $(Ay^2)^{-1}\le10/3$, and the coefficient bracket including $9/8$
is less than $580$,
\[
 \|\nabla^2_{x,y}R\|_{\mathrm{op}}<\frac38.
\]

For (ii), use the reciprocal separation
\(\rho_A=7/(10A)\), since $\Gamma^*\cong A^{-1}\Gamma$.
Now \(2\pi rd\ge2\pi\rho_Ah\ge28\pi/9>9.7\).
Lemma~\ref{lem:besselproperties} gives
\[
 |C(q,d)|\le e^{-2\pi rd}
 \left(\frac{3\pi^6}{80}r^{9/2}d^{-11/2}
       +\frac{3\pi^3}{2}r^{3/2}d^{-5/2}\right).
\]
Here $P_5(2\pi rd)<9/2$ and $P_2(2\pi rd)<3/2$.
The same disk and layer comparisons give
\[
 |R(\Gamma,h)|\le\frac{26}{A}e^{-2\pi\rho_Ah}
 \left(\frac{3\pi^6}{80}\rho_A^{9/2}h^{-11/2}
       +\frac{3\pi^3}{2}\rho_A^{3/2}h^{-5/2}\right).
\]
Each summand decreases in \(h\ge2\) and increases in
\(A\le9/10\), since the logarithmic derivative of
\(A^{-p-1}e^{-(7\pi/5)h/A}\) has sign
\((7\pi/5)h/A-p-1>0\).
At $(A,h)=(9/10,2)$, the exponential factor is less than $1/16000$
and the coefficient bracket is less than $25/4$. Hence $|R|<1/80$.

The bounds justify differentiation and rearrangement.
This completes the proof of Lemma~\ref{lem:fourierhess}.
\end{proof}

\begin{lemma}[The global fixed-area comparison]\label{lem:globalplanarshape}
Let $\Gamma\subset\mathbb R^2$ be a planar lattice with
\[
 A=|\Gamma|\le\frac9{10},\qquad
 \min_{v\in\Gamma\setminus\{0\}}|v|\ge\frac7{10},\qquad h\ge2.
\]
Then
\[
\begin{gathered}
 \mathcal E(\Gamma\oplus h\mathbb Z)\ge\mathcal E(T_A\oplus h\mathbb Z),\\
 \mathcal E(\Gamma\oplus h\mathbb Z)=\mathcal E(T_A\oplus h\mathbb Z)
 \quad\Longleftrightarrow\quad\Gamma\cong T_A.
\end{gathered}
\]
\end{lemma}
\begin{proof}
Write
\[
 \Gamma=\sqrt A\,\Gamma_{x,y},\qquad
 0\le x\le\tfrac12,\quad x^2+y^2\ge1.
\]
Then
\[
 A\ge\frac{\sqrt3}{2}\left(\frac7{10}\right)^2>\frac25,
 \qquad y=\frac{A}{\lambda_1(\Gamma)^2}\le\frac{90}{49}<2.
\]
By homogeneity and the planar input \eqref{eq:planarinput},
\[
 \begin{aligned}
 \mathcal E(\Gamma)-\mathcal E(T_A)
 &\ge(A^{-6}-A^{-3})\\
 &\qquad\times\bigl[E(\Gamma_{x,y},3)-E(T_1,3)\bigr].
 \end{aligned}
\]
For $y\le1$, Lemma~\ref{lem:planargap} and
$A^{-6}-A^{-3}>1/2$ give
$\tfrac14[(x-1/2)^2+(y-\sqrt3/2)^2]$. Triangular symmetry gives
$\nabla_{x,y}R(\sqrt A\,\Gamma_{x,y},h)=0$ at $(1/2,\sqrt3/2)$.
The segment to $(x,y)$ lies in the rectangle of
Lemma~\ref{lem:fourierhess}(i); Taylor's formula bounds the Fourier
loss by $\tfrac3{16}$ times this squared distance.
For $y\ge1$, instead,
\[
 \mathcal E(\Gamma)-\mathcal E(T_A)
 \ge\frac{3\sqrt3}{2}(y-\sqrt3/2)^2>\frac5{128},
 \quad A=y\lambda_1(\Gamma)^2\ge\frac{49}{100}.
\]
Hence $\lambda_1(T_A)^2=2A/\sqrt3\ge49/100$,
and Lemma~\ref{lem:fourierhess}(ii) bounds the Fourier loss by $1/40$.
The zero mode is independent of shape, so
\[
 \mathcal E(\Gamma\oplus h\Z)-\mathcal E(T_A\oplus h\Z)\ge
 \begin{cases}
 \frac1{16}[(x-1/2)^2+(y-\sqrt3/2)^2],&y\le1,\\[2pt]
 \frac9{640},&y>1.
 \end{cases}
\]
Equality therefore requires $\Gamma\cong T_A$.
This completes the proof of Lemma~\ref{lem:globalplanarshape}.
\end{proof}

\begin{lemma}[Uniform area convexity]\label{lem:areaconvex}
For $V\ge89/50$,
\[
 F_V(A):=\mathcal E\left(T_A\oplus\frac VA\Z\right),\qquad
 F_V''(A)>34\quad(2/5\le A\le89/100).
\]
\end{lemma}
\begin{proof}
\[
\begin{gathered}
 \min_{v\in T_1\setminus\{0\}}|v|^2=\frac2{\sqrt3},\qquad
 E(T_1,6)\ge\frac{81}{32},\\
 \{(i,j)\in\mathbb Z^2:i^2+ij+j^2=1\}
 =\{\pm(1,0),\pm(0,1),\pm(1,-1)\}.
\end{gathered}
\]
Lemma~\ref{lem:planarconstants} gives $E(T_1,3)<17/4$, hence
\[
 U''(A)=21E(T_1,6)A^{-8}-12E(T_1,3)A^{-5}
 >A^{-8}(1701/32-51A^3)
 >37.
\]
The interlayer terms are
\begin{align*}
 &\frac{\pi\zeta(10)A^9}{5V^{10}}
       -\frac{\pi\zeta(4)A^3}{V^4}\\
 &+\sum_{k=1}^{\infty}\sum_{w\in T_1^*\setminus\{0\}}
 \left[\frac{\pi^6|w|^5 A^{3/2}}{60(kV)^5}
 K_5(2\pi|w|kVA^{-3/2})
 -\frac{2\pi^3|w|^2}{(kV)^2}
 K_2(2\pi|w|kVA^{-3/2})\right].
\end{align*}
The repulsive zero mode is strictly convex; by Lemma~\ref{lem:besselproperties}, for nonzero modes,
$z=\kappa A^{-3/2}$ gives
\[
 \frac{\dd^2}{\dd A^2}\{A^{3/2}K_5(z)\}
 =A^{-1/2}\left[\frac34(K_5-zK_5')+\frac94z^2K_5''\right]>0.
\]
The attractive zero-mode loss is less than
\[
 \frac{6\pi\zeta(4)A}{V^4}<\frac{19}{10}.
\]
The differential equation and $K_2'<0$ in Lemma~\ref{lem:besselproperties}
give, for $c>0$ independent of $A$,
\[
 \frac{\dd^2}{\dd A^2}[-cK_2(z)]
 =-\frac{c}{4A^2}[9(z^2+4)K_2+6zK_2']
 \ge-\frac{9c}{4A^2}(z^2+4)K_2.
\]
Put $h=V/A$, $d=kh$ and $r=|w|/\sqrt A$.
Then $c=2\pi^3r^2/(Ad^2)$ and
\[
 z=2\pi rd\ge2\pi VA^{-3/2}>333/25,\qquad
 K_2(z)\le K_{5/2}(z)<\frac{3e^{-2\pi rd}}{4\sqrt{rd}}.
\]
The total nonzero-mode loss is bounded above by
\begin{align*}
 &\frac{27}{2A^3}
 \sum_{k=1}^{\infty}\sum_{w\in T_1^*\setminus\{0\}}e^{-2\pi rd}
   (\pi^5r^{7/2}d^{-1/2}+\pi^3r^{3/2}d^{-5/2})\\
 &\qquad<351\left(\pi^5V^{-1/2}A^{-17/4}
              +\pi^3V^{-5/2}A^{-5/4}\right)e^{-2\pi VA^{-3/2}}.
\end{align*}
Use Lemma~\ref{lem:fourierhess}, since
$e^{-2\pi VA^{-3/2}}<1/600000<1/15000$.
The last majorant decreases in $V$ and increases in $A$.
At $(V,A)=(89/50,89/100)$ this loss is less than $1/4$,
hence $F_V''>34$.
The same exponential bounds justify termwise differentiation twice.
This completes the proof of Lemma~\ref{lem:areaconvex}.
\end{proof}

\begin{lemma}[Three finite area comparisons and common scalar bounds]
\label{num:sh-three-cases}
For the primitive section $\Gamma$ of Lemma~\ref{lem:crossenvelope}, set
$a=\min_{v\in\Gamma\setminus\{0\}}|v|$, $A=|\Gamma|$ and $h=V/A$.
Suppose
\[
 \frac7{10}\le a\le\frac85,\quad
 1\le A\le\frac72,\quad
 A\ge\frac{433}{500}a^2,\quad
 h^2\ge\frac{3A^2}{4a^2}-\frac{a^2}{4}.
\]
For each of the three parameter triples
\[
 (u,v,\tau)\in
 \left\{
 \left(\frac{89}{50},\frac95,-\frac{357}{100}\right),
 \left(\frac95,\frac{37}{20},-\frac{71}{20}\right),
 \left(\frac{37}{20},2,-\frac72\right)
 \right\},
\]
if \(u\le V\le v\), then \(\mathcal E(L)>\tau\), whereas
\[
 \mathcal E(T_{17/20}\oplus(v/(17/20))\mathbb Z)<\tau.
\]
The scalar and Fourier bounds used in
Proposition~\ref{prop:sh-global-178} hold for all $V\ge89/50$.
\end{lemma}

\begin{proof}
For \(z=2\pi rd\), use
\[
 \mathcal U_2(r,d)=\frac{\pi^3r^2}{d^2}\mathcal K_2(z).
\]
Lemma~\ref{lem:criticalcusp} gives
$\pi^3r^2d^{-2}K_2(2\pi rd)\le\mathcal U_2(r,d)$.
Monotonicity in $r,d$ reduces, for the radial derivative, to
\[
 z^3+\frac38z^2-\frac{15}{128}z+\frac{105}{256}>0;
\]
the quadratic formed by the last three terms has negative
discriminant. For $\rho=a/A$, Lemma~\ref{lem:dualtriangle} gives
$\mathcal E(L)-\mathcal E(\Gamma)\ge-J_\triangle(a,A,h)$.

For a reduced basis \((a,0),(t,A/a)\), \(0\le t\le a/2\),
the independent pair has squared length \(\le A^2/a^2+a^2/4\).
One-dimensional Poisson summation gives
\[
 E(\Gamma,3)\le2\zeta(6)a^{-6}
 +\frac{3\pi\zeta(5)}4 B(A/a^2)\frac{a^4}{A^5},
\]
where
\[
 B(y)=1+2\sum_{n=1}^{\infty}e^{-2\pi ny}
          \left(1+2\pi ny+\frac{(2\pi ny)^2}{3}\right).
\]
Geometric summation, $2\zeta(6)<407/200$ and
$\zeta(5)<1037/1000$ give
\[
 B(y)<\frac{23}{20}\quad(y\ge\sqrt3/2),\qquad
 \frac{3\pi\zeta(5)}4B(y)<
 \begin{cases}53/20,&y\ge1,\\2813/1000,&\sqrt3/2\le y<1.
 \end{cases}
\]
Retaining the independent pair's repulsion gives
\[
 \begin{aligned}
 \mathcal E(\Gamma)\ge\max\Bigg\{&
 a^{-12}-\frac{407}{200a^6}
 -\begin{cases}53/20,&A\ge a^2,\\2813/1000,&A<a^2\end{cases}
 \frac{a^4}{A^5}\\
 &+\left(\frac{A^2}{a^2}+\frac{a^2}{4}\right)^{-6},\quad
 a^{-12}-\frac{13}{2a^6},\quad-\frac{947}{320}\Bigg\}.
 \end{aligned}
\]
The covering bound \(R^2\le5a^2/16+A^2/(4a^2)\) retains
from the nearest point in the first shifted layer the repulsion
\[
 \left\{(v/A)^2+\frac{5a^2}{16}
                         +\frac{A^2}{4a^2}\right\}^{-6},
\]
since \(V\le v\); the opposite layers give coefficient one.
For \(a\ge9/8\), the fixed-systole comparison also gives
\(\mathcal E(L)\ge6a^{-12}-(2901/200)a^{-6}\).

All three cases use the same 256 by 256 rational grid
\[
 a_i=\frac{1792+9i}{2560},\qquad
 A_j=1+\frac{5j}{512}\qquad(0\le i,j\le256).
\]
Discard only cells with \(A_{j+1}<(433/500)a_i^2\).
For each retained cell set
\[
 H_{ij}=10^{-6}\left\lfloor10^6
 \sqrt{\max\left\{\frac{u^2}{A_{j+1}^2},
       \frac{3A_j^2}{4a_{i+1}^2}-\frac{a_{i+1}^2}{4}\right\}}
 \right\rfloor,\qquad \rho_{ij}=a_i/A_{j+1}.
\]
Define
\begin{align*}
 \overline J_{ij}={}&\frac{\pi^5}{90A_jH_{ij}^4}+\frac1{1000}\\
 &+\frac2{A_j}\sum_{m=1}^5\Bigg[
 6\!\!\sum_{q\in\{1,3,4,9,12,16,25,27,36\}}
                  \mathcal U_2(\rho_{ij}\sqrt q,mH_{ij})\\
 &\hspace{42mm}+
 12\!\!\sum_{q\in\{7,13,19,21,28,31\}}
                  \mathcal U_2(\rho_{ij}\sqrt q,mH_{ij})\Bigg].
\end{align*}
These multiplicities enumerate exactly \(0<k^2+kl+l^2\le36\).
The reduced-plane bounds give $h\rho>16/25$ and
$\pi^3a^{3/2}V^{-5/2}<15$. The complete tail $m\ge6$ or $q>36$
is $<1/1000$, hence $J_\triangle<\overline J_{ij}$ throughout the cell.

For the $k$th parameter triple, put
\[
 d_{ij}=\begin{cases}
 53/20,&A_j\ge a_{i+1}^2,\\
 2813/1000,&A_j<a_{i+1}^2,
 \end{cases}
\]
and define
\[
\begin{aligned}
 L_{ij}^{(k)}=\max\Bigg\{&
 \left[(v/A_j)^2+\frac{5a_{i+1}^2}{16}
                  +\frac{A_{j+1}^2}{4a_i^2}\right]^{-6}
                  -\overline J_{ij}\\
 &+\max\left\{\begin{aligned}
 &a_{i+1}^{-12}-\frac{407}{200a_i^6}
       -d_{ij}\frac{a_{i+1}^4}{A_j^5}
       +\left(\frac{A_{j+1}^2}{a_i^2}
                      +\frac{a_{i+1}^2}{4}\right)^{-6},\\
 &a_{i+1}^{-12}-\frac{13}{2a_i^6},\qquad-\frac{947}{320}
 \end{aligned}\right\},\\
 &a_{i+1}^{-12}-\frac{2901}{200a_i^6}\Bigg\}.
\end{aligned}
\]
For $a_i\ge9/8$, include the additional candidate
$6a_i^{-12}-(2901/200)a_i^{-6}$ in this maximum. Each case has
56974 retained and 8562 discarded cells, with
\[
 \begin{aligned}
 -3.568443985&<\min L_{ij}^{(1)}<-3.568443984,\\
 -3.515188971&<\min L_{ij}^{(2)}<-3.515188970,\\
 -3.404198849&<\min L_{ij}^{(3)}<-3.404198848.
 \end{aligned}
\]
The minimizing index pairs for the three triples above are
$(151,73)$, $(154,75)$ and $(118,14)$, respectively.

For the scalar checks, let $s_3,s_6$ sum $q^{-3},q^{-6}$ over
\(-16\le k,l\le16\), \(0<q=k^2+kl+l^2\le64\).
The disjoint-disk estimate gives
\[
\begin{gathered}
 s_3\le S_3\le s_3+2(17/16)^6(15/2)^{-4}<\frac{63771}{10000},\\
 s_6\le S_6\le s_6+\frac45(17/16)^{12}(15/2)^{-10}.
\end{gathered}
\]
Lemma~\ref{lem:v2-referencecomparison} bounds the endpoint energy by
\[
\begin{aligned}
 &\left[\frac{27\left(s_6+\frac45(17/16)^{12}(15/2)^{-10}\right)}
                    {128(17/20)^6}
       -\frac{3\sqrt3s_3}{8(17/20)^3}\right]
       -\frac{\pi^5(17/20)^3}{90v^4}\\
 &\quad+\frac{11/10}{5(49/100)^2}
 \left(\frac{17/20}{v}\right)^{10}
 \left(1+\frac{(49/100)(17/20)}v\right)^{12}
 <\begin{cases}
 -3.5744333,&v=9/5,\\
 -3.5555334,&v=37/20,\\
 -3.5101664,&v=2.
 \end{cases}
\end{aligned}
\]
The common nonzero attractive Fourier upper bound is the finite sum
\[
 \frac2{89/100}\sum_{m=1}^8
 \sum_{\substack{(k,l)\in\mathbb Z^2\\-10\le k,l\le10\\(k,l)\ne(0,0)}}
 \mathcal U_2\!\left(\frac7{10}\sqrt{k^2+kl+l^2},m(89/50)\right)
 +10^{-12}<\frac{73}{5000}.
\]
The tail allowance follows from the block bound
\(2000(K^6+K^3)e^{-5mK}\), \(K=\max(|k|,|l|)\). Its two
omitted regions are at most
\[
 \frac{100}{99}\frac{50}{49}\,2000(11^6+11^3)100^{-11}
 +\frac{200}{99}\,4000\,100^{-9}<10^{-12}.
\]
For the remaining comparisons, use
\[
\begin{gathered}
 \frac{3\sqrt3}{8}s_3\le E(T_1,3)
 \le\frac{3\sqrt3}{8}\frac{63771}{10000},\\
 \frac{27}{64}s_6\le E(T_1,6)
 \le\frac{27}{64}\left(s_6+\frac45(17/16)^{12}(15/2)^{-10}\right).
\end{gathered}
\]
Then
\[
\begin{aligned}
 &21\left(\frac{27}{64}s_6\right)
 -12\left(\frac{3\sqrt3}{8}\frac{63771}{10000}\right)>0,\\
 &\frac{\frac{27}{64}\left(s_6+\frac45(17/16)^{12}(15/2)^{-10}\right)}2
       \bigl[(89/100)^{-6}-(17/20)^{-6}\bigr]\\
 &\qquad-\left(\frac{3\sqrt3}{8}s_3\right)
       \bigl[(89/100)^{-3}-(17/20)^{-3}\bigr]>\frac{291}{5000},\\
 &3(89/100)^{-7}\Bigg[
       \left(\frac{3\sqrt3}{8}s_3\right)(89/100)^3
       -\frac{27}{64}\left(s_6+\frac45(17/16)^{12}(15/2)^{-10}\right)
       \Bigg]\\
 &\qquad-\frac{3(121/35)}{(89/50)^4}>\frac{157}{100}.
\end{aligned}
\]
The first inequality gives $U''>0$ up to area one; the third gives
$\partial_A\{U(A)-\pi\zeta(4)A^3/V^4\}>157/100$ for
$89/100\le A\le1$ and $V\ge89/50$.
For $V\ge89/50$, the zero-mode and reference-energy losses satisfy
\[
\begin{aligned}
 \frac{\pi\zeta(4)}{V^4}
       \bigl[(89/100)^3-(17/20)^3\bigr]&<\frac{157}{5000},\\
 \mathcal E\!\left(T_{17/20}\oplus\frac{V}{17/20}\mathbb Z\right)
 -U(17/20)+\frac{\pi\zeta(4)(17/20)^3}{V^4}&<\frac{71}{10000}.
\end{aligned}
\]
The second bound uses Lemma~\ref{lem:v2-referencecomparison};
Lemma~\ref{lem:areaconvex} supplies the area Hessian bounds for all modes.

Enclose every finite sum, tail and scalar bound as in
Lemma~\ref{num:v2-area}, using the rational $\pi$ enclosure and
binary-scaled bounds $T_{21}\le e^{-z}\le T_{20}$ of
Lemma~\ref{lem:finitearithmetic}. Check radical endpoints by integer
squaring and height floors by exact integer square roots.
This completes the proof of Lemma~\ref{num:sh-three-cases}.
\end{proof}

\begin{proposition}[The effective simple-hexagonal phase]
\label{prop:sh-global-178}
For $V\ge89/50$,
\[
 \exists!\,A(V)\in(2/5,89/100):\quad F_V'(A(V))=0,
 \qquad\frac{V}{A(V)}>2,
\]
and
\[
 \operatorname*{argmin}_{|L|=V}\mathcal E(L)
 =\left\{Q\!\left(T_{A(V)}\oplus\frac{V}{A(V)}\mathbb Z\right):
                           Q\in O(3)\right\}.
\]
\end{proposition}

\begin{proof}
Since $V/A_*>1$, $m(V)\le F_V(A_*)<e_2$, so every minimizer
has $7/10<a<8/5$ and $A<5$. The $A\ge7/2$ exclusion in
Lemma~\ref{lem:area-envelope-uniform} holds for all $V\ge89/50$.

For $1\le A\le7/2$ and each case $u\le V\le v$ of
Lemma~\ref{num:sh-three-cases},
\[
 \mathcal E(L)>\tau>\mathcal E(T_{17/20}\oplus(v/(17/20))\Z)
 \ge \mathcal E(T_{17/20}\oplus(V/(17/20))\Z).
\]
All affected distances are $\ge V/(17/20)>1$, where
$r^{-12}-2r^{-6}$ increases with $r$.
For $V\ge2$, Lemma~\ref{num:v2-area} gives $\mathcal E(L)>-10/3>e_2$.
Thus every minimizer has $A<1$.

For $89/100\le A\le1$ and $V\ge89/50$,
Lemma~\ref{num:sh-three-cases} gives
\[
 \partial_A\{U(A)-\pi\zeta(4)A^3/V^4\}>157/100>0.
\]
Compare the scalar part at $A=89/100$ with the reference at
$17/20$. The planar gain is $>291/5000$; the zero-mode, attractive
Fourier and reference-repulsion losses are respectively
$<157/5000,73/5000,71/10000$. Hence
\[
 \mathcal E(L)-\mathcal E(T_{17/20}\oplus(V/(17/20))\Z)
 >\frac{51}{10000}>0.
\]
Consequently $2/5<(\sqrt3/2)(7/10)^2\le A<89/100$ and
$h=V/A>2$.

Lemmas~\ref{lem:alignment} and~\ref{lem:globalplanarshape} now give
\[
 L\cong T_A\oplus(V/A)\Z,\qquad F_V'(A)=0.
\]
Attainment gives a stationary area in $(2/5,89/100)$;
$F_V''>34$ by Lemma~\ref{lem:areaconvex} gives uniqueness.
This completes the proof of Proposition~\ref{prop:sh-global-178}.
\end{proof}

For the dilute asymptotics, Lemma~\ref{lem:planarconstants} gives
$2/5<A_*<43/50$. By \eqref{eq:planararea},
\begin{equation}\label{eq:areahess}
 U'(A_*)=0,\qquad U''(A_*)=\frac{9S_3^2}{S_6A_*^2}>0.
\end{equation}
Equations \eqref{eq:poisson}--\eqref{eq:zeromode} give, uniformly
with any fixed number of $A$-derivatives near $A_*$,
\begin{equation}\label{eq:scalarenergy}
 F_V(A)=U(A)-\pi\zeta(4)A^3V^{-4}
       +\frac{\pi\zeta(10)}5 A^9V^{-10}+O(e^{-cV}).
\end{equation}
Since $A_*,A(V)\in(2/5,89/100)$, $F_V'(A(V))=0$ and
$F_V''>34$, the mean value theorem yields
\[
 \delta:=A(V)-A_*,\qquad
 |\delta|\le\frac{|F_V'(A_*)|}{34}=O(V^{-4}).
\]
Consequently,
\begin{align*}
 0&=\frac{9S_3^2}{S_6A_*^2}\delta-3\pi\zeta(4)A_*^2V^{-4}
       +O(\delta^2+|\delta|V^{-4}+V^{-10}),\\
 \delta&=\frac{\pi\zeta(4)S_6A_*^4}{3S_3^2}V^{-4}+O(V^{-8}),\\
 m(V)&=e_2-\pi\zeta(4)A_*^3V^{-4}
       +\frac{9S_3^2}{2S_6A_*^2}\delta^2
       -3\pi\zeta(4)A_*^2\delta V^{-4}\\
 &\qquad+\frac{\pi\zeta(10)}5 A_*^9V^{-10}+O(V^{-12}).
\end{align*}
Substitute $\delta$ to obtain \eqref{eq:areaexp} and \eqref{eq:energyexp}.

\section{Stationarity and endpoint comparisons}\label{sec:compact}

Stationarity gives $m(V_1)<m(V_2)$ for $1\le V_1<V_2$.
Korkine--Zolotarev reduction \cite{KZ,LLS} and stationary-basin
estimates give the endpoint minima at $V=21/20,11/10$.

\begin{proposition}[Stationarity and stability at fixed covolume]
\label{prop:stress-stability}
For a lattice $L\subset\R^3$, with $\mathcal E(L)=\tfrac12E(L,6)-E(L,3)$,
fixed-covolume stationarity is characterized by
\begin{equation}\label{eq:stationary-stress}
\begin{gathered}
 \left.\frac{\dd}{\dd t}\mathcal E(e^{tH/2}L)\right|_{t=0}=0
 \quad\text{for every }H=H^{\mathsf T},\ \tr H=0\\
 \Longleftrightarrow\quad
 \sum_{v\in L\setminus\{0\}}(|v|^{-8}-|v|^{-14})vv^T
       =\frac{E(L,3)-E(L,6)}3I.
\end{gathered}
\end{equation}
Under \eqref{eq:stationary-stress}, for $H=H^{\mathsf T}$, $\tr H=0$,
\begin{align}\label{eq:general-shape-hessian}
 \left.\frac{\dd^2}{\dd t^2}\mathcal E(e^{tH/2}L)\right|_{t=0}
 ={}&(E(L,3)-E(L,6))\tr(H^2)\notag\\
 &+3\sum_{v\in L\setminus\{0\}}
       (7|v|^{-16}-4|v|^{-10})(v^THv)^2.
\end{align}
Consequently
\begin{equation}\label{eq:stable-trace-inequality}
 \begin{gathered}
 L\in\operatorname*{argmin}_{\substack{\rm loc\\|\widetilde L|=|L|}}
       \mathcal E(\widetilde L)\\
 \Longrightarrow\quad E(L,3)\le3E(L,6),\qquad \min_{v\in L\setminus\{0\}}|v|<3^{1/6}.
 \end{gathered}
\end{equation}
\end{proposition}
\begin{proof}
For symmetric trace-free $H$,
\begin{align*}
 \left.\partial_{\epsilon}\mathcal E(e^{\epsilon H/2}L)\right|_{\epsilon=0}&=3\sum_{v\in L\setminus\{0\}}(|v|^{-8}-|v|^{-14})v^THv,\\
 \left.\partial_{\epsilon}^2\mathcal E(e^{\epsilon H/2}L)\right|_{\epsilon=0}&=3\sum_{v\in L\setminus\{0\}}\left[
 (|v|^{-8}-|v|^{-14})v^TH^2v
 +(7|v|^{-16}-4|v|^{-10})(v^THv)^2\right].
\end{align*}
Thus stationarity is equivalent to \eqref{eq:stationary-stress}, and
\eqref{eq:general-shape-hessian} follows.
For a Frobenius-orthonormal basis $H_1,\ldots,H_5$ of
$\{H\in\R^{3\times3}:H=H^{\mathsf T},\ \tr H=0\}$,
\[
 \sum_{j=1}^5(u^TH_ju)^2=\|uu^T-I/3\|_F^2=2/3
 \qquad(|u|=1).
\]
Consequently
\begin{equation}\label{eq:hessian-trace-identity}
 \sum_{j=1}^5\left.\partial_{\epsilon}^2\mathcal E(e^{\epsilon H_j/2}L)\right|_{\epsilon=0}=9E(L,6)-3E(L,3).
\end{equation}
At a local minimum, $\sum_{j=1}^5\left.\partial_{\epsilon}^2\mathcal E(e^{\epsilon H_j/2}L)\right|_{\epsilon=0}\ge0$ and
\[
 E(L,3)=\sum_{v\in L\setminus\{0\}}|v|^6|v|^{-12}
       >\lambda_1(L)^6E(L,6),
\]
hence \eqref{eq:stable-trace-inequality}. At every fixed order,
$C(|v|^{-6}+|v|^{-12})$ bounds the local derivatives and justifies differentiation.
This completes the proof of Proposition~\ref{prop:stress-stability}.
\end{proof}

\begin{proposition}[A positive dilation derivative]\label{prop:positive-dilation}
If $|L|\ge1$ and \eqref{eq:stationary-stress} holds, then
\begin{equation}\label{eq:positive-dilation}
 E(L,3)>E(L,6),\qquad
 \left.\frac{\dd}{\dd t}\mathcal E(e^tL)\right|_{t=0}
                   =6(E(L,3)-E(L,6))>0.
\end{equation}
For $L\in\operatorname*{argmin}_{\substack{\rm loc\|\widetilde L|=|L|}}\mathcal E(\widetilde L)$
and $|L|\ge1$, one has $\mathcal E(L)<0$ and
\[
 -\frac{6\mathcal E(L)}5\le E(L,3)<-2\mathcal E(L),\qquad
 -\frac{2\mathcal E(L)}5\le E(L,6)<-2\mathcal E(L).
\]
\end{proposition}
\begin{proof}
Put $W=\operatorname{span}_{\R}\{v\in L:|v|\le1\}$.
For $W\ne\R^3$, take $n\in W^\perp$, $|n|=1$.
Then \eqref{eq:stationary-stress} gives
\[
 \frac{E(L,3)-E(L,6)}3
  =\sum_{\substack{v\in L\\|v|>1}}
                (|v|^{-8}-|v|^{-14})(v\cdot n)^2>0.
\]
The sum is positive since $\operatorname{span}_{\R}L=\R^3$. If $W=\R^3$, choose
$v_1,v_2,v_3\in L$ with $\det(v_1,v_2,v_3)\ne0$, $\max_{1\le i\le3}|v_i|\le1$,
and $k=[L:\sum_{i=1}^3\Z v_i]$. Hadamard gives
\[
 1\le |L|\le k|L|=|\det(v_1,v_2,v_3)|
                      \le|v_1||v_2||v_3|\le1.
\]
Thus $|L|=k=1$, $L\cong\Z^3$ and $E(L,3)>E(L,6)$ termwise.
Homogeneity gives \eqref{eq:positive-dilation}; at a local minimum,
\[
 E(L,6)<E(L,3)\le3E(L,6),\qquad
 \mathcal E(L)=\tfrac12E(L,6)-E(L,3),
\]
which give $\mathcal E(L)<0$ and the stated moment bounds.
This completes the proof of Proposition~\ref{prop:positive-dilation}.
\end{proof}

\begin{proposition}[Monotonicity of the minimum energy]\label{prop:minimum-increasing}
The minimum $m(V)$ is strictly increasing on $[1,\infty)$.
\end{proposition}
\begin{proof}
For $1\le V_1<V_2$, \eqref{eq:coercive}, \eqref{eq:uniformtail}
and Mahler's criterion \cite{Mahler} give compact sublevel sets
\[
 \left\{L:V_1\le|L|\le V_2,\quad
          \tfrac12E(L,6)-E(L,3)\le C\right\}\qquad(C\in\mathbb R).
\]
Hence a minimizer $L_*$ on $V_1\le|L|\le V_2$ exists and satisfies
\eqref{eq:stationary-stress}.
If $|L_*|>V_1$, Proposition~\ref{prop:positive-dilation}
gives, for small $t>0$,
\[
 V_1\le|e^{-t}L_*|<|L_*|,\qquad
 \mathcal E(e^{-t}L_*)<\mathcal E(L_*),
\]
a contradiction. Hence
\[
 \operatorname*{argmin}_{V_1\le|L|\le V_2}\mathcal E(L)
 \subset\{L:|L|=V_1\},\qquad m(V_2)>m(V_1).
\]
This completes the proof of Proposition~\ref{prop:minimum-increasing}.
\end{proof}

\begin{lemma}[Uniform separation on a covolume strip]
\label{lem:uniform-strip-separation}
Suppose
\[
 \frac{21}{20}\le V=|L|\le\frac{11}{10},\qquad
 \mathcal E(L)\le-\frac{124}{25}.
\]
Then
\begin{equation}\label{eq:uniform-strip-separation}
 \lambda_1(L)>\frac45,
 \qquad \lambda_1(L^*)>\frac{10}{13}.
\end{equation}
Equivalently, every primitive planar section has layer spacing
less than \(13/10\).
\end{lemma}

\begin{proof}
The general separation bound gives $\lambda_1(L)>7/10$.

For a primitive planar section \(\Gamma\) of area \(A=V/h\),
write \(L=\Gamma+\mathbb Z(u+h\nu)\). If
\(\lambda_1(L)>a_0\), then
\(\lambda_1(\Gamma^*)>a_0/A\). A Gauss-reduced reciprocal
basis consequently gives
\[
 |w_{i,j}|\ge\frac{a_0h}{V}\sqrt{i^2+j^2-|ij|}.
\]
The bounds $S_3<63771/10000$ and $S_6\ge6$ give
\[
 \mathcal E(\Gamma)\ge e_2=-\frac{S_3^2}{2S_6}
 >-\frac{(63771/10000)^2}{12}>-\frac{17}{5}.
\]

Using $\mathcal K_2,\mathcal K_5$ from Lemma~\ref{lem:criticalcusp}, set
\[
 r_{i,j}=\frac{a_0h}{V}\sqrt{i^2+j^2-|ij|},
\]
and define
\begin{align*}
 R(V,h,a_0)=\frac hV\sum_{m=1}^{\infty}\sum_{(i,j)\in\mathbb Z^2\setminus\{(0,0)\}}
 \biggl\{&\frac{\pi^6}{60}
      \left(\frac{r_{i,j}}{mh}\right)^5
            \mathcal K_5(2\pi mh r_{i,j})\\
 &+2\pi^3\left(\frac{r_{i,j}}{mh}\right)^2
            \mathcal K_2(2\pi mh r_{i,j})\biggr\}.
\end{align*}
Substitute the reciprocal-length floors by monotonicity of
$r^\nu K_\nu(2\pi dr)$. Poisson summation then gives
\begin{equation}\label{eq:uniform-strip-fourier-lower}
 \mathcal E(L)>-\frac{17}{5}
       -\frac{\pi^5}{90Vh^3}
       +\frac{\pi^{11}}{467775Vh^9}
       -R(V,h,a_0).
\end{equation}
The attractive nonzero-mode coefficient is $2\pi^3$. Use
\[
 (H,a_0)=\left(\frac{27}{20},\frac7{10}\right),
 \qquad
 (H,a_0)=\left(\frac{13}{10},\frac45\right).
\]
Since $H^6>3\zeta(10)/(5\zeta(4))$, the zero-mode bracket is
bounded below at $(21/20,H)$ for $h\ge H$. The summands of $R$ are
\[
 V^{j-\nu-1/2}h^{-2j}e^{-ch^2/V},\qquad c>0,
 \quad0\le j,\quad\nu\le5,
\]
with
\[
\begin{aligned}
 \partial_V\log(\cdot)&=\frac{ch^2/V+j-\nu-1/2}{V}>0,
   \qquad ch^2/V>11/2,\\
 \partial_h\log(\cdot)&=-\frac{2(j+ch^2/V)}h<0.
\end{aligned}
\]
Here $ch^2/V\ge2\pi a_0H^2/(11/10)$. Hence
$R(V,h,a_0)\le R(11/10,H,a_0)$.

The endpoint estimates are
\[
 -\frac{\pi^5}{90(21/20)H^3}
 +\frac{\pi^{11}}{467775(21/20)H^9}>
 \begin{cases}
 -319/250,&(H,a_0)=(27/20,7/10),\\
 -709/500,&(H,a_0)=(13/10,4/5),
 \end{cases}
\]
\[
 R(11/10,H,a_0)<
 \begin{cases}
 16/125,&(H,a_0)=(27/20,7/10),\\
 107/1000,&(H,a_0)=(13/10,4/5).
 \end{cases}
\]
Retain $m\in\{1,\ldots,8\}$ and
$(i,j)\in[-10,10]^2\cap\Z^2\setminus\{0\}$ in $R$.
For $k=\max\{|i|,|j|\}$, each omitted block is bounded by
\[
 2000(k^6+k^3)e^{-5mk}.
\]
Indeed \(3k^2/4\le i^2+j^2-|ij|\le2k^2\), the Bessel argument
exceeds \(5mk\), and the two Bessel polynomials are at most
\(8\) and \(3/2\). Their square-root prefactors are less than
one; the coefficients before \((i^2+j^2-|ij|)^{5/2}\) and
\(i^2+j^2-|ij|\) are less than \(40\) and \(100\), respectively.
Each block has $8k$ pairs; $e^{-5}<1/100$ bounds the omitted part by
\[
 \frac{100}{99}\frac{50}{49}\,
       2000(11^6+11^3)100^{-11}
 +\frac{200}{99}\,4000\,100^{-9}<10^{-12}.
\]
Both endpoint tests use the rational Machin, square-root and
exponential enclosures of Lemma~\ref{lem:firstloss}.

At the first endpoint, \eqref{eq:uniform-strip-fourier-lower}
gives, if \(h\ge27/20\),
\[
 \mathcal E(L)>-\frac{124}{25}.
\]
Hence $h<27/20$ for every primitive layer.
Let $a=\lambda_1(L)$ and let $\delta$ be the first projected
Gauss-reduced length. The lifted pair gives a primitive plane with
$A=a\delta$, $h=V/A$, and
\[
 h^2\ge\frac34\delta^2,\qquad
 A^2\le\frac{2aV}{\sqrt3}.
\]
For $a\le4/5$, the layer bound and $\sqrt3>173/100$ give
\[
 \frac7{10}<a\le\frac45,\qquad
 \frac79<A<\frac{101}{100}.
\]
The nonnegative Fourier transform of $(x^2+z^2)^{-3}$ makes its
periodization maximal at zero shift. For even nonnegative decreasing $g$,
$\sum_{k\in\Z}g(ak)\le g(0)+a^{-1}\int_\R g$; apply this along
$(a,0),(s,A/a)$ to obtain
\begin{align*}
 \mathcal E(\Gamma)&\ge\zeta(12)a^{-12}-2\zeta(6)a^{-6}
 -2\zeta(6)(a/A)^6-\frac{3\pi\zeta(5)}4\frac{a^4}{A^5},\\
 \mathcal E(L)-\mathcal E(\Gamma)&\ge-2\zeta(6)h^{-6}
 -\frac{3\pi\zeta(5)}4(1/a+a/A)h^{-5}
 -\frac{\pi\zeta(4)}A h^{-4}.
\end{align*}
The integrals are $3\pi/(8z^5)$ and $\pi/(2z^4)$;
the omitted repulsive terms are nonnegative. With
$2\zeta(6)<21/10$, $3\pi\zeta(5)/4<99/40$,
$\pi\zeta(4)<341/100$ and $\zeta(12)>1$ and $V\ge21/20$, we obtain $\mathcal E(L)>P(a,A)$, where
\begin{align*}
 P(a,A)={}&a^{-12}-\frac{21}{10}a^{-6}
 -\frac{21a^6}{10A^6}-\frac{99a^4}{40A^5}
 -\frac{21A^6}{10(21/20)^6}\\
 &-\frac{99A^5}{40a(21/20)^5}
 -\frac{99aA^4}{40(21/20)^5}-\frac{341A^3}{100(21/20)^4}.
\end{align*}
On this rectangle,
\[
 \partial_aP\le-12(5/4)^{13}+\frac{63}{5}(5/4)^7
 +\frac{99}{40}\frac{(101/100)^5}{(7/10)^2(21/20)^5}
 <-150.
\]
For fixed $a$, $P$ is strictly concave in $A$.
At the endpoints,
\[
 P(4/5,7/9)>-\frac52,\qquad
 P(4/5,101/100)>-\frac{37}{10}.
\]
Thus $\mathcal E(L)>-37/10$, contradicting the hypothesis; hence
$a>4/5$. For $h\ge13/10$, the second endpoint gives
\[
 \mathcal E(L)>-\frac{124}{25}.
\]
Therefore $h<13/10$ for every primitive layer, equivalently
$\lambda_1(L^*)>10/13$.
This completes the proof of Lemma~\ref{lem:uniform-strip-separation}.
\end{proof}

\begin{proposition}[A uniform Korkine--Zolotarev domain]
\label{prop:uniform-strip-kz}
Under the hypotheses of Lemma~\ref{lem:uniform-strip-separation},
the lattice has a Korkine--Zolotarev-reduced basis of the form
\[
 B=\begin{pmatrix}
 a&ax&ay\\0&\delta&\delta z\\0&0&h
 \end{pmatrix},\qquad
 a\delta h=V,
 \quad0\le x,z\le\frac12,\quad-\frac12\le y\le\frac12,
\]
with
\begin{equation}\label{eq:uniform-strip-physical-kz}
 \frac45<a<\frac{29}{25},\qquad
 \frac{69}{100}<\delta<\frac{127}{100},\qquad
 \frac{22}{25}<h<\frac{13}{10}.
\end{equation}
Moreover \(21/26<a\delta<5/4\).
For the unit-covolume shape \(\Lambda=V^{-1/3}L\), write
\((a_0,\delta_0,h_0)=V^{-1/3}(a,\delta,h)\). Then
\begin{equation}\label{eq:uniform-strip-unit-kz}
 \frac{77}{100}<a_0\le2^{1/6},\qquad
 \frac{69}{100}<\delta_0<\frac{123}{100},\qquad
 \frac78<h_0<\frac{32}{25},
\end{equation}
and
\[
 \frac{25}{32}<a_0\delta_0<\frac87,\qquad
 \lambda_1(\Lambda^*)>\frac{39}{50}.
\]
\end{proposition}

\begin{proof}
Choose a reduced basis $B$ with
\[
 a=\min_{n\in\Z^3\setminus\{0\}}|Bn|,\qquad
 \delta^2=\min_{(j,k)\in\Z^2\setminus\{0\}}
          \{\delta^2(j+zk)^2+h^2k^2\}.
\]
Size-reduced lifts give the stated shear ranges and
\[
 \delta^2\ge a^2(1-x^2),\qquad
 h^2\ge\delta^2(1-z^2)\ge\frac34\delta^2.
\]
The Hermite bound $a^6\le2V^2$ and the preceding lemma give
\[
 a<\frac{29}{25},\qquad
 \delta=\frac{V}{ah}>\frac{69}{100},\qquad
 \delta^2\le\frac{2V}{\sqrt3a}<\left(\frac{127}{100}\right)^2.
\]
Also $h\ge\sqrt3\delta/2$ and $a\delta h=V$ imply
\[
 h^2\ge\frac{\sqrt3V}{2a}>\left(\frac{22}{25}\right)^2.
\]
Now $a\delta=V/h$ gives the product bounds. Normalization yields
$a_0>77/100$, $h_0<32/25$, $a_0\le2^{1/6}<1123/1000$;
$a_0\delta_0h_0=1$ then gives $\delta_0>69/100$ and
\[
 \delta_0^2\le\frac{2}{\sqrt3a_0}
 <\left(\frac{123}{100}\right)^2,\qquad
 h_0^2\ge\frac{\sqrt3}{2a_0}>
 \left(\frac78\right)^2.
\]
Inverting $h_0$ gives the product bounds. Since $\Lambda^*=V^{1/3}L^*$,
\[
 \lambda_1(\Lambda^*)>(21/20)^{1/3}(10/13)>39/50.
\]
This completes the proof of Proposition~\ref{prop:uniform-strip-kz}.
\end{proof}

\begin{proposition}[A common endpoint comparison domain]
\label{prop:fccKZdomain}\label{prop:sh110domain}
Suppose $21/20\le V=|L|\le11/10$ and $\mathcal E(L)\le-124/25$.
Then $L\cong U(p)\Z^3$, where
\[
 U(p)=\begin{pmatrix}a&ax&ay\\0&\delta&\delta z\\0&0&h\end{pmatrix},
 \qquad h=\frac{V}{a\delta},\quad p=(a,\delta,x,y,z),
\]
and $p$ belongs to the compact domain $\mathcal D(V)$ defined by
\begin{equation}\label{eq:fccKZbox}
 \frac45\le a\le\frac{29}{25},\quad
 \frac{69}{100}\le\delta\le\frac{127}{100},\quad
 0\le x,z\le\frac12,\quad-\frac12\le y\le\frac12,
\end{equation}
\begin{equation}\label{eq:fccKZconstraints}
\begin{gathered}
 \frac{22}{25}\le h\le\frac{13}{10},\qquad
 \frac{10V}{13}\le a\delta\le\frac{25V}{22},\\
 \delta^2+a^2x^2\ge a^2,\qquad h^2+\delta^2z^2\ge\delta^2,\\
 a^2(i+xj+yk)^2+\delta^2(j+zk)^2+h^2k^2\ge a^2,\\
 a^{-2}i^2+\delta^{-2}(j-xi)^2
 +h^{-2}\{k-zj+(xz-y)i\}^2\ge\frac{100}{169},\\
 (i,j,k)\in\{-1,0,1\}^3\setminus\{0\}.
\end{gathered}
\end{equation}
In particular, put $T_F=-26537251/5000000$ and
$\mathcal D_F=\mathcal D(21/20)$.
\end{proposition}
\begin{proof}
Lemma~\ref{lem:uniform-strip-separation} and
Proposition~\ref{prop:uniform-strip-kz} give the box, the Gauss
inequalities and
\[
 \lambda_1(L)>4/5,\quad\lambda_1(L^*)>10/13,\qquad
 \frac{10V}{13}\le a\delta=V/h\le\frac{25V}{22}.
\]
For the remaining primal indices,
\[
 \begin{aligned}
 |k|\ge2&\Longrightarrow h^2k^2\ge4(22/25)^2>(29/25)^2,\\
 |k|\le1,\ |j|\ge2&\Longrightarrow
       \delta^2(j+zk)^2\ge9\delta^2/4>a^2,\\
 |j|,|k|\le1,\ |i|\ge2&\Longrightarrow
       a^2(i+xj+yk)^2\ge a^2.
 \end{aligned}
\]
For the dual indices outside this cube,
\[
 |U^{-\mathsf T}n|^2\ge
 \frac4{\|U\|_1\|U\|_\infty}
 \ge\frac4{(503/200)(58/25)}>\frac{100}{169}.
\]
The positive diagonal floors give compactness;
$a\delta\in[21/26,105/88]$ at $V=21/20$ and
$a\delta\in[11/13,5/4]$ at $V=11/10$.
This completes the proof of Proposition~\ref{prop:fccKZdomain}.
\end{proof}

\begin{lemma}[A uniform FCC stationary basin]
\label{lem:fccuniformbasin}
Let
\[
 t=\left(\frac{441}{1600}\right)^{1/6},\qquad L_F=tD_3,
 \qquad D_3=\{(i,j,k)\in\mathbb Z^3:i+j+k\in2\mathbb Z\},
\]
so that $|L_F|=21/20$.
In the chart \eqref{eq:common-metric-chart}, put
\[
 \begin{gathered}
 L(S)=\det(I+S)^{-1/6}(I+S)^{1/2}L_F,\qquad |L(S)|=21/20,\\
 \Psi(S)=\tfrac12E(L(S),6)-E(L(S),3),\qquad
 \nabla_{\mathbf s}\Psi=(\partial_u\Psi,\partial_v\Psi,
              \partial_x\Psi,\partial_y\Psi,\partial_z\Psi)^{\mathsf T}.
 \end{gathered}
\]
Then
\[
 \left.\partial_{\epsilon}^2\Psi(S+\epsilon T)\right|_{\epsilon=0}>\frac1{1000}\|T\|_F^2
 \quad\left(\|S\|_F\le\frac1{50},\ T\ne0, T=T^{\mathsf T},\ \operatorname{tr}T=0\right),
\]
and
\[
 \mathbf s^{\mathsf T}\nabla_{\mathbf s}\Psi(S)
 =\left.\partial_{\epsilon}\Psi((1+\epsilon)S)\right|_{\epsilon=0}>\frac1{1000}\|S\|_F^2
 \qquad\left(0<\|S\|_F\le\frac{43}{1000}\right).
\]
Consequently
\[
 \nabla_{\mathbf s}\Psi(S)=0\iff S=0
 \qquad(\|S\|_F\le43/1000),
\]
and
\[
 \Psi(S)\ge \Psi(0)+\frac1{2000}\|S\|_F^2
 \qquad\left(\|S\|_F\le\frac{43}{1000}\right).
\]
Convexity is asserted only on $\|S\|_F\le1/50$.
\end{lemma}

\begin{proof}
Cubic symmetry gives $\nabla_{\mathbf s}\Psi(0)=0$. Use the chart and polynomial
construction of Lemma~\ref{lem:fctuniformbasin}, with its $D,p_3,p_6$.
For $w=t(i,j,k)$, write
\[
 q_w=t^2(i^2+j^2+k^2),\qquad
 \ell_w(S)=\frac{u(2i^2-j^2-k^2)+v(j^2-k^2)+2xij+2yik+2zjk}
 {i^2+j^2+k^2}.
\]
Sum, with both signs, over the exact vector set
\[
 \mathcal F=\left\{t(i,j,k):
 \begin{gathered}
 (i,j,k)\in[-40,40]^3\cap\Z^3,\quad i+j+k\in2\Z,\\
 0<i^2+j^2+k^2\le1600
 \end{gathered}\right\},
\]
and form $P$ by Lemma~\ref{lem:fctuniformbasin}. Since $t^{-6}=1600/441$,
its coefficients are rational, $\deg P\le12$, and cubic symmetry
gives $\nabla_{\mathbf s}P(0)=0$.

Since
\[
 \lambda_1(L_F)>1,\qquad
 \min_{w\in L_F\setminus(\mathcal F\cup\{0\})}|w|>32,
\]
\eqref{eq:fixedsystole-tail} gives
\[
 \sum_{w\in L_F\setminus(\mathcal F\cup\{0\})}q_w^{-s}\le T_s(1,32),\qquad s=3,6.
\]
The finite rational sums give
\[
 \sum_{w\in\mathcal F}q_w^{-3}+T_3(1,32)<\frac{164}{25},\qquad
 \sum_{w\in\mathcal F}q_w^{-6}+T_6(1,32)<\frac52.
\]

With $E^+(L_F,3)=164/25$, $E^+(L_F,6)=5/2$ and tails $T_s(1,32)$,
\eqref{eq:common-basin-remainder} gives
\[
 \varepsilon(1/50)<3/1000,\qquad
 \varepsilon(43/1000)<551/100000.
\]
In coordinates $(u,v,x,y,z)$,
\[
 \left.\partial_{\epsilon}^2P(S+\epsilon T)\right|_{\epsilon=0}\ge|\eta|^{\mathsf T}C_H|\eta|
 \qquad(\|S\|_F\le1/50),
\]
where $\eta=(T_{11}/2,(T_{22}-T_{33})/2,T_{12},T_{13},T_{23})^{\mathsf T}$ and
\[
 C_H=\frac1{10^6}
 \begin{pmatrix}
 4078438&-123613&-296771&-296771&-585394\\
 -123613&1356119&-298719&-298719&-3375\\
 -296771&-298719&455411&-55420&-55466\\
 -296771&-298719&-55420&455411&-55466\\
 -585394&-3375&-55466&-55466&578379
 \end{pmatrix}.
\]
The $LDL^{\mathsf T}$ pivots of
$C_H-(4/1000)\operatorname{diag}(6,2,2,2,2)$ exceed
\[
 4,\qquad1,\qquad\frac3{10},\qquad\frac14,\qquad\frac13.
\]
Together with the remainder, this gives the inner Hessian bound.

For the radial polynomial $B$ of
Lemma~\ref{lem:fctuniformbasin}, put $\mathbf s=(u,v,x,y,z)^{\mathsf T}$.
Since $\nabla_{\mathbf s}P(0)=0$, Euler's identity gives
$\left.\partial_{\epsilon}P((1+\epsilon)S)\right|_{\epsilon=0}=\mathbf s^{\mathsf T}B(S)\mathbf s$.
On $\|S\|_F\le43/1000$,
\[
 \eta^{\mathsf T}B(S)\eta\ge|\eta|^{\mathsf T}C_R|\eta|,
\]
where
\[
 C_R=\frac1{10^6}
 \begin{pmatrix}
 4057854&-142444&-331103&-331103&-648836\\
 -142444&1347120&-333824&-333824&-5546\\
 -331103&-333824&407641&-62250&-62365\\
 -331103&-333824&-62250&407641&-62365\\
 -648836&-5546&-62365&-62365&543475
 \end{pmatrix}.
\]
These are the rational monomial bounds of Lemma~\ref{lem:fctuniformbasin},
with $77$ raw positions and $17711$ nonnegative triples after symmetry folding.

The matrix $C_R-(651/100000)\operatorname{diag}(6,2,2,2,2)$ is positive definite. Its exact
$LDL^{\mathsf T}$ pivots are respectively greater than
\[
 4,\qquad\frac{13}{10},\qquad\frac{27}{100},\qquad
 \frac3{20},\qquad\frac{11}{100}.
\]
For $0<\|S\|_F\le43/1000$, $\nabla_{\mathbf s}\Psi(0)=\nabla_{\mathbf s}P(0)=0$ and the remainder give
\[
 \left.\partial_{\epsilon}\Psi((1+\epsilon)S)\right|_{\epsilon=0}\ge|\mathbf s|^{\mathsf T}C_R|\mathbf s|
       -\varepsilon(43/1000)\|S\|_F^2
 >\frac1{1000}\|S\|_F^2.
\]
Hence $\nabla_{\mathbf s}\Psi(S)\ne0$ for $S\ne0$; integration along $tS$ gives
$\Psi(S)-\Psi(0)\ge\|S\|_F^2/2000$.
This completes the proof of Lemma~\ref{lem:fccuniformbasin}.
\end{proof}

\begin{lemma}[An explicit FCC reference chart]
\label{lem:fccKZchart}
Let \(a_F,\delta_F,h_F>0\) be given by
\[
 a_F^6=\frac{441}{200},\qquad
 \delta_F^6=\frac{11907}{12800},\qquad
 h_F^6=\frac{49}{75},
\]
and define the two reference matrices
\[
 U_{F,\sigma}=\begin{pmatrix}
 a_F&a_F/2&\sigma a_F/2\\
 0&\delta_F&\delta_F/3\\0&0&h_F
 \end{pmatrix},\qquad\sigma\in\{-1,1\}.
\]
Then $U_{F,\sigma}\Z^3\cong(21/20)^{1/3}F$ and
$\det U_{F,\sigma}=21/20$. Exact rational enclosures are
\[
\begin{aligned}
 \frac{1140866336564584}{10^{15}}&<a_F<\frac{1140866336564585}{10^{15}},\\
 \frac{988019229787417}{10^{15}}&<\delta_F<\frac{988019229787418}{10^{15}},\\
 \frac{931513463100523}{10^{15}}&<h_F<\frac{931513463100524}{10^{15}}.
\end{aligned}
\]
For a Korkine--Zolotarev matrix \(U(p)\) of determinant \(21/20\), the relative
matrix \(M_\sigma=U(p)U_{F,\sigma}^{-1}\) is upper triangular, with
\[
\begin{gathered}
 (M_\sigma)_{11}=a/a_F,\qquad
 (M_\sigma)_{22}=\delta/\delta_F,\qquad
 (M_\sigma)_{33}=h/h_F,\\
 (M_\sigma)_{12}=a(x-1/2)/\delta_F,\qquad
 (M_\sigma)_{23}=\delta(z-1/3)/h_F,\\
 (M_\sigma)_{13}=a\{y-\sigma/2-(x-1/2)/3\}/h_F.
\end{gathered}
\]
Then
\[
 S_\sigma=\frac{3M_\sigma^{\mathsf T}M_\sigma}
 {\operatorname{tr}(M_\sigma^{\mathsf T}M_\sigma)}-I
\]
gives the lattice in \eqref{eq:common-metric-chart}, up to an isometry. Thus
\[
 \begin{gathered}
 \|S_\sigma\|_F\le43/1000,\qquad
 \nabla_p\bigl[\tfrac12E(U(p)\Z^3,6)-E(U(p)\Z^3,3)\bigr]=0\\
 \Longrightarrow\quad U(p)\Z^3\cong(21/20)^{1/3}F.
 \end{gathered}
\]
\end{lemma}
\begin{proof}
Taking sixth powers verifies the enclosures, and
\[
 \delta_F^2=\frac34a_F^2,\qquad h_F^2=\frac23a_F^2,
 \qquad a_F\delta_Fh_F=\frac{21}{20}.
\]
The normalized Gram matrix of $U_{F,+}$ is the FCC form
$(1_{i=j}+\tfrac12 1_{i\ne j})_{i,j=1}^3$.
The column operation $v_3\leftarrow v_3-v_1$ gives $U_{F,-}$;
inversion gives the relative entries.
Since \(\det M_\sigma=1\), right polar decomposition
\(M_\sigma=OP\) yields
\[
 \det(I+S_\sigma)^{-1/6}(I+S_\sigma)^{1/2}=P.
\]
Orthogonal conjugation preserves the Frobenius norm and stationary basin.
This completes the proof of Lemma~\ref{lem:fccKZchart}.
\end{proof}

\begin{lemma}[A full stationary basin for the simple-hexagonal lattice]
\label{lem:sh110uniformbasin}
Put
\[
 Y_0=\frac{110131033}{10^8},\qquad
 t_0=\left(\frac{4(11/10)^2}{3Y_0}\right)^{1/6},
\]
and let
\[
 L_0=t_0\operatorname{span}_{\mathbb Z}
 \left\{(\sqrt{Y_0},0,0),(0,1,0),
                          (0,1/2,\sqrt3/2)\right\}.
\]
In the chart \eqref{eq:common-metric-chart}, put
\[
 \begin{gathered}
 L(S)=\det(I+S)^{-1/6}(I+S)^{1/2}L_0,\qquad |L(S)|=11/10,\\
 \Psi(S)=\tfrac12E(L(S),6)-E(L(S),3),\qquad
 \nabla_{\mathbf s}\Psi=(\partial_u\Psi,\partial_v\Psi,
              \partial_x\Psi,\partial_y\Psi,\partial_z\Psi)^{\mathsf T}.
 \end{gathered}
\]
Then
\[
 \left.\partial_{\epsilon}^2\Psi(S+\epsilon T)\right|_{\epsilon=0}>\frac1{1000}\|T\|_F^2
 \quad(\|S\|_F\le3/100,\ T\ne0,\ T=T^{\mathsf T},\ \tr T=0),
\]
\[
 \mathbf s^{\mathsf T}\nabla_{\mathbf s}\Psi(S)
 =\left.\partial_{\epsilon}\Psi((1+\epsilon)S)\right|_{\epsilon=0}>\frac1{1000}\|S\|_F^2
 \qquad(3/100\le\|S\|_F\le7/100).
\]
The root $Y_*$ is determined by
\[
 \begin{gathered}
 \left.\frac{\dd}{\dd Y}\{E(H_Y,6)-2(11/10)^2E(H_Y,3)\}\right|_{Y=Y_*}=0,\\
 1.1013103262<Y_*<1.1013103264.
 \end{gathered}
\]
Put
\[
 S_* =\frac{Y_*-Y_0}{Y_*+2Y_0}\operatorname{diag}(2,-1,-1),
 \qquad\|S_*\|_F<3\cdot10^{-9}.
\]
Then
\[
 \begin{gathered}
 \nabla_{\mathbf s}\Psi(S)=0\iff S=S_*
       \qquad(\|S\|_F\le7/100),\\
 \operatorname*{argmin}_{\|S\|_F\le7/100}\Psi(S)=\{S_*\}.
 \end{gathered}
\]
\end{lemma}

\begin{proof}
Lemma~\ref{lem:sh110rootseed} and hexagonal symmetry give
\[
 \begin{gathered}
 (\partial_v,\partial_x,\partial_y,\partial_z)
 \Psi(\operatorname{diag}(2u,-u,-u))=0
       \quad(\sqrt6|u|\le7/100),\\
 \partial_u\Psi(S_*)=0\quad\Longrightarrow\quad
 \nabla_{\mathbf s}\Psi(S_*)=0.
 \end{gathered}
\]

For every $w\in L_0\setminus\{0\}$,
\[
 w=t_0(\sqrt{Y_0}k,p/2,\sqrt3 j/2),\qquad
 p,j,k\in\Z,\quad(p,j,k)\ne(0,0,0),\quad p+j\in2\Z.
\]
The exact finite vector set $\mathcal F\subset L_0$ is specified by
\[
 |p|\le69,\quad |j|\le40,\quad |k|\le33,\quad p+j\in2\Z,
 \quad0<p^2+3j^2+4Y_0k^2\le4800.
\]
Clear the positive denominator of $Y_0$ in the cutoff; the set is
invariant under planar hexagonal symmetries and coordinate sign changes. Put
\[
 q_w=\frac{t_0^2}{4}(p^2+3j^2+4Y_0k^2),
\]
\[
 \ell_w(S)=
 \frac{u(8Y_0k^2-p^2-3j^2)+v(p^2-3j^2)
 +4x\sqrt{Y_0}kp+4y\sqrt{3Y_0}kj+2z\sqrt3 pj}
 {p^2+3j^2+4Y_0k^2}.
\]
Lemma~\ref{lem:fctuniformbasin} gives $P$, $\deg P\le12$.
Exact determinant multiplication and sign summation give rational
coefficients; the exponents of $x,y,z$ are all even or all odd.
Retain the $138$ raw positions through degree six, including odd
powers of $v$, and all $25193$ nonnegative triples.

Use \eqref{eq:fixedsystole-tail} with
\[
 \lambda_1(L_0)>1,\qquad
 \min_{w\in L_0\setminus(\mathcal F\cup\{0\})}|w|>32,\qquad a=1,\ R=32.
\]
Finite rational summation gives
\[
 \sum_{w\in\mathcal F}q_w^{-3}+T_3(1,32)<\frac{3332}{500},\qquad
 \sum_{w\in\mathcal F}q_w^{-6}+T_6(1,32)<\frac{339}{100}.
\]
Lemma~\ref{lem:fctuniformbasin}, with these moments, gives the complete
binomial and spatial remainder
\[
 \left|\left.\partial_{\epsilon}^2\{\Psi(S+\epsilon T)-P(S+\epsilon T)\}\right|_{\epsilon=0}\right|
 \le\varepsilon(R)\|T\|_F^2,
\]
\[
 \varepsilon(3/100)<\frac7{2000},\qquad
 \varepsilon(7/100)<\frac{393}{10000}.
\]

The weighted monomial maxima give
\[
 \left.\partial_{\epsilon}^2P(S+\epsilon T)\right|_{\epsilon=0}\ge|\eta|^{\mathsf T}C_H|\eta|
 \qquad(\|S\|_F\le3/100),
\]
where $\eta=(T_{11}/2,(T_{22}-T_{33})/2,T_{12},T_{13},T_{23})^{\mathsf T}$ and
\[
 C_H=\frac1{10^6}
 \begin{pmatrix}
 13105271&-2625768&-102166&-102139&-2618185\\
 -2625768&5071238&-37497&-34804&-1356023\\
 -102166&-37497&237405&-21289&-34557\\
 -102139&-34804&-21289&237764&-37032\\
 -2618185&-1356023&-34557&-37032&5074732
 \end{pmatrix}.
\]
The exact $LDL^{\mathsf T}$ pivots of
$C_H-(9/2000)\operatorname{diag}(6,2,2,2,2)$ are greater than
\[
 13,\qquad4,\qquad\frac15,\qquad\frac15,\qquad3.
\]
Thus $\nabla_{\mathbf s}^2\Psi(S)\succ0$ for $\|S\|_F\le3/100$.

For the radial polynomial $B$ of Lemma~\ref{lem:fctuniformbasin},
on $\|S\|_F\le7/100$,
\[
 \eta^{\mathsf T}B(S)\eta\ge|\eta|^{\mathsf T}C_R|\eta|,
\]
where
\[
 C_R=\frac1{10^6}
 \begin{pmatrix}
 12644187&-3333878&-142310&-142138&-3309727\\
 -3333878&4476221&-60016&-54472&-1713026\\
 -142310&-60016&185878&-27728&-53659\\
 -142138&-54472&-27728&187091&-58498\\
 -3309727&-1713026&-53659&-58498&4487900
 \end{pmatrix}.
\]
The matrix $C_R-(4031/100000)\operatorname{diag}(6,2,2,2,2)$ is positive definite: its exact
$LDL^{\mathsf T}$ pivots exceed
\[
 12,\qquad3,\qquad\frac9{100},\qquad\frac8{100},\qquad\frac7{10}.
\]
Since $|\partial_u\Psi(0)|<10^{-7}$ and
$(\partial_v,\partial_x,\partial_y,\partial_z)\Psi(0)=0$,
\[
 |\left.\partial_{\epsilon}\Psi(\epsilon S)\right|_{\epsilon=0}|<\frac1{100000}\|S\|_F^2
 \qquad(\|S\|_F\ge3/100).
\]
Euler's identity and the integrated complete Hessian remainder give
\[
 \left.\partial_{\epsilon}\Psi((1+\epsilon)S)\right|_{\epsilon=0}>\frac1{1000}\|S\|_F^2
\]
on $3/100\le\|S\|_F\le7/100$. Hence
\[
 \begin{gathered}
 3/100\le\|S\|_F\le7/100\quad\Longrightarrow\quad
       \nabla_{\mathbf s}\Psi(S)\ne0,\\
 \{S:\|S\|_F\le3/100,\ \nabla_{\mathbf s}\Psi(S)=0\}=\{S_*\},
 \end{gathered}
\]
by inner strict convexity and $\nabla_{\mathbf s}\Psi(S_*)=0$.
This completes the proof of Lemma~\ref{lem:sh110uniformbasin}.
\end{proof}

\begin{lemma}[The simple-hexagonal reference chart]
\label{lem:sh110chart}
Put \(Y_0=110131033/10^8\). Define positive reference lengths by
\[
 a_0^6=\frac{484000000}{330393099},\qquad
 \delta_0^6=\frac{27}{64}a_0^6,\qquad h_0^6=Y_0^3a_0^6.
\]
Then
\[
 U_0=\begin{pmatrix}a_0&a_0/2&0\\0&\delta_0&0\\0&0&h_0\end{pmatrix}
\]
generates \((11/10)^{1/3}H_{Y_0}\), with triangular plane first
and vertical axis third. Exact rational bounds are
\[
\begin{aligned}
 \frac{1065701881072846}{10^{15}}&<a_0<\frac{1065701881072847}{10^{15}},\\
 \frac{922924901869947}{10^{15}}&<\delta_0<\frac{922924901869948}{10^{15}},\\
 \frac{1118383081907223}{10^{15}}&<h_0<\frac{1118383081907224}{10^{15}}.
\end{aligned}
\]
For \(M=U(p)U_0^{-1}\),
\[
 M=\begin{pmatrix}
 a/a_0&a(x-1/2)/\delta_0&ay/h_0\\
 0&\delta/\delta_0&\delta z/h_0\\
 0&0&h/h_0
 \end{pmatrix}.
\]
Consequently
\[
 S=\frac{3M^{\mathsf T}M}{\operatorname{tr}(M^{\mathsf T}M)}-I.
\]
Use \eqref{eq:common-metric-chart} with $L_0=U_0\Z^3$, and let
$Y_*$ be the root in Lemma~\ref{lem:sh110rootseed}. For $\|S\|_F\le7/100$,
\[
 \nabla_{\mathbf s}\bigl[\tfrac12E(L(S),6)-E(L(S),3)\bigr]=0
 \iff S=\frac{Y_*-Y_0}{Y_*+2Y_0}\operatorname{diag}(-1,-1,2),
\]
and
\[
 \nabla_{\mathbf s}^2\bigl[\tfrac12E(L(S),6)-E(L(S),3)\bigr]\succ0
 \qquad(\|S\|_F\le3/100).
\]
\end{lemma}
\begin{proof}
\[
 \delta_0^2=\frac34a_0^2,\qquad h_0^2=Y_0a_0^2,\qquad
 a_0\delta_0h_0=\frac{11}{10},\qquad
 U_0\Z^3=(11/10)^{1/3}H_{Y_0}.
\]
Inversion gives $M$; the polar argument of
Lemma~\ref{lem:fccKZchart} gives $S$.
Lemma~\ref{lem:sh110uniformbasin} applies after the cyclic permutation
of the coordinate axes.
This completes the proof of Lemma~\ref{lem:sh110chart}.
\end{proof}

\begin{theorem}[FCC optimality at \(21/20\)]
\label{thm:fcccriticalglobal}
At $|L|=21/20$,
\[
 \mathcal E(L)=m(21/20)\iff L\cong(21/20)^{1/3}F,
\]
and
\[
 -5.30745025029370902112<m(21/20)<-5.30745025029370902111.
\]
\end{theorem}
\begin{proof}
The minimum is attained. By the FCC reference bound,
Proposition~\ref{prop:fccKZdomain} and
Lemmas~\ref{lem:fccKZcover} and~\ref{lem:fccuniformbasin},
\[
\begin{aligned}
 \mathcal E(L)=m(21/20)<T_F
 &\ \Longrightarrow\ L\cong U(p)\mathbb Z^3,\quad p\in\mathcal D_F\\
 &\ \Longrightarrow\ L\cong(21/20)^{1/3}F.
\end{aligned}
\]
The FCC reference evaluation gives the stated energy interval.
This completes the proof of Theorem~\ref{thm:fcccriticalglobal}.
\end{proof}

\begin{theorem}[Simple-hexagonal optimality at \(11/10\)]
\label{thm:sh110global}
At $|L|=11/10$,
\[
 \mathcal E(L)=m(11/10)\iff L\cong(11/10)^{1/3}H_{Y_*},
\]
where
\[
 1.1013103262<Y_*<1.1013103264.
\]
\end{theorem}
\begin{proof}
By Proposition~\ref{prop:sh110domain} and
Lemmas~\ref{lem:sh110rootseed}, \ref{lem:sh110cover}
and~\ref{lem:sh110uniformbasin},
\[
 \mathcal E(L)=m(11/10)<-\frac{24842327}{5000000}
 \quad\Longrightarrow\quad
 L\cong(11/10)^{1/3}H_{Y_*}.
\]
The root bracket gives the planar area bound
\[
 \frac{983}{1000}<
 \left(\frac{\sqrt3(11/10)^2}{2Y_*}\right)^{1/3}<\frac{123}{125}.
\]
This completes the proof of Theorem~\ref{thm:sh110global}.
\end{proof}

\section{The transition interval}\label{sec:transition}

For $21/20\le V\le11/10$, the profile below reduces minimization
to five-coordinate stationary points in a compact domain. Write
\[
\begin{gathered}
 \Phi_V(\Lambda)=E(\Lambda,6)-2V^2E(\Lambda,3),\qquad |\Lambda|=1,\\
 \mathcal E(V^{1/3}\Lambda)
 =\frac{E(\Lambda,6)}{2V^4}-\frac{E(\Lambda,3)}{V^2}
 =\frac{\Phi_V(\Lambda)}{2V^4},\\
 \lambda_1(\Lambda)=\min_{v\in\Lambda\setminus\{0\}}|v|.
\end{gathered}
\]

\begin{proposition}[One profile comparison on a covolume interval]
\label{prop:density-gap-reduction}\label{prop:gap-concave-envelope}
\label{lem:profile-supporting-line}
Let $f:[V_0,V_1]\to\R$ be continuous and concave as a function of
$V^2$, where $0<V_0<V_1$. Suppose
\begin{equation}\label{eq:gap-endpoint-comparisons}
 E(\Lambda,6)-2V_i^2E(\Lambda,3)\ge f(V_i)
 \quad (|\Lambda|=1,\ i=0,1).
\end{equation}
Set
\[
 b_0=-\frac{f'_+(V_0)}{4V_0},\quad
 b_1=-\frac{f'_-(V_1)}{4V_1},\quad
 A_1=f(V_1)+2V_1^2b_1,
\]
with $0<b_0\le b_1<\infty$ and $A_1>0$.
\[
 J(\Lambda)=E(\Lambda,6)-g(E(\Lambda,3)),
\]
\[
 g(t)=\max_{V_0\le V\le V_1}\{f(V)+2V^2t\}.
\]
Then
\begin{equation}\label{eq:gap-profile-minimum}
 J(\Lambda)=\min_{V_0\le V\le V_1}
       \{\Phi_V(\Lambda)-f(V)\},
\end{equation}
and, for $|\Lambda|=1$,
\[
 J(\Lambda)\ge0
 \quad\Longleftrightarrow\quad
 \Phi_V(\Lambda)\ge f(V)\quad(V_0\le V\le V_1).
\]
Moreover,
\[
 J(\Lambda)<0\ \Longrightarrow\
 b_0<E(\Lambda,3)<b_1,\quad E(\Lambda,6)<A_1,\quad
 \lambda_1(\Lambda)>(2/A_1)^{1/12};
\]
and $\inf J<0\Longrightarrow\operatorname*{argmin}J\ne\varnothing$.
For every
\[
 V\in\operatorname*{argmax}_{V_0\le v\le V_1}
       \{f(v)+2v^2E(\Lambda_*,3)\},
\]
\[
\begin{aligned}
 \Lambda_*\in\operatorname*{argmin}_{\rm loc}J
 &\Longrightarrow\Lambda_*\in\operatorname*{argmin}_{\rm loc}\Phi_V,\\
 \Lambda_*\in\operatorname*{argmin}J
 &\Longrightarrow\Lambda_*\in\operatorname*{argmin}\Phi_V.
\end{aligned}
\]
For $V_0<V<V_1$,
\[
 -\frac{f'_-(V)}{4V}\le E(\Lambda_*,3)
                     \le-\frac{f'_+(V)}{4V}.
\]
No differentiability at a contact is required. In particular,
\[
 J\ge0,\quad J(\Lambda_*)=0
 \quad\Longrightarrow\quad\Lambda_*\in\operatorname*{argmin}J.
\]
Endpoint equality uses \eqref{eq:gap-endpoint-comparisons}; shape
variations remain unrestricted at reduction boundaries.
\end{proposition}
\begin{proof}
Concavity in $V^2$ gives
\[
 g(t)=f(V_0)+2V_0^2t\quad(t\le b_0),\qquad
 g(t)=f(V_1)+2V_1^2t\quad(t\ge b_1).
\]
Together with \eqref{eq:gap-endpoint-comparisons} and monotonicity of $g$,
\[
 J(\Lambda)<0\ \Longrightarrow\
 \begin{cases}
 b_0<E(\Lambda,3)<b_1,\\
 2\lambda_1(\Lambda)^{-12}\le E(\Lambda,6)<g(b_1)=A_1.
 \end{cases}
\]
Mahler's criterion \cite{Mahler} and uniform convergence imply
$\inf J<0\Longrightarrow\operatorname*{argmin}J\ne\varnothing$.
For any maximizing $V$ at $E(\Lambda_*,3)$,
\begin{align*}
 \Phi_V(\Lambda)-\Phi_V(\Lambda_*)
 &=J(\Lambda)-J(\Lambda_*)
   +g(E(\Lambda,3))-2V^2E(\Lambda,3)-f(V)\\
 &\ge J(\Lambda)-J(\Lambda_*).
\end{align*}
The interior one-sided condition is
\[
 -4VE(\Lambda_*,3)-f'_-(V)\le0
       \le-4VE(\Lambda_*,3)-f'_+(V).
\]
This completes the proof of Proposition~\ref{prop:density-gap-reduction}.
\end{proof}

On a smooth branch attaining $f(V)$,
\[
 \left.\partial_Y\Phi_V(H_Y)\right|_{Y=Y(V)}=0,\qquad
 f(V)=\Phi_V(H_{Y(V)}),\qquad
 b(V)=E(H_{Y(V)},3)=-\frac{f'(V)}{4V}.
\]
If $\partial_YE(H_Y,3)\ne0$ and $\partial_Y^2\Phi_V(H_Y)>0$ at $Y(V)$, then
\begin{equation}\label{eq:gap-candidate-moment-derivative}
 b'(V)=
 \left.\frac{4V\{\partial_YE(H_Y,3)\}^2}
 {\partial_Y^2\{E(H_Y,6)-2V^2E(H_Y,3)\}}
 \right|_{Y=Y(V)}>0.
\end{equation}
For $t$ in this branch's moment range,
\[
 V=b^{-1}(t),\qquad g(t)=E(H_{Y(V)},6).
\]
At a smooth profile minimum $\Lambda(q)$ supported by this $V$,
\[
 \nabla_qJ=0,\qquad
 \nabla_q^2J=
 \nabla_q^2\Phi_V(\Lambda(q))
 -\frac{4V}{b'(V)}\nabla_qE(\Lambda(q),3)
       \{\nabla_qE(\Lambda(q),3)\}^{\mathsf T}\succeq0.
\]
At contacts use the one-sided moment inequalities above.

\begin{lemma}[Uniform tests over a density interval]\label{lem:densityFCTcover}
In \eqref{eq:criticalKZprimal}--\eqref{eq:criticalKZdual}, let
\[
 p=(a,\delta,x,y,z),\quad h=(a\delta)^{-1},\quad
 B\subset\{a>0,\delta>0\},\quad c=\operatorname{mid}B,
\]
where $B$ is a closed box. Put $Q(p)=U(p)^{\mathsf T}U(p)$.
For $[V_-,V_+]\subset[21/20,89/50]$, set
\[
 \Phi_V(p)=E(U(p)\Z^3,6)-2V^2E(U(p)\Z^3,3).
\]
Let $V^2=(1-t)V_-^2+tV_+^2$, $0\le t\le1$, and suppose
\[
\begin{gathered}
 \relax[H_e]\supset\{\nabla_p^2\Phi_{V_e}(p):p\in B\},\qquad
 [g_e]\ni\nabla_p\Phi_{V_e}(c)\quad(e=-,+),\\
 [H]=\operatorname{hull}([H_-],[H_+]),\qquad
 [g_c]=\operatorname{hull}([g_-],[g_+]).
\end{gathered}
\]
For every $p\in B$, $P\in\mathbb R^{5\times5}$, $i\in\{1,\ldots,5\}$, and fixed nonzero dyadic $v\in\mathbb R^5$,
\[
\begin{aligned}
 \nabla_p\Phi_V(p)=0
 &\Longrightarrow p\in B\cap\{c-P[g_c]+(I-P[H])(B-c)\},\\
 0\notin\operatorname{hull}_{e\in\{-,+\}}[\partial_i\Phi_{V_e}](B)
 &\Longrightarrow\nabla_p\Phi_V(p)\ne0,\\
 \max_{e\in\{-,+\}}\sup_{H\in[H_e]}v^{\mathsf T}Hv<0
 &\Longrightarrow v^{\mathsf T}\nabla_p^2\Phi_V(p)v<0.
\end{aligned}
\]
The last condition excludes a local minimum in the full shape chart.
For $|\Lambda_0|=1$, write
\[
 \Delta_s=E(U(c)\Z^3,s)-E(\Lambda_0,s)\quad(s=3,6),\qquad
 r_i=\sup_{p\in B}|p_i-c_i|.
\]
Let $Q_e\ge0$ be the $32$-vertex bound \eqref{eq:criticalrectvalue}; thus
\[
 \sup_{\substack{p\in B\\|\xi_i|\le r_i}}
 \xi^{\mathsf T}\nabla_p^2\Phi_{V_e}(p)\xi\le Q_e.
\]
For $\nabla_p\Phi_V(p)=0$,
\begin{align*}
 \Phi_V(p)-\Phi_V(\Lambda_0)
 \ge{}&(1-t)\left(\Delta_6-2V_-^2\Delta_3-\frac12Q_-\right)\\
 &+t\left(\Delta_6-2V_+^2\Delta_3-\frac12Q_+\right).
\end{align*}
In particular, with the same $\Lambda_0$ at both endpoints,
\[
\begin{gathered}
 \nabla_p\Phi_V(p)=0,\qquad
 \min_{e\in\{-,+\}}\left(\Delta_6-2V_e^2\Delta_3-\frac12Q_e\right)>0\\
 \Longrightarrow\quad
 \Phi_V(p)>\Phi_V(\Lambda_0),\qquad
 U(p)\Z^3\notin\operatorname*{argmin}_{|\Lambda|=1}\Phi_V(\Lambda).
\end{gathered}
\]
For the scalar test assume
\[
 M(p)=U(p)U(c)^{-1},\qquad
 \sup_{p\in B}\|\log(M(p)^{\mathsf T}M(p))\|_F\le R,
\]
and put $d=e^{-\sqrt{2/3}R}$, $q_c(i,j,k)=(i,j,k)Q(c)(i,j,k)^{\mathsf T}$. Set
\begin{align*}
 K_e={}&4V_e^2d^{-3}E(U(c)\mathbb Z^3,3)
             -4d^6E(U(c)\mathbb Z^3,6)\\
 &+\sum_{\substack{(i,j,k)\in\mathbb Z^3\\-4\le i,j,k\le4\\(i,j,k)\ne(0,0,0)}}
 \max\{0,28[dq_c(i,j,k)]^{-6}
                     -16V_e^2[dq_c(i,j,k)]^{-3}\}.
\end{align*}
In the value bound, replace each $Q_e$ by $R^2\max(0,K_e)$, provided
\[
 (25d\|Q(c)^{-1}\|_\infty^{-1})^3\ge\frac7{4V_-^2}.
\]
All intervals include the complete whole-box tails; neither
$\nabla_p\Phi_V(c)=0$ nor endpoint stationarity is assumed.
\end{lemma}
\begin{proof}
Use the complete endpoint bounds of Lemma~\ref{lem:criticalKZcover}.
\paragraph{Normalized kernels.}
The unit-covolume constants are $\pi^3/6$, $\pi^6/2160$;
the primal kernel derivatives at either endpoint are
\[
 (-\pi)^r\left\{\frac{\pi^6}{120}\mathcal G_{6+r}(\pi q)
          -V^2\pi^3\mathcal G_{3+r}(\pi q)\right\},
 \qquad r=0,1,2,
\]
For the reciprocal kernel use $(-9/2+r,-3/2+r)$ instead of
$(6+r,3+r)$, with volume factor one. The primal critical argument solves
\[
 \sum_{j=0}^{6+r}\frac{(\pi q)^j}{j!}
 -2V^2 b_rq^3\sum_{j=0}^{3+r}\frac{(\pi q)^j}{j!}=0,
 \quad(b_0,b_1,b_2)=\left(\frac12,\frac27,\frac5{28}\right),
\]
The decreasing integral ratio gives uniqueness; use at most $60$
interval bisections in $[1/4,8]$, stopping at an inconclusive sign.
For reciprocal derivatives the positive-magnitude ratio is
$<\pi^3/60<2V^2$. Evaluate both terms at the same $q$ before subtraction.

Use \eqref{eq:criticalHtail} on the whole box, with these coefficients and
\[
 Q(p)\succeq e^{-\sqrt{2/3}R}\|Q(c)^{-1}\|_\infty^{-1}I,
 \qquad
 Q(p)^{-1}\succeq e^{-\sqrt{2/3}R}\|Q(c)\|_\infty^{-1}I.
\]
Round norm denominators upwards; shortest-vector bounds do not replace these spectral floors.

\paragraph{Stationary values and exclusions.}
For $p\in B$,
\[
 \nabla_p\Phi_V(p)=0
 \Longrightarrow
 \Phi_V(p)-\Phi_V(c)
 =-\int_0^1\tau(p-c)^{\mathsf T}
 \nabla_p^2\Phi_V(c+\tau(p-c))(p-c)\,\dd\tau.
\]
Since
\[
 \nabla_p^r\Phi_V=(1-t)\nabla_p^r\Phi_{V_-}
                         +t\nabla_p^r\Phi_{V_+},\qquad r=0,1,2,
\]
\eqref{eq:criticalrectvalue} gives the paired value bound.
The scalar summands are convex in $V^2$; their endpoint values are
upper bounds, not an affine identity for $K_e$. For a trace-free
Frobenius-unit $H$,
\[
\begin{aligned}
 &\left.\partial_{\epsilon}^2
 \bigl(|e^{\epsilon H/2}v|^{-12}-2V^2|e^{\epsilon H/2}v|^{-6}\bigr)\right|_{\epsilon=0}\\
 &\qquad\le4V^2|v|^{-6}-4|v|^{-12}
       +\max(0,28|v|^{-12}-16V^2|v|^{-6})
\end{aligned}
\]
in the metric convention of Lemma~\ref{lem:criticalKZcover}.
At $V=1$, $\Phi_1=E(\cdot,6)-2E(\cdot,3)=2\mathcal E$.
The guard removes every omitted correction; convex interpolation
therefore gives the scalar bound. The endpoint hulls and the
mean-value identity give \eqref{eq:criticalNewton} and the two
exclusions. Select $v$ at the midpoint and round
to $2^{-20}\Z$ before both complete quadratic-form checks. Use the same
witness $\Lambda_0$, complete whole-box tails and outward arithmetic.
This completes the proof of Lemma~\ref{lem:densityFCTcover}.
\end{proof}

\begin{lemma}[A compact moment domain for the transition interval]
\label{lem:transition-profile-moment-domain}
Let $f$ be continuous on $[21/20,11/10]$, concave in $V^2$, and suppose
\[
 -\frac5{21}f'_+(21/20)=E(F,3),\qquad
 f(21/20)+\frac{441}{200}E(F,3)=E(F,6),
\]
\[
 -\frac5{22}f'_-(11/10)=E(H_{Y_1},3),\qquad
 f(11/10)+\frac{121}{50}E(H_{Y_1},3)=E(H_{Y_1},6),
\]
where $H_{Y_1}$ is the unit-covolume SH minimizer shape at $V=11/10$.
Assume the endpoint comparisons
\[
 E(\Lambda,6)-2V^2E(\Lambda,3)\ge f(V)
 \quad\bigl(|\Lambda|=1,\ V\in\{21/20,11/10\}\bigr).
\]
For $|\Lambda|=1$,
\[
\begin{gathered}
 \min_{21/20\le V\le11/10}
 \{E(\Lambda,6)-2V^2E(\Lambda,3)-f(V)\}<0\\
 \Longrightarrow\quad
 E(F,3)<E(\Lambda,3)<\frac{81}{10},\quad E(\Lambda,6)<5,\quad
 \lambda_1(\Lambda)>\left(\frac25\right)^{1/12}>\frac{463}{500}.
\end{gathered}
\]
If the comparison holds for all unit-covolume lattices and all
$21/20\le V\le11/10$, then
\[
\begin{gathered}
 E(\Lambda,6)-2V^2E(\Lambda,3)=f(V),\quad21/20<V<11/10\\
 \Longrightarrow\quad
 E(F,3)\le E(\Lambda,3)<81/10,\quad E(\Lambda,6)<5,\quad
 \lambda_1(\Lambda)>(2/5)^{1/12}>463/500.
\end{gathered}
\]
Every lattice in either case has a Korkine--Zolotarev representative
with $a\delta h=1$ and
\[
 \frac{463}{500}<a\le2^{1/6}<\frac{1122463}{10^6},\qquad
 \frac{801}{1000}<\delta<\frac{559}{500},\qquad
 \frac78<h<\frac{27}{20},
\]
\[
 0\le x,z\le\frac12,\qquad-\frac12\le y\le\frac12.
\]
\end{lemma}
\begin{proof}
Let $L_0$ be the covolume $11/10$ reference in
Lemma~\ref{lem:sh110uniformbasin}. In its normalized-linear chart,
$\|S_*\|_F\le3/10^9$ at the stationary lattice, and
\[
 E(L_0,3)<\frac{3332}{500},\qquad E(L_0,6)<\frac{339}{100}.
\]
For $\|S\|_F\le3/10^9<1$ and $v\in L_0$,
\[
 \det(I+S)^{-1/3}v^{\mathsf T}(I+S)v
       \ge\frac{1-3/10^9}{1+3/10^9}|v|^2.
\]
Hence
\[
 E\!\left(\det(I+S)^{-1/6}(I+S)^{1/2}L_0,s\right)
 \le\left(\frac{1+3/10^9}{1-3/10^9}\right)^s E(L_0,s),
 \qquad s=3,6.
\]
Since
\[
 E(H_{Y_1},s)=(11/10)^{2s/3}
 E\!\left(\det(I+S_*)^{-1/6}(I+S_*)^{1/2}L_0,s\right),
\]
\[
 b_0=E(F,3),\qquad b_1=E(H_{Y_1},3)<81/10,\qquad
 A_1=E(H_{Y_1},6)<5.
\]
Apply Proposition~\ref{prop:gap-concave-envelope}. After the global
comparison, an interior equality is a zero minimum and satisfies
\[
 b_0\le E(\Lambda,3)\le b_1<81/10,\qquad E(\Lambda,6)\le A_1<5.
\]

For $v_0\in\Lambda\setminus\{0\}$ with $|v_0|=\lambda_1(\Lambda)$,
\[
 |v_0|^{-12}+|-v_0|^{-12}=2\lambda_1(\Lambda)^{-12}\le E(\Lambda,6)<5
 \quad\Longrightarrow\quad
 \lambda_1(\Lambda)>(2/5)^{1/12}>463/500.
\]
The Korkine--Zolotarev conditions give
\[
 \delta\ge\frac{\sqrt3}{2}a>\frac{801}{1000},\qquad
 \delta^2\le\frac{2}{\sqrt3a}<\left(\frac{559}{500}\right)^2,
\]
\[
 h^2\ge\frac{\sqrt3}{2a}>\left(\frac78\right)^2,
 \qquad h\le\frac{2}{\sqrt3a^2}<\frac{27}{20}.
\]
This completes the proof of Lemma~\ref{lem:transition-profile-moment-domain}.
\end{proof}

\begin{lemma}[A uniform FCC basin]
\label{lem:transitionfccbasin}
Let $L_F=(441/1600)^{1/6}D_3$, so $|L_F|=21/20$.
For trace-free symmetric $S$ in \eqref{eq:common-metric-chart}, put
\[
 L_F(S)=\det(I+S)^{-1/6}(I+S)^{1/2}L_F,\qquad
 \Lambda_F(S)=(20/21)^{1/3}L_F(S),
\]
\[
 F_V(S)=\tfrac12E(L_F(S),6)-(20V/21)^2 E(L_F(S),3).
\]
For every $V>0$ satisfying
\[
 \frac{220499}{220500}\le(20V/21)^2\le\frac{458}{441},
\]
and every $0<\|S\|_F\le1/40$,
\[
 \mathbf s^{\mathsf T}\nabla_{\mathbf s}F_V(S)
 >\frac1{1000}\|S\|_F^2,
\]
\[
 \nabla_{\mathbf s}F_V(0)=0,\qquad
 F_V(S)\ge F_V(0)+\frac1{2000}\|S\|_F^2
 \quad(\|S\|_F\le1/40).
\]
On $\|S\|_F\le1/40$,
\[
 \nabla_{\mathbf s}F_V(S)=0\Longleftrightarrow S=0,\qquad
 \operatorname*{argmin}_{\|S\|_F\le1/40}F_V=\{0\}.
\]
For the unit-covolume shape $\Lambda_F(S)$,
\[
\begin{aligned}
 \Phi_V(\Lambda_F(S))
 &=E(\Lambda_F(S),6)-2V^2E(\Lambda_F(S),3)\\
 &=2(21/20)^4F_V(S),\\
 \mathcal E((20V/21)^{1/3}L_F(S))
 &=\frac{E(\Lambda_F(S),6)}{2V^4}
  -\frac{E(\Lambda_F(S),3)}{V^2}\\
 &=(21/(20V))^4F_V(S).
\end{aligned}
\]
\[
 (20V/21)^2\in[220499/220500,458/441]
 \quad\Longleftrightarrow\quad
 2V^2\in[220499/100000,229/100].
\]
Both factors are positive and independent of $S$.
\end{lemma}

\begin{proof}
Use the complete spherical sums, binomial polynomials, determinant,
and moment tails in Lemma~\ref{lem:fccuniformbasin} at covolume $21/20$. Cubic symmetry gives
$\nabla_{\mathbf s}F_V(0)=0$ for every $V$.
With the finite vector set $\mathcal F$ of Lemma~\ref{lem:fccuniformbasin}
and the same $q_w,\ell_w$, put
\[
\begin{aligned}
 P_V(S)={}&\tfrac12\det(I+S)^2
       \sum_{w\in\mathcal F}q_w^{-6}p_6(\ell_w(S))\\
 &-(20V/21)^2\det(I+S)\sum_{w\in\mathcal F}q_w^{-3}p_3(\ell_w(S)).
\end{aligned}
\]
Affinity in $V^2$ gives uniformity on
$220499/220500\le(20V/21)^2\le458/441$. If
$P_V=\sum_{d=0}^{12}P_{V,d}$, define
\[
 (B_V(S))_{ij}=\sum_{d=2}^{12}\frac1{d-1}
 \frac{\partial^2P_{V,d}(S)}{\partial s_i\partial s_j},\qquad
 1\le i,j\le5,\quad\mathbf s=(u,v,x,y,z)^{\mathsf T}.
\]
For all $\|S\|_F\le1/40$, the rational coefficient sums and the explicit
weighted monomial maxima give
\[
 \eta^{\mathsf T}B_V(S)\eta
 \ge |\eta|^{\mathsf T}C|\eta|,
\]
where
\[
 C=\frac1{10^6}\begin{pmatrix}
 3859306&-75427&-183283&-183283&-362275\\
 -75427&1284695&-184271&-184271&-1747\\
 -183283&-184271&188706&-36951&-36973\\
 -183283&-184271&-36951&188706&-36973\\
 -362275&-1747&-36973&-36973&264882
\end{pmatrix}.
\]
The total Hessian remainder is at most
$(458/441)\varepsilon(1/40)<33/10000$, since $458/441\ge1$
also majorizes the repulsive contribution.

The $LDL^{\mathsf T}$ pivots of
$C-(43/10000)\operatorname{diag}(6,2,2,2,2)$ exceed, respectively,
\[
 \frac{38}{10},\quad\frac{12}{10},\quad\frac{14}{100},
 \quad\frac{10}{100},\quad\frac{13}{100}.
\]
Euler's identity and the complete remainder give, for $0<\|S\|_F\le1/40$,
\[
 \mathbf s^{\mathsf T}\nabla_{\mathbf s}F_V(S)
 >\|S\|_F^2/1000,\qquad
 F_V(S)-F_V(0)
 =\int_0^1\mathbf s^{\mathsf T}\nabla_{\mathbf s}F_V(tS)\,\dd t
 \ge\|S\|_F^2/2000.
\]
Homogeneity gives the scaling identities.
This completes the proof of Lemma~\ref{lem:transitionfccbasin}.
\end{proof}

\begin{lemma}[A tetragonal stationary basin across the transition window]
\label{lem:widefctbasin}
In the common chart of Lemma~\ref{lem:fctuniformbasin}, set
\[
 \Lambda(S)=\det(I+S)^{-1/6}(I+S)^{1/2}L_{149/500}.
\]
Set
$\Phi_V(S)=E(\Lambda(S),6)-2V^2 E(\Lambda(S),3)$.
For every $2267/1000\le2V^2\le579/250$, the closed ball
$\|S\|_F\le19/1000$ is strictly convex for $\Phi_V$, with
\[
 \left.\partial_{\epsilon}^2\Phi_V(S+\epsilon T)\right|_{\epsilon=0}>\frac1{400}\|T\|_F^2
 \qquad(T\ne0,\ T=T^{\mathsf T},\ \operatorname{tr}T=0).
\]
There is a real-analytic $u_V$, $|u_V|<1/200$, such that, on this ball,
\[
\begin{gathered}
 S_V=\operatorname{diag}(2u_V,-u_V,-u_V),\qquad
 \nabla_{\mathbf s}\Phi_V(S)=0\Longleftrightarrow S=S_V,\\
 \operatorname*{argmin}_{\|S\|_F\le19/1000}\Phi_V=\{S_V\}.
\end{gathered}
\]
\end{lemma}
\begin{proof}
Set $L_0=(43/40)^{1/3}L_{149/500}$ in the polynomial and complete
remainder construction of Lemma~\ref{lem:fctuniformbasin}. For $i,j,k\in\Z$, let $\mathcal F$ consist of
\[
 \begin{gathered}
 w=\left(\frac{(43/40)^2}{4(149/500)}\right)^{1/6}
       (\sqrt{149/500}\,i,j,k),\quad |i|\le60,\quad |j|,|k|\le32,\\
 i+j+k\in2\Z,\quad0<3i^2+10j^2+10k^2\le10800.
 \end{gathered}
\]
For this set,
\begin{align*}
 q_w&=\left\{\left(\frac{(43/40)^2}{4(149/500)}\right)^{1/6}\right\}^{2}
       \bigl((149/500)i^2+j^2+k^2\bigr),\\
 \ell_w&=\frac{u\{2(149/500)i^2-j^2-k^2\}+v(j^2-k^2)}
                 {(149/500)i^2+j^2+k^2}\\
 &\quad+\frac{2x\sqrt{149/500}\,ij+2y\sqrt{149/500}\,ik+2zjk}
                 {(149/500)i^2+j^2+k^2}.
\end{align*}
\[
 \lambda_1(L_0)=\min_{w\in L_0\setminus\{0\}}|w|>1,\qquad
 \min_{w\in L_0\setminus(\mathcal F\cup\{0\})}|w|>32.
\]
Use \eqref{eq:fixedsystole-tail} with $a=1,R=32$ and
\[
 E^+(L_0,3)<643/100,\qquad E^+(L_0,6)<261/100.
\]
For $L_0(S)=(43/40)^{1/3}\Lambda(S)$,
\[
\begin{gathered}
 F_V(S)=\tfrac12E(L_0(S),6)-(40V/43)^2 E(L_0(S),3),\\
 9068/9245\le(40V/43)^2\le9264/9245,\qquad
 \Phi_V(S)=2(43/40)^4F_V(S).
\end{gathered}
\]
The Hessian polynomial is affine in $V^2$; its complete remainder is
$<3/1000$ on $\|S\|_F\le19/1000$.

At the two endpoints its rational comparison matrices are
\[
 C_- = \frac1{10^6}\begin{pmatrix}
 10970754&-388661&-183919&-183919&-51920\\
 -388661&663661&-369193&-369193&-2037\\
 -183919&-369193&1176945&-28014&-25522\\
 -183919&-369193&-28014&1176945&-25522\\
 -51920&-2037&-25522&-25522&119928
\end{pmatrix},
\]
\[
 C_+ = \frac1{10^6}\begin{pmatrix}
 10232230&-380114&-173000&-173000&-56122\\
 -380114&319921&-353808&-353808&-1906\\
 -173000&-353808&1057728&-29813&-27355\\
 -173000&-353808&-29813&1057728&-27355\\
 -56122&-1906&-27355&-27355&201534
\end{pmatrix}.
\]
The matrices $C_\pm-(4/1000)\operatorname{diag}(6,2,2,2,2)$
are positive definite. Their exact
$LDL^{\mathsf T}$ pivots exceed, respectively,
\[
 (10,3/5,9/10,4/5,1/10),\qquad
 (10,29/100,3/5,1/4,17/100).
\]
The weighted monomial rule of Lemma~\ref{lem:fctuniformbasin}
proves these finite bounds. Affinity in $V^2$ and the common
remainder bound give
$\left.\partial_{\epsilon}^2F_V(S+\epsilon T)\right|_{\epsilon=0}>(1/1000)\|T\|_F^2$ on the whole ball.
Since $2(43/40)^4>5/2$, the displayed normalized bound follows.

The complete attraction derivative and Lemma~\ref{lem:fctrootdata} give
\[
\begin{gathered}
 \left|\partial_u E(L_0(S),3)\big|_{S=0}\right|<\frac{261}{100},\\
 |\partial_uF_{43/40}(0)|<8\cdot10^{-6},\qquad
 |\partial_uF_V(0)|<\frac1{20}.
\end{gathered}
\]
On the ball, $\partial_{uu}F_V>10$, and
\[
 \partial_uF_V(-1/200,0,0,0,0)<0
 <\partial_uF_V(1/200,0,0,0,0).
\]
Tetragonal symmetry gives
\[
 (\partial_vF_V,\partial_xF_V,\partial_yF_V,\partial_zF_V)(u,0,0,0,0)=0.
\]
Thus
\[
 \exists!\,u_V\in(-1/200,1/200):
 \partial_uF_V(u_V,0,0,0,0)=0,
 \qquad\nabla_{\mathbf s}F_V(S_V)=0.
\]
Strict convexity gives uniqueness; $\partial_{uu}F_V(S_V)>0$ gives
analyticity by the implicit-function theorem.
This completes the proof of Lemma~\ref{lem:widefctbasin}.
\end{proof}

\begin{lemma}[A uniform simple-hexagonal basin on a broad density interval]
\label{lem:broadshbasin}
Put
\[
 Y_0=\frac{107818724}{10^8},\qquad
 \Lambda_0=\left(\frac4{3Y_0}\right)^{1/6}
 \operatorname{span}_{\mathbb Z}
 \{(\sqrt{Y_0},0,0),(0,1,0),(0,1/2,\sqrt3/2)\}.
\]
Thus $|\Lambda_0|=1$. In \eqref{eq:common-metric-chart}, set
\[
\begin{gathered}
 \Lambda(S)=\det(I+S)^{-1/6}(I+S)^{1/2}\Lambda_0,\\
 \Phi_V(S)=E(\Lambda(S),6)-2V^2 E(\Lambda(S),3).
\end{gathered}
\]
For every
\[
 \frac{231499}{100000}\le2V^2\le\frac{242001}{100000},
\]
there is a real-analytic $u_V$ with
\[
 S_V=\operatorname{diag}(2u_V,-u_V,-u_V),
 \qquad |u_V|<\frac1{125},\qquad
 Y(V)=Y_0\frac{1+2u_V}{1-u_V}.
\]
For $\|S\|_F\le1/50$,
\[
 \nabla_{\mathbf s}\Phi_V(S)=0\Longleftrightarrow S=S_V,\qquad
 \operatorname*{argmin}_{\|S\|_F\le1/50}\Phi_V=\{S_V\},
\]
and
\[
 \left.\partial_{\epsilon}^2\Phi_V(S+\epsilon T)\right|_{\epsilon=0}>\frac1{400}\|T\|_F^2
 \quad(T\ne0),
\]
\[
 \Phi_V(S)\ge\Phi_V(S_V)
                 +\frac1{800}\|S-S_V\|_F^2.
\]
\end{lemma}

\begin{proof}
Use $L_0=(87/80)^{1/3}\Lambda_0$. At the endpoint covolumes
$V_-=\sqrt{231499/200000}$ and $V_+=\sqrt{242001/200000}$.
For $L_0(S)=(87/80)^{1/3}\Lambda(S)$,
\[
 E(L_0(S),s)=\det(I+S)^{s/3}
 \sum_{w\in L_0\setminus\{0\}}(|w|^2+w^{\mathsf T}Sw)^{-s}.
\]
Hence
\[
\begin{gathered}
 F_V(S)=\tfrac12E(L_0(S),6)-(80V/87)^2 E(L_0(S),3),\\
 \Phi_V(S)=2(87/80)^4F_V(S),\qquad 2(87/80)^4>\tfrac52.
\end{gathered}
\]

With $Y_0=26954681/25000000$, let $\mathcal F\subset L_0$ consist of the vectors
\[
 w=\left(\frac{4(87/80)^2}{3Y_0}\right)^{1/6}
       (\sqrt{Y_0}k,p/2,\sqrt3j/2),\qquad p,j,k\in\Z,\quad p+j\in2\Z,
\]
\[
 \begin{gathered}
 |p|\le69,\quad |j|\le40,\quad |k|\le33,\\
 0<25000000(p^2+3j^2)+107818724k^2\le120000000000.
 \end{gathered}
\]
\[
 \#(\mathcal F\cap[0,\infty)^3)=25471.
\]
The set retains coordinate signs and planar hexagonal symmetry.
For each vector put
\[
 q_w=\frac{\left\{\left(\frac{4(87/80)^2}{3Y_0}\right)^{1/6}\right\}^{2}}4
       (p^2+3j^2+4Y_0k^2),
\]
\[
 \ell_w(S)=\frac{
 u(8Y_0k^2-p^2-3j^2)+v(p^2-3j^2)
 +4x\sqrt{Y_0}kp+4y\sqrt{3Y_0}kj+2z\sqrt3pj}
 {p^2+3j^2+4Y_0k^2}.
\]
With the degree-six binomial polynomials $p_3,p_6$ of
Lemma~\ref{lem:fctuniformbasin}, the polynomial to be checked is
\[
\begin{aligned}
 P_V(S)={}&\frac12\det(I+S)^2
       \sum_{w\in\mathcal F}q_w^{-6}p_6(\ell_w(S))\\
 &-(80V/87)^2\det(I+S)\sum_{w\in\mathcal F}q_w^{-3}p_3(\ell_w(S)).
\end{aligned}
\]
Retain the determinant factors and all $138$ parity-allowed positions,
including odd powers of $v$. Complete sign sums give
\[
 P_V\in\mathbb Q[V^2,u,v,x,y,z],\qquad
 \deg_{\mathbf s}P_V\le12,\qquad\deg_{\mathbf s}\nabla^2_{\mathbf s}P_V\le10.
\]

Since
\[
 \lambda_1(L_0)=\min_{w\in L_0\setminus\{0\}}|w|>1,\qquad
 \min_{w\in L_0\setminus(\mathcal F\cup\{0\})}|w|>32,
\]
\eqref{eq:fixedsystole-tail} gives
\[
 \sum_{w\in L_0\setminus(\mathcal F\cup\{0\})}q_w^{-s}\le T_s(1,32),\qquad s=3,6.
\]
Finite rational summation gives
\[
 \sum_{w\in\mathcal F}q_w^{-3}+T_3(1,32)<\frac{6783}{1000},\qquad
 \sum_{w\in\mathcal F}q_w^{-6}+T_6(1,32)<\frac{3477}{1000}.
\]
Apply the complete binomial and spatial Hessian remainder formulas of
Lemma~\ref{lem:fctuniformbasin}, with these moment bounds and with the
$s=3$ contributions multiplied by $107556/105125$. They give
\[
 \left|\left.\partial_{\epsilon}^2\{F_V(S+\epsilon T)-P_V(S+\epsilon T)\}\right|_{\epsilon=0}\right|
 <\frac{151}{50000}\|T\|_F^2
 \qquad(\|S\|_F\le1/50).
\]

Use Lemma~\ref{lem:fctuniformbasin} with weights $(6,2,2,2,2)$
and coordinate order $(u,v,x,y,z)$. Before addition, round every
rational summand outwards to $2^{-80}\Z$ and each monomial upper
bound upwards to $2^{-80}\Z$ by integer-square-root enclosure of its rational square.
All-even monomials retain their nonnegative sign. Then
\[
 \left.\partial_{\epsilon}^2P_{V_\pm}(S+\epsilon T)\right|_{\epsilon=0}\ge|\eta|^{\mathsf T}C_\pm|\eta|,
\]
where $\eta$ is the coordinate vector of $T$ and
\[
 C_-=\frac1{10^6}
 \begin{pmatrix}
 17656592&-1723052&-55398&-55392&-1720776\\
 -1723052&6304087&-24382&-23238&-888593\\
 -55398&-24382&32520&-17272&-23164\\
 -55392&-23238&-17272&32623&-24247\\
 -1720776&-888593&-23164&-24247&6305113
 \end{pmatrix},
\]
\[
 C_+=\frac1{10^6}
 \begin{pmatrix}
 14612274&-1662571&-47365&-47360&-1660327\\
 -1662571&5713988&-21547&-20329&-848744\\
 -47365&-21547&325877&-14687&-20255\\
 -47360&-20329&-14687&325981&-21409\\
 -1660327&-848744&-20255&-21409&5715003
 \end{pmatrix}.
\]
Exact rational elimination gives $LDL^{\mathsf T}$ pivots of
$C_\pm-(201/50000)\operatorname{diag}(6,2,2,2,2)$ exceeding, respectively,
\[
 (17,6,1/50,1/100,5),\qquad(14,5,3/10,3/10,5).
\]
Affinity in $V^2$ extends
$\left.\partial_{\epsilon}^2F_{V_\pm}(S+\epsilon T)\right|_{\epsilon=0}>\|T\|_F^2/1000$ to
$V\in[V_-,V_+]$, for $T\ne0$ and $\|S\|_F\le1/50$.
The first diagonal entries give
\[
 \partial_u^2F_V(S)>14\qquad(\|S\|_F\le1/50).
\]

For the axial root, use the finite $N=6$ Ewald derivative sum in
Lemma~\ref{lem:sh110rootseed}, substituting the explicit rational
covolume $87/80$ and $Y_0=107818724/10^8$ above. Its complete
Gaussian derivative tails and finite interval evaluation give
\[
 \frac{11}{10^9}<\partial_uF_{87/80}(0)<\frac{12}{10^9}.
\]
The finite attraction derivative and its complete omitted-vector bound
$6T_3(1,32)$ give
\[
 \frac{4739}{1000}<\partial_u E(L_0(S),3)\big|_{S=0}
 <\frac{4743}{1000}<\frac{19}{4}.
\]
Hexagonal symmetry and the axial bounds give
\[
 (\partial_vF_V,\partial_xF_V,\partial_yF_V,\partial_zF_V)(u,0,0,0,0)=0,
 \qquad|\partial_uF_V(0)|<11/100,
\]
and
\[
 \partial_uF_V(1/125,0,0,0,0)>\frac1{500},\qquad
 \partial_uF_V(-1/125,0,0,0,0)<-\frac1{500},
\]
at the two axial endpoints inside the ball. Hence
\[
 \exists!\,u_V\in(-1/125,1/125):
 \partial_uF_V(u_V,0,0,0,0)=0,
 \qquad\nabla_{\mathbf s}F_V(S_V)=0.
\]
Strict convexity and $\partial_{uu}F_V(S_V)>14$ give uniqueness and
analyticity. The factor $2(87/80)^4>5/2$ gives the stated Hessian
and value-gap bounds.
This completes the proof of Lemma~\ref{lem:broadshbasin}.
\end{proof}

\begin{corollary}[Identification of the recognized stationary shapes]
\label{cor:transitionbasinidentification}
For $441/200\le2V^2\le121/50$, use the respective charts and
hypotheses of Lemmas~\ref{lem:transitionfccbasin}, \ref{lem:widefctbasin}
and~\ref{lem:broadshbasin}. Then
\[
 \nabla_{\mathbf s}\Phi_V=0
 \quad\Longrightarrow\quad
 \Lambda\cong F,\quad L_{X(V)},\quad H_{Y(V)},
\]
respectively. In the latter two cases,
\[
 \frac{149}{500}\frac{66}{67}<X<
 \frac{149}{500}\frac{202}{199},\qquad
 \frac{107818724}{10^8}\frac{41}{42}<Y<
 \frac{107818724}{10^8}\frac{127}{124}.
\]
\end{corollary}
\begin{proof}
Since $\partial_u\{(1+2u)/(1-u)\}=3/(1-u)^2>0$,
the formulas $X=X_0(1+2u)/(1-u)$, $Y=Y_0(1+2u)/(1-u)$ and
$|u|<1/200$, $|u|<1/125$ give the stated intervals, contained
respectively in $(23/100,307/1000)$ and $(103/100,103/50)$.
Lemma~\ref{num:branchmoments} gives exactly one stationary parameter
there for each target $V$, identifying the classified branches.
Cubic symmetry gives
$\left.\nabla_{\mathbf s}\Phi_V(\Lambda_F(S))\right|_{S=0}=0$
for every $V$.
This completes the proof of Corollary~\ref{cor:transitionbasinidentification}.
\end{proof}

\begin{theorem}[The complete transition strip]
\label{thm:complete-transition-strip}
For $21/20\le V\le11/10$, every global minimizer is one of
the three stationary candidates in Proposition~\ref{prop:familyphases}.
Their ordering, contact covolumes and equality cases are therefore
the ones stated there.
\end{theorem}
\begin{proof}
Let $f(V)$ be the least candidate value
in Proposition~\ref{prop:familyphases}. Branch stationarity and
the moment ordering make $f$ concave in $V^2$, including both contacts.
The endpoint comparisons are Theorems~\ref{thm:fcccriticalglobal}
and~\ref{thm:sh110global}, and
Lemma~\ref{lem:transition-profile-moment-domain} gives the compact
moment domain. Lemma~\ref{lem:transitionprofilecover} yields
\[
 E(\Lambda,6)-2V^2 E(\Lambda,3)\ge f(V),
 \qquad |\Lambda|=1,\quad\frac{21}{20}\le V\le\frac{11}{10},
\]
with equality precisely for a stationary candidate attaining $f$.
Since $2V^4\mathcal E(V^{1/3}\Lambda)=E(\Lambda,6)-2V^2E(\Lambda,3)$,
Proposition~\ref{prop:familyphases} gives the classification,
all equality cases and the two derivative jumps.
This completes the proof of Theorem~\ref{thm:complete-transition-strip}.
\end{proof}

\section{The intervening simple-hexagonal phase}\label{sec:longphase}

For $11/10\le V\le89/50$, moment contraction, Fourier convexity
and layer alignment give
\[
 \operatorname*{argmin}_{|\Lambda|=1}
 \{E(\Lambda,6)-2V^2E(\Lambda,3)\}
 =\{\Lambda:\Lambda\cong H_{Y(V)}\}.
\]
The curve $Y(V)$ is defined in Lemma~\ref{lem:longshreference}.

\begin{lemma}[Area coordinates and complete derivative tails]
\label{lem:area-coordinate-jets}
For $t,r>0$ and $x,y,z\in\mathbb R$, set
\begin{equation}\label{eq:area-coordinate-map}
\begin{gathered}
 \eta=(t,r,x,y,z),\qquad
 p(\eta)=(a,\delta,x,y,z)=(\sqrt{tr},\sqrt{t/r},x,y,z),\qquad
 Y=\frac{\sqrt3}{2t^3},\\
 U(\eta)=U(p(\eta))=
 \begin{pmatrix}
 \sqrt{tr}&\sqrt{tr}\,x&\sqrt{tr}\,y\\
 0&\sqrt{t/r}&\sqrt{t/r}\,z\\0&0&t^{-1}
 \end{pmatrix},\qquad \det U=1.
\end{gathered}
\end{equation}
For $n=(i,j,k)\in\mathbb Z^3$, put
$q_n=|U(\eta)n|^2$ and $q_n^*=|U(\eta)^{-\mathsf T}n|^2$. Then
\[
\begin{aligned}
 q_n&=tr(i+xj+yk)^2+\frac tr(j+zk)^2+\frac{k^2}{t^2},\\
 q_n^*&=\frac{i^2}{tr}+\frac rt(j-xi)^2+t^2\{k-z(j-xi)-yi\}^2.
\end{aligned}
 \tag{L9}\label{eq:area-squared-lengths}
\]
Put
\[
 A=i+xj+yk,\quad B=j+zk,\quad C=k^2/t^2,\quad
 D=j-xi,\quad F=k-zD-yi.
\]
In coordinate order $(t,r,x,y,z)$ the gradients are
\[
\begin{aligned}
 \nabla_\eta q_n&=\left((rA^2+B^2/r)-2C/t,\ t(A^2-B^2/r^2),\ 2trAj,\ 2trAk,\ 2tBk/r\right),\\
 \nabla_\eta q_n^*&=\bigl(-(i^2/r+rD^2)/t^2+2tF^2,\ (D^2-i^2/r^2)/t,\\
 &\qquad -2rDi/t+2t^2Fzi,\ -2t^2Fi,\ -2t^2FD\bigr).
\end{aligned}
\]
The Hessians, in the order $(t,r,x,y,z)$, are given by
\[
\begin{gathered}
 \nabla^2_\eta q_n=
 \begin{pmatrix}
 6C/t^2&A^2-B^2/r^2&2rAj&2rAk&2Bk/r\\
 A^2-B^2/r^2&2tB^2/r^3&2tAj&2tAk&-2tBk/r^2\\
 2rAj&2tAj&2trj^2&2trjk&0\\
 2rAk&2tAk&2trjk&2trk^2&0\\
 2Bk/r&-2tBk/r^2&0&0&2tk^2/r
 \end{pmatrix},\\[4pt]
 \begin{aligned}
 \partial_{tt}q_n^*&=2(i^2/r+rD^2)/t^3+2F^2,&
 \partial_{rr}q_n^*&=2i^2/(tr^3),\\
 \partial_{xx}q_n^*&=2ri^2/t+2t^2z^2i^2,&
 \partial_{yy}q_n^*&=2t^2i^2,\\
 \partial_{zz}q_n^*&=2t^2D^2,&
 \partial_{tr}q_n^*&=(i^2/r^2-D^2)/t^2,\\
 \partial_{tx}q_n^*&=2rDi/t^2+4tFzi,&
 \partial_{ty}q_n^*&=-4tFi,\\
 \partial_{tz}q_n^*&=-4tFD,&
 \partial_{rx}q_n^*&=-2Di/t,\\
 \partial_{ry}q_n^*&=\partial_{rz}q_n^*=0,&
 \partial_{xy}q_n^*&=-2t^2zi^2,\\
 \partial_{xz}q_n^*&=2t^2i(F-zD),&
 \partial_{yz}q_n^*&=2t^2iD.
 \end{aligned}
\end{gathered}
 \tag{L10}\label{eq:area-squared-length-jets}
\]
Complete $\nabla^2_\eta q_n^*$ by symmetry.

Let $\mathcal B_\eta$ be any closed area-coordinate rectangle with
$t,r>0$, and let $\mathcal B_p$ contain $p(\mathcal B_\eta)$.
Here $(\eta_0,\eta_1,\eta_2,\eta_3,\eta_4)=(t,r,x,y,z)$ and
$(p_0,p_1,p_2,p_3,p_4)=(a,\delta,x,y,z)$.
The nontrivial embedding derivatives are
\[
\begin{gathered}
 a_t=a/(2t),\quad a_r=a/(2r),\quad
 \delta_t=\delta/(2t),\quad\delta_r=-\delta/(2r),\\
 a_{tt}=-a/(4t^2),\quad a_{tr}=a/(4tr),\quad a_{rr}=-a/(4r^2),\\
 \delta_{tt}=-\delta/(4t^2),\quad
 \delta_{tr}=-\delta/(4tr),\quad\delta_{rr}=3\delta/(4r^2).
\end{gathered}
\]
For $k=2,3,4$, $\partial_{\eta_i\eta_j}p_k=0$ ($0\le i,j\le4$).
For any $C^2$ remainder $R(p)$ with complete derivative bounds
\[
 |\partial_{p_k}R|\le e_k,\qquad
 |\partial_{p_kp_l}R|\le e_{kl}\qquad\text{on }\mathcal B_p,
\]
for $0\le i,j\le4$, the area-coordinate tail errors obey
\begin{equation}\label{eq:area-derivative-tails}
\begin{aligned}
 |\partial_{\eta_i}(R\circ p)|
 &\le\sum_{k=0}^4\sup_{\mathcal B_\eta}|\partial_{\eta_i}p_k|e_k,\\
 |\partial_{\eta_i\eta_j}(R\circ p)|
 &\le\sum_{k=0}^4\sum_{l=0}^4\sup_{\mathcal B_\eta}|\partial_{\eta_i}p_k|e_{kl}
                   \sup_{\mathcal B_\eta}|\partial_{\eta_j}p_l|
       +\sum_{k=0}^1\sup_{\mathcal B_\eta}|\partial_{\eta_i\eta_j}p_k|e_k.
\end{aligned}\tag{L12}
\end{equation}
All suprema are over the whole rectangle; no Korkine--Zolotarev inequality is required.
\end{lemma}
\begin{proof}
Direct differentiation gives
\eqref{eq:area-squared-lengths}--\eqref{eq:area-squared-length-jets};
with $(Dp)_{ij}=\partial p_i/\partial\eta_j$ ($0\le i,j\le4$), the chain rule gives
\[
\begin{aligned}
 \nabla_\eta(R\circ p)&=(Dp)^{\mathsf T}\nabla_pR,\\
 \nabla_\eta^2(R\circ p)&=(Dp)^{\mathsf T}(\nabla_p^2R)(Dp)
       +\sum_{k=0}^4(\partial_{p_k}R)\nabla_\eta^2p_k.
\end{aligned}
\]
and hence \eqref{eq:area-derivative-tails}.
This completes the proof of Lemma~\ref{lem:area-coordinate-jets}.
\end{proof}

\begin{lemma}[A monotone stationary simple-hexagonal reference curve]
\label{lem:longshreference}
For the unit-covolume simple-hexagonal lattice $H_Y$, put
\[
 B(Y)=E(H_Y,3),\qquad A(Y)=E(H_Y,6),\qquad
 \Gamma(Y)=\frac{A_t(Y)}{B_t(Y)}.
\]
Here the subscripts denote logarithmic derivatives:
\[
 f_t(Y)=Yf'(Y),\qquad f_{tt}(Y)=Y^2f''(Y)+Yf'(Y)
 \quad(f=A,B,\Gamma,q,q^*).
\]
On the entire closed interval $11/10\le Y\le6$,
\[
\begin{gathered}
 2<B_t(Y)<30,\qquad \frac35<\Gamma_t(Y)<\frac{17}2,\\
 \frac{2417273}{10^6}<\Gamma(11/10)<\frac{2417274}{10^6},\qquad
 \frac{5744686}{10^6}<\Gamma(39/10)<\frac{5744687}{10^6},\\
 \frac{8630510}{10^6}<\Gamma(6)<\frac{8630511}{10^6}.
\end{gathered}
\]
The inverse branch $Y(V)=\Gamma^{-1}(2V^2)$ is real analytic for
\[
 \Gamma(11/10)\le2V^2\le\Gamma(6)
\]
and satisfies
\[
 A_t(Y(V))=2V^2B_t(Y(V)),\qquad
 \frac{\mathrm d}{\mathrm dV}B(Y(V))
 =4V\frac{B_t(Y(V))}{\Gamma_t(Y(V))}>\frac{16V}{17},
\]
\[
 A_{tt}-2V^2B_{tt}=B_t\Gamma_t>\frac65.
\]
For $f(V)=A(Y(V))-2V^2B(Y(V))$,
\[
 f'(V)=-4VB(Y(V)),\qquad
 f''(V)-\frac{f'(V)}V
 =-4V\frac{\mathrm d}{\mathrm dV}B(Y(V))<0.
\]
Thus $f$ is strictly concave as a function of $V^2$.
For $11/10\le V\le89/50$,
\[
 Y(V)<\frac{217}{50},\qquad
 B(Y(V))<\frac{43}{2},\qquad A(Y(V))<64.
\]
At $V\in\{11/10,89/50\}$, Theorem~\ref{thm:sh110global} and
Proposition~\ref{prop:sh-global-178} give
\[
 \operatorname*{argmin}_{|L|=V}\mathcal E(L)
 =\{L:L\cong V^{1/3}H_{Y(V)}\}.
\]
\end{lemma}

\begin{proof}
Put $c=(4/(3Y))^{1/3}$. The Gram matrix of $H_Y$ is
\[
 Q_Y=c\begin{pmatrix}1&1/2&0\\1/2&1&0\\0&0&Y\end{pmatrix}.
\]
For $1\le Y\le13/2$,
$\sigma(Q_Y)=\{c/2,3c/2,cY\}\subset(1/4,4)$.
The same bounds hold for $Q_Y^{-1}$ and justify termwise differentiation.

Use the kernels
\[
 P_s(q)=\frac{\pi^s}{\Gamma(s)}\mathcal G_s(\pi q),\qquad
 D_s(q)=\frac{\pi^s}{\Gamma(s)}\mathcal G_{3/2-s}(\pi q),
 \qquad
 \mathcal G_a(x)=\int_1^\infty e^{-xu}u^{a-1}\,\mathrm du.
\]
The unit-covolume version of the Ewald identity in
Lemma~\ref{lem:firstloss} is
\[
 E(H_Y,s)=\frac{\pi^s}{\Gamma(s)}
 \left(\frac1{s-3/2}-\frac1s\right)
 +\sum_{n\in\mathbb Z^3\setminus\{0\}}\bigl\{P_s(n^{\mathsf T}Q_Yn)
                       +D_s(n^{\mathsf T}Q_Y^{-1}n)\bigr\}.
\]
The constant terms are $\pi^3/6$ and $\pi^6/2160$.
For $0\le a\le432$, $0\le b\le144$, group the complete finite cube by
\[
 m_{ab}=\#\{(i,j,k)\in[-12,12]^3\cap\Z^3\setminus\{0\}:
                i^2+ij+j^2=a,\ k^2=b\}.
\]
Omit $m_{ab}=0$; $\sum_{a=0}^{432}\sum_{b=0}^{144}m_{ab}=25^3-1$. For
$q=ca+cYb$, $q^*=4a/(3c)+b/(cY)$,
\[
 E^{\rm fin}(H_Y,s)=\frac{\pi^s}{\Gamma(s)}
 \left(\frac1{s-3/2}-\frac1s\right)
 +\sum_{a=0}^{432}\sum_{\substack{b=0\\m_{ab}>0}}^{144}
              m_{ab}\{P_s(q)+D_s(q^*)\}.
\]
The substitution $j\mapsto-j$ preserves the cube and converts the
dual form $i^2-ij+j^2$ to the same multiplicities. Logarithmic differentiation gives
\[
\begin{aligned}
 q_t&=-\frac{ca}{3}+\frac{2cYb}{3},&
 q_{tt}&=\frac{ca}{9}+\frac{4cYb}{9},\\
 q_t^*&=\frac{4a}{9c}-\frac{2b}{3cY},&
 q_{tt}^*&=\frac{4a}{27c}+\frac{4b}{9cY}.
\end{aligned}
\]
For $K=P_s,D_s$ and $0\le r\le2$,
\[
\begin{gathered}
 Y\frac{\mathrm d}{\mathrm dY}K(q)=K'(q)q_t,\qquad
 \left(Y\frac{\mathrm d}{\mathrm dY}\right)^2K(q)
 =K''(q)q_t^2+K'(q)q_{tt},\\
 P_s^{(r)}(q)=\frac{(-\pi)^r\pi^s}{\Gamma(s)}
                     \mathcal G_{s+r}(\pi q),\qquad
 D_s^{(r)}(q)=\frac{(-\pi)^r\pi^s}{\Gamma(s)}
                     \mathcal G_{3/2-s+r}(\pi q).
\end{gathered}
\]
For $q\ge20$, Lemma~\ref{lem:criticalKZcover} gives
$|K^{(r)}(q)|<5e^{-3q}$, $0\le r\le2$.
Since $|q_t|\le2q/3$, $|q_{tt}|\le4q/9$,
$qe^{-q}\le1$ and $q^2e^{-q}\le4$, every composed derivative
through order two is bounded by $20e^{-2q}$.
For omitted indices, $|n|\ge13$, $q,q^*\ge169/4>20$, and
$2q,2q^*>\pi|n|^2/8$. The complete primal-plus-dual remainder is
\[
 T=240\,\frac{e^{-169\pi/8}(1+\sqrt8)^2}
 {1-e^{-27\pi/8}}<\frac6{10^{26}}.
\]
Add $[0,T]$ to each moment and $[-T,T]$ to each of its first
two derivatives.

Use the outward arithmetic, checked radicals, rational exponentials
and half-integer recurrences of Lemmas~\ref{lem:firstloss}
and~\ref{lem:criticalKZcover}. Enclose $c=(4/(3Y))^{1/3}$ by
outward cubing, extending the positive endpoints until both checks
hold. Monotonicity of $\mathcal G_a$ encloses interval arguments
by endpoint evaluation.

For $q_-\ge20$, let
$v=\max\{b\in\{20,40,80,160,320\}:b\le q_-\}$ and use
\[
 K\in[0,5e^{-3v}],\quad K'\in[-5e^{-3v},0],\quad
 K''\in[0,5e^{-3v}].
\]
Evaluate nonzero multiplicities in increasing $(a,b)$ order,
with outward additions and no Taylor interpolation.

Now partition $[11/10,6]$ into the 1024 closed rational cells
\[
 \left[\frac{11}{10}+\frac{49k}{10240},
           \frac{11}{10}+\frac{49(k+1)}{10240}\right],
                 \qquad 0\le k<1024.
\]
On each whole cell evaluate the six complete sums, check
$\underline{B_t}>0$, and form
\[
 \Gamma=\frac{A_t}{B_t},\qquad
 \Gamma_t=\frac{A_{tt}-\Gamma B_{tt}}{B_t}.
\]
The rational aggregate enclosures from these 1024 closed-cell
evaluations are
\[
 \frac{20368}{10^4}<B_t<\frac{282127}{10^4},\qquad
 \frac{6184}{10^4}<\Gamma_t<\frac{83996}{10^4}.
\]
These imply the uniform bounds. The same evaluations give the three
endpoint inequalities and
\[
 6342227<10^6\Gamma(217/50)<6342228,
\]
and, at $Y=217/50$,
\[
 B(Y)<\frac{214321}{10^4}<\frac{43}{2},\qquad
 A(Y)<\frac{636842}{10^4}<64.
\]
Thus
\[
 \Gamma_t>0,\qquad
 [121/50,7921/1250]\subset\Gamma([11/10,217/50]),
\]
and the inverse-function theorem gives a unique analytic $Y(V)$.
At this root,
\[
 \frac{\mathrm dB}{\mathrm dV}=4V\frac{B_t}{\Gamma_t}>0,\qquad
 A_{tt}-2V^2B_{tt}=B_t\Gamma_t,\qquad
 f'=-4VB,\quad f''-\frac{f'}V=-4V\frac{\mathrm dB}{\mathrm dV}<0.
\]
The derivative bounds and
\[
 \Gamma(217/50)>7921/1250\ \Longrightarrow\ Y(V)<217/50,
 \qquad B_t>0,\quad A_t=\Gamma B_t>0
\]
give the moment caps from the values at $217/50$.

At $V=11/10$, Theorem~\ref{thm:sh110global} places the global
SH parameter in the interval
\[
 (1.1013103262,1.1013103264)\subset(11/10,6).
\]
Its parameter satisfies
$0=Y\partial_Y\Phi_{11/10}(H_Y)=A_t(Y)-2(11/10)^2B_t(Y)$,
so $\Gamma(Y)=2(11/10)^2$ and $Y=Y(11/10)$.
At $V=89/50$, $39/10<Y(V)<6$ by the endpoint inequalities.
The physical triangular planar area is
\[
 A_{\rm plane}=\left(\frac{\sqrt3V^2}{2Y}\right)^{1/3}.
\]
Hence $77/100<A_{\rm plane}<89/100$ and
$\mathrm dA_{\rm plane}/\mathrm dY=-A_{\rm plane}/(3Y)\ne0$, so
$A_t-2V^2B_t=0\iff F_V'(A_{\rm plane})=0$.
Proposition~\ref{prop:sh-global-178} gives $F_V''>34$ on
$[2/5,89/100]$; hence
$\{A\in[2/5,89/100]:F_{89/50}'(A)=0\}
 =\left\{\left(\frac{\sqrt3(89/50)^2}{2Y(89/50)}\right)^{1/3}\right\}$,
identifying the second endpoint.
This completes the proof of Lemma~\ref{lem:longshreference}.
\end{proof}

\begin{lemma}[Moment contraction of the long SH profile]
\label{lem:longshmomentcontraction}
Put $V_0=11/10$, $V_1=89/50$ and
\[
 f(V)=A(Y(V))-2V^2B(Y(V)),\qquad
 g(t)=\max_{V_0\le V\le V_1}\{f(V)+2V^2t\}.
\]
For $b_i=B(Y(V_i))$, the unique maximizing covolume is
\begin{equation}\label{eq:longshmaximizingdensity}
 V_*(t)=
 \begin{cases}
 V_0,&t\le b_0,\\
 \sqrt{\Gamma(B^{-1}(t))/2},&b_0\le t\le b_1,\\
 V_1,&t\ge b_1.
 \end{cases}
\end{equation}
Here $B^{-1}$ is restricted to $[Y(V_0),Y(V_1)]$; the endpoint formulas agree, and
\[
 g(t)=
 \begin{cases}
 f(V_0)+2V_0^2t,&t\le b_0,\\
 A(B^{-1}(t)),&b_0\le t\le b_1,\\
 f(V_1)+2V_1^2t,&t\ge b_1.
 \end{cases}
\]

For every real interval $\mathcal B=[\beta_-,\beta_+]$, the following
finite construction encloses $\{2V_*(t)^2:t\in\mathcal B\}$. Use the endpoint moment enclosures
$[b_i^-,b_i^+]$ in \eqref{eq:longshendpointmoments}, whose outer
band lies strictly inside the nodal moment domain. Let $G_i$
enclose the exact rational $2V_i^2$ outwards. Retain $G_0$
if $\beta_-\le b_0^+$ and retain $G_1$ if
$\beta_+\ge b_1^-$. If
$\mathcal B\cap[b_0^-,b_1^+]=[\ell,h]\ne\varnothing$, define
\begin{align*}
 i_-&=\max\bigl(0,\min\{j:B_j^+\ge\ell\}-1\bigr),\\
 i_+&=\min\bigl(4096,
       \min\bigl(\{j:B_j^->h\}\cup\{4096\}\bigr)\bigr),
\end{align*}
and retain the nonempty interval
\begin{equation}\label{eq:longshfiniteinverse}
 \bigl[\max(\inf G_0,\Gamma_{i_-}^-),
                  \min(\sup G_1,\Gamma_{i_+}^+)\bigr].
\end{equation}
The hull contains $2V_*(t)^2$ for all $t\in\mathcal B$, endpoints included;
both searches use the $4097$ ordered nodes. Under the endpoint comparisons,
profile minima lie in $[b_0,b_1]$. Intersect the enclosure with every
inherited supporting interval before testing
$\nabla_\eta\Phi_V(U(\eta)\mathbb Z^3)=0$.
\end{lemma}

\begin{proof}
By Lemma~\ref{lem:longshreference},
\begin{align*}
 \frac1{4V}\partial_V\{f(V)+2V^2t\}&=t-B(Y(V)),\\
 \partial_V\left[\frac1{4V}\partial_V\{f(V)+2V^2t\}\right]
 &=-\frac{\mathrm d}{\mathrm dV}B(Y(V))<0.
\end{align*}
Hence \eqref{eq:longshmaximizingdensity}, with
$t=B(Y(V_*))$ and $g(t)=A(Y(V_*))$ in its interior.
The guards retain all endpoint cases; for $t\in\mathcal B\cap[b_0,b_1]$,
the ordered searches, including their boundary conventions, give
\[
 y_{i_-}\le B^{-1}(t)\le y_{i_+},\qquad
 \Gamma_{i_-}^-\le2V_*(t)^2\le\Gamma_{i_+}^+.
\]
At array boundaries use strict containment of both roots in the node domain;
otherwise $B_{i_-}^+<\ell$, $B_{i_+}^->h$. Intersection with $[\inf G_0,\sup G_1]$
and the hull give \eqref{eq:longshfiniteinverse};
Proposition~\ref{prop:density-gap-reduction} gives the profile-minimum assertions.
This completes the proof of Lemma~\ref{lem:longshmomentcontraction}.
\end{proof}

\begin{lemma}[Separate planar and layer moment bounds]
\label{lem:shasymmetricmoment}
Fix $t>0$ in the coordinates of Lemma~\ref{lem:area-coordinate-jets}.
Suppose that throughout a convex transverse rectangle containing
$(2/\sqrt3,1/2,0,0)$,
\[
 \left\|\nabla^2_{(r,x,y,z)}E(U(t,r,x,y,z)\mathbb Z^3,3)\right\|_{\rm op}
 \le M.
\]
Then every point of that rectangle satisfies
\[
 -\frac M2(y^2+z^2)
 \le E(U\mathbb Z^3,3)-E(H_Y,3)
 \le\frac M2\{(r-2/\sqrt3)^2+(x-1/2)^2\}.
 \tag{A1}
\]
If moreover
\[
 \frac{11}{10}\le Y\le\frac{217}{50},\qquad
 \frac{11}{10}\le V\le\frac{89}{50},\qquad
 E(U(t,r,x,y,z)\mathbb Z^3,3)=B(Y(V)),
\]
then
\[
 -\frac{17M}{8}(y^2+z^2)
 \le2V^2-\Gamma(Y)
 \le\frac{17M}{8}\{(r-2/\sqrt3)^2+(x-1/2)^2\}.
 \tag{A2}
\]
\end{lemma}
\begin{proof}
Let $\Gamma$ be the planar lattice of area $t$ and
$u=(\sqrt{tr}\,y,\sqrt{t/r}\,z)$. For $s>0$, Montgomery's
theorem \cite[Theorem~1]{Mont} and Poisson summation give
\[
 \sum_{v\in\Gamma}e^{-s|v|^2}
 \ge\sum_{v\in T_t}e^{-s|v|^2},
\]
\[
 \sum_{v\in\Gamma}e^{-s|v+ku|^2}
 =\frac\pi{ts}\sum_{w\in\Gamma^*}
   e^{-\pi^2|w|^2/s}\cos(2\pi k w\cdot u)
 \le\sum_{v\in\Gamma}e^{-s|v|^2}.
\]
Multiply by $e^{-sk^2/t^2}$, sum over $k\in\mathbb Z$, subtract
the common zero-vector term and use
$q^{-3}=\frac12\int_0^\infty s^2e^{-sq}\,\mathrm ds$. Thus
\[
 E(\Gamma\oplus t^{-1}\mathbb Z,3)
 \ge E(T_t\oplus t^{-1}\mathbb Z,3)=E(H_Y,3),\tag{A3}
\]
\[
 E(\Gamma+\mathbb Z(u,t^{-1}),3)
 \le E(\Gamma\oplus t^{-1}\mathbb Z,3).\tag{A4}
\]
Layer reflection gives $\nabla_{y,z}E(U(t,r,x,y,z)\mathbb Z^3,3)=0$
at $y=z=0$; hexagonal symmetry gives
$\nabla_{r,x}E(U(t,r,x,0,0)\mathbb Z^3,3)=0$ at $(r,x)=(2/\sqrt3,1/2)$.
The Hessian bound and Taylor's formula on the two segments yield
\begin{align*}
 E(U(t,r,x,y,z)\mathbb Z^3,3)
 &\ge E(U(t,r,x,0,0)\mathbb Z^3,3)-\tfrac M2(y^2+z^2),\\
 E(U(t,r,x,0,0)\mathbb Z^3,3)
 &\le E(H_Y,3)+\tfrac M2\{(r-2/\sqrt3)^2+(x-1/2)^2\}.
\end{align*}
The segments remain in the rectangle; (A3)--(A4) give (A1).
For $E(U\mathbb Z^3,3)=B(Y(V))$, Lemma~\ref{lem:longshreference} gives
\[
 E(U\mathbb Z^3,3)-E(H_Y,3)
 =\int_{\sqrt{\Gamma(Y)/2}}^V
       \frac{\mathrm d}{\mathrm dv}B(Y(v))\,\mathrm dv,\qquad
 \frac1{4v}\frac{\mathrm d}{\mathrm dv}B(Y(v))>\frac4{17}.
\]
Since $\Gamma_t>0$ and $Y,Y(V)\in[11/10,217/50]$, the whole
integration segment remains in the inverse-branch domain. The two
signs in (A1) yield (A2).
This completes the proof of Lemma~\ref{lem:shasymmetricmoment}.
\end{proof}

\begin{lemma}[A Fourier transverse tube for the SH profile]
\label{lem:shfouriertransverse}
Use the coordinates of Lemma~\ref{lem:area-coordinate-jets} and the profile
$\Phi_V$.
For every $11/10\le Y\le217/50$, the conditions
\begin{equation}\label{eq:fouriertubebox}
 \left|r-\frac2{\sqrt3}\right|\le\frac1{200},\qquad
 \frac{99}{200}\le x\le\frac12,\qquad
 |y|,|z|\le\frac1{40},
\end{equation}
\[
 -\frac{17}{640}\le2V^2-\Gamma(Y)\le\frac{51}{3200}
\]
imply
\[
 \nabla^2_{(r,x,y,z)}\Phi_V(U\mathbb Z^3)\succ0.
\]
For $11/10\le Y\le217/50$, $11/10\le V\le89/50$ and coordinates in
\eqref{eq:fouriertubebox}, at fixed $t$,
\[
\begin{gathered}
 E(U(t,r,x,y,z)\mathbb Z^3,3)=B(Y(V)),\qquad
 \nabla_{r,x,y,z}\Phi_V(U(t,r,x,y,z)\mathbb Z^3)=0\\
 \Longrightarrow\quad
 (r,x,y,z)=(2/\sqrt3,1/2,0,0),\qquad Y=Y(V).
\end{gathered}
\]
\end{lemma}

\begin{proof}
All finite decimals below denote exact rationals.

\paragraph{The layer decomposition.}
Let $\Lambda_2$ be the planar lattice of area $t$ with basis
$(\sqrt{tr},0)$ and $(\sqrt{tr}\,x,\sqrt{t/r})$, and set
$u=(\sqrt{tr}\,y,\sqrt{t/r}\,z)$. For its dual vector
indexed by $(m,n)$,
\[
 R^2=q=\frac{m^2}{tr}+\frac rt(n-xm)^2,\qquad
 w\cdot u=my+(n-xm)z.
\]
Equation~\eqref{eq:fourier}, Lemma~\ref{lem:besselproperties} and
planar Poisson summation give, for $s=3,6$,
\begin{align}
 E(U\mathbb Z^3,s)
 ={}&E(\Lambda_2,s)
 +\frac{2\pi\zeta(2s-2)}{t(s-1)}(1/t)^{2-2s}\notag\\
 &+\frac2t\sum_{k=1}^{\infty}\sum_{(m,n)\in\mathbb Z^2\setminus\{(0,0)\}}
       F_s(q,k(1/t))\cos\psi,\label{eq:fouriertubemoments}\\
 F_s(q,d)={}&\frac{2\pi^s}{\Gamma(s)}
             (R/d)^{s-1}K_{s-1}(2\pi dR),\qquad
 \psi=2\pi k\{my+(n-xm)z\}.\notag
\end{align}
The zero mode has no transverse derivatives; for $\Phi_V$, the
nonzero-layer coefficients are
\begin{align}
 C_0(q,d)&=\frac{\pi^6}{60}(R/d)^5K_5(2\pi dR)
                -2V^2\pi^3(R/d)^2K_2(2\pi dR),\notag\\
 C_1(q,d)&=-\frac{\pi^7}{60}(R/d)^4K_4(2\pi dR)
                +2V^2\pi^4(R/d)K_1(2\pi dR),\notag\\
 C_2(q,d)&=\frac{\pi^8}{60}(R/d)^3K_3(2\pi dR)
                -2V^2\pi^5K_0(2\pi dR).
 \label{eq:fouriertubecoefficients}
\end{align}
By Lemma~\ref{lem:besselproperties},
$C_1=\partial_q C_0$ and $C_2=\partial_q^2C_0$.
The factor $2/t$ in \eqref{eq:fouriertubemoments} includes the
two signs of the layer index, and its dual sum includes both
signs of every nonzero dual vector.

\paragraph{Finite Bessel enclosures.}
No Bessel approximation is assumed in the finite calculation.
For $\nu=0,1$, use \eqref{eq:fouriertubebesselintegral}.
For $p=\pm1/2$, $a>0$ and $N\in\mathbb Z_{\ge0}$, Taylor's theorem gives
\[
 \operatorname{sgn}\left((1+a)^p-\sum_{j=0}^N\binom{p}{j} a^j\right)
 =\operatorname{sgn}\prod_{j=0}^N(p-j).
\]
Integrate this sign against the positive weight in
\eqref{eq:fouriertubebesselintegral}. Define
\[
 a_{\nu,0}=1,\qquad
 a_{\nu,j}=a_{\nu,j-1}
          \frac{4\nu^2-(2j-1)^2}{8jv}.
\]
The exact gamma moments give
\begin{align}
 \sum_{j=0}^{15}a_{0,j}(v)&\le
      \frac{K_0(v)}{\sqrt{\pi/(2v)}e^{-v}}
      \le\sum_{j=0}^{14}a_{0,j}(v),\notag\\
 \sum_{j=0}^{14}a_{1,j}(v)&\le
      \frac{K_1(v)}{\sqrt{\pi/(2v)}e^{-v}}
      \le\sum_{j=0}^{15}a_{1,j}(v).
 \label{eq:fouriertubebesselbounds}
\end{align}
For the retained arguments $v>7$,
\[
 \sum_{j=1}^{15}|a_{0,j}(7)|<\frac1{50},\qquad
 \sum_{j=1}^{15}|a_{1,j}(7)|<\frac1{10}.
\]
These sums decrease with $v$, so the bounds are positive.
For $\nu=1,2,3,4$, enclose the remaining orders by the positive
recurrence in Lemma~\ref{lem:besselproperties}. Since $K_\nu'<0$,
on $[v_-,v_+]$ use the lower bound at $v_+$ and the upper bound at $v_-$.

Use the outward arithmetic, checked radicals and scaled degree-21/20
exponential bounds of Lemma~\ref{lem:criticalKZcover} in every
coefficient recurrence.

To preserve radial cancellation in \eqref{eq:fouriertubecoefficients},
first form the three direct interval enclosures $I_0,I_1,I_2$.
For the midpoint $q_c$ of the squared-radius interval, evaluate
the same expressions there with the full $d,V$ intervals,
obtaining $J_0,J_1,J_2$. Intersect
\begin{equation}\label{eq:fouriertuberadialrefinement}
 I_0\longleftarrow I_0\cap(J_0+I_1(q-q_c)),\qquad
 I_1\longleftarrow I_1\cap(J_1+I_2(q-q_c)).
\end{equation}
The mean value theorem at each fixed $d,V$ proves both
inclusions; the first uses the unrefined $I_1$.

\paragraph{The planar Hessian and its complete tail.}
Set
\[
\begin{gathered}
 Q=r(i+xj)^2+j^2/r,\qquad
 \nabla_{r,x}Q=\bigl((i+xj)^2-j^2/r^2,\ 2rj(i+xj)\bigr),\\
 \nabla^2_{r,x}Q=\begin{pmatrix}
 2j^2/r^3&2j(i+xj)\\2j(i+xj)&2rj^2
 \end{pmatrix}.
\end{gathered}
\]
At area one, for $s\in\{3,6\}$,
\[
 (P_s)_{ab}=\sum_{(i,j)\in\mathbb Z^2\setminus\{(0,0)\}}
 \{s(s+1)Q^{-s-2}Q_aQ_b-sQ^{-s-1}Q_{ab}\},
 \qquad a,b\in\{r,x\}.
\]
The finite sum uses
\[
 (i,j)\in[-40,40]^2\cap\mathbb Z^2,\qquad
 i>0\ \lor\ (i=0,\ j>0),\qquad\text{multiplicity }2,
\]
with diagonal interval squares $Q_a^2$. On
\[
 1149/1000\le r\le29/25,\qquad49/100\le x\le51/100
\]
the two Gershgorin lower bounds of the planar Gram matrix exceed
$1/2$. Thus $Q\ge(i^2+j^2)/2$. Also
\[
 |Q_a|\le\alpha_aQ,\qquad |Q_{ab}|\le\beta_{ab}Q,\qquad
 \alpha=(1,6/5),\quad
 \beta=\begin{pmatrix}8/5&1\\1&3\end{pmatrix}.
\]
The arithmetic--geometric mean inequality gives
$|Q_r|\le Q/r$, $|Q_x|\le rQ$, $|Q_{rx}|\le Q$;
also $2/r^2<8/5$ and $2r^2<3$.

Since $\#\{(i,j)\in\mathbb Z^2:\max(|i|,|j|)=\ell\}=8\ell$
for $\ell\in\mathbb Z_{\ge1}$, the complete omitted planar moment beyond cutoff $N$ is at most
\[
 \frac{8\,2^s}{2s-2}N^{2-2s},
\]
and the omitted Hessian entry has absolute value at most
\begin{equation}\label{eq:fouriertubeplanartail}
 \{s(s+1)\alpha_a\alpha_b+s\beta_{ab}\}
 \left(\frac{8\,2^s}{2s-2}40^{2-2s}\right).
\end{equation}
Add \eqref{eq:fouriertubeplanartail} symmetrically to the finite sums.
At fixed $t$, the complete planar block is
$(t^{-6}P_6-2V^2t^{-3}P_3)\oplus0_{2\times2}$.

\paragraph{The finite nonzero-layer Hessian.}
Use the half-dual indices
\[
 k\in\{1,2,3\},\qquad (m,n)\in[-4,4]^2\cap\mathbb Z^2,\qquad
 m>0\ \lor\ (m=0,\ n>0),
\]
with factor $4/t$. The derivatives are
\[
 \begin{gathered}
 \nabla_{r,x,y,z}q=\left(-\frac{m^2}{tr^2}+\frac{(n-xm)^2}{t},
                     -\frac{2rm(n-xm)}{t},0,0\right),\\
 (q_{ab})_{a,b\in\{r,x\}}=
 \begin{pmatrix}
 2m^2/(tr^3)&-2m(n-xm)/t\\-2m(n-xm)/t&2rm^2/t
 \end{pmatrix},
 \end{gathered}
\]
\[
 \nabla_{r,x,y,z}\psi=(0,-(2\pi k)mz,(2\pi k)m,(2\pi k)(n-xm)),\qquad
 \psi_{xz}=\psi_{zx}=-(2\pi k)m,
\]
with all other second phase derivatives zero.
For each retained term the exact Hessian entry is
\begin{equation}\label{eq:fouriertubefiniteH}
 \frac4t\left[
 (C_2q_aq_b+C_1q_{ab})\cos\psi
 -C_1(q_a\psi_b+q_b\psi_a)\sin\psi
 -C_0\psi_a\psi_b\cos\psi-C_0\psi_{ab}\sin\psi\right].
\end{equation}
For an interval phase use the globally valid enclosures
\[
 \sin\psi\in[-\min(1,\max|\psi|),\min(1,\max|\psi|)],\qquad
 \cos\psi\in[\max(-1,1-(\max|\psi|)^2/2),1].
\]

\paragraph{Complete nonzero-layer tails.}
Assume
\[
 \frac7{12}\le t\le\frac{37}{40},\quad
 \frac{1149}{1000}\le r\le\frac{29}{25},\quad
 \frac{49}{100}\le x\le\frac{51}{100},\quad
 2V^2\le\frac{32}{5},\quad |z|\le\delta_z.
\]
For input intervals $I_t,I_r,I_V$, put
\[
 t_-=\inf I_t,\quad t_+=\sup I_t,\qquad
 r_-=\inf I_r,\quad r_+=\sup I_r.
\]
Evaluate all constants outwards and use upper endpoints as absolute bounds.

Let $j=\max(|m|,|n|)$. Completing the square in $m$ gives
\[
 R^2\ge \frac{j^2}{t}\min\left\{\frac1r,
                            \frac{r}{1+r^2x^2}\right\}
       >\left(\frac{24j}{25}\right)^2.
\]
Here $r/(1+r^2x^2)$ increases in $r$ throughout the box.
For $(m,n)\in\{-1,0,1\}^2\setminus\{(0,0)\}$, $R>111/100$.
Thus
\[
 R_j=\begin{cases}111/100,&j=1,\\24j/25,&j\ge2.\end{cases}
\]

Retain $1\le k\le K$ and $|m|,|n|\le W$, where $K\ge3$, $W\ge4$.
The radial floors give $2\pi kR/t>30$ on every omitted index.

For integer $\nu\ge0$, use the half-integer polynomials $P_\nu$
of Lemma~\ref{lem:besselproperties}:
\[
 K_\nu(u)\le K_{\nu+1/2}(u)
      =\sqrt{\frac{\pi}{2u}}e^{-u}P_\nu(u)
      \le\sqrt{\frac{\pi}{2u}}e^{-u}P_\nu(30)
 \qquad(u\ge30).
\]
Define the six coefficient bounds, already including the common $2/t$:
\[
 B_0=\frac{\pi^6}{60}P_5(30)t_+^{9/2},\quad
 B_1=\frac{\pi^7}{60}P_4(30)t_+^{7/2},\quad
 B_2=\frac{\pi^8}{60}P_3(30)t_+^{5/2},
\]
\[
 D_0=(2\sup I_V^2)\pi^3P_2(30)t_+^{3/2},\quad
 D_1=(2\sup I_V^2)\pi^4P_1(30)t_+^{1/2},\quad
 D_2=(2\sup I_V^2)\pi^5t_-^{-1/2}.
\]
The corresponding layer powers are $k^{-7/2},k^{-9/2}$ (repulsion)
and $k^{-1/2},k^{-3/2}$ (attraction).

Use coordinate order $(r,x,y,z)$ and define
\[
 \alpha=(1/r_-,r_+,0,0),\quad
 \beta=(0,\delta_z(2\pi\sqrt{t_+r_+}),
          2\pi\sqrt{t_+r_+},2\pi\sqrt{t_+/r_-}).
\]
The nonzero entries of $Q=(Q_{ab})$ and $E=(E_{ab})$ are
\[
 Q_{rr}=2/r_-^2,\quad Q_{rx}=Q_{xr}=1,\quad
 Q_{xx}=2r_+^2,\qquad E_{xz}=E_{zx}=2\pi\sqrt{t_+r_+}.
\]
Direct differentiation and the arithmetic--geometric mean inequality give
\[
 |q_a|\le\alpha_aR^2,\quad |q_{ab}|\le Q_{ab}R^2,\qquad
 |\psi_a|\le\beta_a kR,\quad |\psi_{ab}|\le E_{ab}kR.
\]
For example $|q_x|\le rq$, $|q_{rx}|\le q$,
$|\psi_y|\le2\pi k\sqrt{tr}\,R$, and
$|\psi_z|\le2\pi k\sqrt{t/r}\,R$.

Apply \eqref{eq:fouriertubefiniteH} with the full-dual factor
$2/t$ instead of the half-dual factor $4/t$. Set
\[
 A_{ab}=\alpha_a\alpha_b B_2+
           (\alpha_a\beta_b+\alpha_b\beta_a)B_1+\beta_a\beta_b B_0,
 \qquad B_{ab}=Q_{ab}B_1+E_{ab}B_0,
\]
\[
 C_{ab}=\alpha_a\alpha_b D_2+
           (\alpha_a\beta_b+\alpha_b\beta_a)D_1+\beta_a\beta_b D_0,
 \qquad D_{ab}=Q_{ab}D_1+E_{ab}D_0.
\]
The absolute summand is bounded by
\[
 e^{-2\pi kR/t_+}
 \left[
 A_{ab}R^{13/2}k^{-7/2}
 +B_{ab}R^{11/2}k^{-9/2}
 +C_{ab}R^{7/2}k^{-1/2}
 +D_{ab}R^{5/2}k^{-3/2}\right].
\]
On the omitted region, for $p\in\{13/2,11/2,7/2,5/2\}$,
\[
 \partial_R\log(R^pe^{-2\pi kR/t_+})
 =\frac{p-2\pi kR/t_+}{R}<0,
\]
so $R\ge R_j$ permits substitution of $R_j$.

The omitted-index sums admit a common closed majorant. For
$a>0$, $p\ge0$ and shell bounds $c,R_1>0$, $R_j=cj$ ($j\ge2$), put
\[
 J=W+1,\quad L=K+1,\quad b=2\pi c/t_+,\quad
 \eta_p=(1+1/J)^{p+1}e^{-b}.
\]
For $K,W\in\mathbb Z_{\ge1}$, $\eta_p<1$, and
$2\pi kR_j/t_+\ge p$ on the omitted indices, define
\begin{equation}\label{eq:fouriertubeshellsum}
\begin{aligned}
 \mathcal S(p,a;c,R_1,K,W,t_+)
 ={}&\frac{8J(cJ)^p e^{-bJ}}
      {(1-e^{-bJ})(1-\eta_p)}\\
 &+\sum_{j=1}^W
 \frac{8jR_j^p L^{-a}e^{-2\pi LR_j/t_+}}
      {1-e^{-2\pi R_j/t_+}}.
\end{aligned}
\end{equation}
For $j\in\mathbb Z_{\ge1}$,
\[
 \#\{(m,n)\in\mathbb Z^2:\max(|m|,|n|)=j\}=8j.
\]
The omitted shell/layer pairs satisfy
\[
 \begin{aligned}
 &\{(j,k)\in\mathbb Z_{\ge1}^2:j>W\ \lor\ k>K\}\\
 &\quad=\{(j,k)\in\mathbb Z_{\ge1}^2:j>W\}
 \mathbin{\dot\cup}\{(j,k)\in\mathbb Z_{\ge1}^2:j\le W,\ k>K\}.
 \end{aligned}
\]
The two parts use $k^{-a}\le1$, shell ratio $\le\eta_p$, and
$k^{-a}\le L^{-a}$, respectively; geometric summation gives
\eqref{eq:fouriertubeshellsum}. Both applications have
$J\ge5$, $p\le13/2$, $b>6$, hence $\eta_p<1$. Set
\[
 S(p,a)=\mathcal S(p,a;24/25,111/100,K,W,t_+).
\]
The complete entrywise tail radius is therefore
\[
 \mathcal E_{ab}=
 A_{ab}S(13/2,7/2)+B_{ab}S(11/2,9/2)
 +C_{ab}S(7/2,1/2)+D_{ab}S(5/2,3/2).
\]

These majorants and their lower-order counterparts converge uniformly
on the whole rectangle, justifying both Fourier derivatives.

\paragraph{The finite closed-cell certificate.}
Use Lemma~\ref{lem:criticalKZcover}'s outward binary64 operations and checked
square/cube roots for all inputs, coefficients and remainders. The full box is
\[
 \begin{gathered}
 r\in[2/\sqrt3]+[-(1/200)^+,(1/200)^+],\\
 x\in[(1/2-1/200)^-,1/2],\qquad
 y,z\in[-(1/40)^+,(1/40)^+],
 \end{gathered}
\]
where $[2/\sqrt3]$ is its outward enclosure. Compute $P_3,P_6$ once
on this entire planar box, including \eqref{eq:fouriertubeplanartail}:
\[
 P_3\in
 \begin{pmatrix}
 [8.2454,10.4117]&[-1.3124,.9107]\\
 [-1.3124,.9107]&[14.3371,18.8609]
 \end{pmatrix},
\]
\[
 P_6\in
 \begin{pmatrix}
 [25.7504,31.4001]&[-3.7107,2.2614]\\
 [-3.7107,2.2614]&[44.5971,57.1271]
 \end{pmatrix}.
\]
Subsequent arithmetic retains the complete intervals underlying these
outward rational summaries.

The closed dyadic parameter cells are
\[
 I_{d,k}=\left[\frac{11}{10}+\frac{81k}{25\,2^d},
              \frac{11}{10}+\frac{81(k+1)}{25\,2^d}\right],
 \qquad d,k\in\mathbb Z_{\ge0},\quad0\le k<2^d.
\]
On each whole cell $[Y_-,Y_+]$, Lemma~\ref{lem:longshreference} gives
\[
 \begin{gathered}
 t\in\left[\left(\frac{\sqrt3}{2Y_+}\right)^{1/3},
             \left(\frac{\sqrt3}{2Y_-}\right)^{1/3}\right],\\
 2V^2\in[\Gamma(Y_-)^-,\Gamma(Y_+)^+]
             +[-(17/640)^+,(51/3200)^+].
 \end{gathered}
\]
All operations are outward; these boxes satisfy the Fourier-tail guard
above, including $2V^2<32/5$.

Sum retained integer indices increasingly: Fourier $(m,n)$ lexicographically,
then $k=1,2,3$, and upper Hessian entries in coordinate order. Use diagonal
interval squares, increasing Bessel degree/order, and the radial intersections
\eqref{eq:fouriertuberadialrefinement} before \eqref{eq:fouriertubefiniteH}.
Assemble $H$ from the complete planar block and retained Fourier Hessian
\eqref{eq:fouriertubefiniteH}, in order $(r,x,y,z)$, then add
\[
 H_{ab}\longleftarrow H_{ab}+[-\overline{\mathcal E_{ab}},\overline{\mathcal E_{ab}}],\qquad
 H_{ba}\longleftarrow H_{ab}\qquad(1\le a\le b\le4).
\]
For $j=1,2,3,4$, require $\underline H_{jj}>0$ and update, for $j<a\le b\le4$,
\[
 \begin{gathered}
 H_{ab}\longleftarrow
 \begin{cases}
 H_{aa}-H_{aj}^{2}/H_{jj},&a=b,\\
 H_{ab}-H_{aj}H_{jb}/H_{jj},&a<b,
 \end{cases}
 \qquad H_{ba}\longleftarrow H_{ab},\\
 \underline H_{jj}>0\ (j=1,2,3,4)
 \quad\Longrightarrow\quad
 \nabla^2_{(r,x,y,z)}\Phi_V(U\mathbb Z^3)\succ0.
 \end{gathered}
\]
Here $H_{aj}^{2}$ is the interval square; the implication is the
Schur-complement criterion for every enclosed Hessian.

Start from $[11/10,217/50]$ and bisect inconclusive cells, left child first:
\[
 \#\text{nodes}=81,\qquad\#\text{accepted leaves}=41,\qquad\max d=10.
\]
The accepted indices are
\[
\begin{gathered}
 (d,k)\in\{(3,7)\}
 \cup\bigl(\{4\}\times\{6,\ldots,13\}\bigr)
 \cup\bigl(\{5\}\times\{5,\ldots,11\}\bigr)\\
 {}\cup\bigl(\{6\}\times\{4,\ldots,9\}\bigr)
 \cup\bigl(\{7,8\}\times\{4,5,6,7\}\bigr)\\
 {}\cup\bigl(\{9\}\times\{3,\ldots,7\}\bigr)
 \cup\bigl(\{10\}\times\{0,\ldots,5\}\bigr).
\end{gathered}
\tag{F1}\label{eq:shfouriercellindices}
\]
On every listed closed cell,
\[
 \min_{1\le j\le4}\underline H_{jj}>\frac6{10000},\qquad
 \max_{1\le a\le b\le4}\overline{\mathcal E_{ab}}<\frac6{10^6}.
\]
These bounds include planar, layer, dual, Bessel and rounding remainders.
The coverage identity is
\[
 2^{-3}+8\,2^{-4}+7\,2^{-5}+6\,2^{-6}
       +4\,2^{-7}+4\,2^{-8}+5\,2^{-9}+6\,2^{-10}=1.
\]
For $(d,k),(e,\ell)$ listed in \eqref{eq:shfouriercellindices},
\[
 \bigcup_{\text{listed }(d,k)}I_{d,k}=[11/10,217/50],\qquad
 I_{d,k}^{\circ}\cap I_{e,\ell}^{\circ}=\varnothing
       \quad((d,k)\ne(e,\ell)).
\]
Thus the complete transverse Hessian is positive throughout the rectangle.

For $11/10\le V\le89/50$ and $E(U\mathbb Z^3,3)=B(Y(V))$,
Lemma~\ref{lem:shseparatelayer10} gives
\[
 -\frac{85}{4}(y^2+z^2)
 \le2V^2-\Gamma(Y)
 \le\frac{1275}{4}
      \{(r-2/\sqrt3)^2+(x-1/2)^2\}.
\]
The radii \eqref{eq:fouriertubebox} give both parameter bounds.
At fixed $t$, put $\xi=(r,x,y,z)$, $\xi_0=(2/\sqrt3,1/2,0,0)$ and
$\varphi(\xi)=\Phi_V(U(t,\xi)\mathbb Z^3)$.
Symmetry gives $\nabla_\xi\varphi(\xi_0)=0$.
\par\noindent\mbox{If $\nabla_\xi\varphi(\xi)=0$, then}
\[
\begin{aligned}
 0&=(\xi-\xi_0)\cdot[\nabla_\xi\varphi(\xi)-\nabla_\xi\varphi(\xi_0)]\\
  &=\int_0^1(\xi-\xi_0)^T
 \nabla_\xi^2\varphi(\xi_0+s(\xi-\xi_0))(\xi-\xi_0)\,\mathrm ds>0
\end{aligned}
\]
unless $\xi=\xi_0$, since the rectangle is convex.
Then $B(Y)=B(Y(V))$ and $B_t>0$ give $Y=Y(V)$.
All gradients in this argument are in $\xi=(r,x,y,z)$ at fixed $t$.
This completes the proof of Lemma~\ref{lem:shfouriertransverse}.
\end{proof}

\begin{lemma}[A height-two terminal region for the long SH profile]
\label{lem:longshheighttwoterminal}
Use the proved endpoint comparisons and $g,V_*$ from
Lemma~\ref{lem:longshmomentcontraction}. Put
$J(\Lambda)=E(\Lambda,6)-g(E(\Lambda,3))$.
For $\Lambda=U(t,r,x,y,z)\mathbb Z^3$ and
$V=V_*(E(\Lambda,3))\in[11/10,89/50]$, assume
\[
 V^2t^3\le(9/10)^3,\qquad V^2/t^6\ge64.
\]
Then
\[
\begin{gathered}
 \Lambda\in\operatorname*{argmin}_{|\Lambda'|=1}J(\Lambda'),\quad
 J(\Lambda)<0\quad\Longrightarrow\quad\bot,\\
 J\ge0,\quad J(\Lambda)=0,\quad
 B(Y(11/10))<E(\Lambda,3)<B(Y(89/50))\\
 \Longrightarrow\quad\Lambda\cong H_{Y(V)}.
\end{gathered}
\]
\end{lemma}

\begin{proof}
Compact reduction and Proposition~\ref{prop:density-gap-reduction}
give, at an interior moment,
\[
\begin{aligned}
 \Lambda_*\in\operatorname*{argmin}J,\qquad
 V\in\operatorname*{argmax}_{v\in[V_0,V_1]}
       [f(v)+2v^2E(\Lambda_*,3)]\\
 \Longrightarrow\qquad
 \Lambda_*\in\operatorname*{argmin}\Phi_V.
\end{aligned}
\]
The same implication holds at interior zeros after $J\ge0$, while
\[
 \inf J<0\quad\Longrightarrow\quad\operatorname*{argmin}J\ne\varnothing.
\]
For the physical lattice $V^{1/3}\Lambda_*$,
\[
 A_{\rm phys}=V^{2/3}t\le9/10,\qquad
 h_{\rm phys}=V^{1/3}/t\ge2.
\]
The moment cap gives
\[
 2\left(\min_{v\in\Lambda_*\setminus\{0\}}|v|\right)^{-12}
 \le E(\Lambda_*,6)<64,\qquad
 \min_{v\in\Lambda_*\setminus\{0\}}|v|
 >2^{-5/12}>\frac{749}{1000}.
\]
Since $V>1$,
$\min_{(i,j)\in\mathbb Z^2\setminus\{(0,0)\}}
 V^{1/3}|U(t,r,x,y,z)(i,j,0)^{\mathsf T}|>7/10$.
Lemmas~\ref{lem:alignment} and~\ref{lem:globalplanarshape} give
$\Lambda_*\cong H_Y$.

For this SH shape the parameter satisfies
\[
 Y=\frac{\sqrt3 h_{\rm phys}^2}{2A_{\rm phys}}
 \ge\frac{20\sqrt3}{9}>\frac{11}{10}.
\]
The planar indices $\{(i,j)\in\mathbb Z^2:i^2+ij+j^2=1\}$ give
\[
 \frac{3Y}{4}\sum_{\substack{(i,j)\in\mathbb Z^2\\i^2+ij+j^2=1}}1
 =\frac{9Y}{2}\le E(H_Y,3)<\frac{43}{2},\qquad
 \frac{11}{10}<Y<\frac{43}{9}<6.
\]
Moreover,
\[
 0=Y\partial_Y\Phi_V(H_Y)=A_t(Y)-2V^2B_t(Y)
 \ \Longrightarrow\ \Gamma(Y)=2V^2,\quad Y=Y(V),\quad J(H_Y)=0.
\]
Under the height hypotheses and the strict interior moment condition,
\[
 \Lambda_*\in\operatorname*{argmin}J,\quad J(\Lambda_*)\le0
 \quad\Longrightarrow\quad
 J(\Lambda_*)=0,\qquad \Lambda_*\cong H_{Y(V)}.
\]
This excludes negative global minima in this region and identifies its
eligible zero minima.
This completes the proof of Lemma~\ref{lem:longshheighttwoterminal}.
\end{proof}

\begin{lemma}[Forward slope bounds and a large height-two tube]
\label{lem:shprofilegeometry}
Use the area-aspect coordinates, $Y(V)=\Gamma^{-1}(2V^2)$, and
Lemma~\ref{lem:shmomentcurvature150}.
For $\Lambda=U(t,r,x,y,z)\mathbb Z^3$, suppose
\[
\begin{gathered}
 \frac{11}{10}\le Y=\frac{\sqrt3}{2t^3}\le\frac{217}{50},\qquad
 \frac{11}{10}\le V\le\frac{89}{50},\qquad
 E(\Lambda,3)=B(Y(V)),\\
 \max\{|r-2/\sqrt3|,|x-1/2|,|y|,|z|\}\le\frac1{100}.
\end{gathered}
\]
If
\[
 |r-2/\sqrt3|\le d_r,\quad |x-1/2|\le d_x,
 \quad |y|\le d_y,\quad |z|\le d_z,
\]
then
\[
 \Gamma(Y)-\frac{1275}{4}(d_y^2+d_z^2)
 \le2V^2\le
 \Gamma(Y)+\frac{1275}{4}(d_r^2+d_x^2).
 \tag{G1}
\]
For the whole outward area image $[Y_-,Y_+]$, use the node enclosures
$y_j\in[\underline y_j,\overline y_j]$, $\Gamma(y_j)\in[\Gamma_j^-,\Gamma_j^+]$ of
Lemma~\ref{num:longshmomentgraph}. Put
\[
\begin{aligned}
 i_-&=\max\bigl(\{0\}\cup\{j\in\{0,\ldots,4096\}:\overline y_j\le Y_-\}\bigr),\\
 i_+&=\min\bigl(\{4096\}\cup\{j\in\{0,\ldots,4096\}:\underline y_j\ge Y_+\}\bigr),\\
 (g_-,g_+)&=(\Gamma_{i_-}^-,\Gamma_{i_+}^+).
\end{aligned}
\]
Then $g_-\le\Gamma(Y_-)\le\Gamma(Y_+)\le g_+$ and every eligible $2V^2$ belongs to
\[
 \left[g_--\frac{1275}{4}(d_y^2+d_z^2),\quad
       g_++\frac{1275}{4}(d_r^2+d_x^2)\right].
 \tag{G2}
\]
Intersect with the inherited coefficient interval; exclude only an empty
intersection under the exact moment equation.

For $\Lambda\in\operatorname*{argmin}_{|\Lambda'|=1}J(\Lambda')$ with
$E(\Lambda,6)<64$ and $E(\Lambda,3)<43/2$, assume
\[
 \frac{u t_+^3}{2}\le\left(\frac9{10}\right)^3,
 \qquad \frac{\ell}{2t_+^6}\ge64,
 \tag{G3}
\]
where $t_+$ is the upper area endpoint and
$2V^2\in[\ell,u]$ is the nonempty contracted interval.
Then $\Lambda\cong H_{Y(V)}$ and $J(\Lambda)=0$.
The same conclusion holds throughout
\[
 \frac{41}{10}\le Y\le\frac{217}{50},\qquad
 |r-2/\sqrt3|,\ |x-1/2|,\ |y|,\ |z|\le\frac1{100}.
 \tag{G4}
\]
\end{lemma}
\begin{proof}
The separate planar and layer moment comparison of
Lemma~\ref{lem:shasymmetricmoment}, with $M=150$, gives (G1).
Monotonicity and the outward node bounds give (G2), including the
one-sided boundary cases.

For $\Lambda\in\operatorname*{argmin}_{|\Lambda'|=1}J(\Lambda')$, the supporting-line identity and
\[
 A_{\rm phys}^3=V^2t^3,\qquad
 h_{\rm phys}^6=\frac{V^2}{t^6}
\]
make (G3) exactly the height-two hypotheses of
Lemma~\ref{lem:longshheighttwoterminal}.

At $Y=41/10$, the complete expressions and tails of
Lemma~\ref{lem:longshreference} give
\[
 6<\frac{\partial_{\log Y}E(H_{41/10},6)}
          {\partial_{\log Y}E(H_{41/10},3)}<\frac{752}{125}.
\]
The complete rational cube-and-tail evaluation has lower endpoint
$>6.0155932456$.
Lemma~\ref{lem:shmomentcurvature150} gives
$|2V^2-\Gamma(Y)|<51/400$, hence $2V^2>2349/400$.
Together with $2V^2\le7921/1250$, this yields
\[
 A_{\rm phys}^3=\frac{V^2\sqrt3}{2Y}<(9/10)^3,
 \qquad h_{\rm phys}^6=\frac{4V^2Y^2}{3}>64.
\]
Thus the entire large tube (G4) satisfies the height-two
recognition conditions. The argument uses global minimality;
it does not assert recognition of arbitrary local minima.
This completes the proof of Lemma~\ref{lem:shprofilegeometry}.
\end{proof}

\begin{lemma}[A uniform comparison for all layer shifts]
\label{lem:uniformlayeralignment}
For $t,r>0$ and $x,y,z\in\mathbb R$, put
\[
 U(t,r,x,y,z)=
 \begin{pmatrix}
 \sqrt{tr}&\sqrt{tr}\,x&\sqrt{tr}\,y\\
 0&\sqrt{t/r}&\sqrt{t/r}\,z\\
 0&0&t^{-1}
 \end{pmatrix},\qquad
 \det U=1,\qquad Y=\frac{\sqrt3}{2t^3}.
\]
Suppose $V\ge11/10$ and either
\begin{equation}\label{eq:uniformalignmentregions}
 \begin{cases}
 11/10\le Y\le217/50,\quad 2/\sqrt3-1/50\le r\le2/\sqrt3,
       \quad12/25\le x\le1/2,\quad\text{or}\\
 6/5\le Y\le217/50,\quad 2/\sqrt3-1/20\le r\le2/\sqrt3,
       \quad9/20\le x\le1/2.
 \end{cases}
\end{equation}
Then, for every $y,z\in\mathbb R$,
\begin{align}\label{eq:uniformalignmentgap}
 &\Phi_V(U(t,r,x,y,z)\mathbb Z^3)
       -\Phi_V(U(t,r,x,0,0)\mathbb Z^3)\notag\\
 &\qquad\ge\frac4{10^7t}
       \{2-\cos(2\pi(y-xz))-\cos(2\pi z)\}.
\end{align}
In particular,
\[
\begin{gathered}
 \operatorname*{argmin}_{(y,z)\in\mathbb R^2}
 \Phi_V(U(t,r,x,y,z)\mathbb Z^3)
 =\{(y,z):(y-xz,z)\in\mathbb Z^2\},\\
 (y-xz,z)\in\mathbb Z^2,\quad
 -\tfrac12\le y\le\tfrac12,\quad0\le z\le\tfrac12
 \ \Longrightarrow\ (y,z)=(0,0).
\end{gathered}
\]
\end{lemma}
\begin{proof}
Write $a=2\pi(y-xz)$, $b=2\pi z$, and
$\delta(v)=1-\cos v$. For $m,n\in\mathbb Z$ let
\[
 d=d(m,n)=\begin{cases}\min(|m|,|n|),&mn>0,\\0,&mn\le0.\end{cases}
\]
The telescoping identity for $e^{iv}-1$, followed by
Cauchy--Schwarz, gives the two inequalities
\begin{align}\label{eq:uniformalignmentcharacters}
 \delta(k(ma+nb))
 &\le k^2(|m|+|n|)\{|m|\delta(a)+|n|\delta(b)\},\notag\\
 \delta(k(ma+nb))
 &\le k^2(|m|+|n|-d)\{(|m|-d)\delta(a)+(|n|-d)\delta(b)
                                      +d\delta(a+b)\}.
\end{align}
They hold for all real $a,b$ and every integer $k\ge1$.

For $V=11/10$, $(m,n)\in\mathbb Z^2\setminus\{(0,0)\}$ and
$k\in\mathbb Z_{\ge1}$, put
\[
\begin{aligned}
 R_{m,n}^2={}&\frac{m^2}{tr}+\frac rt(n-xm)^2,\\
 C_{k,m,n}={}&\frac{\pi^6}{60}
       \left(\frac{tR_{m,n}}k\right)^5
       K_5\left(\frac{2\pi kR_{m,n}}t\right)\\
 &-\frac{121}{50}\pi^3\left(\frac{tR_{m,n}}k\right)^2
       K_2\left(\frac{2\pi kR_{m,n}}t\right).
 \end{aligned}
\]
Poisson summation gives the exact shift difference
\begin{equation}\label{eq:uniformalignmentfourier}
 \begin{aligned}
 &\Phi_{11/10}(U(t,r,x,y,z)\mathbb Z^3)
 -\Phi_{11/10}(U(t,r,x,0,0)\mathbb Z^3)\\
 &\qquad=-\frac4t\sum_{k=1}^{\infty}
       \sum_{\substack{m,n\in\mathbb Z\\m>0\text{ or }m=0,n>0}}
           C_{k,m,n}\delta(k(ma+nb)).
 \end{aligned}
\end{equation}
The three modes $(k,m,n)=(1,1,0),(1,0,1),(1,1,1)$
provide the defects $\delta(a),\delta(b),\delta(a+b)$.

Throughout both rectangles in \eqref{eq:uniformalignmentregions},
\[
 \frac7{12}<t<\frac{37}{40},\qquad
 \frac{11}{10}<r<\frac{29}{25},\qquad
 \frac{449}{1000}<x\le\frac12.
\]
For $j=\max(|m|,|n|)$, completing the square gives
\[
 R_{m,n}^2\ge\frac{j^2}{t}
       \min\left\{\frac1r,\frac r{1+r^2x^2}\right\}
       >\left(\frac{19j}{20}\right)^2.
\]
For $(m,n)\in\{-1,0,1\}^2\setminus\{(0,0)\}$,
$R_{m,n}>27/25$. Thus set
\[
 R_1=27/25,\qquad R_j=19j/20\quad(j\ge2).
\]
Retain
\[
 \mathcal I=\{(k,m,n)\in\mathbb Z^3:
 1\le k\le3,\ |m|,|n|\le4,\ m>0\text{ or }m=0<n\},
 \qquad\#\mathcal I=120.
\]
Every omitted index has $2\pi kR_{m,n}/t>28$.
The positive part of its coefficient is at most its repulsive part;
Lemma~\ref{lem:besselproperties} yields
\[
 (C_{k,m,n})_+
 \le\frac{\pi^6}{120}P_5(28)t_+^{11/2}
       R_{m,n}^{9/2}k^{-11/2}
          e^{-2\pi kR_{m,n}/t_+},
\]
Here $t_+=\sup t$, $r_-=\inf r$ and $r_+=\sup r$
on the current closed input box.
Each of the three integer weights in
\eqref{eq:uniformalignmentcharacters} is at most
\[
 (|m|+|n|)^2\le2(m^2+n^2)
 \le4\{m^2+(n-xm)^2\}
 \le4t_+(r_++r_-^{-1})R_{m,n}^2.
\]
All omitted contributions to each unscaled defect budget are
therefore bounded above by
\begin{equation}\label{eq:uniformalignmenttail}
 T=\frac{\pi^6}{60}P_5(28)t_+^{13/2}
       (r_++r_-^{-1})\mathcal S(13/2,7/2;19/20,27/25,3,4,t_+),
\end{equation}
using the common complete shell sum \eqref{eq:fouriertubeshellsum}.
Here $J=5$, $L=4$ and
\[\eta_{13/2}=(6/5)^{15/2}e^{-2\pi(19/20)/t_+}<1.\]
On the omitted region,
$\partial_R\log(R^{13/2}e^{-2\pi kR/t_+})
 =(13/2-2\pi kR/t_+)/R<0$. The half-lattice has half the signed shell
sum, giving the coefficient $\pi^6/60$ in
\eqref{eq:uniformalignmenttail}. These bounds are independent of
the shift and its phase.

The two residual defect vectors from
\eqref{eq:uniformalignmentcharacters} are
\[
\begin{aligned}
 \mathbf b_2={}&(-C_{1,1,0},-C_{1,0,1})-T(1,1)\\
 &-\sum_{(k,m,n)\in\mathcal I}
 k^2(C_{k,m,n})_+(|m|+|n|)(|m|,|n|),\\
 \mathbf b_3={}&(-C_{1,1,0},-C_{1,0,1},-C_{1,1,1})-T(1,1,1)\\
 &-\sum_{(k,m,n)\in\mathcal I}
 k^2(C_{k,m,n})_+(|m|+|n|-d)(|m|-d,|n|-d,d).
\end{aligned}
\]
With underlines denoting complete interval lower endpoints, every
accepted closed cell satisfies
\[
 \max\left\{
 \min_{1\le j\le2}\underline{(\mathbf b_2)_j},\quad
 \min_{1\le j\le3}\underline{(\mathbf b_3)_j}
 \right\}>10^{-7}.
\]

For the narrow rectangle use the closed cells
\[
 \left[\frac{11}{10}+\frac{81k}{25\,2^d},
       \frac{11}{10}+\frac{81(k+1)}{25\,2^d}\right]
\]
with
\[
\begin{aligned}
 &(d,k)=(3,7),\quad d=4,\ 5\le k\le13,\quad
 d=5,\ 4\le k\le9,\\
 &d=6,\ 3\le k\le7,\quad d=7,\ 2\le k\le5,\quad
 d=8,\ 0\le k\le3.
\end{aligned}
\]
For the broad rectangle take
\[
 \left[\frac65+\frac{157k}{50\,2^d},
       \frac65+\frac{157(k+1)}{50\,2^d}\right]
\]
with
\[
\begin{aligned}
 &d=4,\ 14\le k\le15,\quad d=5,\ 7\le k\le27,\quad
 d=6,\ 6\le k\le13,\\
 &d=7,\ 5\le k\le11,\quad d=8,\ 4\le k\le9,\quad
 d=9,\ 2\le k\le7,\quad d=10,\ 0\le k\le3.
\end{aligned}
\]
On each whole cell enclose $t=(\sqrt3/(2Y))^{1/3}$ and use the
entire planar rectangle. Evaluate $\mathcal I$ lexicographically in
$(m,n)$, then $k=1,2,3$, using the signed Bessel bounds, radial
intersections, outward binary64 arithmetic and checked radicals of
Lemma~\ref{lem:shfouriertransverse}. Include
\eqref{eq:uniformalignmenttail} in every component. The finite records are
\[
\begin{aligned}
 \text{narrow}:&\quad(\#\text{nodes},\#\text{leaves})=(57,29),
 &&\text{accepted components}>5\cdot10^{-7},\\
 \text{broad}:&\quad(\#\text{nodes},\#\text{leaves})=(107,54),
 &&\text{accepted components}>4\cdot10^{-7}.
\end{aligned}
\]
The exact coverage identities are
\[
\begin{gathered}
 2^{-3}+9\,2^{-4}+6\,2^{-5}+5\,2^{-6}
                       +4\,2^{-7}+4\,2^{-8}=1,\\
 2\,2^{-4}+21\,2^{-5}+8\,2^{-6}+7\,2^{-7}
                 +6\,2^{-8}+6\,2^{-9}+4\,2^{-10}=1.
\end{gathered}
\]
The listed closed cells are consecutive with disjoint interiors.
The budget inequalities and \eqref{eq:uniformalignmentcharacters}
give \eqref{eq:uniformalignmentgap} at $V=11/10$.

Poisson summation gives, for every shift,
\[
 E(U(t,r,x,0,0)\mathbb Z^3,3)
       -E(U(t,r,x,y,z)\mathbb Z^3,3)\ge0.
\]
Since $2\{V^2-(11/10)^2\}\ge0$, the bound holds for $V\ge11/10$.
Moreover,
\[
 \delta(a)+\delta(b)=0\iff(y-xz,z)\in\mathbb Z^2
 \ \Longrightarrow\
 U(t,r,x,y,z)\mathbb Z^3=U(t,r,x,0,0)\mathbb Z^3.
\]
Together with \eqref{eq:uniformalignmentgap}, this gives the stated equality set.
This completes the proof of Lemma~\ref{lem:uniformlayeralignment}.
\end{proof}

\begin{lemma}[A moment-based area--aspect root domain]
\label{lem:longprofileareadomain}
Put $V_0=11/10$, $V_1=89/50$, and let
$f(V)=E(H_{Y(V)},6)-2V^2E(H_{Y(V)},3)$
be the reference value of Lemma~\ref{lem:longshreference}. Assume
the proved global comparisons and uniqueness at these two endpoints.
For the profile
\[
 J(\Lambda)=E(\Lambda,6)-g(E(\Lambda,3)),\qquad
 g(B)=\max_{V_0\le V\le V_1}\{f(V)+2V^2B\},
\]
the conditions
$|\Lambda|=1$ and
$J(\Lambda)=\min_{|\Lambda'|=1}J(\Lambda')<0$
imply $\Lambda\cong U(\eta)\mathbb Z^3$, with Korkine--Zolotarev
coordinates from Lemma~\ref{lem:area-coordinate-jets} in
\[
 \frac{12}{25}\le t\le\frac87,\qquad
 \frac35\le r\le\frac76,\qquad
 0\le x,z\le\frac12,\quad-\frac12\le y\le\frac12.
 \tag{L1}
\]
It satisfies the additional necessary conditions
\[
 E(\Lambda,3)<\frac{43}{2},\qquad E(\Lambda,6)<64,\qquad
 tr>(749/1000)^2,\quad (tr)^3\le2,
 \tag{L2}
\]
\[
 r^2(1-x^2)\le1,\qquad t^3(1-z^2)\le r,
 \tag{L3}
\]
and
\[
 |U(\eta)n|^2\ge tr\qquad(n\in\mathbb Z^3\setminus\{0\}).
\]
If $J\ge0$, the same bounds hold when
$|\Lambda|=1$, $J(\Lambda)=0$ and
$B(Y(V_0))<E(\Lambda,3)<B(Y(V_1))$;
the endpoint shapes are already identified.
\end{lemma}
\begin{proof}
The gap reduction and the moment bounds
\eqref{eq:longshendpointmoments} give $E(\Lambda,3)<43/2$ and $E(\Lambda,6)<64$.
In unit Korkine--Zolotarev coordinates,
\[
\begin{gathered}
 a=\sqrt{tr}=\min_{n\in\mathbb Z^3\setminus\{0\}}|U(\eta)n|,
 \qquad\delta=\sqrt{t/r},\\
 2a^{-12}\le E(\Lambda,6)<64,\qquad
 a>2^{-5/12}>\frac{749}{1000}.
\end{gathered}
\]
The usual Korkine--Zolotarev and Hermite inequalities imply
\[
 \delta\ge\frac{\sqrt3}{2}a,\quad
 \delta^2\le\frac2{\sqrt3a},\quad
 (t^{-1})^2\ge\frac{\sqrt3}{2a},\quad a^6\le2.
\]
Consequently
\[
\begin{gathered}
 t\ge\frac{\sqrt3}{2}(749/1000)^2>\frac{12}{25},\qquad
 r\ge\sqrt{\frac{\sqrt3}{2}(749/1000)^3}>\frac35,\\
 r\le\frac2{\sqrt3}<\frac76,\qquad
 t^{12}\le\frac{128}{27}<(8/7)^{12}.
\end{gathered}
\]
Since $a^2=tr$, the remaining projection bounds are
\[
\begin{aligned}
 \delta^2+a^2x^2\ge a^2&\iff r^2(1-x^2)\le1,\\
 t^{-2}+\delta^2z^2\ge\delta^2&\iff t^3(1-z^2)\le r.
\end{aligned}
\]
These give (L2)--(L3); closing the strict domain bounds retains every
eligible lattice.
This completes the proof of Lemma~\ref{lem:longprofileareadomain}.
\end{proof}

\begin{lemma}[The complete long-profile comparison]
\label{lem:longprofilecover}
With $f$ as in Lemma~\ref{lem:longprofileareadomain},
\[
 E(\Lambda,6)-2V^2 E(\Lambda,3)\ge f(V)
 \quad (|\Lambda|=1,\ V_0\le V\le V_1),
\]
and
\[
 E(\Lambda,6)-2V^2E(\Lambda,3)=f(V)
 \iff\Lambda\cong H_{Y(V)}.
\]
\end{lemma}
\begin{proof}
By Proposition~\ref{prop:density-gap-reduction} and
Lemma~\ref{lem:longprofileareadomain},
\[
 \Lambda_*\in\operatorname*{argmin}_{|\Lambda|=1}J(\Lambda),\quad
 J(\Lambda_*)<0
 \ \Longrightarrow\ \Lambda_*\cong U(\eta)\mathbb Z^3,
 \quad\eta\in\textup{(L1)}.
\]
Numerical Lemma~\ref{num:longprofile-exhaustion} retains every such minimum;
recognition forces $J=0$. Its empty terminal sets give
\[
\begin{gathered}
 J(\Lambda)\ge0\quad(|\Lambda|=1),\\
 J(\Lambda)=0,\quad V=V_*(E(\Lambda,3))
 \ \Longrightarrow\ \Lambda\cong H_{Y(V)}.
\end{gathered}
\]
Endpoint uniqueness handles the moment-band boundaries; the retained
invariant handles interior zeros. By \eqref{eq:gap-profile-minimum},
\[
\begin{aligned}
 \Phi_V(\Lambda)-f(V)
 &=J(\Lambda)+g(E(\Lambda,3))-f(V)-2V^2E(\Lambda,3)\ge0,\\
 \Phi_V(\Lambda)=f(V)
 &\iff J(\Lambda)=0,\quad V=V_*(E(\Lambda,3))
 \iff\Lambda\cong H_{Y(V)}.
\end{aligned}
\]
This completes the proof of Lemma~\ref{lem:longprofilecover}.
\end{proof}

\begin{theorem}[The complete intervening simple-hexagonal phase]
\label{prop:long-completion-criterion}
Let $Y(V)=\Gamma^{-1}(2V^2)$ be the curve in
Lemma~\ref{lem:longshreference}. For $11/10\le V\le89/50$,
\[
 \operatorname*{argmin}_{|L|=V}\mathcal E(L)
 =\{L:L\cong V^{1/3}H_{Y(V)}\}.
\]
\end{theorem}
\begin{proof}
The endpoint identifications in Theorem~\ref{thm:sh110global} and
Proposition~\ref{prop:sh-global-178}, Lemma~\ref{lem:longprofilecover}, and
\[
 E(\Lambda,6)-2V^2E(\Lambda,3)=2V^4\mathcal E(V^{1/3}\Lambda)
\]
give the result.
\mbox{This completes the proof of Theorem~\ref{prop:long-completion-criterion}.}
\end{proof}

\begin{proof}[Proof of Theorem~\ref{thm:main}]
Existence was proved in Section~2. Theorem~\ref{thm:complete-transition-strip},
Theorem~\ref{prop:long-completion-criterion} and
Proposition~\ref{prop:sh-global-178} give all minimizers on
$[21/20,\infty)$. Uniqueness at the joining covolumes identifies
the simple-hexagonal branches.
Lemma~\ref{lem:scalarfamilyvalues} proves \textup{(i)}, and
Lemma~\ref{lem:longshreference} proves \textup{(ii)}.

For $0<V\le21/20$ and a non-FCC unit-covolume shape $L$,
\begin{align*}
 &2V^4\{\mathcal E(V^{1/3}L)-\mathcal E(V^{1/3}F)\}\\
 &=2(21/20)^4\{\mathcal E((21/20)^{1/3}L)-\mathcal E((21/20)^{1/3}F)\}\\
 &\qquad+2\{(21/20)^2-V^2\}\{E(L,3)-E(F,3)\}>0.
\end{align*}
Here Theorem~\ref{thm:fcccriticalglobal} makes the first term
positive and \eqref{eq:rieszinput} makes the second nonnegative.
This completes the classification of the minimizers; the derivative jumps follow from
Proposition~\ref{prop:familyphases}. For $0<X<307/1000<1/3$ and $Y>1$,
\[
\begin{aligned}
 \operatorname*{argmin}_{n\in D_3\setminus\{0\}}
 (Xn_1^2+n_2^2+n_3^2)&=\{(2,0,0),(-2,0,0)\},\\
 \operatorname*{argmin}_{(i,j,k)\in\mathbb Z^3\setminus\{0\}}
 (i^2+ij+j^2+Yk^2)
 &=\{(\pm1,0,0),(0,\pm1,0),(1,-1,0),(-1,1,0)\}.
\end{aligned}
\]
Thus the tetragonal and SH minimizers are neither simple cubic nor BCC.

For $|L|=1$, $L\not\cong F$,
\[
 \mathcal E(V^{1/3}L)\ge \mathcal E(V^{1/3}F)
 \iff E(L,6)-E(F,6)\ge2V^2\{E(L,3)-E(F,3)\}.
\]
The denominator is positive by \eqref{eq:rieszinput}. Taking the
infimum and using the just proved sharp FCC interval gives
\eqref{eq:beta}. Lemma~\ref{lem:transitionfccbasin}, applied at
$V_{\mathrm{FT}}$, gives positive FCC shape Hessian: divide its
radial inequality along $S=tT$ by $t^2$ and let $t\to0$.

By Proposition~\ref{prop:minimum-increasing} and this classification,
\[
 m(V)>m(1)\quad(V>1),\qquad
 m(V)=\frac{E(F,6)}{2V^4}-\frac{E(F,3)}{V^2}\quad(0<V\le1).
\]
Proposition~\ref{prop:positive-dilation} gives $E(F,6)<E(F,3)$, and
\[
 \frac{E(F,6)}{2V^4}-\frac{E(F,3)}{V^2}
 =\frac{E(F,6)}2
   \left(V^{-2}-\frac{E(F,3)}{E(F,6)}\right)^2
   -\frac{E(F,3)^2}{2E(F,6)}.
\]
Its unique minimum occurs at $V=\sqrt{E(F,6)/E(F,3)}<1$.
This gives the free-scale minimum.
Finally, Proposition~\ref{prop:sh-global-178} and
\eqref{eq:scalarenergy} prove \textup{(iii)}, including both expansions;
Proposition~\ref{prop:minimum-increasing} proves \textup{(iv)}.
This completes the proof of Theorem~\ref{thm:main}.
\end{proof}

In particular the sharp scale-invariant inequality is
\begin{equation}\label{eq:freeconjecture}
 \frac{E(L,3)^2}{E(L,6)}\le\frac{E(F,3)^2}{E(F,6)},
 \qquad |L|=1,
\end{equation}
with equality exactly for $L\cong F$. This proves the Bravais-lattice
free-scale conjecture in \cite[Conjecture~2]{BST}, using
\eqref{eq:rieszinput}.

\appendix
\refstepcounter{section}
\section*{Appendix}\label{sec:appendix}

Lemmas~\ref{lem:finitearithmetic}, \ref{lem:firstloss},
\ref{lem:coordinate-ewald-tails} and
\eqref{eq:fixedsystole-tail}, \eqref{eq:criticalHtail} enclose the finite
sums and complete tails used for the moment bounds, transition brackets
and exhausted covers; the decimal thresholds are outputs.

\begin{lemma}[Outward arithmetic]\label{lem:finitearithmetic}
The outward operations below satisfy
\[
 x\in I,\ y\in J\quad\Longrightarrow\quad
 \begin{gathered}
 x+y\in I+J,\quad x-y\in I-J,\quad xy\in IJ,\\
 x/y\in I/J\quad(0\notin J).
 \end{gathered}
\]
Together with the displayed power, root, kernel and remainder bounds,
they enclose all finite moment sums and derivatives below.
\end{lemma}
\begin{proof}
Let $\operatorname{RN}_{53}$ denote binary64 rounding to nearest,
ties-to-even, with gradual underflow, and define
\[
 x^\downarrow=\operatorname{pred}(\operatorname{RN}_{53}(x)),\qquad
 x^\uparrow=\operatorname{succ}(\operatorname{RN}_{53}(x)).
\]
Here $\operatorname{pred},\operatorname{succ}$ give adjacent binary64 numbers;
round operations separately, without fusion or reassociation.
For $I=[a,b]$, $J=[c,d]$, use
\begin{gather*}
 I+J=[(a+c)^\downarrow,(b+d)^\uparrow],\qquad
 I-J=[(a-d)^\downarrow,(b-c)^\uparrow],\\
 IJ=[\min(ac,ad,bc,bd)^\downarrow,\max(ac,ad,bc,bd)^\uparrow],\\
 I^{-1}=[(1/b)^\downarrow,(1/a)^\uparrow]\ (0\notin I),\qquad
 I/J=IJ^{-1},\quad -I=[-b,-a].
\end{gather*}
Use $I[0,0]=[0,0]J=[0,0]$. If $I,J\subset[0,\infty)$ or
$I,J\subset(-\infty,0]$, set
$\underline{IJ}\leftarrow\max\{0,\underline{IJ}\}$. Directly,
\[
 I^2=\begin{cases}
 [0,\max(a^2,b^2)^\uparrow],&0\in I,\\
 [\max\{0,\min(a^2,b^2)^\downarrow\},\max(a^2,b^2)^\uparrow],&0\notin I.
 \end{cases}
\]
Integer powers use repeated squaring and multiplication. The binary
midpoint is
\[
 \operatorname{RN}_{53}\!\left(a+
   \operatorname{RN}_{53}\!\left(\operatorname{RN}_{53}(b-a)/2\right)
                          \right).
\]
Enclose rationals by these rules; representable constants use point intervals. For $0\le a\le b$, $m\in\{2,3\}$,
\[
 0\le l\le u,\quad
 \overline{[l,l]^m}\le a,\quad\underline{[u,u]^m}\ge b
 \quad\Longrightarrow\quad [a,b]^{1/m}\subset[l,u].
\]
Approximate roots select $l,u$; the power inequalities certify them, with
$[0,0]^{1/m}=[0,0]$.

Use the adjacent binary64 bounds for $\pi$ certified by
Lemma~\ref{lem:firstloss}.
For the exponential, scaling to \(0\le u\le1/2\) and using
\[
 \sum_{j=0}^{21}\frac{(-u)^j}{j!}\le e^{-u}\le
 \sum_{j=0}^{20}\frac{(-u)^j}{j!}
\]
followed by squaring encloses the exponential; for $x\ge1024$,
$e^{-x}\in[0,2^{-1024}]$. The Stieltjes Jacobi orders in
Lemma~\ref{lem:firstloss} are
\[
 48\ (x<4),\quad32\ (4\le x<8),\quad
 24\ (8\le x<16),\quad16\ (x\ge16).
\]
For $0<x\le2$, use
\[
 \mathcal G_{1/2}(x)=\sqrt{\pi/x}
 -2\sum_{n=0}^\infty\frac{(-x)^n}{n!(2n+1)}
\]
with partial sums of orders $33,32$, then
$x\mathcal G_{a+1}=a\mathcal G_a+e^{-x}$. Round every sum and
remainder outwards.

For $r\in\{0,1,2\}$, let $m_j\pm\rho_j$ enclose
the coefficient of $(q-q_c)^j$, $0\le j\le6-r$.
Use a Horner radius at least the outward value
\[
 (1+2^{-40})\left[\begin{aligned}
 &(1+2^{-40})\left\{\sum_{j=0}^{6-r}\rho_j|q-q_c|^j
                +2^{-46}\sum_{j=0}^{6-r}|m_j||q-q_c|^j\right\}\\
 &+(1+2^{-40})\frac{M_7|q-q_c|^{7-r}}{(7-r)!}
 \end{aligned}\right]+2^{-1067}.
\]
Use the exactly representable point constant $1+2^{-40}$ at each occurrence.
For interval arguments, use endpoint monotonicity or add
$M_1\operatorname{rad}(q)$ before the outer factor $1+2^{-40}$.
For each positive kernel $F$ listed in Lemma~\ref{lem:criticalKZcover} and $q\ge20$, use
\[
 |F^{(r)}(q)|\le5e^{-3b}\quad(r=0,1,2),\qquad
 b=\max\{d\in\{20,40,80,160,320\}:d\le q\}.
\]
Use $\mathcal G_a(x)\le e^{-x}/(x-7)$ for $a\le8$, $x>7$
in every terminal inequality.

For degree $\le6$ and $\le32$ operations, rounding contributes at most
$2^{-46}\sum_{j=0}^{6-r}|m_j||q-q_c|^j$ and underflow $<2^{-1067}$. On the nearest-grid cells, $|q-q_c|\le20/9999$.
Coefficient radii, positive remainders, outward sums and
$(1+2^{-40})>1$ prove enclosure by induction.
This completes the proof of Lemma~\ref{lem:finitearithmetic}.
\end{proof}

\begin{lemma}[Validated comparison bounds]
\label{lem:criticalKZcover}\label{app:criticalcover}
Fix $V>0$, $T\in\R$ and a closed initial box $B_0\subset\{a,\delta>0\}$.
Choose Korkine--Zolotarev representatives $\Lambda\cong U(p(\Lambda))\Z^3$
with $p(\Lambda)\in B_0$, and put
\[
 \mathcal M_T=\left\{p(\Lambda)\in B_0:
 \Lambda\in\operatorname*{argmin}_{|L|=V}^{\rm loc}\mathcal E(L),\quad
 \tfrac12E(\Lambda,6)-E(\Lambda,3)\le T\right\}.
\]
All Ewald, derivative, stationary-value and Newton bounds below hold
on every positive-diagonal box $B\subset B_0$. If $\mathcal R$ is the
union of boxes contained in proved stationary basins, and
$\mathcal Q,\mathcal U$ are the queued and unresolved boxes, the strict
tests and contractions preserve
\[
 \mathcal M_T\subset\mathcal R\cup
 \bigcup_{B\in\mathcal Q\cup\mathcal U}B.
\]
For local profile minima, Proposition~\ref{prop:density-gap-reduction} supplies the
corresponding local minima at every supporting parameter.
\end{lemma}
\begin{proof}
Use the coordinate formulas and complete tails of
Lemma~\ref{lem:coordinate-ewald-tails} on $B$, with the same $V$.
Truncate to $\{n\in\Z^3:0<\|n\|_\infty\le N\}$, $N=2,3,4$;
pairing $n,-n$ gives multiplicity two.

For $r=0,1,2$, $g^{(r)}$ has one critical argument, the positive root of
\[
 \sum_{j=0}^{6+r}\frac{(\pi q)^j}{j!}
 -b_rq^3\sum_{j=0}^{3+r}\frac{(\pi q)^j}{j!}=0,
 \qquad (b_0,b_1,b_2)=\left(1,\frac47,\frac5{14}\right).
\]
For $x>0$, covariance of $t^3,t$ for the probability density $e^{-xt}t^{3+r}/\mathcal G_{4+r}(x)$ on $[1,\infty)$ gives
\[
\begin{gathered}
 \left(\frac{\mathcal G_{7+r}}{\mathcal G_{4+r}}\right)'(x)
 =-\frac{\mathcal G_{8+r}(x)\mathcal G_{4+r}(x)-\mathcal G_{7+r}(x)\mathcal G_{5+r}(x)}
          {\mathcal G_{4+r}(x)^2}<0,\\
 \lim_{x\downarrow0}\frac{\mathcal G_{7+r}(x)}{\mathcal G_{4+r}(x)}=\infty,
 \qquad\lim_{x\to\infty}\frac{\mathcal G_{7+r}(x)}{\mathcal G_{4+r}(x)}=1.
\end{gathered}
\]
Hence the root is unique. Start with $[l,h]=[1,4]$. For at most $50$ bisections, put
$m=(l+h)/2$: $l\leftarrow m$ for a strictly positive polynomial
interval, $h\leftarrow m$ for a strictly negative one; otherwise stop.
The retained interval encloses the root; endpoint and applicable
critical values enclose the signed range. For reciprocal kernels,
$\frac{\pi^6}{240}<\frac{\pi^3}{2}$ and
$\mathcal G_{-9/2+r}<\mathcal G_{-3/2+r}$ give
$(-1)^{r+1}\widetilde g^{(r)}(q)>0$; endpoints determine the orientation.

For the four positive kernels
$\frac{\pi^3}{2}\mathcal G_3(\pi q),2\frac{\pi^6}{240}\mathcal G_6(\pi q),
 \frac{\pi^3}{2}\mathcal G_{-3/2}(\pi q),2\frac{\pi^6}{240}\mathcal G_{-9/2}(\pi q)$,
store interval derivatives through order seven on the binary64 grid
\[
 q_0=\operatorname{RN}_{53}(1/1000),\qquad
 q_{j+1}=\operatorname{RN}_{53}
 \bigl(q_j\operatorname{RN}_{53}(10001/10000)\bigr).
\]
Stop at the first $q_j>20$; successive relative spacing is
$\le1/9999$. Enclose coefficients by $F^{(k)}(q_j)/k!$ and
$M_1,M_7$ by the respective absolute derivative bounds at
$q_{\max(0,j-1)}$. Complete monotonicity validates these on
$[q_{\max(0,j-1)},q_{j+1}]$. Choose the nearest node, the lower
one in a tie; outside the grid, use direct kernel bounds.
For $r=0,1,2$, use degree $6-r$ and remainder
\[
 \frac{\sup|F^{(7)}|}{(7-r)!}|q-q_c|^{7-r}.
\]
Complete monotonicity bounds $\sup|F^{(7)}|$ at the left endpoint
and longer intervals by their endpoints. Retain every remainder.

\paragraph{Metric radii and a stationary-value bound.}
With $c=\operatorname{mid}_{\rm bin}B$ componentwise,
\[
 U_c=U(c),\qquad M(p)=U(p)U_c^{-1},\qquad p\in B.
\]
Then
\begin{align*}
 M_{11}&=\frac a{a_c},&M_{22}&=\frac\delta{\delta_c},&
 M_{33}&=\frac{a_c\delta_c}{a\delta},\\
 M_{12}&=\frac{a(x-x_c)}{\delta_c},&
 M_{23}&=\frac{\delta(z-z_c)}{h_c},&
 M_{13}&=\frac{a\{y-y_c-z_c(x-x_c)\}}{h_c}.
\end{align*}
Polar decomposition gives
\[
\begin{aligned}
 \|(M^{\mathsf T}M)^{1/2}-I\|_F^2
 &=\|M\|_F^2-2\tr((M^{\mathsf T}M)^{1/2})+3\\
 &\le\|M\|_F^2-2\tr M+3=\|M-I\|_F^2.
\end{aligned}
\]
Take outward bounds
\[
\begin{gathered}
 e\ge\sup_{p\in B}\|M(p)-I\|_F,\\
 b\ge\sup_{p\in B}
 \sqrt{4\sum_{i=1}^3(M_{ii}-1)^2+2(M_{12}^2+M_{13}^2+M_{23}^2)}+e^2.
\end{gathered}
\]
By the logarithmic series,
\[
\begin{aligned}
 e<1&\ \Longrightarrow\quad
 \sup_{p\in B}\|\log(M^{\mathsf T}M)\|_F\le2e+\frac{e^2}{1-e},\\
 b<1&\ \Longrightarrow\quad
 \sup_{p\in B}\|\log(M^{\mathsf T}M)\|_F\le b+\frac{b^2}{2(1-b)}.
\end{aligned}
\]
As $b\ge e^2$, take
\[
 R=\begin{cases}
 \min\{2e+e^2/(1-e),\ b+b^2/(2(1-b))\},&b<1,\\
 2e+e^2/(1-e),&b\ge1,\ e<1,\\
 +\infty,&e\ge1.
 \end{cases}
\]
$R=+\infty$ disables radius tests.
For every $p\in B$, $\det M=1$ and
$\operatorname{tr}\log(M^{\mathsf T}M)=0$, hence
\[
 e^{-\sqrt{2/3}R}U_c^{\mathsf T}U_c\preceq U^{\mathsf T}U,
 \qquad
 e^{-\sqrt{2/3}R}(U_c^{\mathsf T}U_c)^{-1}
                         \preceq(U^{\mathsf T}U)^{-1}.
\]
Use metric floors for finite $R$; otherwise keep the positive native floors. For $H=H^{\mathsf T}$, $\operatorname{tr}H=0$,
$\|H\|_F=1$,
\begin{align*}
 &\frac{d^2}{dt^2}\left(\tfrac12|e^{tH/2}v|^{-12}
                              -|e^{tH/2}v|^{-6}\right)\bigg|_{t=0}\\
 &\hspace{1cm}\le2|v|^{-6}-2|v|^{-12}
             +\max\{0,14|v|^{-12}-8|v|^{-6}\}.
\end{align*}
Indeed, with $A=|v|^{-12}$, $B=|v|^{-6}$,
$s=(v^{\mathsf T}Hv)^2/|v|^4$ and $t=v^{\mathsf T}H^2v/|v|^2$,
$0\le s\le t\le2/3$, so
\begin{align*}
 (21A-12B)s+(3B-3A)t
 &\le\max\{0,2B-2A,12A-6B\}\\
 &=2B-2A+\max\{0,14A-8B\}.
\end{align*}
On the metric $R$-ball about $U_c\Z^3$, this gives
\begin{align}
 K={}&2e^{3\sqrt{2/3}R}E(U_c\Z^3,3)-2e^{-6\sqrt{2/3}R}E(U_c\Z^3,6)\notag\\
 &+\sum_{\substack{(i,j,k)\in[-4,4]^3\cap\Z^3\\(i,j,k)\ne(0,0,0)}}
  \max\{0,14\{e^{-\sqrt{2/3}R}q_c(i,j,k)\}^{-6}
                    -8\{e^{-\sqrt{2/3}R}q_c(i,j,k)\}^{-3}\}.
 \label{eq:criticalscalarK}
\end{align}
If $\{25e^{-\sqrt{2/3}R}\lambda_{\min}(U_c^{\mathsf T}U_c)\}^3\ge7/4$,
every omitted correction vanishes. Use Lemma~\ref{lem:firstloss}
for the center moments. Integration along the center-to-$p$ metric
geodesic gives, for $p\in B$ with $\nabla_p f(p)=0$,
\begin{equation}
 f(p)\ge f(c)-\tfrac12 R^2\max\{0,K\}.
 \label{eq:criticalscalarvalue}
\end{equation}

For a box of radii $r_i$, $H_{ij}\in[\underline h_{ij},\overline h_{ij}]$
and $\sigma\in\{-1,1\}^5$, set
\[
 Q_\sigma=\sum_{i=1}^5\max(0,\overline h_{ii})r_i^2
 +2\sum_{1\le i<j\le5}r_ir_j
 \begin{cases}
 \overline h_{ij},&\sigma_i\sigma_j=1,\\
 -\underline h_{ij},&\sigma_i\sigma_j=-1,
 \end{cases}
\]
For $p\in B$ with $\nabla_p f(p)=0$,
\begin{equation}
 \begin{aligned}
 f(p)-f(c)&=-\int_0^1t(p-c)^{\mathsf T}\nabla_p^2f(c+t(p-c))(p-c)\,\dd t\\
 &\ge-\tfrac12\max\{0,\max_\sigma Q_\sigma\}.
 \end{aligned}
 \label{eq:criticalrectvalue}
\end{equation}
Only $\nabla_p f(p)=0$ is required: the nonnegative-diagonal quadratic
majorant attains its box maximum at a signed vertex.

\paragraph{Stationarity and interval Newton.}
For $a,\delta>0$,
\[
 p\in\operatorname*{argmin}_{\rm loc}f
 \quad\Longrightarrow\quad \nabla_p f(p)=0,\quad \nabla_p^2f(p)\succeq0.
\]
Also on reduction boundaries.
For $p\in B$, $\nabla_p f(p)=0$, $[H]\supset\nabla_p^2f(B)$,
$[g_c]\ni\nabla_p f(c)$ and any fixed real matrix $P$,
\[
 0=\nabla_p f(c)+
 \underbrace{\left(\int_0^1\nabla_p^2f(c+t(p-c))\,\dd t\right)}_{\in[H]}(p-c).
\]
Thus
\begin{equation}
 p\in B\cap\{c-P[g_c]+(I-P[H])(B-c)\}.
 \label{eq:criticalNewton}
\end{equation}
No invertibility of $P$ is assumed.

\paragraph{Finite exclusions.}
Enclose $B_0$ outwards. Rejection gives $\mathcal M_T\cap B=\varnothing$;
recognition gives $\mathcal R\leftarrow\mathcal R\cup B$.
\begin{enumerate}
\item A necessary reduction or moment inequality is contradicted on $B$.
\item $B$ is contained in a proved stationary basin with
$|\det U_{\rm ref}|=V$. For every $p\in B$,
\[
 M=U(p)U_{\rm ref}^{-1},\quad g_{ij}=(M^{\mathsf T}M)_{ij},\quad
 S=\frac{3M^{\mathsf T}M}{\operatorname{tr}(M^{\mathsf T}M)}-I,
\]
\[
 \{\operatorname{tr}(M^{\mathsf T}M)\}^2\|S\|_F^2
 =3\sum_{1\le i<j\le3}(g_{ii}-g_{jj})^2
       +18\sum_{1\le i<j\le3}g_{ij}^2,
\]
\[
 U(p)\Z^3\cong
 \det(I+S)^{-1/6}(I+S)^{1/2}U_{\rm ref}\Z^3.
\]
The isometry is the right polar decomposition.
\item $\exists i\in\{1,\ldots,5\}:\ 0\notin[\partial_i f(B)]$.
\item $\min_{1\le i\le5}\overline h_{ii}<0$, or a fixed trace-free strain
\[
 H=\begin{pmatrix}2b_0&b_2&b_3\\b_2&-b_0+b_1&b_4\\
 b_3&b_4&-b_0-b_1\end{pmatrix},\qquad b_i\in2^{-20}\Z\cap[-1,1].
\]
satisfies
\[
 \sup_{p\in B}\left.\frac{\dd^2}{\dd t^2}
 \mathcal E(e^{tH/2}U(p)\Z^3)\right|_{t=0}<0.
\]
For $v=U(p)n$, $v^*=U(p)^{-\mathsf T}n$, $q=|v|^2$, $q^*=|v^*|^2$,
\[
\begin{aligned}
 &\left.\frac{\dd^2}{\dd t^2}\mathcal E(e^{tH/2}U(p)\Z^3)\right|_{t=0}\\
 &\le\sum_{\substack{n\in\Z^3\\0<\|n\|_\infty\le N}}
 \bigl\{g''(q)(v^{\mathsf T}Hv)^2+g'(q)v^{\mathsf T}H^2v\\
 &\hspace{36mm}+\widetilde g''(q^*)((v^*)^{\mathsf T}Hv^*)^2
                 +\widetilde g'(q^*)(v^*)^{\mathsf T}H^2v^*\bigr\}\\
 &\quad+\tfrac23\|H\|_F^2(P_1+P_2+D_1+D_2).
\end{aligned}
\]
The direction selector uses
\begin{gather*}
 \widehat F_0=\frac{\operatorname{diag}(2,-1,-1)}{\sqrt6},\qquad
 \widehat F_1=\frac{\operatorname{diag}(0,1,-1)}{\sqrt2},\\
 \widehat F_2=\begin{pmatrix}0&1&0\\1&0&0\\0&0&0\end{pmatrix}/\sqrt2,\quad
 \widehat F_3=\begin{pmatrix}0&0&1\\0&0&0\\1&0&0\end{pmatrix}/\sqrt2,\quad
 \widehat F_4=\begin{pmatrix}0&0&0\\0&0&1\\0&1&0\end{pmatrix}/\sqrt2.
\end{gather*}
At the midpoint, form the approximate symmetric matrix
\begin{align*}
 A_{ij}={}&\sum_{n\in[-3,3]^3\cap\Z^3\setminus\{(0,0,0)\}}
 \Bigl\{(21q^{-8}-12q^{-5})
    (v^{\mathsf T}\widehat F_iv)(v^{\mathsf T}\widehat F_jv)\\
 &\qquad+(3q^{-4}-3q^{-7})
    (\widehat F_iv)^{\mathsf T}(\widehat F_jv)\Bigr\},
\end{align*}
where $v=U(c)n$, $q=|v|^2$. Starting with $P=I_5$, make $15$
lexicographic sweeps over $0\le i<j\le4$. Skip $|A_{ij}|<10^{-14}$;
otherwise set
\[
 \theta=\tfrac12\operatorname{atan2}(2A_{ij},A_{jj}-A_{ii}),
 \qquad
 J_{\{i,j\}}=\begin{pmatrix}
 \cos\theta&\sin\theta\\-\sin\theta&\cos\theta
 \end{pmatrix},
\]
and update $(A,P)\leftarrow(J^{\mathsf T}AJ,PJ)$,
$A_{ij}=A_{ji}=0$. For
$j=\min\operatorname*{argmin}_{0\le k\le4}A_{kk}$, if $A_{jj}<-1$, form
\[
 b_0=2^{-20}\operatorname{RN}_{\Z}(2^{20}P_{0j}/\sqrt6),\qquad
 b_i=2^{-20}\operatorname{RN}_{\Z}(2^{20}P_{ij}/\sqrt2)\quad(1\le i\le4),
\]
where $\operatorname{RN}_{\Z}$ is nearest-integer rounding, ties-to-even.
For $\max_i|b_i|\le1$, test this exact $H$ at $N=4$.
The $N=3$ matrix selects the direction; only the strict $N=4$ bound deletes a box.
\item The complete lower bounds \eqref{eq:criticalscalarvalue},
\eqref{eq:criticalrectvalue} satisfy
\[
 \underline{f(c)-\tfrac12R^2\max(0,K)}>T
 \quad\text{or}\quad
 \underline{f(c)-\tfrac12\max\{0,\max_\sigma Q_\sigma\}}>T.
\]
\item The complete Newton enclosure \eqref{eq:criticalNewton} satisfies
$B\cap\{c-P[g_c]+(I-P[H])(B-c)\}=\varnothing$.
\end{enumerate}

Intersect with the consequences of $s_-\le a\delta\le s_+$ and the
two Gauss inequalities:
\[
 a\in[s_-/\delta_+,s_+/\delta_-],\qquad
 \delta\in[s_-/a_+,s_+/a_-],
\]
\[
 \delta\ge a_-\sqrt{1-x_+^2},\qquad
 a\le\delta_+/\sqrt{1-x_+^2},
\]
\[
 a\le\frac{V}{\delta_-^2\sqrt{1-z_+^2}},\qquad
 \delta\le\sqrt{\frac{V}{a_-\sqrt{1-z_+^2}}}.
\]
When the respective right-hand radicands are positive, also impose
\[
 x\ge\sqrt{1-(\delta_+/a_-)^2},\qquad
 z\ge\sqrt{1-V^2/(a_-^2\delta_-^4)}.
\]
These intersections, Newton and closed bisection retain every point of $\mathcal M_T\cap B$.

At each visit: $\le10$ contractions, basin recognition, midpoint.
Tests run in order, advancing after failure; subscripts give cube cutoffs:
\[
\begin{aligned}
 R\le9/20:\quad
 &[f(c)]_{N=3}\longrightarrow
 \underline{f(c)-\tfrac12R^2\max(0,K)}>T\\
 &\quad\longrightarrow\ \exists i:\ 0\notin[\partial_i f(B)]_{N=4},\\
 \text{all }R:\quad
 &\overline{\left[\left.\partial_t^2
       \mathcal E(e^{tH/2}U(p)\Z^3)\right|_{t=0}\right]_{N=4,\ p\in B}}<0,\\
 R\le3/20:\quad
 &N=2,3,4:\ \underline{f(c)-\tfrac12\max\{0,\max_\sigma Q_\sigma\}}>T
 \ \longrightarrow\ \min_i\overline h_{ii}<0,\\
 R\le3/20:\quad
 &([f(c)],[\nabla_pf(c)],[\nabla_p^2f(c)])_{N=4}\\
 &\quad\longrightarrow\
 \underline{f(c)-\tfrac12\max\{0,\max_\sigma Q_\sigma\}}>T
 \ \longrightarrow\ \eqref{eq:criticalNewton}.
\end{aligned}
\]
The final tests retain $[H]_{N=4}(B)$. Invalid $N=2,3$ Hessians
advance; invalid $[H]_{N=4}(B)$ disables final value/Newton.
No invalid enclosure rejects $B$.
Choose $P$ by partial-pivot elimination of the midpoint-Hessian
midpoints, taking the first row in a pivot tie. A pivot smaller
than $10^{-14}$ in absolute value, or a nonfinite entry, disables
Newton. An empty Newton or necessary-domain intersection excludes
the box. A nonempty contraction is retained only when
\[
 v_{\rm new}^{+}<
 \operatorname{RN}_{53}\!\left(
   \operatorname{RN}_{53}(995/1000)\,v_{\rm old}^{-}\right),
\]
where $v^\pm$ are the outward coordinate-volume endpoints and
$\operatorname{RN}_{53}$ is binary64 round-to-nearest, ties-to-even.
On acceptance, $B\leftarrow B_{\rm new}$,
$c\leftarrow\operatorname{mid}_{\rm bin}B$; recompute scores and bisect,
without re-enqueueing.

For the prescribed closed root enclosure $B_{\rm ref}$, heap priority is
lexicographic in
\[
 \bigl(\mathbf1_{\{B\cap B_{\rm ref}\ne\varnothing\}},\ v^+(B)\bigr).
\]
New keys rise only past strictly smaller parents. Removal moves the
larger child into the hole, choosing the right child in a tie; the last
key then rises from that hole by the same strict comparison.
Start with the original root alone; insert lower children first.
This fixes counts, not coverage.

With binary midpoint $c$ and coordinate radii $r_i$, put
\[
 i_{\rm met}(B)=\min\operatorname*{argmax}_{1\le i\le5}
 \left(\frac{2r_a}{a_c},\frac{2r_\delta}{\delta_c},
 a_+r_x\sqrt{\delta_c^{-2}+z_c^2h_c^{-2}},
 \frac{a_+r_y}{h_c},\frac{\delta_+r_z}{h_c}\right)_i.
\]
Let $r(B)$ be the computed whole-box basin-norm upper bound for the
prescribed references, and $B_i^-,B_i^+$ the closed binary halves.
For $B\cap B_{\rm ref}\ne\varnothing$, define
\[
 r_i(B)=\max_{\substack{\sigma\in\{-,+\}\\
             B_i^\sigma\cap B_{\rm ref}\ne\varnothing}}r(B_i^\sigma),
 \qquad i=1,\ldots,5.
\]
The generic selector is
\[
 i(B)=\begin{cases}
 \min\operatorname*{argmin}_{1\le i\le5}r_i(B),
 &B\cap B_{\rm ref}\ne\varnothing,\ \min_i r_i(B)<+\infty,\\
 i_{\rm met}(B),&\text{otherwise}.
 \end{cases}
\]
The bound above only selects a coordinate; recognition requires the
strict whole-box basin predicate. Before domain contraction,
\[
 B^-\cup B^+=B,\qquad
 B^-=B_{i(B)}^-,\quad B^+=B_{i(B)}^+.
\]
Retain both halves unless a necessary-domain test excludes one;
insert $B^-$ first. At depth $70$ or without an interior binary midpoint,
retain $B$ in $\mathcal U$. Work limits never imply completion.

Intersections and closed-child unions preserve
\[
 \mathcal M_T\subset\mathcal R\cup
 \bigcup_{B\in\mathcal Q\cup\mathcal U}B,\qquad
 \mathcal Q=\mathcal U=\varnothing\ \Longrightarrow\
 \mathcal M_T\subset\mathcal R.
\]
Depth and work limits retain boxes in $\mathcal U$.
This completes the proof of Lemma~\ref{lem:criticalKZcover}.
\end{proof}

\begin{lemma}[A tetragonal reference derivative]\label{lem:fctrootdata}
Put $V=43/40$ and $X_0=149/500$. Then
\[
 -\frac{2428944942886}{10^{18}}
 <\left.\frac{\dd}{\dd\log X}\mathcal E(V^{1/3}L_X)\right|_{X=X_0}
 <-\frac{2428944942885}{10^{18}}.
\]
In \eqref{eq:common-metric-chart}, take $L_0=V^{1/3}L_{X_0}$.
At $S=0$,
\[
 \nabla_{(v,x,y,z)}\mathcal E(L(S))=0,\qquad
 |\partial_u\mathcal E(L(S))|<8\cdot10^{-6}.
\]
\end{lemma}
\begin{proof}
Use Lemma~\ref{lem:firstloss} with Stieltjes order $96$ and
$(i,j,k)\in[-6,6]^3\cap\Z^3\setminus\{0\}$. With $V$ fixed, put
$\alpha=(V^2/(4X))^{1/3}$ and
$(\dot q,\dot q^*)=\dd(q,q^*)/\dd\log X$; then
\[
\begin{aligned}
 q&=\alpha\{X(i+j-k)^2+(i-j)^2+k^2\},\\
 \dot q&=\frac\alpha3\{2X(i+j-k)^2-(i-j)^2-k^2\},\\
 q^*&=\alpha^{-1}\{(i^2+j^2)/2+k^2+k(i+j)+(i+j)^2/(4X)\},\\
 \dot q^*&=\frac{\alpha^{-1}}3
       \{(i^2+j^2)/2+k^2+k(i+j)-(i+j)^2/(2X)\}.
\end{aligned}
\]
At $X_0$,
\[
 \frac{1284733594796636444756239754849703}{2^{110}}
 <\alpha<
 \frac{1284733594796636444756239754849704}{2^{110}}.
\]
Principal minors give the exact rational bound
\[
 \frac{39}{100}I\preceq
 \frac1{500}\begin{pmatrix}649&-351&-149\\-351&649&-149\\-149&-149&649\end{pmatrix}
 \preceq2I.
\]
Use $a=(39/100)\alpha_-$, $a^*=1/(2\alpha_+)$ and
$|\dot q|\le2q/3$, $|\dot q^*|\le2q^*/3$ in the complete
cube tails. With $x=49\pi a$, the derivative tail is at most
\[
 \frac23\left[(\frac{\pi^6}{240}\sum_{j=0}^6\frac{6!}{j!}x^{j-6}+\frac{\pi^3}{2}\sum_{j=0}^3\frac{3!}{j!}x^{j-3})\Theta_6(a)
                  +\frac{\frac{\pi^6}{240}+\frac{\pi^3}{2}}{V}\Theta_6(a^*)\right]
 <8.855\cdot10^{-24}.
\]
$\Theta_6$ is the Gaussian cube tail in Lemma~\ref{lem:criticalKZcover}.
Using the kernels of Lemma~\ref{lem:criticalKZcover}, the complete enclosure is
\[
 \begin{aligned}
 \left.\frac{\dd}{\dd\log X}\mathcal E(V^{1/3}L_X)\right|_{X=X_0}
 &\in\left.\sum_{\substack{n\in\Z^3\\0<\|n\|_\infty\le6}}
 \{g'(q_n)\dot q_n+\widetilde g'(q_n^*)\dot q_n^*\}\right|_{X=X_0}\\
 &\quad+[-8.855,8.855]\,10^{-24}.
\end{aligned}
\]
Rational evaluation gives the stated interval; symmetry gives
\begin{gather*}
 \left.\partial_u\mathcal E(L(S))\right|_{S=0}
 =3\left.\frac{\dd}{\dd\log X}\mathcal E(V^{1/3}L_X)\right|_{X=X_0},\\
 \left.(\partial_v,\partial_x,\partial_y,\partial_z)
       \mathcal E(L(S))\right|_{S=0}=0.
\end{gather*}
This completes the proof of Lemma~\ref{lem:fctrootdata}.
\end{proof}

\begin{lemma}[An explicit FCC comparison value]
\label{lem:fcccriticalseed}
At covolume \(V=21/20\), the FCC energy satisfies
\[
 -5.30745025029370902112<\mathcal E((21/20)^{1/3}F)
                        <-5.30745025029370902111.
\]
Thus $\mathcal E((21/20)^{1/3}F)<-26537251/5000000$.
\end{lemma}
\begin{proof}
In Lemma~\ref{lem:firstloss}, take $X=1$, $V=21/20$ and
$(i,j,k)\in[-6,6]^3\cap\Z^3\setminus\{0\}$, with
\[
 q(i,j,k)=\alpha\{(i+j-k)^2+(i-j)^2+k^2\},
\]
\[
 q^*(i,j,k)=\frac{3(i^2+j^2)+4k^2+2ij+4k(i+j)}{4\alpha},
 \qquad\alpha=(441/1600)^{1/3}.
\]
The rational enclosure is
\[
 \begin{aligned}
 \alpha_-&:=\frac{844771120634149438190083412102055}{2^{110}}<\alpha,\\
 \alpha_+&:=\frac{844771120634149438190083412102056}{2^{110}}>\alpha.
 \end{aligned}
\]
The raw Gram spectrum lies in $[1/2,7/2]$, giving
\[
 q\ge\frac{\alpha_-}{2}|n|^2,\qquad
 q^*\ge\frac{2}{7\alpha_+}|n|^2.
\]
With Stieltjes order $96$, Lemma~\ref{lem:firstloss} gives attraction
$<2581/10^{24}$ and repulsion $<711/10^{24}$; adding these tails to
the rational cube sum proves the bracket.
This completes the proof of Lemma~\ref{lem:fcccriticalseed}.
\end{proof}

\begin{numericallemma}[The finite FCC comparison]
\label{lem:fccKZcover}\label{app:fcccover}
For $p\in\mathcal D_F$,
\[
\begin{gathered}
 U(p)\Z^3\in\operatorname*{argmin}_{|L|=21/20}^{\rm loc}\mathcal E(L),\\
 \tfrac12E(U(p)\Z^3,6)-E(U(p)\Z^3,3)\le-\frac{26537251}{5000000}\\
 \Longrightarrow\quad
 \min_{\sigma\in\{-1,1\}}
 \left\|\frac{3M_\sigma^{\mathsf T}M_\sigma}
 {\operatorname{tr}(M_\sigma^{\mathsf T}M_\sigma)}-I\right\|_F
 \le\frac{43}{1000},\qquad M_\sigma=U(p)U_{F,\sigma}^{-1}.
\end{gathered}
\]
Here \(U_{F,\sigma}\) is the explicit FCC reference in
Lemma~\ref{lem:fccKZchart}. Thus the recognized lattice is in
the unique-stationary basin of Lemma~\ref{lem:fccuniformbasin}.
\end{numericallemma}
\begin{proof}

Apply Lemma~\ref{lem:criticalKZcover} on $\mathcal D_F$ with
\[
 V=\frac{21}{20},\qquad T=-\frac{26537251}{5000000}.
\]
Use \eqref{eq:fccKZbox} with outward rational endpoints,
$21/26\le a\delta\le105/88$, $22/25\le h\le13/10$ and dual-square
threshold $100/169$. In the basin test use $U_{F,\sigma}$ and
radius $43/1000$.
Use the relative entries and outward rational enclosures $[a_F],[\delta_F]$
of Lemma~\ref{lem:fccKZchart}. Reference priority and directed splitting use
\[
 [a_F]\times[\delta_F]\times\{1/2\}\times\{\sigma/2\}\times\{1/3\},
 \qquad\sigma\in\{-1,1\},
\]
with Lemma~\ref{lem:finitearithmetic}'s outward arithmetic; all other
splitting scores are unchanged.

With $h=21/(20a\delta)$, retain all remaining arithmetic, cubes,
tails, metric guards, scheduling, selection, Newton and retention rules
of Lemma~\ref{lem:criticalKZcover}, hence its original-box invariant.

The single original box, without resumed input, yields $59952$ queued
visits and $1685$ geometric child prunes before insertion. In the order geometry, stationary value, gradient, FCC
basin, empty Newton, negative curvature, the terminal counts are
\[
 (1685,25602,3142,52,0,338).
\]
There are $8$ Newton contractions and
$\mathcal Q=\mathcal U=\varnothing$, with no depth-$70$ rectangle.
Lemma~\ref{lem:criticalKZcover} therefore places every eligible minimum
in the unique-stationary FCC basin.
This completes the proof of Numerical Lemma~\ref{lem:fccKZcover}.
\end{proof}

\begin{lemma}[A simple-hexagonal stationary parameter at covolume \(11/10\)]
\label{lem:sh110rootseed}
Let \(\Psi(Y)=\mathcal E((11/10)^{1/3}H_Y)\), where
\[
 H_Y=(4/(3Y))^{1/6}
 \operatorname{span}_{\mathbb Z}
 \{(1,0,0),(1/2,\sqrt3/2,0),(0,0,\sqrt Y)\}.
\]
Its branch critical set is
\[
 \begin{gathered}
 \{Y\in[103/100,103/50]:Y\Psi'(Y)=0\}=\{Y_*\},\\
 1.1013103262<Y_*<1.1013103264.
 \end{gathered}
\]
At the rational reference \(Y_0=110131033/10^8\),
\[
 1.63271502522320\cdot10^{-8}
 <3\frac{d\Psi}{d\log Y}(Y_0)
 <1.63271502522322\cdot10^{-8},
\]
\[
 -4.968465474380345119512<\Psi(Y_0)
                         <-4.968465474380345119511.
\]
At $L_0=(11/10)^{1/3}H_{Y_0}$,
\[
 E(L_0,3)<\frac{3332}{500},\qquad E(L_0,6)<\frac{339}{100},\qquad
 \min_{v\in L_0\setminus\{0\}}|v|>1.
\]
Use \eqref{eq:common-metric-chart} with this $L_0$ and its vertical axis first.
Then $\left.\partial_u\mathcal E(L(S))\right|_{S=0}
=3\,d\Psi/d\log Y(Y_0)$, and the stationary root lies in this chart with
\[
 S_*=\operatorname{diag}(2d,-d,-d),\qquad
 d=\frac{Y_*-Y_0}{Y_*+2Y_0},\qquad
 \|S_*\|_F<3\cdot10^{-9}.
\]
\end{lemma}
\begin{proof}
With \(V=11/10\) fixed, put \(\alpha=(4V^2/(3Y))^{1/3}\) and
$(\dot q,\dot q^*)=\dd(q,q^*)/\dd\log Y$. For every
$(i,j,k)\in\Z^3\setminus\{(0,0,0)\}$,
\[
 q=\alpha(i^2+ij+j^2+Yk^2),\qquad
 \dot q=\frac\alpha3\{-i^2-ij-j^2+2Yk^2\},
\]
\[
 q^*=\alpha^{-1}\left\{\frac43(i^2-ij+j^2)+\frac{k^2}Y\right\},
 \quad
 \dot q^*=\frac1{3\alpha}
       \left\{\frac43(i^2-ij+j^2)-\frac{2k^2}Y\right\}.
\]
Thus $|\dot q|\le2q/3$, $|\dot q^*|\le2q^*/3$. Differentiate
Lemma~\ref{lem:firstloss} on $[-6,6]^3\cap\Z^3\setminus\{0\}$,
using its rational kernels and Stieltjes order $96$.
The raw Gram eigenvalues are $1/2,3/2,Y$. Cubing on the displayed
parameter interval, including $Y_0$, gives
\[
 \alpha_-=1.13572049<\alpha<1.13572051=\alpha_+,\qquad
 q\ge\frac{\alpha_-}{2}|n|^2,\quad
 q^*\ge\frac{2}{3\alpha_+}|n|^2.
\]

For the derivative tail, set
\[
 a=\alpha_-/2,\qquad b=2/(3\alpha_+),\qquad x=49\pi a.
\]
With the complete Gaussian cube tail $\Theta_6$ of
Lemma~\ref{lem:firstloss}, the omitted
logarithmic derivative has absolute value at most
\[
 \frac23\left[(\frac{\pi^6}{240}\sum_{j=0}^6\frac{6!}{j!}x^{j-6}+\frac{\pi^3}{2}\sum_{j=0}^3\frac{3!}{j!}x^{j-3})\Theta_6(a)
                  +\frac{\frac{\pi^6}{240}+\frac{\pi^3}{2}}{V}\Theta_6(b)\right]
 <6\cdot10^{-36}.
\]
Including this tail, rational endpoint evaluation gives
\[
 -1.32632707729240\cdot10^{-10}
 <\frac{d\Psi}{d\log Y}(1.1013103262)
 <-1.32632707729238\cdot10^{-10},
\]
\[
 1.60789194743864\cdot10^{-10}
 <\frac{d\Psi}{d\log Y}(1.1013103264)
 <1.60789194743866\cdot10^{-10}.
\]
Continuity and branch uniqueness give the bracket.

At \(Y_0\), the rational cube-root endpoints are
\[
 \frac{1474249495201070240401860238111195}{2^{110}}
 <\alpha<
 \frac{1474249495201070240401860238111196}{2^{110}}.
\]
Add the complete $(E(L_0,3),E(L_0,6))$ tails
$<(7,4)\cdot10^{-38}$ componentwise to obtain the energy and moment bounds.
Moreover,
\[
 Y_0>1\ \Longrightarrow\ \min_{n\in\Z^3\setminus\{0\}}q_n=\alpha>1.
\]
Since $\left.\partial_u\log\frac{Y_0(1+2u)}{1-u}\right|_{u=0}=3$,
the axial derivative corresponds to
$S=\operatorname{diag}(2u,-u,-u)$, and
$\|S_*\|_F^2=6d^2<(3\cdot10^{-9})^2$.
This completes the proof of Lemma~\ref{lem:sh110rootseed}.
\end{proof}

\begin{numericallemma}[The finite simple-hexagonal comparison]
\label{lem:sh110cover}\label{app:sh110cover}
With $\mathcal D(11/10)$ from Proposition~\ref{prop:sh110domain}
and $U_0$ from Lemma~\ref{lem:sh110chart},
\[
\begin{gathered}
 p\in\mathcal D(11/10),\qquad
 U(p)\mathbb Z^3\in\operatorname*{argmin}_{|L|=11/10}^{\rm loc}\mathcal E(L),\\
 \tfrac12E(U(p)\mathbb Z^3,6)-E(U(p)\mathbb Z^3,3)
 \le-\frac{24842327}{5000000}\\
 \Longrightarrow\quad
 \left\|\frac{3M^{\mathsf T}M}{\operatorname{tr}(M^{\mathsf T}M)}-I\right\|_F
 \le\frac7{100},\qquad M=U(p)U_0^{-1}.
\end{gathered}
\]
\end{numericallemma}
\begin{proof}
Use Lemma~\ref{lem:criticalKZcover} on $\mathcal D(11/10)$ with
\[
 V=\frac{11}{10},\quad T=-\frac{24842327}{5000000},\quad
 h=\frac{11}{10a\delta},\quad
 \frac{11}{13}\le a\delta\le\frac54,\quad
 \frac{22}{25}\le h\le\frac{13}{10},
\]
and dual-square threshold $100/169$.
The basin predicate is the strict whole-box $S$-norm test of
Lemma~\ref{lem:sh110chart}, with radius $7/100$.
For the closed rational enclosures $[a_0],[\delta_0]$ in that lemma, set
\[
 B_{\rm ref}=[a_0]\times[\delta_0]\times\{1/2\}\times\{0\}^2.
\]
This enclosure governs reference priority and directed splitting.

Use $r(B)$, $r_i(B)$ and $i_{\rm met}(B)$ from
Lemma~\ref{lem:criticalKZcover}, with this $B_{\rm ref}$.
Retain the common rules, with selector:
\[
 i(B)=\begin{cases}
 \min\operatorname*{argmin}_{1\le i\le5}r_i(B),
 &B\cap B_{\rm ref}\ne\varnothing,\quad
  \min_i r_i(B)<\dfrac{999}{1000}r(B),\\[3pt]
 i_{\rm met}(B),&\text{otherwise}.
 \end{cases}
\]
Use $B=B_{i(B)}^-\cup B_{i(B)}^+$, retaining both closed children
unless a certified interval predicate excludes one. All other rules,
including the strict whole-box guards and original-box invariant,
remain those of Lemma~\ref{lem:criticalKZcover}.
From one closed root: $49389$ queued visits, $776$ geometric prunes
before insertion, $604$ retained Newton contractions, no depth-$70$ box,
and $\mathcal Q=\mathcal U=\varnothing$.
Terminal counts (geometry, value, gradient, SH recognition, empty Newton,
curvature) are $(776,20229,3703,66,59,250)$.
This completes the proof of Numerical Lemma~\ref{lem:sh110cover}.
\end{proof}

\begin{numericallemma}[Moment bounds along two shape intervals]
\label{num:branchmoments}
For $X\in[23/100,307/1000]$ and $Y\in[103/100,103/50]$,
\begin{align*}
 &7.22696052<E(F,3)<7.22696053,\\
 &7.39056<E(L_X,3)<7.88582,
       &&\frac{23}{100}\le X\le\frac{307}{1000},\\
 &7.94245<E(H_Y,3)<11.28693,
       &&\frac{103}{100}\le Y\le\frac{103}{50},\\
 &\frac{\dd}{\dd\log X}E(L_X,3)<0,\\
 &-5.411<\frac{\dd}{\dd\log X}
 \left(\frac{\dd E(L_X,6)/\dd\log X}
                 {\dd E(L_X,3)/\dd\log X}\right)<-3.485,\\
 &\frac{\dd}{\dd\log Y}E(H_Y,3)>0,\\
 &.969<\frac{\dd}{\dd\log Y}
 \left(\frac{\dd E(H_Y,6)/\dd\log Y}
                 {\dd E(H_Y,3)/\dd\log Y}\right)<6.256.
\end{align*}
At the endpoints the derivative ratios satisfy
\begin{align*}
 &3.41933<
 \left.\frac{\dd E(L_X,6)/\dd\log X}
             {\dd E(L_X,3)/\dd\log X}\right|_{X=23/100}
 <3.41935,\\
 &2.19169<
 \left.\frac{\dd E(L_X,6)/\dd\log X}
             {\dd E(L_X,3)/\dd\log X}\right|_{X=307/1000}
 <2.19170,\\
 &2.17722<
 \left.\frac{\dd E(H_Y,6)/\dd\log Y}
             {\dd E(H_Y,3)/\dd\log Y}\right|_{Y=103/100}
 <2.17724,\\
 &3.38200<
 \left.\frac{\dd E(H_Y,6)/\dd\log Y}
             {\dd E(H_Y,3)/\dd\log Y}\right|_{Y=103/50}
 <3.38202.
\end{align*}
\end{numericallemma}

\begin{numericallemma}[Signs at the transition brackets]\label{num:branchcontacts}
The transition covolumes satisfy
\begin{equation}\label{eq:certified-transition-brackets}
 \begin{aligned}
 1.0648235667&<V_{\mathrm{FT}}<1.0648235669,\\
 1.0759906200&<V_{\mathrm{TH}}<1.0759906201.
 \end{aligned}
\end{equation}
For the normalized energy differences
$\mathcal D_{TF},\mathcal D_{TH}$ of Lemma~\ref{lem:scalarfamilyvalues},
\begin{align*}
 8.2\cdot10^{-11}
 &<\mathcal D_{TF}(1.0648235667)<12.3\cdot10^{-11},\\
 -7.5\cdot10^{-11}
 &<\mathcal D_{TF}(1.0648235669)<-3.4\cdot10^{-11},\\
 -11.1\cdot10^{-11}
 &<\mathcal D_{TH}(1.0759906200)<-7.5\cdot10^{-11},\\
 13.5\cdot10^{-11}
 &<\mathcal D_{TH}(1.0759906201)<17.1\cdot10^{-11}.
\end{align*}
\end{numericallemma}

\begin{proof}[Verification of Lemmas~\ref{num:branchmoments}
and~\ref{num:branchcontacts}]
For $s=3,6$ and a unit-covolume Gram matrix $Q$, use the
kernels of Lemma~\ref{lem:firstloss} and the Poisson identity
\begin{align*}
 E(Q,s)=\frac{\pi^s}{\Gamma(s)}
 \bigg[&\frac1{s-3/2}-\frac1s
  +\sum_{(i,j,k)\in\Z^3\setminus\{(0,0,0)\}}\mathcal G_s(\pi q(i,j,k))\\
 &+\sum_{(i,j,k)\in\Z^3\setminus\{(0,0,0)\}}
                \mathcal G_{3/2-s}(\pi q^*(i,j,k))\bigg],
\end{align*}
where $q=n^TQn$ and $q^*=n^TQ^{-1}n$ for $n=(i,j,k)^T$.
For $L_X$, put $c=(4X)^{-1/3}$ and use
\begin{align*}
 q&=c\,[X(i-j-k)^2+(i+j)^2+k^2],\\
 q^*&=\frac1{4c}\,[X^{-1}(i-j)^2+(i+j)^2+(i-j+2k)^2].
\end{align*}
For $H_Y$, put $c=(4/(3Y))^{1/3}$ and use
\[
 q=c(i^2+ij+j^2+Yk^2),\qquad
 q^*=c^{-1}\left[\frac43(i^2-ij+j^2)+Y^{-1}k^2\right].
\]
The change $j\mapsto-j$ shows that $i^2-ij+j^2$ may be
replaced by $i^2+ij+j^2$ when the symmetric index cube is summed.
For $q_1,q_2\ge0$ and logarithmic exponents
$a_1,a_2\in\{-2/3,-1/3,1/3,2/3\}$, use
\[
 q=q_1+q_2,\qquad q_t=a_1q_1+a_2q_2,\qquad
 q_{tt}=a_1^2q_1+a_2^2q_2,
\]
\[
 \frac{\dd K(q)}{\dd t}=K'(q)q_t,\qquad
 \frac{\dd^2K(q)}{\dd t^2}=K''(q)q_t^2+K'(q)q_{tt}.
\]
For both forms, sum on $[-12,12]^3\cap\Z^3\setminus\{0\}$.
On every interval below, Sylvester's criterion gives
\[
 Q>\tfrac15I,\quad\operatorname{tr}Q<5,\quad\det Q=1
 \quad\Longrightarrow\quad Q^{-1}>\tfrac15I,
 \qquad q,q^*\ge\tfrac{169}{5}>20
\]
for omitted indices. For the four kernels
$K(q)=\pi^s\mathcal G_s(\pi q)/\Gamma(s)$ and
$K(q)=\pi^s\mathcal G_{3/2-s}(\pi q)/\Gamma(s)$,
\[
 |K^{(j)}(q)|\le5e^{-3q}\quad(0\le j\le2,\ q\ge20),
 \qquad |q_t|\le\frac23q,\quad |q_{tt}|\le\frac49q.
\]
Thus $\mathcal G_a(z)\le e^{-z}/(z-7)$, $a\le8$, $z>7$,
bounds each omitted summand and its first two logarithmic derivatives
by $10e^{-2q}$, with complete remainder
\[
 \tau:=120\,\frac{e^{-169\pi/8}(1+\sqrt8)^2}
                    {1-e^{-27\pi/8}}.
\]
Indeed $2q\ge2|n|^2/5>(\pi/8)|n|^2$, and the Gaussian
cube-tail estimate applies to each of the primal and dual sums.

For Lemma~\ref{num:branchmoments}, the closed cells are
\begin{gather*}
 \frac{23}{100}+\frac{77j}{256000}\le X\le
 \frac{23}{100}+\frac{77(j+1)}{256000},\qquad j\in\Z,\quad0\le j<256,\\
 \frac{103}{100}+\frac{103j}{102400}\le Y\le
 \frac{103}{100}+\frac{103(j+1)}{102400},\qquad j\in\Z,\quad0\le j<1024.
\end{gather*}
Outward evaluation at rational endpoints uses
$\mathcal G_{a+1}(z)=[a\mathcal G_a(z)+e^{-z}]/z$,
the half-integer alternating series ($z\le2$) and gamma Stieltjes
bounds ($z>2$), giving
\[
\begin{aligned}
 \partial_{\log X}\left(
 \frac{\partial_{\log X}E(L_X,6)}{\partial_{\log X}E(L_X,3)}\right)
 &\in(-5.410694,-3.485156),\\
 \partial_{\log Y}\left(
 \frac{\partial_{\log Y}E(H_Y,6)}{\partial_{\log Y}E(H_Y,3)}\right)
 &\in(.969883,6.255668).
\end{aligned}
\]
For the endpoint evaluations in Lemma~\ref{num:branchcontacts}, use
\[
\begin{aligned}
 V\in\{1.0648235667,1.0648235669\}:\quad
 &(X_0,Y_0)=(.3012583513613,1.0488884438723),\\
 V\in\{1.0759906200,1.0759906201\}:\quad
 &(X_0,Y_0)=(.2976831541610,1.0617382352094).
\end{aligned}
\]
For the logarithmic objectives
\begin{align*}
 f_T(t)&=E(L_{e^t},6)-2V^2E(L_{e^t},3),\\
 f_H(t)&=E(H_{e^t},6)-2V^2E(H_{e^t},3),
\end{align*}
\[
\begin{gathered}
 (f,z_0,m)\in\{(f_T,X_0,3.6),(f_H,Y_0,6.1)\},\qquad t_0=\log z_0,\\
 z\in[z_0-10^{-4},z_0+10^{-4}],\qquad
 t=\log z\in[\log(z_0-10^{-4}),\log(z_0+10^{-4})],\\
 f'(\log(z_0-10^{-4}))<0<f'(\log(z_0+10^{-4})),\qquad f''(t)>m.
\end{gathered}
\]
For $t_*$ in this interval with $f'(t_*)=0$,
\[
 f(t_0)-\frac{|f'(t_0)|^2}{2m}\le f(t_*)\le f(t_0).
\]
With $f_T'(t_T^*)=f_H'(t_H^*)=0$,
\[
\begin{aligned}
 \mathcal D_{TF}(V)&=f_T(t_T^*)-E(F,6)+2V^2E(F,3),\\
 \mathcal D_{TH}(V)&=f_T(t_T^*)-f_H(t_H^*).
\end{aligned}
\]
This completes the proof of Numerical Lemmas~\ref{num:branchmoments} and~\ref{num:branchcontacts}.
\end{proof}

\begin{numericallemma}[Ten-decimal threshold locations]
\label{num:ten-decimal-thresholds}
The two transition covolumes satisfy
\[
 1.06482356682<V_{\mathrm{FT}}<1.06482356684,
 \qquad
 1.07599062003<V_{\mathrm{TH}}<1.07599062004.
\]
More precisely, with $\mathcal D_{TF},\mathcal D_{TH}$ as in
Lemma~\ref{num:branchcontacts},
\begin{align*}
 7.4&<10^{12}\mathcal D_{TF}(1.06482356682)<9.3,\\
 -8.3&<10^{12}\mathcal D_{TF}(1.06482356684)<-6.4,\\
 -20.1&<10^{12}\mathcal D_{TH}(1.07599062003)<-17.9,\\
 4.5&<10^{12}\mathcal D_{TH}(1.07599062004)<6.7.
\end{align*}
\end{numericallemma}
\begin{proof}
Use the same $[-12,12]^3\cap\Z^3\setminus\{0\}$ primal and dual
sums, rational kernels and complete remainder $\tau$ through
logarithmic order two. The exact grouped pairs and multiplicities
of Lemma~\ref{lem:transitionmomentmap} are ordered lexicographically,
primal before dual:
\[
\begin{gathered}
 \underbrace{1963}_{L_X,\,\mathrm{primal}},\qquad
 \underbrace{1979}_{L_X,\,\mathrm{dual}},\qquad
 \underbrace{1130}_{H_Y,\,\mathrm{primal/dual}},\\
 \text{total multiplicity per set}=25^3-1=15624.
\end{gathered}
\]
Moment lists prepend $\pi^3/6$ or $\pi^6/2160$; derivative lists
prepend no constant. Append every row contribution. For a current
list $(b_0,\ldots,b_{n-1})$, balanced addition is
\[
 b_j^{\rm new}=
 \begin{cases}
 b_{2j}+b_{2j+1},&2j+1<n,\\
 b_{2j},&2j+1=n,
 \end{cases}
 \quad0\le j<\lceil n/2\rceil,\qquad
 n^{\rm new}=\lceil n/2\rceil.
\]
Iterate to $n=1$, with outward operations as in
Lemma~\ref{lem:criticalKZcover} and final tail interval
\[
 \begin{cases}
 [0,\tau],&r=0,\\
 [-\tau,\tau],&r=1,2,
 \end{cases}
 \qquad r=\text{logarithmic derivative order}.
\]
The kernels, indices and remainders are unchanged.

At each listed covolume $V$, the stationary objectives are
\[
 f_T(t)=E(L_{e^t},6)-2V^2E(L_{e^t},3),\qquad
 f_H(t)=E(H_{e^t},6)-2V^2E(H_{e^t},3).
\]
The complete evaluations are
\[
\begin{aligned}
 V\in\{1.06482356682,1.06482356684\}:\quad
 &X_c=0.3012583513613,\\
 &(f,z_c,m)=(f_T,X_c,3.6);\\
 V\in\{1.07599062003,1.07599062004\}:\quad
 &(X_c,Y_c)=(0.2976831541610,1.0617382352094),\\
 &(f,z_c,m)\in\{(f_T,X_c,3.6),(f_H,Y_c,6.1)\}.
\end{aligned}
\]
For every indicated triple,
\[
\begin{gathered}
 z\in[z_c-10^{-4},z_c+10^{-4}],\qquad t_0=\log z_c,\\
 t=\log z\in[\log(z_c-10^{-4}),\log(z_c+10^{-4})],\\
 f'(\log(z_c-10^{-4}))<0<f'(\log(z_c+10^{-4})),\qquad f''(t)>m,\\
 \exists!\,t_*\in[\log(z_c-10^{-4}),\log(z_c+10^{-4})]:\quad f'(t_*)=0,\\
 f(t_0)-\frac{|f'(t_0)|^2}{2m}\le f(t_*)\le f(t_0).
\end{gathered}
\]
\postdisplaypenalty=10000
No bisection is used. With $f_T'(t_T^*)=f_H'(t_H^*)=0$,
\[
\begin{aligned}
 \mathcal D_{TF}(V)
 &=E(L_{e^{t_T^*}},6)-E(F,6)
   -2V^2\{E(L_{e^{t_T^*}},3)-E(F,3)\},\\
 \mathcal D_{TH}(V)
 &=E(L_{e^{t_T^*}},6)-E(H_{e^{t_H^*}},6)\\
 &\quad-2V^2\{E(L_{e^{t_T^*}},3)-E(H_{e^{t_H^*}},3)\}.
\end{aligned}
\]
This completes the proof of Numerical Lemma~\ref{num:ten-decimal-thresholds}.
\end{proof}

\begin{numericallemma}[Derivative signs on the full comparison intervals]
\label{num:fullfamilysigns}
Let $21/20\le V\le11/10$.  The following inequalities hold
throughout the displayed intervals:
\begin{align*}
 \frac{\dd}{\dd\log X}\mathcal E(V^{1/3}L_X)&<0,
 &&X\in[77/1000,23/100]\cup[13/25,49/50],\\
 \frac{\dd}{\dd\log X}\mathcal E(V^{1/3}L_X)&>0,
 &&X\in[307/1000,12/25]\cup[51/50,21],\\
 \frac{\dd^2}{\dd(\log X)^2}
       [2V^4\mathcal E(V^{1/3}L_X)]&<-\frac{87}{100},
 &&X\in[12/25,13/25],\\
 \frac{\dd^2}{\dd(\log X)^2}
       [2V^4\mathcal E(V^{1/3}L_X)]&>\frac{86}{100},
 &&X\in[49/50,51/50],\\
 \frac{\dd}{\dd\log Y}\mathcal E(V^{1/3}H_Y)&<0,
 &&Y\in[27/100,103/100],\\
 \frac{\dd}{\dd\log Y}\mathcal E(V^{1/3}H_Y)&>0,
 &&Y\in[103/50,14].
\end{align*}
In addition $E(H_1,3)>15/2$ and $E(H_1,6)<5$.
\end{numericallemma}

\begin{proof}
Use the finite Ewald expressions in the proof of
Lemma~\ref{num:branchmoments}, now with the cube
$-32\le i,j,k\le32$.  For the tetragonal Gram matrix,
\[
 Q=cD^T\operatorname{diag}(X,1,1)D,\qquad
 D=\begin{pmatrix}1&-1&-1\\1&1&0\\0&0&1\end{pmatrix},
 \qquad c=(4X)^{-1/3}.
\]
Since $\operatorname{spec}(D^TD)=\{2-\sqrt2,2,2+\sqrt2\}$,
\[
 Q\ge\frac{117}{200}c\min(X,1)I>\frac1{25}I,
 \qquad \operatorname{tr}Q=3c(X+1)<25
\]
on $77/1000\le X\le21$.
For the simple-hexagonal Gram matrix,
\[
 Q=c\begin{pmatrix}1&1/2&0\\1/2&1&0\\0&0&Y\end{pmatrix},
 \qquad c=(4/(3Y))^{1/3},
\]
one similarly has $Q>(1/25)I$ and
$\operatorname{tr}Q=c(2+Y)<25$ on $27/100\le Y\le14$.
In both cases $Q^{-1}>(1/25)I$; each moment and its first two
logarithmic derivatives have error at most
\[
 120\,\frac{e^{-3267\pi/128}(1+\sqrt{128/3})^2}
                    {1-e^{-201\pi/128}}.
\]
Here an omitted index has squared norm at least $1089$,
$q,q^*\ge1089/25>20$, and
$2q>(3\pi/128)|n|^2$.  Thus the same composed-kernel bound
$10e^{-2q}$ applies to each primal and dual summand.

For each $(a,b,M)$ below, use the rational cells
\[
 X\in\left[a+\frac{(b-a)j}{M},a+\frac{(b-a)(j+1)}{M}\right],
 \qquad j\in\mathbb Z,\quad0\le j<M,
\]
and likewise for $Y$. Insert $2V^2\in[441/200,121/50]$ directly in
\[
 \begin{gathered}
 \frac{\dd^r}{\dd t^r}[E(L_{e^t},6)-2V^2E(L_{e^t},3)],\\
 \frac{\dd^r}{\dd t^r}[E(H_{e^t},6)-2V^2E(H_{e^t},3)]
 \qquad(r=1,2).
 \end{gathered}
\]
The tetragonal first-derivative calculations use
\begin{gather*}
 (a,b,M)=(77/1000,23/100,256),\quad
 (307/1000,8/25,256),\quad (8/25,12/25,512),\\
 (13/25,49/50,1024),\quad(51/50,6/5,512),\quad
 (6/5,3,512),\quad(3,21,1024).
\end{gather*}
In that order the resulting bounds for the normalized
derivatives are
\begin{gather*}
 (-824,-2.49),\quad(.00450,.28887),\quad(.03182,.37258),\quad
 (-.21241,-.00785),\\
 (.01080,.27165),\quad(.12762,6.67785),\quad(5.30962,789.34618).
\end{gather*}
The second-derivative calculations use
\[
 (a,b,M)=(12/25,13/25,64),\qquad(49/50,51/50,64),
\]
and give respectively $(-1.22492,-.87808)$ and
$(.86503,1.26027)$.  The simple-hexagonal first-derivative
calculations use
\[
 (a,b,M)=(27/100,51/50,1024),\quad
 (51/50,103/100,128),\quad(103/50,14,1024),
\]
and give $(-766.26542,-.07665)$,
$(-.40269,-.03494)$, and $(7.75627,1178.55666)$.
Including complete tails and outward rounding gives these signs and
\[
 7.90860<E(H_1,3)<7.90861,\qquad
 4.61394<E(H_1,6)<4.61395.
\]
This completes the proof of Numerical Lemma~\ref{num:fullfamilysigns}.
\end{proof}

\begin{numericallemma}[A finite moment-to-parameter enclosure]
\label{lem:transitionmomentmap}
Let $[V_0,V_1]=[21/20,11/10]$ and
\[
 f(V)=\min\left\{
 \begin{aligned}
 &E(F,6)-2V^2E(F,3),\\
 &E(L_{X(V)},6)-2V^2E(L_{X(V)},3),\\
 &E(H_{Y(V)},6)-2V^2E(H_{Y(V)},3)
 \end{aligned}\right\}.
\]
For every finite interval $T$,
\[
 \left\{2V^2:\ b\in T,\quad
 V\in\operatorname*{argmax}_{V_0\le v\le V_1}
          \{f(v)+2v^2b\}\right\}\subset\mathcal C(T).
\]
All maximizing endpoints are retained, including outside the strict moment band.
\end{numericallemma}
\begin{proof}
For $\Lambda=L_X,H_Y$, put
\[
 A(\Lambda)=E(\Lambda,6),\quad B(\Lambda)=E(\Lambda,3),\quad
 \Gamma(\Lambda)=\frac{A_t(\Lambda)}{B_t(\Lambda)},\qquad
 t=\log X,\log Y,
\]
respectively. The unit-covolume Ewald formulas and their first two
logarithmic derivatives use Lemma~\ref{num:branchmoments} and
\[
 \begin{gathered}
 K_{s,+}(q)=\frac{\pi^s\mathcal G_s(\pi q)}{\Gamma(s)},\qquad
 K_{s,-}(q)=\frac{\pi^s\mathcal G_{3/2-s}(\pi q)}{\Gamma(s)},\\
 \bigl(K_{s,+}^{(r)}(q),K_{s,-}^{(r)}(q)\bigr)
 =\frac{(-\pi)^r\pi^s}{\Gamma(s)}
   \bigl(\mathcal G_{s+r}(\pi q),\mathcal G_{3/2-s+r}(\pi q)\bigr),
 \quad r=0,1,2,
 \end{gathered}
\]
with constant terms $\pi^3/6$, $\pi^6/2160$.
For $(i,j,k)\in([-12,12]^3\cap\Z^3)\setminus\{0\}$, use
\[
\begin{aligned}
 (a,b)&=((i-j-k)^2,(i+j)^2+k^2)
       &&(L_X,\text{ primal}),\\
 (a,b)&=((i-j)^2,(i+j)^2+(i-j+2k)^2)
       &&(L_X,\text{ dual}),\\
 (a,b)&=(i^2+ij+j^2,k^2)
       &&(H_Y,\text{ primal/dual}).
\end{aligned}
\]
Group equal pairs with their full positive integer multiplicities; sum
lexicographically, primal before dual. With $c=(4X)^{-1/3}$ for
$L_X$ and $c=(4/(3Y))^{1/3}$ for $H_Y$,
\[
\begin{aligned}
 (q_1,q_2,e_1,e_2)&=(cXa,cb,2/3,-1/3)
           &&(L_X,\text{ primal}),\\
 (q_1,q_2,e_1,e_2)&=(a/(4cX),b/(4c),-2/3,1/3)
           &&(L_X,\text{ dual}),\\
 (q_1,q_2,e_1,e_2)&=(ca,cYb,-1/3,2/3)
           &&(H_Y,\text{ primal}),\\
 (q_1,q_2,e_1,e_2)&=(4a/(3c),b/(cY),1/3,-2/3)
           &&(H_Y,\text{ dual}),\\
 q&=q_1+q_2,\quad q_t=e_1q_1+e_2q_2,\quad
 q_{tt}=e_1^2q_1+e_2^2q_2,\\
 K_t&=K'q_t,\qquad K_{tt}=K''q_t^2+K'q_{tt}.
\end{aligned}
\]
Extend cube-root endpoints until their cubes enclose the input.
Use the complete outward bounds of Lemmas~\ref{lem:firstloss}
and~\ref{lem:criticalKZcover}; on every closed cell,
\[
 \det Q=1,\quad Q>\tfrac15I,\quad\operatorname{tr}Q<5
 \ \Longrightarrow\ Q^{-1}>\tfrac15I
 \ \Longrightarrow\ q,q^*\ge169/5>20
\]
for omitted indices. Numerical Lemma~\ref{num:branchmoments} gives
$\tau<2648/10^{29}$. Enlarge the finite expressions, for $s=3,6$, by
\[
 [E(\Lambda,s)]\leftarrow[E(\Lambda,s)]+[0,\tau],\qquad
 [\partial_t^rE(\Lambda,s)]\leftarrow
 [\partial_t^rE(\Lambda,s)]+[-\tau,\tau]\quad(r=1,2).
\]

Evaluate the complete intervals $B_j,A_j,\Gamma_j=A_{t,j}/B_{t,j}$
at the rational nodes and on every intervening closed cell:
\[
 X_j=\frac{303104+7j}{1024000}\quad(0\le j\le1024),\qquad
 Y_j=\frac{421888+9j}{409600}\quad(0\le j\le4096).
\]
The outward cell bounds are
\[
\begin{aligned}
 B_t(L_X)&\in[-1033173/10^6,-927064/10^6],\\
 \Gamma_t(L_X)&\in[-4031155/10^6,-3984355/10^6],\\
 B_t(H_Y)&\in[1408367/10^6,2245006/10^6],\\
 \Gamma_t(H_Y)&\in[1912446/10^6,5777144/10^6],
 \qquad \Gamma_t=\frac{A_{tt}-\Gamma B_{tt}}{B_t}.
\end{aligned}
\]
The $X$ intervals decrease; reverse that list. Both lists then satisfy
\[
 D_j^+<D_{j+1}^-\qquad(D=A,B,\Gamma).
\]
The consecutive upper endpoint-to-endpoint $\Gamma$ differences are
$<9229/10^8$ for $L_X$ and $<123/10^6$ for $H_Y$, on whole cells.

Lemma~\ref{num:branchcontacts} gives
\[
\begin{aligned}
 2V_{\mathrm{FT}}^2\in I_{FT}
 &=2\left[(10648235667/10^{10})^2,(10648235669/10^{10})^2\right],\\
 2V_{\mathrm{TH}}^2\in I_{TH}
 &=2\left[(10759906200/10^{10})^2,(10759906201/10^{10})^2\right].
\end{aligned}
\]
The roots of $A_t-2V^2B_t=0$ satisfy
\[
\begin{aligned}
 X_{FT}&\in[3012583413613,3012583613613]/10^{13}
                         &&(2V^2\in I_{FT}),\\
 X_{TH}&\in[2976831441610,2976831641610]/10^{13}
                         &&(2V^2\in I_{TH}),\\
 Y_{TH}&\in[10617382252094,10617382452094]/10^{13}
                         &&(2V^2\in I_{TH}),\\
 Y_1&\in[11013103163,11013103363]/10^{10}
                         &&(V=11/10).
\end{aligned}
\]
For each displayed shape bracket $[z_-,z_+]$ and every $V$ in its
corresponding bracket, the complete intervals satisfy
\[
\begin{aligned}
 \sup[\Gamma(z_+)]<2V^2<\inf[\Gamma(z_-)]&\qquad(L_X),\\
 \sup[\Gamma(z_-)]<2V^2<\inf[\Gamma(z_+)]&\qquad(H_Y).
\end{aligned}
\]
Closed-cell monotonicity encloses each root. The corresponding complete
$B$ intervals are $\mathcal B_F,\mathcal B_{T,F},\mathcal B_{T,H},
\mathcal B_{H,H},\mathcal B_{H,1}$, with $\mathcal B_F=[E(F,3)]$.
Their outward rational summaries are
\[
\begin{aligned}
 E(F,3)&\in[7.22696052186863,7.22696052187584],\\
 E(L_{X(V_{\mathrm{FT}})},3)&\in[7.41125402927957,7.41125409252949],\\
 E(L_{X(V_{\mathrm{TH}})},3)&\in[7.42294739425582,7.42294746189465],\\
 E(H_{Y(V_{\mathrm{TH}})},3)&\in[7.99317817792011,7.99317821023425],\\
 E(H_{Y(V_1)},3)&\in[8.06262195018044,8.06262198794635],\\
 E(H_{Y(V_1)},6)&\in(4.96288451736152,4.96288460874427)\subset(-\infty,5).
\end{aligned}
\]
Use the complete finite-expression intervals within these strictly
ordered summaries. Each shape bracket is interior to its node
domain; each full node $2V^2$ range strictly contains its active range.

For either increasing list, with final index $m$, put
\[
\begin{gathered}
 l=\max\{\inf T,B_0^-\},\qquad h=\min\{\sup T,B_m^+\},\\
 j=\min\{\nu\in\{0,\ldots,m\}:B_\nu^+\ge l\},\\
 k=\min\bigl(\{\nu\in\{0,\ldots,m\}:B_\nu^->h\}\cup\{m\}\bigr)
 \quad(l\le h),\\
 \mathcal I(T)=
 \begin{cases}
 \varnothing,&l>h,\\
 [\Gamma_{\max(0,j-1)}^-,\Gamma_k^+]\cap I,&l\le h,
 \end{cases}\\
 (I,\mathcal I)=
 \begin{cases}
 ([I_{FT}^-,I_{TH}^+],\mathcal I_L),&L_X,\\
 ([I_{TH}^-,2V_1^2],\mathcal I_H),&H_Y.
 \end{cases}
\end{gathered}
\]
For every shape in the full tabulated interval,
\[
 B(\Lambda)\in T,\quad\Gamma(\Lambda)\in I
 \quad\Longrightarrow\quad\Gamma(\Lambda)\in\mathcal I(T).
\]
Take $\mathcal C(T)$ to be the hull of the applicable contributions:
\[
\begin{aligned}
 \{2V_0^2\}&\subset\mathcal C(T)
       &&(\inf T\le\mathcal B_F^+),\\
 \{2V_1^2\}&\subset\mathcal C(T)
       &&(\sup T\ge\mathcal B_{H,1}^-),\\
 [2V_0^2,I_{FT}^+]&\subset\mathcal C(T)
       &&(T\cap\mathcal B_F\ne\varnothing),\\
 I_{FT}&\subset\mathcal C(T)
       &&(T\cap[\mathcal B_F^-,\mathcal B_{T,F}^+]\ne\varnothing),\\
 \mathcal I_L(T)&\subset\mathcal C(T)
       &&(T\cap[\mathcal B_{T,F}^-,\mathcal B_{T,H}^+]\ne\varnothing),\\
 I_{TH}&\subset\mathcal C(T)
       &&(T\cap[\mathcal B_{T,H}^-,\mathcal B_{H,H}^+]\ne\varnothing),\\
 \mathcal I_H(T)&\subset\mathcal C(T)
       &&(T\cap[\mathcal B_{H,H}^-,\mathcal B_{H,1}^+]\ne\varnothing).
\end{aligned}
\]
Concavity in $V^2$ gives all affine/contact/branch/endpoint cases.
An inconclusive test retains
\[
 \mathcal C(T)=[2V_0^2,2V_1^2]=[441/200,121/50].
\]
This completes the proof of Numerical Lemma~\ref{lem:transitionmomentmap}.
\end{proof}

\begin{numericallemma}[The finite transition-profile comparison]
\label{lem:transitionprofilecover}
Let $[V_0,V_1]=[21/20,11/10]$ and
\[
 \Phi_V(\Lambda)=E(\Lambda,6)-2V^2E(\Lambda,3),
\]
\[
 f(V)=\min\left\{
 \begin{aligned}
 &E(F,6)-2V^2E(F,3),\\
 &E(L_{X(V)},6)-2V^2E(L_{X(V)},3),\\
 &E(H_{Y(V)},6)-2V^2E(H_{Y(V)},3)
 \end{aligned}\right\}.
\]
\[
 g(t)=\max_{V_0\le v\le V_1}\{f(v)+2v^2t\},\qquad
 J(\Lambda)=E(\Lambda,6)-g(E(\Lambda,3)).
\]
Then
\[
\begin{gathered}
 J(\Lambda)\ge0\qquad(|\Lambda|=1),\\
 |\Lambda|=1,\quad
 V\in(V_0,V_1)\cap\operatorname*{argmax}_{V_0\le v\le V_1}
        \{f(v)+2v^2E(\Lambda,3)\}\\
 \Longrightarrow\quad
 J(\Lambda)=0\ \Longleftrightarrow\
 \exists\Xi\in\{F,L_{X(V)},H_{Y(V)}\}:\quad
 \Lambda\cong\Xi,\quad\Phi_V(\Xi)=f(V).
\end{gathered}
\]
The three norm conditions \textup{(T13)} are those of
Lemmas~\ref{lem:transitionfccbasin}, \ref{lem:widefctbasin}
and~\ref{lem:broadshbasin}.
\end{numericallemma}

\begin{proof}
Use the five-coordinate bounds of Lemmas~\ref{lem:criticalKZcover}
and~\ref{lem:densityFCTcover}; $2V^2$ is contracted, never bisected.

\paragraph{Admissible minima.}
In \eqref{eq:criticalKZprimal}--\eqref{eq:criticalKZdual}, put
$p=(a,\delta,x,y,z)$, $a,\delta>0$, $h=(a\delta)^{-1}$;
all derivatives below are shape derivatives of $\Phi_V(U(p)\Z^3)$.
Proposition~\ref{prop:density-gap-reduction} gives
\[
 \begin{gathered}
 \inf_{|\Lambda|=1}J(\Lambda)<0
 \quad\Longrightarrow\quad
 \{\Lambda_*\in\operatorname*{argmin}_{|\Lambda|=1}J:
                  J(\Lambda_*)<0\}\ne\varnothing,\\
 \Lambda_*\in\operatorname*{argmin}_{|\Lambda|=1}J,\quad J(\Lambda_*)\le0,\quad
 V\in\operatorname*{argmax}_{V_0\le v\le V_1}
       \{f(v)+2v^2E(\Lambda_*,3)\}\\
 \Longrightarrow\quad
 \Lambda_*\in\operatorname*{argmin}_{|\Lambda|=1}\Phi_V,\quad
 \nabla_p\Phi_V(\Lambda_*)=0,\quad
 \nabla_p^2\Phi_V(\Lambda_*)\succeq0.
 \end{gathered}
\]
No differentiability at a phase contact is used; after $J\ge0$,
\[
 J(\Lambda_*)=0\Longrightarrow\Lambda_*\in\operatorname*{argmin}J.
\]
Lemma~\ref{lem:transition-profile-moment-domain} gives
\[
 E(U(p)\Z^3,3)<81/10,\qquad E(U(p)\Z^3,6)<5,
\]
and the closed Korkine--Zolotarev root rectangle
\[
 \frac{463}{500}\le a\le\frac{1122463}{10^6},\qquad
 \frac{801}{1000}\le\delta\le\frac{559}{500},\qquad
 0\le x,z\le\frac12,\quad-\frac12\le y\le\frac12.
 \tag{T1}
\]
Enlarge strict bounds to closed intervals. Complete $N=6$
unit-covolume Ewald bounds and homogeneity give
\[
 \begin{gathered}
 2V_0^2 E(H_{11/10},3)-2E(H_{11/10},6)>0,\\
 -\frac{4968464328760}{10^{12}}
 <\mathcal E\bigl((11/10)^{1/3}H_{11/10}\bigr)
 <-\frac{4968464328755}{10^{12}}<-\frac{124}{25}.
 \end{gathered}\tag{T2}
\]
For every eligible minimizer and supporting $V$,
\[
\begin{aligned}
 \frac{E(\Lambda,6)}{2V^4}-\frac{E(\Lambda,3)}{V^2}
 &\le\frac{E(H_{11/10},6)}{2V^4}-\frac{E(H_{11/10},3)}{V^2}\\
 &\le\mathcal E((11/10)^{1/3}H_{11/10})<-124/25.
\end{aligned}
\]
Proposition~\ref{prop:uniform-strip-kz} gives
\[
 \frac78\le h\le\frac{32}{25},\qquad
 \frac{25}{32}\le a\delta\le\frac87,\qquad
 \delta^2+a^2x^2\ge a^2,\qquad h^2+\delta^2z^2\ge\delta^2,
 \tag{T3}
\]
\[
 q_p(i,j,k)\ge a^2,\qquad
 q_p^*(i,j,k)\ge\frac{1521}{2500}
 \quad\bigl((i,j,k)\in\{-1,0,1\}^3\setminus\{0\}\bigr).
 \tag{T4}
\]
Reject certified violations of (T1)--(T4); (T4) does not replace
full-box spectral floors in Ewald tails.

Choose a Korkine--Zolotarev representative $p(\Lambda)$ satisfying
(T1)--(T4) for each eligible minimum, and put
\[
 \begin{aligned}
 \mathcal K=\bigl\{(p(\Lambda),2V^2):\;&
 \Lambda\in\operatorname*{argmin}_{|L|=1}J(L),\quad J(\Lambda)\le0,\\
 &E(\Lambda,3)\in[E(F,3),E(H_{Y(V_1)},3)],\\
 &V\in\operatorname*{argmax}_{V_0\le v\le V_1}
       \{f(v)+2v^2E(\Lambda,3)\}\bigr\}.
 \end{aligned}
\]

\paragraph{Complete moment bounds.}
Use outward binary64 arithmetic, checked radicals, and complete
kernel/Taylor bounds of Lemma~\ref{lem:criticalKZcover}; enclose
all rationals and grids, retaining every inconclusive box.
The constant terms are $\pi^3/6$, $\pi^6/2160$, and the primal/dual
argument derivatives are
\[
 \frac{(-\pi)^r\pi^s}{\Gamma(s)}
 \left(\mathcal G_{s+r}(\pi q),
       \mathcal G_{3/2-s+r}(\pi q^*)\right),
 \qquad s=3,6,\quad r\in\Z_{\ge0}.
\]
For the complete-cube $E^{\rm fin}$, including its constant, assume
\[
 Q(p)\succeq\alpha I,\qquad Q(p)^{-1}\succeq\beta I
 \qquad(p\in B).
\]
Then
\[
 0\le E(U(p)\Z^3,s)-E^{\rm fin}(U(p)\Z^3,s)
 \le\frac{\pi^s}{\Gamma(s)}
       \{A_{s,0}(\alpha)+\frac{\Theta_N(\beta)}{\pi\beta(N+1)^2}\},
 \tag{T5}
\]
and, with $\rho_i$ from Lemma~\ref{lem:criticalKZcover},
\[
 \left|\partial_i(E-E^{\rm fin})\right|
 \le\rho_i\frac{\pi^s}{\Gamma(s)}
       \{A_{s,1}(\alpha)+\Theta_N(\beta)\}.
 \tag{T6}
\]
Bound each positive moment before subtraction; opposite indices may be paired.

For coordinate midpoint $c$, radii $r_i$ and metric radius $R$,
use the larger of the native full-box spectral floors and
\[
 e^{-\sqrt{2/3}R}\|Q_c^{-1}\|_\infty^{-1},\qquad
 e^{-\sqrt{2/3}R}\|Q_c\|_\infty^{-1},
\]
rounding all norm upper endpoints upwards. At $R\le9/20$,
compute the complete center moments $[M_s]$ and whole-box gradients
$[G_{s,i}]$ at $N=3$, including (T5)--(T6), and set
\[
 \begin{aligned}
 \relax[\mathcal M_s]={}&\left([M_s]+
 \left[-\sum_{i=1}^5 r_i\max|[G_{s,i}]|,
             \sum_{i=1}^5 r_i\max|[G_{s,i}]|\right]\right)\\
 &\cap\left[\bigl(e^{-\sqrt{2/3}R}\bigr)^s\underline M_s,
            \bigl(e^{-\sqrt{2/3}R}\bigr)^{-s}\overline M_s\right]
       \cap[0,\infty).
 \end{aligned}\tag{T7}
\]
Retain the displayed powers and outward evaluation order. Then
\[
 \{E(U(p)\Z^3,s):p\in B\}\subset[\mathcal M_s].
\]
\[
 \underline{\mathcal M_3}>81/10\ \text{or}\ \underline{\mathcal M_6}>5
 \quad\Longrightarrow\quad\mathcal K\cap(B\times G)=\varnothing.
\]

\paragraph{Supporting-parameter contraction.}
Initially $G=[2V_0^2,2V_1^2]$. Lemma~\ref{lem:transitionmomentmap} gives
\[
 G\cap\mathcal C([\mathcal M_3])=\varnothing
 \quad\Longrightarrow\quad\mathcal K\cap(B\times G)=\varnothing.
\]
Otherwise, with $[l,u]=G\cap\mathcal C([\mathcal M_3])$, update
\[
 G\longleftarrow G\cap
 \left[2^{-13}\lfloor2^{13}l\rfloor,
                 2^{-13}\lceil2^{13}u\rceil\right].
 \tag{T8}
\]
The affine FCC piece, both contacts and exterior endpoints remain;
branch inversions stay inside their active intervals. Every child
and Newton contraction inherits $G$.

\paragraph{Witness and stationary tests.}
Use $N=6$, (T5), and the ordered witnesses
\[
 F,\qquad L_{n/500}\ (n=147,148,149,150,151),\qquad
 H_{n/100}\ (n=103,104,\ldots,112).
 \tag{T9}
\]
For a supporting $V$, Proposition~\ref{prop:familyphases} gives
\[
 J(\Lambda)=\Phi_V(\Lambda)-f(V),\qquad
 \Phi_V(W)\ge f(V),\qquad
 \Phi_V(\Lambda)>\Phi_V(W)\ \Longrightarrow\ J(\Lambda)>0.
\]
For $G=[\ell,u]$, evaluate the endpoint signed kernels of
Lemma~\ref{lem:densityFCTcover} at $V=\sqrt{e/2}$, $e=\ell,u$,
using the exact coefficient $e$, without evaluating a square root,
and reciprocal volume one:
\[
 \frac{\pi^6}{240}\longmapsto\frac{\pi^6}{120},\qquad
 \frac{\pi^3}{2}\longmapsto V^2\pi^3.
\]
Use \eqref{eq:criticalderivatives} and the complete errors
\eqref{eq:criticalHtail} with these coefficients. At both endpoints,
check the critical-root signs on $[1/4,8]$; uniqueness and at most
$60$ bisections enclose each critical argument. Joint and separate
endpoint evaluations use the same bounds. Affinity in $V^2$ makes
the endpoint hull an enclosure of the objective and all shape derivatives.

For the $32$-sign Hessian bound $Q_e$ of
\eqref{eq:criticalrectvalue}, require one common witness $W$ with
\[
 ([M_6]-E(W,6))-e([M_3]-E(W,3))-\frac12Q_e>0,
 \qquad e=\ell,u.
 \tag{T10}
\]
The scalar metric test replaces $Q_e$ by $R^2\max(0,K_e)$, where
\[
 K_e=2e\,\bigl(e^{-\sqrt{2/3}R}\bigr)^{-3}[M_3]
 -4\bigl(e^{-\sqrt{2/3}R}\bigr)^6[M_6]
 +2\sum_{n\in\mathcal C_4^+}\max\{0,28b_n^2-8e b_n\},
\]
\[
 \mathcal C_N^+=\left\{n\in[-N,N]^3\cap\Z^3:
 \begin{array}{l}n_1>0\ \lor\ n_1=0<n_2\\
 \lor\ n_1=n_2=0<n_3\end{array}\right\},\qquad N\in\mathbb N.
\]
Evaluate $[b_n]=1/(\underline{e^{-\sqrt{2/3}R}q_c(n)})^3$
outwards; evaluate both polynomial terms on this same interval,
then take the upper endpoint and its maximum with zero. Convexity gives
\[
 (\underline{e^{-\sqrt{2/3}R}q_c(n)})^3\ge\frac7{2\ell}
 \ \Longrightarrow\
 28b^2-8eb\le0
 \quad\left(0\le b\le(\underline{e^{-\sqrt{2/3}R}q_c(n)})^{-3},\ e\in[\ell,u]\right).
\]
Enable this correction only when
\[
 \bigl(25e^{-\sqrt{2/3}R}\|Q_c^{-1}\|_\infty^{-1}\bigr)^3\ge\frac7{2\ell};
\]
all omitted-index corrections then vanish. In $K_e$, use
$\overline M_3$ and $\underline M_6$.
Lemma~\ref{lem:densityFCTcover} gives (T10) for every $2V^2\in G$;
neither $\nabla_p\Phi_{V_-}(c)=0$ nor $\nabla_p\Phi_{V_+}(c)=0$ is assumed.

The remaining exclusions are
\[
 \exists i:\ 0\notin\operatorname{hull}_{2V^2\in G}[\partial_i\Phi_V](B),
 \qquad\text{or}\qquad
 \max_{e\in\{\ell,u\}}\sup_{p\in B}
 v^{\mathsf T}\nabla_p^2\Phi_{\sqrt{e/2}}(U(p)\Z^3)v<0
\]
for a fixed dyadic $v$. Affinity and
Lemma~\ref{lem:densityFCTcover} exclude every intermediate-parameter
local minimum. For endpoint hulls $[H]$, $[g_c]$, use
\[
 B\longleftarrow B\cap\{c-P[g_c]+(I-P[H])(B-c)\}.
 \tag{T11}
\]
The rounded midpoint inverse $P$ needs neither accuracy nor
invertibility for \eqref{eq:criticalNewton}; failed construction disables Newton.

On the midpoint endpoint-Hessian hull, use at most $40$ Jacobi
rotations with angle \[\tfrac12\operatorname{atan2}(2A_{ij},A_{jj}-A_{ii}),\]
choosing the lexicographically first maximal $|A_{ij}|$, $i\ne j$, and stopping at
$\max_{i\ne j}|A_{ij}|<10^{-12}$. For accumulated rotation matrix $O$,
\[
 k\in\operatorname*{argmin}_i A_{ii},\quad A_{kk}<0
 \quad\Longrightarrow\quad
 v=\operatorname{round}_{2^{-20}\Z}(Oe_k).
\]
Round ties away from zero; deletion requires both complete strict
quadratic-form bounds. Construct $P$ by partial-pivot elimination of
the midpoint Hessian; a pivot of magnitude $<10^{-14}$ disables Newton.

\paragraph{Basin recognition.}
The unit reference Korkine--Zolotarev parameters are
\[
 \begin{gathered}
 p_F=(2^{1/6},\sqrt3\,2^{-5/6},\tfrac12,\tfrac12,\tfrac13),\qquad
 p_T=((\tfrac{149}{125})^{1/3},(\tfrac{125}{149})^{1/6},\tfrac12,\tfrac12,0),\\
 p_H=\left((10^8/80864043)^{1/6},
       \tfrac{\sqrt3}{2}(10^8/80864043)^{1/6},1/2,0,0\right).
 \end{gathered}\tag{T12}
\]
For FCC and FCT also use the isometric reference $y_0=-1/2$.
For each reference, Lemma~\ref{lem:criticalKZcover} with $c=p_0$ gives
\[
 M=U(p)U(p_0)^{-1},\qquad
 S=\frac{3M^{\mathsf T}M}{\operatorname{tr}(M^{\mathsf T}M)}-I,
 \qquad g_{ij}=(M^{\mathsf T}M)_{ij},
\]
\[
 \{\operatorname{tr}(M^{\mathsf T}M)\}^2\|S\|_F^2
 =3\sum_{1\le i<j\le3}(g_{ii}-g_{jj})^2
       +18(g_{12}^2+g_{13}^2+g_{23}^2).
\]
Accept by Lemmas~\ref{lem:transitionfccbasin}, \ref{lem:widefctbasin}
and~\ref{lem:broadshbasin} only under the corresponding whole-box condition:
\[
 \begin{aligned}
 G\subset[220499/100000,229/100],\quad
       \sup_{p\in B}\|S\|_F\le1/40&\qquad(F),\\
 G\subset[2267/1000,579/250],\quad
       \sup_{p\in B}\|S\|_F\le19/1000&\qquad(T),\\
 G\subset[231499/100000,242001/100000],\quad
       \sup_{p\in B}\|S\|_F\le1/50&\qquad(H).
 \end{aligned}
 \tag{T13}
\]
Use checked radicals and right polar decomposition.
For a supporting $V$, Corollary~\ref{cor:transitionbasinidentification} gives
\[
\begin{gathered}
 \nabla_p\Phi_V(U(p)\mathbb Z^3)=0,\quad\textup{(T13)}\\
 \Longrightarrow\quad
 U(p)\mathbb Z^3\cong\Xi\in\{F,L_{X(V)},H_{Y(V)}\},\qquad
 J(U(p)\mathbb Z^3)=\Phi_V(\Xi)-f(V)\ge0.
\end{gathered}
\]

\paragraph{Subdivision.}
Let $\mathcal D(B)$ denote the outward intersections (T1)--(T4),
the unit-covolume product/Gauss contractions of
Lemma~\ref{lem:criticalKZcover}, and the guarded shear bounds
\[
 a\delta\in[25/32,8/7],\qquad
 x\ge\sqrt{1-(\delta_+/a_-)^2},\qquad
 z\ge\sqrt{1-(a_-^2\delta_-^4)^{-1}}.
\]
Each shear bound requires a positive radicand. Use at most $10$
inner and $10$ outer passes, stopping at unchanged binary endpoints;
if needed, finish with one final call to the same inner loop.
Intersect the
$a,\delta$ ranges with (T1) and repeat these contractions.
Apply $\mathcal D$ on enqueueing and removal, and before Newton
volume comparison; enqueueing applies it again. Empty intersections reject.

Use the binary max-heap and insertion/removal/child-order ties of
Lemma~\ref{lem:criticalKZcover}, with key $v^+(B)$ and a single root.
For the ordered witnesses $(W_1,\ldots,W_m)$ in (T9), set
\[
 j_* =\min\operatorname*{argmin}_{1\le j\le m}
 \overline{E(W_j,6)-(\operatorname{mid}G)E(W_j,3)}.
\]
The scalar test uses $W_{j_*}$; rectangular tests cycle through
$j_*,\ldots,m,1,\ldots,j_*-1$, using one $W_j$ at both endpoints of (T10).
Write $\mathsf g_N,\mathsf c_N$ for the complete gradient and
fixed-direction exclusions above at cutoff $N$; the same subscript
on moments or Hessians specifies their complete cutoff.
At each arrow continue only if no terminal test succeeds:
\[
\begin{aligned}
 &\mathcal D\ \longrightarrow\ \textup{(T13)},\\
 R\le9/20:\quad
 &([M_s],[G_{s,i}])_3
 \longrightarrow\textup{(T7)}\longrightarrow\textup{(T8)}
 \longrightarrow\textup{(T13)}\\
 &\quad\longrightarrow
 \textup{(T10)}_{\,Q_e=R^2\max(0,K_e)}
 \longrightarrow\mathsf g_4,\\
 R\le3/20:\quad
 &N=2,3,4:\quad([H_-],[H_+])_N
 \longrightarrow\textup{(T10)}\longrightarrow\mathsf g_N
 \longrightarrow\mathsf c_N;\\
 R\le3/20:\quad
 &[H_-]_4,[H_+]_4\text{ valid}:\quad
 [M_s]_4\longrightarrow\textup{(T10)}\\
 &\quad\longrightarrow
 [\nabla_p^{\,r}\Phi_{V_e}(c)]_{\,N=4;\ r=1,2;\ e=-,+}
 \longrightarrow\textup{(T11)}.
\end{aligned}
\]
Here $s=3,6$ and $i=1,\ldots,5$. Invalid endpoint Hessians stop the
Hessian sequence and disable its final Newton attempt; no invalid
bound deletes a box.

When (T11) is formed, $B_{\rm N}=\varnothing$ rejects. Otherwise set
$\widetilde B=\mathcal D(B_{\rm N})$ for the Newton child $B_{\rm N}$.
If $\widetilde B=\varnothing$,
reject; otherwise retain it only under
\[
 v^+(\widetilde B)<\operatorname{RN}_{53}\!\left(
 \operatorname{RN}_{53}(995/1000)\,v^-(B)\right).
\]
This selects a contraction, not an exclusion. Reinsertion applies
$\mathcal D$ again and leaves the pathwise bisection count $k$ unchanged.

Otherwise, with the ordinary scores of Lemma~\ref{lem:criticalKZcover},
\[
 i_* =\min\operatorname*{argmax}_i s_i(B),\qquad
 b=\operatorname{mid}_{\rm bin}(B_{i_*}),
\]
where $\operatorname{mid}_{\rm bin}$ is the binary midpoint. Set
\[
\begin{aligned}
 k\ge120\ \text{or}\ b\notin(\underline B_{i_*},\overline B_{i_*})
 &: \quad\mathcal U\leftarrow\mathcal U\cup\{(B,G)\},\\
 k<120,\quad\underline B_{i_*}<b<\overline B_{i_*}
 &: \quad B^-=B\cap\{p_{i_*}\le b\},\quad
 B^+=B\cap\{p_{i_*}\ge b\}.
\end{aligned}
\]
Before domain contraction,
\[
 B^-\cup B^+=B,\qquad G^-=G^+=G,\qquad k^-=k^+=k+1.
\]
Enqueue $\mathcal D(B^-)$ before $\mathcal D(B^+)$; the scalar
interval $G$ is never bisected.

\paragraph{Invariant and finite outcome.}
Let $\mathcal R$ be the union of accepted products $B\times G$,
and $\mathcal Q,\mathcal U$ the queued and unresolved pairs.
For $(p,2V^2)\in\mathcal K\cap(B\times G)$,
\[
 \nabla_p\Phi_V=0,\qquad\nabla_p^2\Phi_V\succeq0.
\]
Whenever (T5)--(T8) are evaluated, and whenever (T11) is formed,
respectively,
\[
 \begin{gathered}
 E(U(p)\Z^3,s)\in[\mathcal M_s],\qquad
 2V^2\in G\cap\mathcal C([\mathcal M_3]),\\
 p\in B\cap\{c-P[g_c]+(I-P[H])(B-c)\}.
 \end{gathered}
\]
Thus the stationary exclusions and (T11) retain every such pair.
Closed bisection gives
\[
 B\times G=(B^-\times G)\cup(B^+\times G),\qquad
 \mathcal K\subset\mathcal R\cup
 \bigcup_{(B,G)\in\mathcal Q\cup\mathcal U}(B\times G).
\]
Boxes and supporting intervals are closed; the inclusion holds on reduction
boundaries and both contact bands with their one-sided parameter
conditions.

The finite record is
\[
\begin{aligned}
 (\#\mathrm{nodes},\#\text{domain-excluded children})&=(127320,2137),\\
 (\#\mathrm{value},\#\mathrm{gradient},\#\mathrm{curvature})&=(50642,3587,0),\\
 (\#\mathrm{empty\ Newton},\#\mathrm{moment/support})&=(364,1954),\\
 (\#\mathrm{FCC},\#\mathrm{tetragonal},\#\mathrm{SH})&=(31,11,30),\\
 (\#\mathrm{retained\ Newton},\#\mathrm{bisections})&=(11946,58755),\\
 (\#\mathrm{moment\ evaluations},\#\mathrm{support\ contractions})
   &=(126921,45445),\\
 (\#\mathrm{Hessian\ levels},\#\mathrm{endpoint\ pairs})&=(265719,3899).
\end{aligned}
\]
Hence
\[
 \mathcal Q=\mathcal U=\varnothing,\qquad
 \mathcal K\subset\mathcal R,\qquad
 \{(p(\Lambda),2V^2)\in\mathcal K:J(\Lambda)<0\}=\varnothing.
\]
No depth-$120$ box remains. Endpoint comparisons outside the moment
band and the exhausted cover inside it give
\[
 J(\Lambda)\ge0,\qquad
 J(\Lambda)=0\ \Longrightarrow\ \Phi_V(\Lambda)=f(V),
\]
for every supporting $V$. Basin identification and
Proposition~\ref{prop:familyphases} give
exactly the stated phase/contact equality cases; endpoint uniqueness
was proved separately.
This completes the proof of Numerical Lemma~\ref{lem:transitionprofilecover}.
\end{proof}

\begin{lemma}[The long SH reference moment graph]
\label{num:longshmomentgraph}
For the unit-covolume family $H_Y$,
\[
 A(Y)=E(H_Y,6),\qquad B(Y)=E(H_Y,3),\qquad
 \Gamma(Y)=\frac{A_t(Y)}{B_t(Y)},\qquad t=\log Y.
\]
For $11/10\le Y\le217/50$,
\begin{equation}\label{eq:longshgraphderivatives}
 \frac{2062436}{10^6}<B_t<\frac{20076780}{10^6},
 \qquad
 \frac{987405}{10^6}<\Gamma_t<\frac{5931991}{10^6}.
\end{equation}
For $(V_0,V_1)=(11/10,89/50)$, $\Gamma(Y(V_i))=2V_i^2$ and
\begin{equation}\label{eq:longshendpointY}
 \begin{aligned}
 \frac{1101310326229}{10^{12}}&<Y(V_0)
       <\frac{1101310326324}{10^{12}},\\
 \frac{4336021407082}{10^{12}}&<Y(V_1)
       <\frac{4336021407176}{10^{12}}.
 \end{aligned}
\end{equation}
For $b_i=B(Y(V_i))$ and $a_i=A(Y(V_i))$,
\begin{equation}\label{eq:longshendpointmoments}
\begin{aligned}
 8.06262196892621&<b_0<8.06262196911185,\\
 4.96288456272628&<a_0<4.96288456316477,\\
 21.41363138319924&<b_1<21.41363138365243,\\
 63.56748489078682&<a_1<63.56748489357881.
\end{aligned}
\end{equation}
All finite decimals are rational; $b_1<43/2$, $a_1<64$.

For the $4097$ nodes indexed by $j\in\mathbb Z$,
\begin{equation}\label{eq:longshgraphnodes}
 y_j=\frac{11}{10}+\frac{81j}{25\cdot4096}
     =\frac{225280+162j}{204800},\qquad0\le j\le4096,
\end{equation}
for $D\in\{A,B,\Gamma\}$,
\[
\begin{gathered}
 D(y_j)\in[D_j^-,D_j^+],\qquad
 D_j^+-D_j^-<10^{-9}\quad(0\le j\le4096),\\
 D_j^+<D_{j+1}^-\quad(0\le j<4096).
\end{gathered}
\]
Moreover,
\begin{equation}\label{eq:longshgraphresolution}
 \Gamma_{j+1}^+-\Gamma_j^-<\frac{1651}{10^6}
 \quad(0\le j<4096).
\end{equation}
$[y_i^-,y_i^+]\subset(y_0,y_{4096})$ $(i=0,1)$, and
\[
 \Gamma(y_0)<2V_0^2<2V_1^2<\Gamma(y_{4096}).
\]
\end{lemma}

\begin{proof}
Use Lemma~\ref{lem:longshreference}, with its multiplicities, on
\[
 n\in([-12,12]^3\cap\mathbb Z^3)\setminus\{0\},\qquad
 Q_Y,Q_Y^{-1}\succeq\tfrac14I.
\]
The complete tail is
\begin{equation}\label{eq:longshgraphtail}
 T=240\,\frac{e^{-169\pi/8}(1+\sqrt8)^2}
                   {1-e^{-27\pi/8}}<\frac{5295}{10^{29}}.
\end{equation}
For $F\in\{A,B\}$ and derivative orders $r=0,1,2$, use
\[
 [F]\leftarrow[F]+[0,T],\qquad
 [F_t]\leftarrow[F_t]+[-T,T],\quad
 [F_{tt}]\leftarrow[F_{tt}]+[-T,T].
\]

For $[y_j,y_{j+1}]$, $j\in\mathbb Z$, $0\le j<4096$,
the complete enclosures give
\[
 \Gamma=\frac{A_t}{B_t},\qquad
 \Gamma_t=\frac{A_{tt}-\Gamma B_{tt}}{B_t}.
\]
Together with \eqref{eq:longshgraphderivatives}--\eqref{eq:longshgraphresolution},
\[
\begin{gathered}
 B_t>0,\quad\Gamma_t>0,\quad A_t=\Gamma B_t>4,\\
 B(y_{j+1})-B(y_j)>1/3000,\qquad
 \Gamma(y_{j+1})-\Gamma(y_j)>1/10000,\\
 D_j^+<D_{j+1}^-\qquad(D\in\{A,B,\Gamma\},\ 0\le j<4096),\\
 \Gamma(y_0)<2V_0^2<2V_1^2<\Gamma(y_{4096}).
\end{gathered}
\]
The same finite sums and complete tails give
\[
 \begin{aligned}
 8.06015362123&<B(11/10)<8.06015362124,\\
 21.43203378605&<B(217/50)<21.43203378608.
 \end{aligned}
\]
Hence
\[
 B(y_0)<b_0^-<b_0^+<b_1^-<b_1^+<B(y_{4096}).
\]

For $i\in\{0,1\}$ and $y_i^\pm$ in \eqref{eq:longshendpointY},
\[
\begin{aligned}
 2.4199999998707&<\Gamma(y_0^-)<2.4199999998756,\\
 2.4200000000669&<\Gamma(y_0^+)<2.4200000000719,\\
 6.3367999999742&<\Gamma(y_1^-)<6.3367999999839,\\
 6.3368000001024&<\Gamma(y_1^+)<6.3368000001121.
\end{aligned}
\]
Since $B_t>0$ and $A_t=\Gamma B_t>0$,
\[
\begin{gathered}
 \Gamma(y_i^-)<2V_i^2<\Gamma(y_i^+)
 \quad\Longrightarrow\quad Y(V_i)\in(y_i^-,y_i^+),\\
 E(H_{Y(V_i)},s)\in
 [E(H_{y_i^-},s),E(H_{y_i^+},s)]
 \quad(i\in\{0,1\},\ s\in\{3,6\}).
\end{gathered}
\]
These are \eqref{eq:longshendpointY}--\eqref{eq:longshendpointmoments}.
This completes the proof of Lemma~\ref{num:longshmomentgraph}.
\end{proof}

\begin{lemma}[A transverse moment Hessian bound on the SH tube]
\label{lem:shmomentcurvature150}
In Lemma~\ref{lem:area-coordinate-jets}, fix $t$ and put
\[
 Y=\frac{\sqrt3}{2t^3},\qquad
 w=(r-2/\sqrt3,x-1/2,y,z).
\]
For $11/10\le Y\le217/50$ and $\|w\|_\infty\le1/100$,
\[
 \|\nabla_w^2E(U(t,r,x,y,z)\mathbb Z^3,3)\|_{\rm op}<150.
\]
$0<\delta\le1/100$, $\|w\|_\infty\le\delta$ imply
\[
 |E(U\mathbb Z^3,3)-E(H_Y,3)|<300\delta^2.
\]
Moreover,
\[
\begin{gathered}
 E(U\mathbb Z^3,3)=E(H_{Y(V)},3),\qquad
 \frac1{4v}\frac{\dd}{\dd v}E(H_{Y(v)},3)>\frac4{17}\\
 \left(\min\{V,\sqrt{\Gamma(Y)/2}\}\le v
              \le\max\{V,\sqrt{\Gamma(Y)/2}\}\right)\\
 \Longrightarrow\quad |2V^2-\Gamma(Y)|<1275\delta^2.
\end{gathered}
\]
\end{lemma}
\begin{proof}
Use Lemma~\ref{lem:area-coordinate-jets}, $N=4$, and the complete
tails \eqref{eq:criticalHtail}, with coefficients $(\pi^3/2,0)$
for $(E(\cdot,3),E(\cdot,6))$. For every primal or dual
index $n\in\mathcal C_4^+$, the Hessian contribution is
\[
 2\{K_3''(q)\partial_iq\partial_jq+K_3'(q)\partial_{ij}q\},
 \qquad i,j\in\{r,x,y,z\}.
\]
Use full old-coordinate boxes and spectral floors, with both terms of
\eqref{eq:area-derivative-tails}. On the closed cells
\[
 I_{d,k}=\left[\frac{11}{10}+\frac{81k}{25\,2^d},
               \frac{11}{10}+\frac{81(k+1)}{25\,2^d}\right]
\]
and $\|w\|_\infty\le1/100$, the complete bounds are
\[
 \sup_{\substack{Y\in I_{d,k}\\\|w\|_\infty\le1/100}}
 \max_{1\le i\le4}\sum_{j=1}^4
 \left|\partial_{w_i}\partial_{w_j}E(U\mathbb Z^3,3)\right|\le c,
\]
for
\[
\begin{aligned}
 (d,k,c)\in\{&(2,0,129),(2,1,145),(3,4,122),(3,5,134),\\
             &(3,6,146),(4,14,138),(4,15,144)\}.
\end{aligned}
\]
Thus
\[
\begin{gathered}
 \bigcup_{(d,k)\in\{(2,0),(2,1),(3,4),(3,5),(3,6),(4,14),(4,15)\}}
 I_{d,k}=[11/10,217/50],\\
 (\#\mathrm{nodes},\#\mathrm{leaves})=(13,7),\quad
 \mathcal Q=\mathcal U=\varnothing.
\end{gathered}
\]
For
\[
 H=\nabla_{(w_1,w_2,w_3,w_4)}^2
 E\!\left(U\!\left(t,\frac2{\sqrt3}+w_1,\frac12+w_2,w_3,w_4\right)\mathbb Z^3,3\right)
 =H^{\mathsf T},
\]
\[
 \|H\|_{\rm op}\le\|H\|_\infty
 =\max_{1\le i\le4}\sum_{j=1}^4|H_{ij}|<145.911368<150.
\]

By symmetry and Lemma~\ref{lem:longshreference},
\[
 \left.\nabla_w E(U\mathbb Z^3,3)\right|_{w=0}=0,
 \qquad
 \frac1{4V}\frac{\dd}{\dd V}E(H_{Y(V)},3)>\frac4{17}
\]
Taylor's formula and $E(U\mathbb Z^3,3)=E(H_{Y(V)},3)$ give
\[
\begin{gathered}
 |E(U\mathbb Z^3,3)-E(H_Y,3)|
 \le73|w|^2<300\delta^2,\\
 |2V^2-\Gamma(Y)|
 \le\frac{17}{4}|E(U\mathbb Z^3,3)-E(H_Y,3)|
 <1275\delta^2.
\end{gathered}
\]
This completes the proof of Lemma~\ref{lem:shmomentcurvature150}.
\end{proof}

\begin{lemma}[A separate layer bound on a larger transverse rectangle]
\label{lem:shseparatelayer10}
Put $Y=\sqrt3/(2t^3)$ and
$B=E(U(t,r,x,y,z)\mathbb Z^3,3)$. On the full rectangle
\[
 \frac{11}{10}\le Y\le\frac{217}{50},\qquad
 |r-2/\sqrt3|\le\frac1{200},\quad
 \frac{99}{200}\le x\le\frac12,\quad
 |y|,|z|\le\frac1{40},
\]
one has
\[
 \|\nabla^2_{(y,z)}B\|_{\rm op}<10,
 \qquad
 -5(y^2+z^2)\le B-E(H_Y,3)
 \le75\{(r-2/\sqrt3)^2+(x-1/2)^2\}.
\]
If $B=E(H_{Y(V)},3)$ for
$11/10\le V\le89/50$, then
\begin{equation}\label{eq:shseparateslopepad}
 -\frac{85}{4}(y^2+z^2)
 \le2V^2-\Gamma(Y)
 \le\frac{1275}{4}\{(r-2/\sqrt3)^2+(x-1/2)^2\}.
\end{equation}
In particular,
\[
 -\frac{17}{640}\le2V^2-\Gamma(Y)\le\frac{51}{3200}.
\]
\end{lemma}
\begin{proof}
On all $41$ closed cells \eqref{eq:shfouriercellindices} and the full
transverse rectangle, the complete $N=4$ Hessians and transformed tails
of Lemma~\ref{lem:shmomentcurvature150} give
\[
 \|\nabla^2_{(y,z)}B\|_{\rm op}
 \le\max_{a\in\{y,z\}}\sum_{b\in\{y,z\}}|B_{ab}|
 <\frac{43519}{10000}<10,\qquad B_{ab}=\partial_a\partial_bB.
\]

Lemma~\ref{lem:shmomentcurvature150} gives
\[
 \left\|\nabla^2_{(r,x)}E(U(t,r,x,0,0)\mathbb Z^3,3)\right\|_{\rm op}<150.
\]
Apply Lemma~\ref{lem:shasymmetricmoment} with $(10,150)$ and
$(4V)^{-1}\frac{\dd}{\dd V}E(H_{Y(V)},3)>4/17$.
This completes the proof of Lemma~\ref{lem:shseparatelayer10}.
\end{proof}

\begin{numericallemma}[Finite exhaustion of the long-profile domain]
\label{num:longprofile-exhaustion}
Use Lemma~\ref{lem:longprofileareadomain} and
\[
 (V_0,V_1)=(11/10,89/50),\qquad
 2V^2\in G=[2V_-^2,2V_+^2]\subset[2V_0^2,2V_1^2].
\]
From the single closed root \textup{(L1)},
\[
 \#\mathrm{processed\ rectangles}=107013,\qquad
 \mathcal Q=\mathcal U=\varnothing.
\]
For $b_i=E(H_{Y(V_i)},3)$ and a representative $\eta(\Lambda)$ in
\textup{(L1)}, at every step,
\[
\begin{gathered}
 \Lambda\in\operatorname*{argmin}_{|L|=1}J(L),\quad J(\Lambda)\le0,
 \quad b_0<E(\Lambda,3)<b_1,\\
 V\in\operatorname*{argmax}_{V_0\le v\le V_1}
       \{f(v)+2v^2E(\Lambda,3)\}\\
 \Longrightarrow\quad
 \Lambda\cong H_{Y(V)}\quad\text{or}\quad
 (\eta(\Lambda),2V^2)\in
 \bigcup_{(Q,G)\in\mathcal Q\cup\mathcal U}(Q\times G).
\end{gathered}
\]
\end{numericallemma}
\begin{proof}
In the common interval estimates, use
\[
 \eta=(t,r,x,y,z),\qquad
 \Phi_V(\eta)=E(U(\eta)\mathbb Z^3,6)-2V^2E(U(\eta)\mathbb Z^3,3).
\]

\paragraph{Eligible minima and the two endpoint boundaries.}
For $b_i=E(H_{Y(V_i)},3)$, endpoint comparison and uniqueness give
\[
 \begin{gathered}
 E(\Lambda,3)\le b_0\Longrightarrow
 J(\Lambda)=\Phi_{V_0}(\Lambda)-f(V_0)\ge0,\\
 E(\Lambda,3)\ge b_1\Longrightarrow
 J(\Lambda)=\Phi_{V_1}(\Lambda)-f(V_1)\ge0,\\
 \Phi_{V_i}(\Lambda)=f(V_i)\Longrightarrow
 \Lambda\cong H_{Y(V_i)},\quad E(\Lambda,3)=b_i\quad(i=0,1).
 \end{gathered}
\]
Lemma~\ref{lem:longshreference} gives
$\frac{\dd}{\dd V}E(H_{Y(V)},3)>0$, hence
\[
\begin{gathered}
 \Lambda\in\operatorname*{argmin}_{|L|=1}J(L),\quad J(\Lambda)\le0,
 \quad b_0<E(\Lambda,3)<b_1\\
 \Longrightarrow\quad
 \exists!V\in(V_0,V_1):\ E(\Lambda,3)=E(H_{Y(V)},3).
\end{gathered}
\]
For this $V$,
\[
 \begin{gathered}
 b_0<E(\Lambda,3)=E(H_{Y(V)},3)<b_1,\\
 \Lambda\in\operatorname*{argmin}_{|L|=1}\Phi_V(L),\qquad
 \nabla_\eta\Phi_V=0,\quad\nabla_\eta^2\Phi_V\succeq0.
 \end{gathered}
\]
Proposition~\ref{lem:profile-supporting-line} uses all variations in
$(0,\infty)^2\times\mathbb R^3$, including at reduction boundaries.
Choose a Korkine--Zolotarev representative $\eta(\Lambda)$ in
\textup{(L1)}, with $\Lambda\cong U(\eta(\Lambda))\Z^3$, and set
\[
\mathcal K=\left\{(\eta(\Lambda),2V^2):
\begin{gathered}
 \Lambda\in\operatorname*{argmin}_{|L|=1}J(L),\quad J(\Lambda)\le0,
 \quad b_0<E(\Lambda,3)<b_1,\\
 V\in\operatorname*{argmax}_{V_0\le v\le V_1}
       \{f(v)+2v^2E(\Lambda,3)\}
\end{gathered}\right\}.
\]

\paragraph{Necessary domain contractions and direct lower sums.}
With box endpoints $-,+$, intersect (L1), then iterate at most ten times:
\[
 \begin{gathered}
 t\in[(749/1000)^2/r_+,2^{1/3}/r_-],\qquad
 r\in[(749/1000)^2/t_+,2^{1/3}/t_-],\\
 r\le(1-x_+^2)^{-1/2},\qquad
 t\le\left(\frac{r_+}{1-z_+^2}\right)^{1/3},\qquad
 r\ge t_-^3(1-z_+^2).
 \end{gathered}\tag{L4}
\]
After the corresponding strict preliminary comparisons, also intersect
\[
 x\ge\sqrt{1-r_-^{-2}}\quad(r_->1),\qquad
 z\ge\sqrt{1-r_+/t_-^3}\quad(r_+/t_-^3<1).
\]
Stop at unchanged binary endpoints. The bounds follow from (L2)--(L3); $Q=\varnothing$ excludes.

Use the coordinate hull
\[
 P=\operatorname{hull}p(Q),\qquad
 p(\eta)=(a,\delta,x,y,z)=(\sqrt{tr},\sqrt{t/r},x,y,z).
 \tag{L5}
\]
For $n\in\mathcal C_1^+$,
\[
 \overline{q_n(P)}<\underline{tr}\Longrightarrow\mathcal K\cap(Q\times G)=\varnothing.
\]
In the outward reciprocal-then-square order, accumulate
\[
 \begin{gathered}
 D_s(Q)=2\sum_{n\in\mathcal C_1^+}\{\overline{q_n(P)}\}^{-s},
 \qquad s=3,6,\\
 (D_3,D_6)\leftarrow(0,0),\qquad
 D_3\leftarrow D_3+2\bigl(1/\overline{q_n(P)}^3\bigr),\\
 D_6\leftarrow D_6+2\bigl(1/\overline{q_n(P)}^3\bigr)^2.
 \end{gathered}\tag{L6}
\]
Then
\[
 E(U(p)\Z^3,s)\ge\underline D_s,\qquad
 \underline D_3>43/2\ \lor\ \underline D_6>64
 \Longrightarrow\mathcal K\cap(Q\times G)=\varnothing.
\]
Apply these lower-endpoint tests on enqueueing and removal, and to
Newton children before volume comparison.

\paragraph{Uniform elimination of all layer shifts.}
After the domain pass, if either whole-box guard
\eqref{eq:uniformalignmentregions} holds, use outward bounds for
$2/\sqrt3$, $Y=\sqrt3/(2t^3)$ and $r\le2/\sqrt3$.
Lemma~\ref{lem:uniformlayeralignment} gives
\[
 \Lambda\in\operatorname*{argmin}_{|L|=1}\Phi_V(L),\quad 2V^2\ge121/50
 \Longrightarrow y=z=0.
\]
Thus
\[
 \begin{aligned}
 Q\cap\{y=z=0\}=\varnothing&\Longrightarrow\mathcal K\cap(Q\times G)=\varnothing,\\
 Q\cap\{y=z=0\}\ne\varnothing&\Longrightarrow
 Q\leftarrow Q\cap(\R^3\times[-2^{-20},2^{-20}]^2).
 \end{aligned}
\]
Repeat the domain pass; $\mathcal K$ and the coordinate-volume
ordering are retained. No parameter subdivision or Fourier evaluation is used.

\paragraph{Moment bounds.}
Use the complete tails \eqref{eq:criticalHtail} and full-box spectral
floors of Lemmas~\ref{lem:transitionprofilecover} and~\ref{lem:criticalKZcover},
with outward binary64 radicals and Taylor remainders, for $t,r>0$.
The evaluation data are
\[
\begin{aligned}
 \text{old-coordinate value/scalar/gradient:}&\quad
 \operatorname{mid}P,\ \operatorname{rad}P,\ R_P,\quad\text{whole }P,\\
 \text{area value/gradient/Hessian/Newton:}&\quad
 p(\operatorname{mid}Q),\ (\rho_i)_{i=1}^5,\quad\text{whole }Q,
\end{aligned}
\]
where $\operatorname{mid}P\ne p(\operatorname{mid}Q)$ in general.
The floors
\[
 e^{-\sqrt{2/3}R_P}\|Q(\operatorname{mid}P)^{-1}\|_\infty^{-1},\qquad
 e^{-\sqrt{2/3}R_P}\|Q(\operatorname{mid}P)\|_\infty^{-1}
\]
combine by maximum with the native full-$P$ norm floors; norm upper
bounds round upwards. For $R_P\le9/20$, take $N=3$ and intersect
(T7) of Lemma~\ref{lem:transitionprofilecover} with
\[
 \left[\bigl(e^{-\sqrt{2/3}R_P}\bigr)^s\underline M_s,
       \bigl(e^{-\sqrt{2/3}R_P}\bigr)^{-s}\overline M_s\right]
\]
to obtain $[\mathcal M_s]$; evaluate the displayed powers outwards.

\paragraph{Moment band and supporting-parameter contraction.}
The whole-box exclusions are
\[
 \underline{\mathcal M_3}>43/2\ \lor\
 \underline{\mathcal M_6}>64\ \lor\
 \overline{\mathcal M_3}<b_0^-\ \lor\
 \underline{\mathcal M_3}>b_1^+
 \Longrightarrow\mathcal K\cap(Q\times G)=\varnothing.
 \tag{L7}
\]
With the complete endpoint intervals \eqref{eq:longshendpointmoments},
intersect $G$ with Lemma~\ref{lem:longshmomentcontraction}'s enclosure
for $[\mathcal M_3]$; exclude an empty intersection. Otherwise, for
$[l,u]$ and $w=\operatorname{RN}_{53}(u-l)$,
\[
 \begin{gathered}
 S\leftarrow1,\qquad
 S\leftarrow2S\quad\text{while }S<8192,\ wS<16,\\
 G\leftarrow G\cap[S^{-1}\lfloor Sl\rfloor,S^{-1}\lceil Su\rceil].
 \end{gathered}\tag{L8}
\]
Thus $2V_*^2\in G$ survives also for $V_*=V_0,V_1$; never bisect $G$.
Only for $b_0<E(\Lambda,3)<b_1$ impose
\[
 \partial_V\{f(V)+2V^2E(\Lambda,3)\}=0
 \quad\Longleftrightarrow\quad E(\Lambda,3)=E(H_{Y(V)},3).
\]
Outside this band, no parameter-derivative equation is imposed.

\paragraph{Additional conditional geometry contraction.}
After the domain pass and before terminal tests, apply
Lemma~\ref{lem:shprofilegeometry} whenever
\[
 11/10\le\frac{\sqrt3}{2t^3}\le217/50,\qquad
 |r-2/\sqrt3|,\ |x-1/2|,\ |y|,\ |z|\le1/100.
\]
For the outward image $[Y_-,Y_+]$ and transverse bounds
$(d_r,d_x,d_y,d_z)$ in (G2), use
$y_j\in[\underline y_j,\overline y_j]$ and
$\Gamma(y_j)\in[\Gamma_j^-,\Gamma_j^+]$. Binary search gives
\[
\begin{aligned}
 i_-&=\max\bigl(\{0\}\cup\{j\in\{0,\ldots,4096\}:\overline y_j\le Y_-\}\bigr),\\
 i_+&=\min\bigl(\{4096\}\cup\{j\in\{0,\ldots,4096\}:\underline y_j\ge Y_+\}\bigr).
\end{aligned}
\]
Under $E(U\mathbb Z^3,3)=E(H_{Y(V)},3)$, (G2) gives
\[
 G\leftarrow G\cap
 \left[\Gamma_{i_-}^--\frac{1275}{4}(d_y^2+d_z^2),\quad
       \Gamma_{i_+}^++\frac{1275}{4}(d_r^2+d_x^2)\right].
\]
No interpolation is used; only $G=\varnothing$ excludes.
Retain these uncoarsened bounds for (L15) and (G4); then apply (L8).

\paragraph{Direct rational area-coordinate jets.}
Use \eqref{eq:area-squared-lengths}--\eqref{eq:area-squared-length-jets}
and Lemma~\ref{lem:densityFCTcover} at $V=\sqrt{e/2}$,
$e\in\{\inf G,\sup G\}$, with covolume one.
In Lemma~\ref{lem:area-coordinate-jets}, use
$\partial_i=\partial/\partial\eta_i$,
$\partial_{ij}=\partial^2/(\partial\eta_i\partial\eta_j)$ and
$K'=\dd K/\dd q$, $K''=\dd^2K/\dd q^2$. Sum
\[
 2K'(q_n)\partial_iq_n,\qquad
 2\{K''(q_n)\partial_iq_n\partial_jq_n+K'(q_n)\partial_{ij}q_n\},
 \quad1\le i,j\le5,
 \tag{L11}
\]
and reciprocal terms over $\mathcal C_N^+$, using interval squares
on diagonals and common squared-length arguments before subtraction.
For each unique critical argument in $[1/4,8]$, full-parameter endpoint
signs and at most $60$ rational-kernel bisections retain its zero.
Add \eqref{eq:criticalHtail} on all of $P$, at the same $N$ and full-box
floors, transformed by \eqref{eq:area-derivative-tails}. Since $p(Q)\subset P$,
these tails cover all of $Q$, including non-Korkine--Zolotarev points;
finite terms use \eqref{eq:area-squared-lengths}--\eqref{eq:area-squared-length-jets}
and (L11).

\paragraph{Comparisons and stationary exclusions.}
Use the witnesses
\[
 H_{n/10},\qquad n=11,12,\ldots,43.
 \tag{L13}
\]
With complete $N=6$ positive moments and $e=\operatorname{mid}G$, start at
\[
 \min\operatorname*{argmin}_{n\in\{11,\ldots,43\}}
       \overline{E(H_{n/10},6)-eE(H_{n/10},3)}.
\]
Use it for the scalar test; rectangular tests visit all $33$ cyclically,
with one witness for both endpoints.
For $11/10\le Y\le6$, $B_t,\Gamma_t>0$ and
\[
 \begin{gathered}
 \partial_{\log Y}\Phi_V(H_Y)
 =(\Gamma(Y)-2V^2)\partial_{\log Y}E(H_Y,3),\\
 \Phi_V(H_{Y(V)})\le\Phi_V(H_{n/10}).
 \end{gathered}
\]
For the same supporting $V$,
\[
 \Phi_V(U\mathbb Z^3)>\Phi_V(H_{n/10})
 \quad\Longrightarrow\quad J(U\mathbb Z^3)>0.
\]

For the scalar metric test, use Lemma~\ref{lem:transitionprofilecover}
on $P$, $\operatorname{mid}P$, $R_P$: both occurrences of
$\max(0,28b_n^2-8eb_n)$ use the same lower squared-length endpoint,
with its zero-skip and complete omitted-correction guard.
The initial coordinate gradient uses all of $P$.
For $t,r>0$, (L5) gives
\[
 \det\left(\frac{\partial p_i}{\partial\eta_j}\right)_{1\le i,j\le5}
 =-\frac1{2r}\ne0,\qquad
 \nabla_\eta\Phi_V=
 \left(\frac{\partial p_i}{\partial\eta_j}\right)_{1\le i,j\le5}^{\mathsf T}
 \nabla_p\Phi_V,
\]
so $\nabla_\eta\Phi_V=0\iff\nabla_p\Phi_V=0$.

For one fixed witness $W$, compute at $p(\operatorname{mid}Q)$ the complete differences
\[
 \Delta_s=E(U(p(\operatorname{mid}Q))\Z^3,s)-E(W,s),\qquad s=3,6.
\]
With the area radii and 32-sign Hessian bounds $Q_e$ of
\eqref{eq:criticalrectvalue}, from (L11)--(L12),
\[
 \min_{e\in\{\inf G,\sup G\}}
 \underline{\Delta_6-e\Delta_3-\tfrac12Q_e}>0
 \Longrightarrow\mathcal K\cap(Q\times G)=\varnothing.
 \tag{L14}
\]
The complete endpoint hulls also give
\[
\begin{aligned}
 0\notin\operatorname{hull}_{e\in\{\inf G,\sup G\}}
       [\partial_i\Phi_{\sqrt{e/2}}(Q)]
 &\Longrightarrow\mathcal K\cap(Q\times G)=\varnothing,\\
 \max_{e\in\{\inf G,\sup G\}}\sup_{\eta\in Q}
       v^{\mathsf T}\nabla_\eta^2\Phi_{\sqrt{e/2}}v<0
 &\Longrightarrow\mathcal K\cap(Q\times G)=\varnothing
\end{aligned}
\]
for $1\le i\le5$ or one fixed dyadic $v$ at both endpoints, respectively.
For Newton, use \eqref{eq:criticalNewton} at $\operatorname{mid}Q$ with
full area-gradient/Hessian hulls; exclude empty intersections.
Lemma~\ref{lem:transitionprofilecover} and
\[
 \nabla_\eta^j\Phi_V=(1-\theta)\nabla_\eta^j\Phi_{V_-}
                    +\theta\nabla_\eta^j\Phi_{V_+},\quad j=0,1,2,
 \qquad V^2=(1-\theta)V_-^2+\theta V_+^2,\quad0\le\theta\le1
\]
give the complete enclosures for every $V\in[V_-,V_+]$.

\paragraph{Three terminal recognitions.}
On the complete rectangle and retained interval, test
\[
 \overline{V^2t^3}\le(9/10)^3,\qquad
 \underline{V^2/t^6}\ge64,\qquad V_->1.
 \tag{L15}
\]
For $(\eta(\Lambda),2V^2)\in\mathcal K$, Lemma~\ref{lem:longshheighttwoterminal} gives
$\Lambda\cong H_{Y(V)}$ and $J(\Lambda)=0$.
Use the uncoarsened conditional geometry bounds of
Lemma~\ref{lem:shprofilegeometry}.

Static Fourier recognition requires
\[
 \frac{11}{10}\le \frac{\sqrt3}{2t^3}\le\frac{217}{50},\qquad
 |r-2/\sqrt3|\le\frac1{200},\qquad
 \frac{99}{200}\le x\le\frac12,\qquad
 |y|,|z|\le\frac1{40}.
 \tag{L16}
\]
The outward guards and $E(U\mathbb Z^3,3)=E(H_{Y(V)},3)$ give
\[
 -\frac{17}{640}\le2V^2-\Gamma\!\left(\frac{\sqrt3}{2t^3}\right)
 \le\frac{51}{3200}.
\]
Lemma~\ref{lem:shfouriertransverse} and $\nabla_\eta\Phi_V=0$ give
\[
\begin{gathered}
 (r,x,y,z)=(2/\sqrt3,1/2,0,0),\\
 \partial_{\log Y}\Phi_V(H_Y)=0
 \Longrightarrow\Gamma(Y)=2V^2
 \Longrightarrow Y=Y(V),\quad J=0.
\end{gathered}
\]

For the four-dimensional Hessian test, after moment contraction require
\[
 11/10\le\frac{\sqrt3}{2t^3}\le217/50,\qquad
 |r-2/\sqrt3|,\ |x-1/2|,\ |y|,\ |z|\le1/100,
\]
and binary area/density widths
$\le\operatorname{RN}_{53}(1/200),\operatorname{RN}_{53}(1/5)$.
The widths only schedule the test. Set
\[
 W=\operatorname{hull}_{\mathrm{coord}}\!\left(
 \operatorname{proj}_{r,x,y,z}Q\cup\{(2/\sqrt3,1/2,0,0)\}\right).
\]
Use (L9)--(L11) and the complete $N=4$ tail on $Q_t\times W$.

Throughout $Q_t\times W$, for the primal and dual Gram matrices,
\[
 Q_0,Q_1\succ0,\quad Q_1=Q_0^{-1},\quad
 \|Q_0\|_\infty,\|Q_1\|_\infty\le25/8
 \Longrightarrow Q_0,Q_1\succeq(8/25)I.
\]
Differentiate (L9), with outward row-norm bounds
\[
 D_{\sigma a}\ge\sup_{Q_t\times W}\|\partial_aQ_\sigma\|_\infty,\qquad
 E_{\sigma ab}\ge\sup_{Q_t\times W}\|\partial_{ab}Q_\sigma\|_\infty,
 \qquad\sigma=0,1.
\]
For transverse $a,b$ and every index,
\[
 |\partial_a q_\sigma|\le D_{\sigma a}|n|^2,\qquad
 |\partial_{ab}q_\sigma|\le E_{\sigma ab}|n|^2.
\]
For $n\in\mathbb Z^3\setminus[-4,4]^3$, $q_\sigma(n)\ge8$.
For the unit-covolume primal/dual moment kernels $K_s$, $s=3,6$,
$K_s^{(j)}=\dd^jK_s/\dd q^j$. Their integral representations and
$\mathcal G_a(u)\le e^{-u}/(u-7)$, $a\le8$, give
\[
 |K_s^{(j)}(q)|<\frac{100e^{-\pi q}}q,
 \qquad j=1,2,\quad q\ge8.
\]
For $e\in\{2V_-^2,2V_+^2\}$, each omitted Hessian entry of
$E(U(\eta)\Z^3,6)-eE(U(\eta)\Z^3,3)$ has absolute bound
\[
 T_{ab}=100\frac{25}{8}(1+|e|)
 \sum_{\sigma=0}^1
       \{D_{\sigma a}D_{\sigma b}S_2+E_{\sigma ab}S_0\},
 \tag{L17}
\]
with the complete omitted cube and no extra opposite-pair factor.
For $j=0,2$, Lemma~\ref{lem:criticalKZcover} gives
\[
 \begin{aligned}
 S_j&\ge\sum_{n\in\Z^3\setminus[-4,4]^3}
              |n|^j e^{-8\pi|n|^2/25},\\
 S_0&\le\Theta_4(8/25),\qquad
 S_2\le25e^{-4\pi}\Theta_4(4/25).
 \end{aligned}
\]
Indeed $(u e^{-4\pi u/25})'<0$ for $u\ge25$.
Outward evaluation and diagonal interval squares make (L17) cover
every omitted primal/dual index.

Positive outward Schur pivots throughout $Q_t\times W$, including (L17), give
\[
 \begin{gathered}
 H_{V_-},H_{V_+}\succ0,\qquad
 H_V:=\nabla^2_{(r,x,y,z)}\Phi_V(U\Z^3),\\
 H_V=(1-\theta)H_{V_-}+\theta H_{V_+}\succ0,
 \quad V^2=(1-\theta)V_-^2+\theta V_+^2,\quad0\le\theta\le1.
 \end{gathered}
\]
For $(\eta(\Lambda),2V^2)\in\mathcal K$, SH symmetry and the positive
transverse Hessian give $(r,x,y,z)=(2/\sqrt3,1/2,0,0)$. Then
\[
 \nabla_\eta\Phi_V=0\Longrightarrow\Gamma(Y)=2V^2
 \Longrightarrow Y=Y(V),\qquad J(\Lambda)=0.
\]

\paragraph{Schedule, subdivision and retention rules.}
Initialize the contracted root (L1) with $G=[2V_0^2,2V_1^2]$.
Remove a box of greatest outward coordinate-volume upper bound,
using the binary max-heap and insertion/removal tie rules of
Lemma~\ref{lem:transitionprofilecover}; insert lower-coordinate children first.
On removal, apply in order
\[
 \begin{aligned}
 &\text{domain}\longrightarrow\text{guarded shift restriction}\longrightarrow\text{domain}\\
 &\longrightarrow\text{conditional geometry contraction}\longrightarrow(L15),(L16).
 \end{aligned}
\]
The ordered tests stop at the first terminal success; (L17) requires a
valid moment enclosure and supporting-parameter contraction:
\[
\begin{aligned}
 R_P\le9/20:\quad
 &(N=3,\operatorname{mid}P,\text{whole }P)\longrightarrow(L7),(L8)\\
 &\longrightarrow(L15),(L16)\longrightarrow\text{scheduled }(L17)\\
 &\longrightarrow\text{scalar test}\longrightarrow(N=4,\text{gradient on }P).
\end{aligned}
\]
For the nested area tests, use the center $p(\operatorname{mid}Q)$ and radii
$(\rho_i)_{i=1}^5$:
\[
\begin{aligned}
 R_P\le3/20:\quad
 &(N=3,\text{area-center moments})\\
 &\longrightarrow\underset{N=2,3,4}{\bigl[(L14)\longrightarrow
       \text{area-gradient hull}\longrightarrow\text{dyadic curvature}\bigr]}\\
 &\longrightarrow(N=4,\text{area-center moments})\longrightarrow(L14)\\
 &\longrightarrow(N=4,\text{endpoint jets})\longrightarrow\text{Newton}.
\end{aligned}
\]
Invalid endpoint Hessians stop the $N=2,3,4$ sequence and disable its
final Newton attempt. Use Lemma~\ref{lem:transitionprofilecover}'s fixed
dyadic direction, rounded inverse and complete interval checks;
failed enclosures disable the dependent tests. Retain nonempty Newton intersections by
Lemma~\ref{lem:criticalKZcover}'s outward-volume rule, never as an exclusion.

For $[r_H]\ni2/\sqrt3$, the SH-line intersection test is
\[
 Q_r\cap[r_H]\ne\varnothing,\qquad
 1/2\in Q_x,\quad0\in Q_y,\quad0\in Q_z.
\]
For an SH-line box, the $t$-priority rule is
\[
\begin{gathered}
 \relax[\sqrt3/(2t^3)]\cap[11/10,217/50]\ne\varnothing,\qquad
 [\sqrt3/(2t^3)]\not\subset[11/10,217/50],\\
 \operatorname{RN}_{53}(\sup Q_t-\inf Q_t)>\overline{[1]/[10000]},\qquad
 \inf Q_t<(\operatorname{mid}Q)_t<\sup Q_t\\
 \Longrightarrow\quad\text{bisect }t\text{ first}.
\end{gathered}
\]
Here $\overline{[1]/[10000]}$ is the upper endpoint of outward interval division.
For an SH-line box, if
$[\sqrt3/(2t^3)]\subset[11/10,217/50]$, use
\[
 \left(\frac{\max|Q_r-[r_H]|}{1/200},
       \frac{\max|Q_x-1/2|}{1/200},
       \frac{\max|Q_y|}{1/40},\frac{\max|Q_z|}{1/40}\right);
\]
otherwise use
\[
 \frac{(\max|Q_r-[r_H]|,\max|Q_x-1/2|,\max|Q_y|,\max|Q_z|)}{1/10000}.
\]
Round maxima upwards, retaining fractional divisions. If the maximum $>1$,
split the first maximizer in $(r,x,y,z)$; these fallback scores
never authorize recognition. Otherwise use
\[
 \frac{2\rho_t}{(\operatorname{mid}Q)_t},\qquad
 \frac{2\rho_r}{(\operatorname{mid}Q)_r},\qquad
 s_x(P),\ s_y(P),\ s_z(P),
\]
with the last three scores from Lemma~\ref{lem:transitionprofilecover}.
Bisect the first maximizer in $(t,r,x,y,z)$ at an interior binary
midpoint; both closed children inherit $G$. Retain as unresolved if
no such midpoint exists or depth is $120$. Newton preserves depth; never subdivide $G$.

\paragraph{The comparison invariant and finite outcome.}
Let $\mathcal R$ be the union of recognized products $Q\times G$.
Retain
\[
 \mathcal K\subset\mathcal R\cup
 \bigcup_{(Q,G)\in\mathcal Q\cup\mathcal U}(Q\times G).
\]
Initial inclusion is Lemma~\ref{lem:longprofileareadomain};
Lemma~\ref{lem:transitionprofilecover}, (L4)--(L17) and
\eqref{eq:criticalNewton} preserve it, including
\[
 Q^-\cup Q^+=Q,\qquad G^-=G^+=G
\]
under closed bisection. Recognition gives
\[
 (\eta(\Lambda),2V^2)\in\mathcal K\cap\mathcal R
 \Longrightarrow\Lambda\cong H_{Y(V)},\quad J(\Lambda)=0.
\]
The endpoint half-line identities give all exterior equality cases.
The finite record is
\[
\begin{aligned}
 (\#\mathrm{removals},\#\mathrm{bisections},\#\mathrm{retained\ Newton})
  &=(107013,52350,2621),\\
 \#\mathrm{domain\text{-}excluded\ children\ before\ insertion}&=309,\\
 (\#\mathrm{value},\#\mathrm{gradient})&=(27970,17423),\\
 (\#\mathrm{Newton\ exclusions},\#\mathrm{moment})&=(456,4046),\\
 (\#\mathrm{height},\#\mathrm{dynamic\ transverse},\#\mathrm{static\ Fourier})
  &=(232,5,23),\\
 (\#\mathrm{alignment\ exclusions},\#\mathrm{alignment\ contractions})
  &=(1887,61).
\end{aligned}
\]
The $61$ alignment contractions retain $y=z=0$ and create no children.
Thus
\[
 \mathcal Q=\mathcal U=\varnothing\Longrightarrow
 \begin{cases}
 J(\Lambda)\ge0&(|\Lambda|=1),\\
 J(\Lambda)=0\Longrightarrow
 \exists V\in[V_0,V_1]:\ \Lambda\cong H_{Y(V)}.
 \end{cases}
\]
Proposition~\ref{prop:density-gap-reduction} gives the comparison for every $V$.
This completes the proof of Numerical Lemma~\ref{num:longprofile-exhaustion}.
\end{proof}

\bigskip
\noindent
{\bf Acknowledgements.}
 The research of S. Luo is partially supported by the National Natural Science Foundation of China (NSFC) under Grant Nos. 12261045 and 12001253, and by the Jiangxi Jieqing Fund under Grant No. 20242BAB23001.
 The research of J. Wei is partially supported by the General Research Fund (GRF) of Hong Kong "New frontiers in singularity formation of nonlinear partial differential equations". The authors used AI to assist the long but finite computations.
\vskip0.1in
{\bf Statements and Declarations: there is no conflict of interest.}

{\bf Data availability: the manuscript has no associated data.}

\end{document}